\documentclass[11pt,reqno]{amsart}

\usepackage{mathtools,amssymb,amsthm}
\usepackage{enumitem}
\usepackage{xcolor}
\usepackage[margin=0.8in]{geometry}
\usepackage[final]{microtype}
\usepackage{amsrefs}
\usepackage[
  colorlinks=true,
  linkcolor=blue,
  citecolor=blue,
  urlcolor=blue
]{hyperref}

\mathtoolsset{showonlyrefs,showmanualtags}
\numberwithin{equation}{section}

\newtheorem{theorem}[equation]{Theorem}
\newtheorem{lemma}[equation]{Lemma}
\newtheorem{proposition}[equation]{Proposition}
\newtheorem{corollary}[equation]{Corollary}

\theoremstyle{definition}
\newtheorem{definition}[equation]{Definition}

\theoremstyle{remark}

\newcommand{\Log}{\operatorname{Log}}
\newcommand{\one}{\mathbf{1}}

\allowdisplaybreaks
\AtBeginDocument{%
  \setlength{\abovedisplayskip}{4pt plus 2pt minus 1pt}%
  \setlength{\belowdisplayskip}{4pt plus 2pt minus 1pt}%
  \setlength{\abovedisplayshortskip}{2pt plus 1pt minus 1pt}%
  \setlength{\belowdisplayshortskip}{3pt plus 1pt minus 1pt}%
  \setlength{\jot}{1.5pt}%
}

\setlist{
  leftmargin=*,
  topsep=2pt,
  itemsep=0pt,
  parsep=0pt,
  partopsep=0pt
}

\title[Prime Squares and Their Applications]
{On the Distribution of Prime Squares\\
and Its Applications}

\author[H. Mousavi]{Hamed Mousavi}
\address{Department of Mathematics, University of Bristol,
Bristol BS8 1QU, UK}
 \email{gj23799@bristol.ac.uk}
\date{}

\begin{document}

\begin{abstract}
We prove that, for every \(A>0\), all but
\(O_A(X(\log X)^{-A})\) integers \(n<X\) satisfying
\(n\equiv11\pmod{24}\) admit a representation
\[
 n=a^2+b^2+p^2,
\]
where \(p\) and \(a^2+b^2+1\) are prime and \((a,b)=1\).
For ergodic averages along prime squares, we prove a universal
divergence result for integrable functions, while establishing
almost-everywhere convergence for near-integrable functions.
We also quantify how the bounds used to establish convergence behave as one approaches integrability. Moreover, we  prove endpoint improving and sparse bounds for
the corresponding discrete averages, with an explicit Orlicz loss.
\end{abstract}

\maketitle
\tableofcontents

\section{Introduction}\label{s:intro}

The object of this paper is the sequence of prime squares
\[
 \mathbb P^2=\{p^2:p\in\mathbb P\}.
\]
We study it from three points of view: improving and sparse inequalities
for discrete averages, pointwise ergodic theory at or near the $L^1$ endpoint,
and an exceptional-set problem of Goldbach type.

These questions arise in different areas of analysis and number theory,
but three arithmetic features of prime squares recur.  First, the
distribution of primes in arithmetic progressions enters all three
applications.  In particular, for each fixed \(q\), it implies that
prime squares are asymptotically equidistributed among the nonzero
quadratic residues modulo \(q\).  Second, one needs cancellation in
quadratic exponential sums.  For
\[
 S_N(\alpha)=\sum_{p\leq N}(\log p)e(\alpha p^2),
\]
the behaviour near a rational point \(\alpha=a/q\) is described by
quadratic Gauss sums and primes in arithmetic progressions, while away
from such points one needs cancellation.  Finally, along suitable square-free moduli \(q\) for which the density
of the nonzero quadratic residues tends to zero, their normalized gaps
are asymptotically Poisson distributed: for every \(\theta>0\), the
proportion of consecutive gaps exceeding \(\theta\) times the mean gap
tends to \(e^{-\theta}\) \cite{GK08}.    The arguments also use more general arithmetic facts, including
the restrictions on the prime factorization of sums of two squares
and Dirichlet's approximation theorem.

\medskip

The sums \(S_N\) belong to a classical line of work in the
Waring--Goldbach problem going back to Vinogradov and Hua
\cite{Vinogradov,Hua}.  An influential  consequence is Hua's theorem that
every sufficiently large integer congruent to \(5\) modulo \(24\) is a
sum of five prime squares.  Subsequent work strengthened the relevant
exponential-sum estimates and treated prime squares drawn from thinner
sets of primes; see, for example,
\cite{Xiumin2005,B2016,GLZ22,SDP19}.  The estimates used below belong to
this circle-method framework.

\bigskip

We begin with discrete improving inequalities.  Define
\[
 A_{N,\Lambda}f(x)
 =\frac1N\sum_{n\leq N}\Lambda(n)f(x-n^2),
\]
where $\Lambda$ is the von Mangoldt function.  The normalization is
chosen so that $A_{N,\Lambda}$ has mass comparable to one and is
supported at spatial scale $N^2$.

The continuous model explains the expected exponents.  Convolution with
a bounded kernel sends $L^1$ to $L^\infty$, but averaging over a curved
measure also improves integrability even though the measure is
singular.  For arclength measure $d\sigma$ on $\mathbb S^1$, the
Littman--Strichartz estimate gives
\[
 \|f*d\sigma\|_{L^3(\mathbb R^2)}
 \lesssim \|f\|_{L^{3/2}(\mathbb R^2)};
\]
Such improving estimates arise from the curvature of the sphere and
the resulting decay of the Fourier transform of its surface measure.
The same ingredients play a central role in the Stein--Tomas restriction
theorem \cite{Littman,Strichartz,Tomas}.
  The exponents are sharp, as is seen by
taking $f$ to be the indicator of a thin annulus around the circle.

A useful comparison is provided by a Salem measure of full dimension
$\mu$ supported on $\mathbb S^1$.  If $\mu_j$ denotes its \(j\)-th
dyadic frequency piece, Fourier decay gives an \(L^2\) estimate for
convolution with $\mu_j$, while the Frostman estimate gives an
\(L^1\to L^\infty\) bound.  Applied to characteristic functions, these
estimates give
\[
 |\langle \mathbf1_F*\mu_j,\mathbf1_G\rangle|
 \lesssim
 \min\bigl\{2^{-j/2+o(j)}(|F||G|)^{1/2},
             2^{j+o(j)}|F||G|\bigr\}.
\]
Balancing the two estimates yields \(L^p\to L^{p'}\) improvement for
\(p>3/2\).  For arclength measure, curvature removes the \(o(j)\) losses,
so the endpoint \(p=3/2\) is also obtained.

The same critical exponent governs averages along the squares on
\(\mathbb Z\).  Indeed, up to spatial scale \(N^2\) there are about
\(N\) squares, so they form a set of dimension \(1/2\).  Han, Lacey,
and Yang proved the corresponding improving and sparse estimates for
\(p>3/2\) \cite{LHF}, and their argument can be refined to reach the
endpoint.  The endpoint distinction has little effect on the usual
ergodic applications.  In the
prime-square setting, however, establishing the improving and sparse estimates
at the endpoint is technically more delicate because the relevant
exponential-sum estimates carry additional arithmetic losses.  These
losses give rise to the Orlicz factor in the following theorem.

For $t>0$, write
\[
 \Log t=2+|\log t|,
 \qquad \Log\Log t=2+\log(2+|\log t|),
\]
and, for $c>0$, put
\[
 \Phi_c(t)=\exp\left(c\frac{(\Log t)^{1/2}}{\Log\Log t}\right).
\]
Thus \(\Phi_c(t)\) is large when \(t\) is small. Our endpoint estimate is as follows.

\begin{theorem}\label{t:improvingfixed}
Let \(I=[-N^2,2N^2]\cap\mathbb Z\).  For all \(F,G\subseteq I\),
\begin{align}\label{e:introimproving}
 \frac1{|I|}
 \bigl|\langle A_{N,\Lambda}\mathbf1_F,\mathbf1_G\rangle\bigr|
 &\ll
 \left(\frac{|F|}{|I|}\frac{|G|}{|I|}\right)^{2/3}
 \Phi_{4\sqrt3/9+o(1)}
 \left(\frac{|F|}{|I|}\frac{|G|}{|I|}\right).
\end{align}
Consequently, for every \(3/2<p\leq2\) and every
\(f\in\ell^p(\mathbb Z)\) supported in $I$,
\[
 \left(\frac1{|I|}\sum_{x\in\mathbb Z}
 |A_{N,\Lambda}f(x)|^{p'}\right)^{1/p'}
 \lesssim_p
 \left(\frac1{|I|}\sum_{x\in\mathbb Z}|f(x)|^p\right)^{1/p}.
\]
\end{theorem}

Here and below the notation $4\sqrt3/9+o(1)$ means that the constant
may be replaced by $4\sqrt3/9+\varepsilon$ for every fixed
$\varepsilon>0$.

Thus \eqref{e:introimproving} reaches the same critical pair
$(p,p')=(3/2,3)$ as the circle and the unweighted square averages, with
the arithmetic loss confined to the Orlicz factor.  In other words, up
to scale-dependent normalization, convolution with $A_{N,\Lambda}$
softens the ``singularities'' of a function on the integers more strongly
than convolution with a general Salem measure on
$\mathbb S^1\subset\mathbb R^2$.

The exponent $3/2$ cannot be improved.  Indeed, take
\[
 F=\{-p^2:p\leq N\},\qquad G=\{0\}.
\]
Then $|F|\asymp N/\log N$ and the prime number theorem gives
\begin{align*}
 \frac1{|I|}\langle A_{N,\Lambda}\mathbf1_F,\mathbf1_G\rangle
 &\asymp N^{-2}\asymp
 \left(\frac{|F||G|}{|I|^2}\right)^{2/3}
 \left(\Log\frac{|F||G|}{|I|^2}\right)^{2/3}.
\end{align*}
The example also shows that some endpoint loss is necessary.  There is
still a  gap between the logarithmic loss forced here and
the quasi-logarithmic function $\Phi_{4\sqrt3/9+o(1)}$ in the theorem.

Improving inequalities for discrete squares, primes, and polynomial
averages have been developed in
\cite{LHF,HLKF,LMR,HKLMY}; Theorem~\ref{t:improvingfixed} is the
prime-square endpoint of this circle of questions.

\medskip

A stronger and more directly applicable result concerns the associated
maximal operator, which controls the averages simultaneously over all
scales.  Put
\[
 A_\Lambda^*f(x)=\sup_{N\geq2}|A_{N,\Lambda}f(x)|.
\]
Maximal estimates along arithmetic sets originate in Bourgain's work on
polynomial averages and Wierdl's work on the primes
\cite{B88,B88b,B89,W89}; later developments include
\cite{MST,LHF,HLKF,HKLMY}. 
A more recent approach to such maximal estimates is sparse domination,
developed for discrete operators in \cite{KL20,H20,CDY,FN}.  We obtain
the endpoint sparse bound corresponding to \eqref{e:introimproving}.

\begin{theorem}\label{t:sparsebound}
Let \(F,G\subset\mathbb Z\) be finite.  There is a sparse collection of
intervals \(\mathcal S\) such that
\begin{align}\label{e:introsparse}
 \bigl|\langle A_\Lambda^*\mathbf1_F,\mathbf1_G\rangle\bigr|
 &\ll \sum_{I\in\mathcal S}|I|
 \left(\frac{|F\cap I||G\cap I|}{|I|^2}\right)^{2/3}
 \Phi_{4\sqrt3/9+o(1)}
 \left(\frac{|F\cap I||G\cap I|}{|I|^2}\right).
\end{align}
Moreover, if \((1/p,1/q)\) lies in the interior of the triangle with
vertices $(2/3,2/3),$ $(1,0)$, and $(0,1)$, 
then there is a sparse collection \(\mathcal S_{p,q}\) such that
\begin{align}\label{e:introsparseopen}
 \bigl|\langle A_\Lambda^*\mathbf1_F,\mathbf1_G\rangle\bigr|
 &\ll_{p,q}
 \sum_{I\in\mathcal S_{p,q}}|I|
 \left(\frac{|F\cap I|}{|I|}\right)^{1/p}
 \left(\frac{|G\cap I|}{|I|}\right)^{1/q}.
\end{align}
\end{theorem}

Here sparsity means that each \(I\in\mathcal S\) contains a set
\(E_I\subset I\) with \(|E_I|\geq |I|/10\), and the sets \(E_I\) are
pairwise disjoint.

\medskip

A contrasting endpoint phenomenon arises in ergodic theory.
The arithmetic distribution and spacing properties of prime squares
again play a role, but through a different mechanism.
Given a sequence of integers \((n_k)\), call it \(L^r\)-good if,
for every measure-preserving probability system \((X,\mu,T)\)
and every \(f\in L^r(X)\), the averages
\[
 \frac1N\sum_{k\leq N}f(T^{n_k}x)
\]
converge almost everywhere.
Birkhoff's theorem shows that \(n_k=k\) is \(L^r\)-good for every
\(r\geq1\) \cite{B31}.
Bourgain proved \(L^r\)-goodness for the squares and, more generally,
for integer polynomial sequences when \(r>1\), while Wierdl proved
the corresponding result for the primes \cite{B89,W89}.
Nair extended this to \(n_k=P(p_k)\) for every integer-valued
polynomial \(P\), and Mehlhop obtained a multidimensional extension
\cite{nair,melhop}.

The endpoint \(L^1\) behaves differently.
Buczolich and Mauldin proved that the squares are universally
\(L^1\)-bad \cite{BM10}.
LaVictoire proved persistent universal \(L^1\)-badness for the
monomials \(n_k=k^d\), \(d\geq2\), and for the primes
\cite{L09};
see also the weighted divergence result of Buczolich and Eisner
\cite{BuczolichEisner}.
Recent work of Fragkos, the author, and Stokolosa shows a related
endpoint failure for polynomial averages over short intervals
along every infinite subsequence
\cite{FragkosMousaviStokolosa}.
For prime squares, we use the distribution and spacing of the
squares of units modulo suitable square-free integers.
These properties yield a sequence of finite periodic
measure-preserving systems for which the prime-square maximal
operators have arbitrarily large weak \((1,1)\) constants.

There is also a quantitative question on the positive side.
The usual proof of pointwise convergence for \(L^p\) functions,
\(p>1\), combines convergence on a dense class with a maximal
inequality. As \(p\downarrow1\), the growth of the maximal
constant measures the loss in this approach to the endpoint.
For averages along the primes, Lacey, the author, and Rahimi
proved sparse bounds for indicator functions with factors
\((p-1)^{-2}\) unconditionally and \((p-1)^{-1}\) under GRH
\cite[Corollary~1.4]{LMR}.
Such growth alone does not imply \(L^1\)-divergence, but its rate
gives information beyond qualitative badness: it quantifies how
close to \(L^1\) the argument based on maximal estimates remains
effective and helps describe the classes of unbounded functions
for which convergence survives.
We establish both persistent divergence and an explicit maximal
bound for prime squares.

\begin{theorem}\label{t:pet}
Let \(p_k\) denote the \(k\)-th prime, and put \(n_k=p_k^2\).
\begin{enumerate}[label=\textup{(\roman*)}]
\item
The sequence \((n_k)\) is persistently universally \(L^1\)-bad.
More precisely, for every nonatomic ergodic measure-preserving
probability system \((X,\mu,T)\) and every infinite set
\(\mathcal S\subset\mathbb N\), there is an \(f\in L^1(X)\)
for which
\[
 \frac1N\sum_{k\leq N}f(T^{n_k}x),
 \qquad N\in\mathcal S,
\]
diverges on a set of positive measure.

\item
There is an absolute constant \(C>0\) such that, for every
invertible measure-preserving probability system \((X,\mu,T)\),
every \(1<p\leq2\), and every \(f\in L^p(X)\),
\[
 \left\|
 \sup_{N\geq1}
 \left|\frac1N\sum_{k\leq N}f(T^{n_k}x)\right|
 \right\|_{L^p(X)}
 \leq
 \exp\left(
 \frac{C}{p-1}\log\frac{e}{p-1}
 \right)
 \|f\|_{L^p(X)}.
\]
\end{enumerate}
\end{theorem}

We do not expect the upper bound in \textup{(ii)} to be sharp.
Determining the optimal dependence on \(p-1\) remains open,
and the available arguments do not suggest a precise conjecture
for its order of growth.
Section~\ref{s:orlicz} proves a stronger modular maximal estimate
and obtains convergence for an explicit Orlicz class near \(L^1\).

\medskip

We next turn to an additive application.  Let \(E(X)\) denote the
number of even integers \(n\leq X\) that are not sums of two primes.
The classical circle method gives $ E(X)=o(X),$
while Montgomery and Vaughan proved the quantitative bound $E(X)\ll X^{1-\delta}$
for some \(\delta>0\) \cite{MontgomeryVaughan}.
Goldbach problems in which one or more prime variables are restricted
to thin sets have a parallel history.  Linnik proved that the set
\[
 \mathcal P=
 \{a^2+b^2+1\in\mathbb P:(a,b)=1\}
\]
is infinite \cite{L60}.  Its distribution in short intervals and its
additive properties were subsequently studied by Wu, Matom\"aki, Tolev,
and others \cite{W98,Mat07,T99,T09,T09b,motamake}.  In particular,
Matom\"aki proved a binary Goldbach theorem with one prime in
\(\mathcal P\), with the number of exceptions up to \(X\) bounded by $ \ll X(\log X)^{-c}$
for some \(c>0\) \cite{motamake}.
Ter\"av\"ainen later proved an almost-all binary theorem and a ternary
Goldbach theorem with all prime variables in \(\mathcal P\)
\cite{Te16} with exceptional set of size $o(X)$.  Matom\"aki and Shao's work on Vinogradov's theorem with
almost twin primes provides another example of sieve and transference
methods in Goldbach problems with constrained prime variables
\cite{MatomakiShao}.

Our prime-square variant has a forced local condition.  If
$r=a^2+b^2+1\in\mathcal P$ and $r>3$, then $a$ and $b$ are odd, so
$r\equiv3\pmod8$.  Modulo $3$, exactly one of $a,b$ is divisible by
$3$, and hence $r\equiv2\pmod3$.  Thus $r\equiv11\pmod{24}$.  Since
$p^2\equiv1\pmod{24}$ for every prime $p>3$, one must have
$r+p^2\equiv12\pmod{24}$.  This is the only local progression relevant
to the theorem, apart from the harmless identity $12=3+3^2$.

\begin{theorem}[Quantitative exceptional set]\label{t:repmotthm}
Almost all integers congruent to $12$ modulo $24$ can be written as
\[
 a^2+b^2+p^2+1,
 \qquad a^2+b^2+1\in\mathcal P,\quad p\in\mathbb P.
\]
More precisely, if
\[
 E_2(X)=\{1\leq n<X:n\equiv12\pmod{24},\
 n\notin\mathcal P+\mathbb P^2\},
\]
then, for every $A>0$,
\[
 |E_2(X)|\ll_A X(\log X)^{-A}.
\]
\end{theorem}

The logarithmic saving of arbitrary order is an essential part of the
result.  While \(|E_2(X)|=o(X)\) would follow from substantially softer
input, the bound $ |E_2(X)|\ll_A X(\log X)^{-A}$
requires circle-method estimates that are uniform in the moduli and
auxiliary parameters arising from the sieve.

\medskip

Several questions remain open.  For a polynomial \(P\in\mathbb Z[x]\) of
degree at least two, one may ask for endpoint improving inequalities and sparse
bounds for averages along \(P(p)\), with the critical exponent dictated by the
degree and the geometry of \(P\).  Even for \(P(p)=p^2\), it would be
interesting to replace the loss \(\Phi_{4\sqrt3/9+o(1)}\) by a logarithmic
loss closer to the lower-bound scale
\((\Log(1/t))^{2/3}\).
At the \(L^1\) endpoint, a natural next problem is to prove persistent
universal \(L^1\)-badness for nonmonomial polynomial sequences
\(P(n)\).  The corresponding question for \(P(p)\) is the prime
analogue suggested by the present paper. On the
additive side, a power-saving estimate $|E_2(X)|\ll X^{1-\delta},$
or an asymptotic formula for the number of representations, would require
stronger uniform information than is presently used here.

\subsection{Ideas of the proofs}

The harmonic-analytic argument is organized around a general high--low
principle proved in Sections~\ref{s:improving} and \ref{s:sparse}.  In
the prime-square application, the model estimates have the form
\begin{align*}
 \|\operatorname{Low}_{Q,N}\|_{\ell^1\to\ell^\infty}
 &\lesssim N^{-2}Q\log Q,\\
 \|\operatorname{High}_{Q,N}\|_{\ell^2\to\ell^2}
 &\lesssim Q^{-1/2}
 \exp\left((2+o(1))\frac{\sqrt{\log Q}}{\log\log Q}\right).
\end{align*}
The same high estimate holds maximally in $N$.  Balancing these two
bounds on characteristic functions produces the power $2/3$ and the
Orlicz constant in Theorem~\ref{t:improvingfixed}.  This abstract form
of the argument also isolates the local low estimate needed by the
stopping-time proof of Theorem~\ref{t:sparsebound}.

The decomposition is proved in Section~\ref{s:hilo}.  A circle-method
approximation expresses the multiplier of the annular prime-square
average as a sum of terms localized near reduced fractions \(a/q\).
The ordinary main terms contain the normalized incomplete Gauss sums
\(\mathfrak s(a,q)\).  If a real exceptional character exists, Page's
theorem produces an additional term involving
\(\mathfrak s_{\chi_q}(a,q)\) and its exceptional zero\footnote{For the recent exclusion of Landau--Siegel zeros
through AI-assisted work, see
\cite{GreenSiegel2026,TaoSiegel2026} and the discussion in
Section~\ref{ss:notaion}.} \(\beta\).
We then apply the high--low decomposition at denominator \(Q\), using
pointwise kernel bounds for the low part and multifrequency and
square-function estimates for the high part.

For the main terms, our arithmetic input is
\[
 |\mathfrak s(a,q)|
 \leq q^{-1/2}
 \exp\left((2+o(1))
 \frac{\sqrt{\log q}}{\log\log q}\right).
\]
For incomplete quadratic Gauss sums over reduced residues, this
subexponential factor is part of the natural maximal estimate, up to the
value of the constant.  Since the high estimate is uniform in the
denominator, the same loss appears in the endpoint Orlicz factor.  Removing
it would require input beyond pointwise bounds for the individual Gauss
sums.

The exceptional-character term is too large to be absorbed into the
remainder and must itself be divided between the low and high parts.
Write \(q=q_0d\), where \(q_0\) is the exceptional conductor.  Two
complementary arithmetic estimates are used: a bound for the complete
sum over \(a\in\mathbb A_q\), suited to the pointwise kernel estimate,
and a bound for each coefficient
\(\mathfrak s_{\chi_q}(a,q)\), which decays like
\(d^{-1/2 + o(1)}\).
We then assign the exceptional terms to the low or high part according
to the size of \(d\).

The multifrequency maximal estimate requires the smoothing cutoff to be
fixed on each denominator scale, rather than depending on \(q\) or
\(d\).  We therefore replace the varying cutoff by a common one.  The
term with the common cutoff is controlled by the multifrequency maximal
estimate, and the difference between the two cutoffs by a square-function
estimate.
Together these estimates give
\[
 \|E_{Q,N}f\|_2
 \lesssim_A Q^{-1/2}(\log Q)^{O_A(1)}
 (1+\log N)^{-A}\|f\|_2,
\]
which is strong enough for the maximal stopping-time argument.  The
multifrequency and high--low framework follows the work on square and
prime averages in \cite{LHF,HLKF,LMR}; the exceptional-character
decomposition and the scale-dependent cutoffs are specific to the
prime-square setting.

\medskip

For the near-\(L^1\) estimate in Theorem~\ref{t:pet}\textup{(ii)},
we first prove a restricted weak-type bound with an Orlicz loss.
We seek a decomposition into a model with a maximal
\(\ell^p(\mathbb Z)\) bound for \(p\) close to one, and an error
with small maximal \(\ell^2(\mathbb Z)\) norm.
The circle-method approximation supplies the initial model,
while the Gauss-sum bounds and maximal \(\ell^2(\mathbb Z)\)
estimates from Section~\ref{s:hilo} control the error.
To estimate the model, we represent its arithmetic coefficients
by measures on finite cyclic groups.
Tao's quantitative Ionescu--Wainger multiplier and sampling
theorems \cite{Tao2021} then allow us to transfer continuous
estimates to \(\mathbb Z\) without estimating each rational
frequency separately. The formulation for canonical major arcs \cite{mirekcanonical} is unnecessary here,
and its stated bounds leave the dependence on \(p\) unspecified,
making Tao's quantitative formulation better suited to our
Orlicz estimates.

The additional work is to control the supremum over averaging
scales while tracking the constants from Ionescu--Wainger theory
and the arithmetic estimates.
We use dyadic increments and direct estimates for a finite range
of scales. For larger scales, we narrow the frequency cutoff and
apply the sampling theorem to a continuous maximal operator.
Fourier decay controls the errors introduced by changing the
cutoffs. We also divide the exceptional-character contributions
according to their conductors: the possible small-conductor
character is included in the model, while large-conductor terms
are controlled in the \(\ell^2\) error.
All estimates remain uniform in the possible Siegel zero.
Balancing the model and error estimates on indicators and
optimizing the parameters gives the restricted Orlicz bound.
Standard harmonic analysis arguments and Calder\'on's transference
principle then give the modular inequality, the stated dependence
on \(p-1\), and the corresponding ergodic estimates.

The proof of Theorem~\ref{t:pet}\textup{(i)} uses LaVictoire's criterion for
persistent universal \(L^1\)-badness \cite{L09}.  For an odd modulus
\(q\), let
\[
 \Lambda_q=\{b^2\pmod q:(b,q)=1\}.
\]
For each fixed \(q\), the prime number theorem in arithmetic
progressions shows that every element of \(\Lambda_q\) occurs with the
same limiting frequency among the prime squares.  Consequently, their
distribution modulo \(q\) converges to the normalized counting measure
on \(\Lambda_q\) as the averaging scale tends to infinity. By a theorem of Granville and
Kurlberg, along suitable square-free moduli the gaps between consecutive
elements of \(\Lambda_q\), after division by their mean, have the
exponential distribution associated with Poisson spacings
\cite{GK08}.  This gives the density, thickening, and difference
estimates required in LaVictoire's criterion \cite[Section 3]{L09}.  Conze's transference
principle and Sawyer's Banach principle then convert the resulting
finite cyclic maximal estimates into Theorem~\ref{t:pet}\textup{(i)}.

\medskip

For Theorem~\ref{t:repmotthm}, consider representations
\(n=r+p^2\), where \(r\) and \(p\) are prime and \(r\) is required to
belong to \(\mathcal P\).  By the elementary criterion for primitive
representations as a sum of two squares, this amounts to excluding
prime divisors of \(r-1\) that are congruent to \(3\pmod4\).

Fix a sifting threshold \(z\), and let \(A\) count the representations
for which \(r-1\) has no such prime divisor below \(z\).  Every desired
representation is counted by \(A\), but \(A\) may also count values of
\(r\) for which all the offending prime factors of \(r-1\) lie above
\(z\).  Since \(r\equiv11\pmod{24}\), the number of prime factors
congruent to \(3\pmod4\), counted with multiplicity, is even.  Thus any
such unwanted value contains at least a pair of these larger prime
factors.  A second count \(B\), which may overcount, is constructed to
include every representation arising from such a pair.  Consequently,
the desired number of representations is at least \(A-B\).
We apply a lower-bound sieve to \(A\) and an upper-bound sieve to \(B\).
It therefore remains to show that the main term supplied by the lower
sieve is larger than the correction supplied by the upper sieve.

The analytic input is a Bombieri--Vinogradov approximation, in mean
square over the represented integer $n$, for the equation
\[
 m^2+\ell t=n,
 \]
where $m,t$ are prime and $t$ lies in an arithmetic progression. On the minor arcs, a prime-square exponential-sum estimate and
Bessel's inequality reduce the required mean-square bound to a
combinatorial count.  The diagonal and off-diagonal contributions are
then estimated separately.

On the major arcs, Siegel--Walfisz and Bombieri--Vinogradov estimates
approximate both the quadratic prime sum, whose phase is
\(m^2\theta\), and the prime sum in an arithmetic progression, whose
phase is \(\ell t\theta\).  Here \(R=(\log N)^C\) is the major-arc
cutoff, so the initial approximation uses reduced fractions with
denominator \(q\leq R\).  We use the same parameter in a common smooth
cutoff with nonnegative inverse Fourier transform, producing a positive
archimedean factor \(H_{\ell,R}\).  The corresponding arithmetic factor
is the truncated local product \(G_R\). All of these quantities are defined precisely in Sections~\ref{s:app2trial} and \ref{s:sieveGB}.

The main new point is to pass from denominators \(q\leq R\) to all
\(R\)-smooth moduli.  In the corresponding argument of Matom\"aki
\cite{motamake}, stronger decay of the arithmetic coefficients permits
absolute summation.  For prime squares, the Gauss sums have only
\(q^{-1/2+o(1)}\) decay, which is not sufficient for this approach.
Instead, we retain the frequency localization. 
At each dyadic scale \(q\sim2^s\leq N/(100R^{1/2})\), the
localized multipliers centered at distinct reduced fractions
have disjoint supports. Plancherel therefore expresses the
squared global \(\ell^2(\mathbb Z)\) norm as the sum of the
squared norms of the individual frequency terms, using
frequency separation in the same spirit as the large sieve.
Consequently, the Gauss-sum coefficients are summed through their
squares rather than their absolute values, which provides sufficient
decay for the range \(R<q\leq \frac{N}{100R^{1/2}}\).

For the remaining \(R\)-smooth moduli \(q>\frac{N}{100R^{1/2}}\), a Rankin trick argument and
multiplicativity show that their total contribution is negligible.
Using the same parameter \(R\) in \(H_{\ell,R}\) and \(G_R\) preserves
the positive main term and the local factors needed in the sieve.

The sieve arguments of Wu and Matom\"aki proceed with a
half-dimensional lower sieve for \(A\) and a linear upper sieve for
\(B\).  The same direct approach does not work in our setting. For the small primes \(p\leq R\), the sieve-density products associated
with \(A\) and \(B\) do not satisfy the product estimates required for
sieve dimensions \(1/2\) and \(1\), respectively.  Consequently, a
single sieve over all primes \(p<z\) is justified only with the uniform
dimensions \(2/3\) for \(A\) and \(3/2\) for \(B\). These estimates make the lower bound for \(A\)
too small relative to the upper bound for \(B\), and the resulting
numerical lower bound for $A-B$ is not positive.

We therefore separate the small primes \(p\leq R\) from the
intermediate primes \(R<p<z\). The first range is removed by preliminary Rosser--Iwaniec weights of
level $ D_0=\exp\bigl((\log\log N)^3\bigr).
$
Since \(R=(\log N)^C\), one has
\[
 \frac{\log D_0}{\log R}\longrightarrow\infty,
 \qquad
 \log D_0=o(\log D(\ell))
\]
uniformly throughout the retained \(\ell\)-range, where \(D(\ell)\)
denotes the available sieve level.  Indeed, throughout this range
\(D(\ell)\geq N^\vartheta\) for some fixed \(\vartheta>0\).  The
fundamental lemma therefore gives the exact  local product over
\(p\leq R\) with relative error \(o(1)\).

We therefore first sieve the primes \(p\leq R\), taking $R=(\log N)^C$. 
The fundamental lemma of the Rosser--Iwaniec sieve shows that the
preliminary upper and lower sums both equal the product over
\(p\leq R\), up to a relative error \(O_J((\log N)^{-J})\) for every
fixed \(J>0\).  At the same time, division by \(D_0\) reduces the
remaining sieve level by only an \(N^{o(1)}\) factor and therefore
changes the outer sieve parameter by \(o(1)\).  Since the outer weights
have both signs, the upper or lower preliminary weights are selected
according to the sign of the corresponding outer weight.

The outer sieve is then applied only to the intermediate primes
\(R<p<z\).  In this range  the lower sieve for \(A\) has dimension \(1/2\), while
the upper sieve for \(B\) has dimension \(1\), matching the dimensions used in Matomäki's argument \cite{Mat07}.

To compare the resulting main terms, we sum the upper-sieve
contributions over \(\ell=2up_2\).  Here \(p_2\equiv3\pmod4\) is one
of the two large prime factors isolated above, while every prime factor
of \(u\) is congruent to \(1\pmod4\).  The local products over \(p\leq R\),
which have already been included in the preliminary sieve, cancel
explicitly against the corresponding factors in the lower-sieve main
term.  After this cancellation, the dependence on \(u\) is given by a
multiplicative weight.  Its mean value is estimated using the theorem
of Granville--Koukoulopoulos, while Shiu's theorem gives the required
uniformity in the ranges arising from the sum over \(p_2\)
\cite{GranvilleKoukoulopoulos,Shiu}.  Summation over \(p_2\) then
produces the Wu--Matom\"aki integral factor \(\mathcal W(\alpha)\).  The error
from replacing \(H_{\ell,R}\) by a uniform Archimedean multiplier \(H_\ell\) is summable over \(\ell\)
and is absorbed into the \(o(1)\) term.

It follows that the difference between the lower bound for \(A\) and
the upper bound for \(B\) has main coefficient
\[
 f_{1/2}(\alpha)
 -c_\alpha F_1(2)\mathcal W(\alpha),
 \qquad
 c_\alpha=\frac{16e^{-\gamma/2}}{\sqrt{\pi\alpha}}.
\]
Here \(f_\kappa\) and \(F_\kappa\) are the lower and upper
Rosser--Iwaniec sieve functions, respectively, \(\alpha\) is the sieve
parameter, and \(\gamma\) is Euler's constant.  For \(\alpha=9/8\),
\[
 f_{1/2}(9/8)
 -c_{9/8}F_1(2)\mathcal W(9/8)>0.17,
\]
so the lower-sieve contribution $A$ exceeds the correction term $B$. In summary, the
pre-sieving of \(p\leq R\) and the resulting weighted summation over
\(\ell\) are the principal differences from the earlier sieving arguments for
\(\mathcal P\) in \cite{W98,Mat07,motamake,Te16}.

\subsection{Organization of the paper}
Section~\ref{ss:notaion} fixes notation and records the standard
analytic and arithmetic inputs. Section~\ref{s:nt} proves the
quadratic Gauss-sum and correlation estimates used throughout.
Sections~\ref{s:improving} and \ref{s:sparse} establish the abstract
improving and sparse principles, and Section~\ref{s:hilo} proves the
prime-square high--low decomposition.
Section~\ref{s:orlicz} proves the Orlicz maximal inequalities and
quantifies the growth of the maximal \(L^p\) bounds as
\(p\downarrow1\). Section~\ref{s:l1} establishes persistent universal
\(L^1\)-badness for prime squares.
The Goldbach argument begins with the reductions in
Section~\ref{s:reductionsGB}; the circle-method approximation is
proved in Section~\ref{s:app2trial}, and the lower and upper sieves
are completed in Section~\ref{s:sieveGB}.

\subsection*{Acknowledgements}

We thank Ben Krause and Yu-Chen Sun for their help with editing
the manuscript, Michael Lacey for proofreading and edits, and
Joni Ter\"av\"ainen for proofreading, edits, and help with the
literature review. We also thank Yaghoub Rahimi, Anastasios
Fragkos, and Oleksiy Klurman, and acknowledge helpful discussions
with Jan Fornal on typical upper bounds for Gauss sums.

ChatGPT was used as an auxiliary editorial tool for proofreading and preliminary literature searches, and to generate an accompanying Lean development. This development, which can be found in \cite{MousaviPrimeSquaresLean2026}, aims to formalize the deductions within the paper conditional on thirteen explicitly documented families of external analytic and arithmetic inputs, which are not themselves formalized. All mathematical statements, references, arguments, and final wording in the manuscript were checked and approved by the author, who takes full responsibility for its contents. The Lean development is publicly available but has not yet undergone independent expert review.

\section{Notation and preliminaries}\label{ss:notaion}

This section collects the notation and standard analytic and arithmetic facts used throughout the paper.

\subsection{General conventions}

We write $e(t)=e^{2\pi it}$ and $\mathbb T=\mathbb R/\mathbb Z$, equipped
with normalized Lebesgue measure.  For $X\geq1$, denote $ [X]=\{1,\ldots,\lfloor X\rfloor\}.$
Intervals in $\mathbb Z$ are finite sets of consecutive integers. For an odd integer \(c\geq1\), \(cI\) denotes the interval
with the same centre as \(I\) and \(c\) times its length.  The indicator of a set $E$, or of a statement $E$, is
denoted by $\mathbf1_E$.  If $E$ is finite and nonempty, put
\begin{align*}
 \langle f\rangle_{E,p}
 =\left(\frac1{|E|}\sum_{x\in E}|f(x)|^p\right)^{1/p},
 \qquad
 \langle f\rangle_{E,\infty}=\sup_{x\in E}|f(x)|.
\end{align*}
All unnormalized $\ell^p$-norms and inner products are taken with
respect to counting measure:
\begin{align*}
 \|f\|_p=\left(\sum_{x\in\mathbb Z}|f(x)|^p\right)^{1/p},
 \qquad \langle f,g\rangle=\sum_{x\in\mathbb Z}f(x)\overline{g(x)}.
\end{align*}
As usual, $p'$ is determined by $1/p+1/p'=1$.

The relations $X\ll Y$, $Y\gg X$, and $X=O(Y)$ mean that
$|X|\leq CY$.  Dependence on parameters is indicated by subscripts,
as in $\ll_A$.  We write $X\sim Y$ if $X\ll Y\ll X$, and
$X\asymp Y$ if $X/Y\to1$.  The symbol $o(1)$ always refers to the main
asymptotic parameter tending to infinity.  The value of a positive
absolute constant may change from line to line.

For $x>0$, define
\begin{align}\label{e:Lognotation}
 \Log x=2+|\log x|,
 \qquad
 \Log\Log x=2+\log(2+|\log x|).
\end{align}
Thus $\Log x\asymp\log x$ for large $x$, while
$\Log x\asymp\log(x^{-1})$ for $0<x\leq1/2$.  For
$\alpha>0$ and $0<\beta<1$, put
\begin{align}\label{e:defPhi}
 \Phi_{\alpha,\beta}(x)
 =\exp\left(\alpha\frac{(\Log x)^\beta}{\Log\Log x}\right),
 \qquad \Phi_\alpha=\Phi_{\alpha,1/2}.
\end{align}

\subsection{Arithmetic notation}

When used as arithmetic variables, the letters \(p,p_1,p_2,\ldots\) denote primes; \(p\) may also denote an exponent in \(L^p\) or \(\ell^p\).  We write
$\mathbb P$ for the primes and $\mathbb P^2=\{p^2:p\in\mathbb P\}$.
The set
\begin{align}\label{e:defcalP}
 \mathcal P=\{x^2+y^2+1\in\mathbb P:(x,y)=1\}
\end{align}
is fixed throughout.  We use $(m,n)$ and $[m,n]$ for the greatest
common divisor and least common multiple.  The notation
$p^a\parallel n$ means that $p^a\mid n$ but $p^{a+1}\nmid n$.
Further, $P^+(n)$ is the largest prime factor of $n$, with
$P^+(1)=1$, $\omega(n)$ is the number of distinct prime factors, and
$\tau(n)=\sum_{m\mid n}1$. 
The M\"obius, Euler totient, and von Mangoldt functions are denoted by
$\mu$, $\phi$, and $\Lambda$, respectively.  Thus
\begin{align*}
 \phi(q)=q\prod_{p\mid q}\left(1-\frac1p\right),
 \qquad
 \Lambda(n)=\begin{cases}\log p,&n=p^j,\ j\geq1,\\0,&\text{otherwise}.
 \end{cases}
\end{align*}
We use without further comment the standard estimates
\begin{align}\label{e:totient}
 \phi(q)\gg\frac q{\log\log(3q)},\qquad
 \omega(q)\ll\frac{\log(2q)}{\log\log(3q)},\qquad
 \tau(q)\ll_\varepsilon q^\varepsilon.
\end{align}
Given a prime \(p\) and a nonzero integer \(n\), let
\[
 v_p(n)=\max\{\nu\geq0:p^\nu\mid n\}
\]
denote the \(p\)-adic valuation of \(n\).

For an odd prime $p$, $(\frac np)$ is the Legendre symbol; the Jacobi
symbol is used for odd composite moduli.  In particular the symbol is
zero when the numerator is not coprime to the modulus.

Put $\mathbb Z_q=\mathbb Z/q\mathbb Z$ and
\begin{align*}
 \mathbb A_q=\mathbb Z_q^\times
 =\{a\pmod q:(a,q)=1\}.
\end{align*}
Unless another range is displayed,
$\sum_{a\in\mathbb A_q}$ means a sum over $0\leq a<q$ with $(a,q)=1$.
The Ramanujan sum is
\begin{align}\label{e:tau}
 c_q(n)=\sum_{a\in\mathbb A_q}e(an/q).
\end{align}
We shall repeatedly use
\begin{align}\label{e:ramanujanformula}
 c_q(n)=\sum_{d\mid(q,n)}d\mu(q/d)
 =\mu\left(\frac q{(q,n)}\right)
   \frac{\phi(q)}{\phi(q/(q,n))}.
\end{align}
In particular, $c_q(n)=\mu(q)$ when $(q,n)=1$, and for prime powers
\begin{align}\label{e:ramprimepowerprelim}
 c_{p^a}(n)=
 \begin{cases}
  0,&p^{a-1}\nmid n,\\
  -p^{a-1},&p^{a-1}\parallel n,\\
  \phi(p^a),&p^a\mid n.
 \end{cases}
\end{align}
Cohen's identity, in the sign convention \eqref{e:tau}, is
\begin{align}\label{eq:Cohen}
 \sum_{r\in\mathbb A_q}c_q(n-r)=\mu(q)c_q(n);
\end{align}
see \cite[Corollary~7.2]{Cohen}.  The quadratic analogue required in
the paper is proved in Lemma~\ref{l:cohenquadratic}.

If $\chi$ is a Dirichlet character modulo $q$, extended by zero away
from the units, its conductor is denoted by $q^*$ and its inducing
primitive character by $\chi^*$.  We use the normalized linear Gauss
sum
\begin{align}\label{e:lineargauss}
 G(\chi,a)=\frac1{\phi(q)}
 \sum_{r\bmod q}\chi(r)e(ar/q).
\end{align}
If $\chi$ is primitive and $(a,q)=1$, then
\begin{align}\label{e:primitivegauss}
 G(\chi,a)=\overline{\chi}(a)G(\chi,1),\qquad
 |G(\chi,1)|=\frac{\sqrt q}{\phi(q)}.
\end{align}
If $\chi\bmod q$ is induced by a primitive character $\chi^*$ of
conductor $q_\chi$, we use the explicit formulas in
\cite[(3.3)--(3.5)]{trojan}, which express $G(\chi,a)$ in terms of
$\chi^*$, $q_\chi$, and $q/q_\chi$.

\subsection{Fourier analysis and arithmetic averages}

For \(f\in\ell^1(\mathbb Z)\), we use the Fourier transform
\begin{align}\label{e:FourierZ}
 \widehat f(\xi)
 =\sum_{x\in\mathbb Z}f(x)e(-x\xi),
 \qquad
 f(x)=\int_{\mathbb T}\widehat f(\xi)e(x\xi)\,d\xi.
\end{align}
For \(F\in L^1(\mathbb R)\), we write
\begin{align}\label{e:FourierR}
 \widehat F(\xi)
 =\int_{\mathbb R}F(x)e(-x\xi)\,dx,
 \qquad
 F^\vee(x)
 =\int_{\mathbb R}F(\xi)e(x\xi)\,d\xi.
\end{align}
Thus \(F^\vee\) denotes the inverse Euclidean Fourier transform. For \(f:\mathbb Z_q\to\mathbb C\), we use the normalization
\begin{align}\label{e:FTZp}
 \widehat f(a)
 =\frac1q\sum_{x\in\mathbb Z_q}f(x)e(-ax/q),
 \qquad
 f(x)=\sum_{a\in\mathbb Z_q}\widehat f(a)e(ax/q).
\end{align}
These normalizations give
\begin{align}\label{e:ramnujansumfourier}
 \frac1q\sum_{x\in\mathbb Z_q}c_q(x)e(-ax/q)
 =\mathbf1_{(a,q)=1}.
\end{align}
For functions on $\mathbb Z$, convolution is
$(f*g)(x)=\sum_yf(y)g(x-y)$.  We use Plancherel's identity and the
multiplier estimate
\begin{align}\label{e:plancherel}
 \|f\|_2=\|\widehat f\|_{L^2(\mathbb T)},\qquad
 \|T_mf\|_2\leq\|m\|_\infty\|f\|_2,\qquad
 \widehat{T_mf}=m\widehat f.
\end{align}
The two principal averaging operators are
\begin{align}\label{e:basicaverages}
 A_{N,\Lambda}f(x)&=\frac1N\sum_{n\leq N}\Lambda(n)f(x-n^2),\\
 A_{N,\mathbb N^2}f(x)&=\frac1N\sum_{n\leq N}f(x-n^2).
\end{align}
Their kernels are supported in an interval of length $N^2$.  Thus in
abstract improving or sparse statements the spatial scale corresponding
to these operators is $N^2$, not $N$.  For a family of operators
$T_N$, write
\begin{align*}
 T^*f(x)=\sup_N|T_Nf(x)|.
\end{align*}


\subsection{Primes in arithmetic progressions}

We  recall the Landau--Page theorem.  There is an absolute
$c_0>0$ such that, for every $Q\geq3$, among primitive real characters
of conductor at most $Q$ there is at most one character $\chi_0$ for
which $L(s,\chi_0)$ has a zero
\begin{align}\label{e:landaupagezero}
 \beta>1-\frac{c_0}{\log Q}.
\end{align}
If it exists, this zero is real and simple.  We call $\chi_0$ and
$\beta$ the exceptional character and exceptional zero at level $Q$.
This is the Landau--Page theorem; see
\cite[Exercise~54]{Tao2014Notes2}.

Recent work of OpenAI rules out Landau--Siegel zeros; see
Green's reformulation, obtained using ChatGPT,
\cite{GreenSiegel2026}, and Tao's exposition
\cite{TaoSiegel2026}.
For completeness, we retain the classical exceptional-character
and exceptional-zero conventions throughout, although the
exceptional case does not occur once \(c_0\) is chosen
sufficiently small.

We use the following consequence of Page's theorem and the prime number
theorem in arithmetic progressions.  Put
\[
 Q_x=\exp(c\sqrt{\log x}),
\]
where \(c>0\) is sufficiently small.  If an exceptional character
exists at level \(Q_x\), denote it by \(\chi_0\), its conductor by
\(q_0\), and its exceptional zero by \(\beta\).  Then, uniformly for
\(q\leq Q_x\) and \((a,q)=1\),
\begin{align}\label{e:landaupage}
 \sum_{\substack{n\leq x\\n\equiv a\pmod q}}\Lambda(n)
 =
 \frac{x}{\phi(q)}
 -\mathbf1_{q_0\mid q}
 \frac{\chi_0(a)x^\beta}{\beta\phi(q)}
 +O\left(xe^{-C\sqrt{\log x}}\right),
\end{align}
for some absolute \(C>0\).  If no exceptional character exists at
level \(Q_x\), the second term is omitted. This follows from the standard
zero-free region and the Landau--Page theorem; see
\cite[Theorem~5.27]{iwaniec}.  Some constants arising from
\eqref{e:landaupage} are ineffective, and we do not track effectivity.

Write \(\pi(x;q,a)=|\{p\leq x:p\equiv a\pmod q\}|\) and
\(\operatorname{Li}(x)=\int_2^x dt/\log t\). We shall also use the Brun--Titchmarsh inequality.  If \(x>q\) and
\((a,q)=1\), then
\[
 \pi(x;q,a)
 \leq
 \frac{2x}{\phi(q)\log(x/q)}.
\]
In particular, when \(q\leq x^{1/2}(\log x)^{-L}\), this gives
\[
 \pi(x;q,a)\ll\frac{x}{\phi(q)\log x}.
\]
We shall use the corresponding uniform consequence
\[
 \max_{\substack{(a,q)=1\\2\leq y\leq x}}
 \left|
 \pi(y;q,a)-\frac{\operatorname{Li}(y)}{\phi(q)}
 \right|
 \ll
 \frac{x}{\phi(q)\log x},
 \qquad
 q\leq x^{1/2}(\log x)^{-L}.
\]
Here the ranges \(y\leq x(\log x)^{-2}\) and
\(y>x(\log x)^{-2}\) are handled respectively by the trivial estimate
and the Brun--Titchmarsh inequality.

The Bombieri--Vinogradov theorem states that, for every \(A>0\) and \(K\geq0\), there is
\(B=B(A,K)>0\) such that
\begin{align}\label{e:BVstandard}
 \sum_{q\leq x^{1/2}(\log x)^{-B}}
 \tau(q)^K
 \max_{\substack{(a,q)=1\\2\leq y\leq x}}
 \left|
 \sum_{\substack{n\leq y\\n\equiv a\;(\mathrm{mod}\;q)}}\Lambda(n)
 -\frac{y}{\phi(q)}
 \right|
 \ll_{A,K}x(\log x)^{-A};
\end{align}
see \cite[Theorem~17.1]{iwaniec}. 

\subsection{Circle-method conventions}

Dirichlet's approximation theorem ensures that, for every
$\xi\in\mathbb T$ and $H\in \mathbb{N}$, there are
$a\in\mathbb A_q$ and $1\leq q\leq H$ such that
\begin{align}\label{e:Dirichletapprox}
 \left|\xi-\frac aq\right|\leq\frac1{qH}.
\end{align}
For quadratic averages with prime parameter $N$, the natural
major-arc width is $N^{-2}$.  Given $1\leq R\leq N$, we therefore set
\begin{align}\label{e:majorDef}
 \mathfrak M(R;N)
 &=\bigcup_{q\leq R}\bigcup_{a\in\mathbb A_q}
 \left\{\xi\in\mathbb T:
 \left|\xi-\frac aq\right|\leq\frac{R}{qN^2}\right\},\\
 \mathfrak m(R;N)&=\mathbb T\setminus\mathfrak M(R;N).
\end{align}
When a different width is required, it is displayed explicitly.  In
the unconditional high--low decomposition, the auxiliary parameter
$Q$ is restricted by
\begin{align}\label{e:defQ}
 \log Q\leq c(\log N)^{1/2}
\end{align}
with $c>0$ sufficiently small.  Under GRH a polynomial range is
available, but all results stated as unconditional use
\eqref{e:defQ}.

We shall also use the standard quadratic Weyl estimate.  If
$(a,q)=1$, $q\leq T$, and $|\lambda_2-a/q|\leq q^{-2}$, then, uniformly
in $\lambda_1$,
\begin{align}\label{vandercorput}
 \left|\frac1T\sum_{n\leq T}e(\lambda_2n^2+\lambda_1n)\right|
 \ll\left(\frac1q+\frac{\log(2q)}T
              +\frac{q\log(2q)}{T^2}\right)^{1/2},
\end{align}
see
\cite[Section~8.2]{iwaniec} for more detail. 

\subsection{Real-variable estimates}

We shall use two standard real-variable estimates.  The one-sided
Hardy--Littlewood maximal function is
\[
 M_{\mathrm{HL}}^+f(x)
 =
 \sup_{h>0}\frac1h\int_0^h|f(x-y)|\,dy.
\]
It satisfies
\[
 \|M_{\mathrm{HL}}^+f\|_{L^{1,\infty}(\mathbb R)}
 \ll\|f\|_{L^1(\mathbb R)},
 \qquad
 \|M_{\mathrm{HL}}^+f\|_{L^p(\mathbb R)}
 \ll\frac1{p-1}\|f\|_{L^p(\mathbb R)},
 \quad 1<p\leq2.
\]
The same bounds hold for a convolution operator \(T_Kf=K*f\)
whose \(L^2(\mathbb R)\)-operator norm is \(O(1)\) and whose
kernel satisfies
\[
 |K(x)|\ll |x|^{-1},
 \qquad
 \sup_{h\neq0}\int_{|x|>2|h|}
 |K(x-h)-K(x)|\,dx\ll1.
\]
These statements follow from the Hardy--Littlewood maximal theorem,
the Calderón--Zygmund decomposition, and interpolation; see
\cite[Chapter~I, Sections~3--4]{SteinHA}.

We shall use Khintchine's inequality in the following form.  If
\((\epsilon_j)\) are independent Rademacher variables, then, uniformly
for \(1\leq p\leq2\),
\begin{align}\label{e:khintchine}
 \left(\sum_j|z_j|^2\right)^{1/2}
 \ll
 \left(\mathbb E\left|
 \sum_j\epsilon_jz_j\right|^p\right)^{1/p}.
\end{align}
Consequently, for functions \(g_j\) on a measure space,
\[
 \left\|\left(\sum_j|g_j|^2\right)^{1/2}\right\|_{L^p}
 \ll
 \left(
 \mathbb E\left\|\sum_j\epsilon_jg_j\right\|_{L^p}^p
 \right)^{1/p}.
\]

\subsection{Sparse forms and stopping times}

We introduce the sparse forms and stopping-time linearization used in
the proof of Theorem~\ref{t:sparsebound}.

\begin{definition}\label{d:sparse}
A collection $\mathcal S$ of intervals is sparse if, for every
$I\in\mathcal S$, there is $E_I\subset I$ such that
$|E_I|\geq|I|/4$ and the sets $E_I$, $I\in\mathcal S$, are pairwise
disjoint.
\end{definition}

To linearize a dyadic maximal operator on an interval $I_0$, we use the
following stopping-time convention.

\begin{definition}\label{d:admiss}
Let $f=\mathbf1_F$ be supported in $3I_0$.  A map $\tau:I_0\longrightarrow\{2^j:1\leq2^j\leq|I_0|\}$
is admissible if, for every $x\in I_0$,
\begin{align}\label{d:admissible}
 \sup_{M\geq\tau(x)}\frac1M
 \sum_{|y-x|\leq CM}f(y)
 \leq C_{\rm adm}\frac{|F|}{|I_0|}.
\end{align}
where $C>1$ is a fixed radius constant and $C_{\rm adm}>1$ is a
sufficiently large fixed constant depending on $C$.
Changing the fixed radius changes only the required admissibility
constant. One may also replace the symmetric local average by the
maximum of the left and right one-sided local averages; these
quantities are comparable with a factor of at most two.
\end{definition}

Condition~\eqref{d:admissible} bounds the local density
\(|F\cap[x-CM,x+CM]|/M\) by the parent density
\(|F|/|I_0|\) at every scale \(M\geq\tau(x)\).


\subsection{Ionescu--Wainger notation}

We record the notation needed for the Ionescu--Wainger and sampling
theorems.  We use the one-dimensional specialization of
\cite[Definition~1.3]{Tao2021}.
For \(P\geq1\), let
\[
 \mathbb T[P]
 =\{\alpha\in\mathbb T:P\alpha=0\}
 =P^{-1}\mathbb Z/\mathbb Z.
\]
For \(\alpha\in\mathbb Q/\mathbb Z\), write
\[
 \operatorname{den}(\alpha)
 =\min\{q\geq1:q\alpha=0\}.
\]

A major-arc parameter set is a triple
\((k,\mathcal Q,\varepsilon)\), where \(k\geq1\),
\(\varepsilon>0\), and \(\mathcal Q\) is a finite set of pairwise
coprime positive integers.  For \(A\subseteq\mathcal Q\), put
\[
 Q_A=\prod_{q\in A}q,
 \qquad Q_\varnothing=1,
\]
and define
\begin{align}\label{e:IWSigma}
 \Sigma_{\leq k}
 =
 \bigcup_{\substack{A\subseteq\mathcal Q\\|A|\leq k}}
 \mathbb T[Q_A].
\end{align}
Thus, on putting $ P=\prod_{q\in\mathcal Q}q,$
we have \(\Sigma_{\leq k}\subseteq\mathbb T[P]\).  The associated
major-arc set is
\[
 \mathfrak M_{\mathrm{IW}}(k,\mathcal Q,\varepsilon)
 =
 \Sigma_{\leq k}+[-\varepsilon,\varepsilon]
 \subseteq\mathbb T.
\]
If \(r\geq1\), \(0<c<1\), and
\(q_{\max}\geq\max_{q\in\mathcal Q}q\), the parameter set is called
\((r,c)\)-good if
\begin{align}\label{e:IWgood}
 \varepsilon<
 \frac{c}{2r q_{\max}^{\,2rk}}.
\end{align}
We next introduce the adelic notation used in the sampling theorems.
Let
\[
 \widehat{\mathbb Z}
 =\varprojlim_{q}\mathbb Z_q,
 \qquad
 \mathbb A_{\mathbb Z}
 =\mathbb R\times\widehat{\mathbb Z}.
\]
The Haar measure on \(\widehat{\mathbb Z}\) is normalized to have
total mass one.  Its dual group is \(\mathbb Q/\mathbb Z\), and hence
the dual of \(\mathbb A_{\mathbb Z}\) is
\(\mathbb R\times\mathbb Q/\mathbb Z\).

In all applications below, the arithmetic Fourier support is contained
in \(\mathbb T[P]\).  The corresponding functions on
\(\mathbb A_{\mathbb Z}\) depend on the second variable only through
its residue class modulo \(P\), and may therefore be identified with
functions on $Y_P=\mathbb R\times\mathbb Z_P,$
equipped with Lebesgue measure in the first variable and normalized
counting measure in the second.
For \(F\in L^1(Y_P)\), define
\begin{align}\label{e:FourierYP}
 \widehat F(\theta,\alpha)
 &=
 \frac1P\sum_{y\in\mathbb Z_P}
 \int_{\mathbb R}
 F(x,y)e(-x\theta-y\alpha)\,dx,
 \qquad
 (\theta,\alpha)\in\mathbb R\times\mathbb T[P].
\end{align}
The inversion formula is
\[
 F(x,y)
 =
 \sum_{\alpha\in\mathbb T[P]}
 \int_{\mathbb R}
 \widehat F(\theta,\alpha)e(x\theta+y\alpha)\,d\theta.
\]
For a bounded function
\(M:\mathbb R\times\mathbb T[P]\to\mathbb C\), let \(T_M^{Y_P}\)
denote the multiplier operator on \(Y_P\) defined by
\[
 \widehat{T_M^{Y_P}F}(\theta,\alpha)
 =
 M(\theta,\alpha)\widehat F(\theta,\alpha).
\]
Suppose that \(M(\theta,\alpha)=0\) unless
\(\alpha\in\Sigma\subseteq\mathbb T[P]\) and
\(|\theta|\leq\varepsilon\).  The corresponding multiplier on
\(\mathbb Z\) is
\begin{align}\label{e:IWlatticemultiplier}
 M_{\Sigma}^{\mathbb Z}(\xi)
 =
 \sum_{\alpha\in\Sigma}\sum_{j\in\mathbb Z}
 M(\xi-\alpha+j,\alpha),
 \qquad \xi\in\mathbb T,
\end{align}
and \(T_{M;\Sigma}^{\mathbb Z}\) denotes the operator satisfying $ \widehat{T_{M;\Sigma}^{\mathbb Z}f}
 =
 M_{\Sigma}^{\mathbb Z}\widehat f.$
If \(M(\theta,\alpha)=m(\theta)\) is independent of \(\alpha\), we
write \(T_{m;\Sigma}^{\mathbb Z}\).
For a complex measure \(\nu\) on \(\mathbb Z_P\), put
\[
 \widehat\nu(\alpha)
 =
 \sum_{y\in\mathbb Z_P}e(-y\alpha)\nu(\{y\}),
 \qquad
 \|\nu\|_{\mathrm{TV}}
 =
 \sum_{y\in\mathbb Z_P}|\nu(\{y\})|.
\]
The symbols used below have the product form $M(\theta,\alpha)=m(\theta)\widehat\nu(\alpha).$

It remains to describe the sampling maps.  Define
\[
 \pi:\mathbb R\times\mathbb T[P]\longrightarrow\mathbb T,
 \qquad
 \pi(\theta,\alpha)=\theta+\alpha\pmod1.
\]
Suppose that \(\Sigma\subseteq\mathbb T[P]\) satisfies
\[
 \|\alpha-\alpha'\|_{\mathbb T}>2\varepsilon
 \qquad
 (\alpha,\alpha'\in\Sigma,\ \alpha\neq\alpha').
\]
Then the intervals $ \alpha+[-\varepsilon,\varepsilon]$ with $\alpha\in\Sigma$
are pairwise disjoint in \(\mathbb T\).  Consequently, every
\(\xi\in\pi([-\varepsilon,\varepsilon]\times\Sigma)\) has a unique
representation
\[
 \xi=\theta+\alpha\pmod1,
 \qquad
 |\theta|\leq\varepsilon,\quad \alpha\in\Sigma.
\]
We write
\[
 \ell^2(\mathbb Z)_{\pi(\Omega)}
 =
 \left\{f\in\ell^2(\mathbb Z):
 \operatorname{supp}\widehat f\subseteq\pi(\Omega)\right\},
 \qquad
 L^2(Y_P)_\Omega
 =
 \left\{F\in L^2(Y_P):
 \operatorname{supp}\widehat F\subseteq\Omega\right\},
\]
where $\Omega=[-\varepsilon,\varepsilon]\times\Sigma$ and $\pi(\Omega)=
 \bigcup_{\alpha\in\Sigma}
 \bigl(\alpha+[-\varepsilon,\varepsilon]\bigr).$

The inverse sampling map is
\begin{align}\label{e:IWinversesampling}
 \mathcal S_\Omega^{-1}f(x,y)
 =
 \sum_{\alpha\in\Sigma}
 \int_{-\varepsilon}^{\varepsilon}
 \widehat f(\alpha+\theta)
 e(x\theta+y\alpha)\,d\theta.
\end{align}
The diagonal sampling map is
\begin{align}\label{e:IWsampling}
 \mathcal S_PF(n)=F(n,n\bmod P),
 \qquad n\in\mathbb Z.
\end{align}
With these definitions, $ \mathcal S_P\mathcal S_\Omega^{-1}f=f$
for every \(f\) whose Fourier transform is supported in
\(\pi(\Omega)\).  The major-arc sampling theorem compares the
\(\ell^p(\mathbb Z)\)-norm of \(f\) with the
\(L^p(Y_P)\)-norm of \(\mathcal S_\Omega^{-1}f\), while the
quantitative Shannon sampling theorem controls
\(\mathcal S_PF\) when the real Fourier support of \(F\) is sufficiently
small relative to \(P^{-1}\).

\subsection{Dynamical systems and transference}
\label{ss:transference}

A measure-preserving system \((X,\mu,T)\) consists of a probability
space \((X,\mu)\) and a measurable map \(T:X\to X\) satisfying
\(\mu(T^{-1}E)=\mu(E)\) for every measurable \(E\).
We write \(T^nf=f\circ T^n\).
For an increasing convex function \(\Psi\) with \(\Psi(0)=0\)
and \(\Psi(t)\to\infty\), let \(L^\Psi(X)\) consist of the
measurable functions satisfying
\(\int_X\Psi(|f|/a)\,d\mu<\infty\) for some \(a>0\).
We call \(\Psi\) doubling if \(\Psi(2t)\ll\Psi(t)\).

For invertible \(T\), a convolution family \(A_Nf=K_N*f\)
with finitely supported kernels has ergodic counterparts
\(A_{N,T}f=\sum_mK_N(m)f\circ T^m\); write
\(A_T^*f=\sup_{N\geq2}|A_{N,T}f|\).
Calder\'on's transference principle \cite{Calderon1968}
transfers strong \(\ell^p(\mathbb Z)\) maximal estimates to
\(L^p(X)\), with the same constants.
The same argument applies to modular estimates: if
\(|\{A^*g>\lambda\}|\leq C\sum_m\Psi(|g(m)|/\lambda)\)
for every finitely supported \(g\) and every \(\lambda>0\), then
\begin{equation}\label{e:transference-modular}
 \mu\{A_T^*f>\lambda\}
 \leq C\int_X\Psi(|f|/\lambda)\,d\mu.
\end{equation}
Indeed, restrict to finitely many kernels supported in \([-R,R]\)
and test on
\(k\mapsto\mathbf1_{[-R,L-1+R]}(k)f(T^{-k}x)\).
Summing over \(0\leq k<L\) and integrating in \(x\) gives
\eqref{e:transference-modular} with the additional factor
\(1+2R/L\).  Let \(L\to\infty\), then enlarge the finite family.

\subsection{Sieve notation}

For a finite set $\mathcal F\subset\mathbb Z$, a set of primes
$\mathfrak P$, and $z\geq2$, put
\begin{align}\label{e:sievenotation}
 P_{\mathfrak P}(z)&=\prod_{\substack{p<z\\p\in\mathfrak P}}p,\\
 \mathcal S(\mathcal F,\mathfrak P,z)
 &=\{a\in\mathcal F:(a,P_{\mathfrak P}(z))=1\},\\
 S(\mathcal F,\mathfrak P,z)&=|\mathcal S(\mathcal F,\mathfrak P,z)|.
\end{align}
If $\mathcal F_d=\{a\in\mathcal F:d\mid a\}$ admits an approximation
\begin{align*}
 |\mathcal F_d|=X\frac{\omega_{\mathcal F}(d)}d+r_d
\end{align*}
for every square-free \(d\) all of whose prime divisors belong to
\(\mathfrak P\), then
\begin{align}\label{e:sieveproduct}
 V(z)=\prod_{\substack{p<z\\p\in\mathfrak P}}
 \left(1-\frac{\omega_{\mathcal F}(p)}p\right)
\end{align}
is the corresponding sieve product.  The dimension is the exponent
$\kappa$ for which, uniformly in the required range,
\begin{align}\label{e:sievedimension}
 \prod_{\substack{w\leq p<z\\p\in \mathfrak{P}}}\left(1-\frac{\omega_{\mathcal F}(p)}p\right)^{-1}
 \ll\left(\frac{\log z}{\log w}\right)^\kappa.
\end{align}
We denote by \(f_\kappa\) and \(F_\kappa\) the lower and upper
Rosser--Iwaniec sieve functions of dimension \(\kappa\), respectively,
with the normalization of \cite[(1.8)--(1.10)]{I80}.  The specific
functions used below are \(f_{1/2}\) and \(F_1\); their required values
and the relevant fundamental-lemma estimates are recorded in
Section~\ref{s:sieveGB}.

\section{Elementary number-theoretic estimates}\label{s:nt}



This section evaluates the local Ramanujan and quadratic exponential sums that occur in the major-arc and sieve arguments. These calculations determine both the low-frequency kernel in the high--low decomposition and the local densities in the sieve.

We first prove a quadratic analogue of Cohen's identity by reducing the
Ramanujan-sum correlation to prime powers and evaluating the odd and
$2$-adic factors separately. We then combine this identity with congruence restrictions to obtain the Euler products used later. Finally, we prove
pointwise bounds for the untwisted and character-twisted quadratic Gauss
sums. 
All residue-class sums in this section are over complete
systems of residues.

\subsection{Cohen's identity for square shifts}

The usual Cohen identity is
\begin{align*}
 \sum_{a\in\mathbb A_q}c_q(n-a)=\mu(q)c_q(n).
\end{align*}
We need the following quadratic version.

\begin{lemma}[Cohen's identity with square shifts]\label{l:cohenquadratic}
Let $2^a\parallel q$.  Then
\begin{align}\label{e:cohenquadratic}
 \sum_{r\in\mathbb A_q}c_q(n-r^2)
 =|\mu(q/2^a)| T_a(n)
 \prod_{\substack{p\mid q\\p>2}}
 \begin{cases}
  p+1,&\left(\frac np\right)=1,\\
  1-p,&\left(\frac np\right)\ne1,
 \end{cases}
\end{align}
where
\begin{align}\label{e:Ta}
 T_a(n)=
 \begin{cases}
  1,&a=0,\\
  (-1)^{n+1},&a=1,\\
  0,&a=2,3\text{ and }n\not\equiv1\pmod{2^{a-1}},\\
  \phi(2^a)c_{2^a}(n-1)=\pm\phi(2^a)^2,
    &a=2,3\text{ and }n\equiv1\pmod{2^{a-1}},\\
  0,&a\geq4.
 \end{cases}
\end{align}
Consequently,
\begin{align}\label{e:CohenQuadratic}
 \left|\sum_{r\in\mathbb A_q}c_q(n-r^2)\right|
 \ll |\mu(q/2^a)|\prod_{p\mid q}(p+1).
\end{align}
Here $(\frac np)$ is the Legendre symbol, with value zero when $p\mid n$.
\end{lemma}

\begin{proof}
Put
\begin{align*}
 C_q(n)=\sum_{r\in\mathbb A_q}c_q(n-r^2).
\end{align*}
If $(q_0,q_1)=1$, the Chinese remainder theorem implies that
$\mathbb A_{q_0q_1} \simeq \mathbb A_{q_0}\times\mathbb A_{q_1}$.
Since $c_{q_0q_1}(m)=c_{q_0}(m)c_{q_1}(m)$, and each factor depends only
on the corresponding component of $r$, we obtain
\begin{align}\label{e:Cmultiplicative}
 C_{q_0q_1}(n)=C_{q_0}(n)C_{q_1}(n).
\end{align}
It therefore suffices to evaluate $C_{p^k}(n)$.
For every prime $p$ and $k\geq1$,
\begin{align}\label{e:ramprimepower}
 c_{p^k}(m)=
 \begin{cases}
  0,&p^{k-1}\nmid m,\\
  -p^{k-1},&p^{k-1}\mid m, p^k\nmid m,\\
  \phi(p^k),&p^k\mid m.
 \end{cases}
\end{align}
Let first $p$ be odd.  For $k=1$, the congruence
$r^2\equiv n\pmod p$ has two solutions in $\mathbb A_p$ if
$(\frac np)=1$ and no solution otherwise.  Hence
\begin{align}\label{e:Cp}
 C_p(n)=
 \begin{cases}
 2(p-1)-(p-3)=p+1,&(\frac np)=1,\\
 -(p-1)=1-p,&(\frac np)\ne1.
 \end{cases}
\end{align}
Suppose that $k\geq2$.  If $(\frac np)\ne1$, then
$p^{k-1}\nmid n-r^2$ for every $r\in\mathbb A_{p^k}$, so
\eqref{e:ramprimepower} gives $C_{p^k}(n)=0$.  If
$(\frac np)=1$, the two roots modulo $p$ lift, by Hensel's lemma, to
exactly two roots modulo every $p^j$.  There are consequently two
classes modulo $p^k$ for which $p^k\mid n-r^2$, and
$2(p-1)$ further classes modulo $p^k$ for which
$p^{k-1}\parallel n-r^2$.  Thus
\begin{align*}
 C_{p^k}(n)=2\phi(p^k)-2(p-1)p^{k-1}=0.
\end{align*}
Together with \eqref{e:Cp}, this gives the odd part of
\eqref{e:cohenquadratic}.

For $a\leq3$, every odd square is $1$ modulo $2^a$, whence
\begin{align*}
 C_{2^a}(n)=\phi(2^a)c_{2^a}(n-1).
\end{align*}
Equation \eqref{e:ramprimepower} gives the first four cases in
\eqref{e:Ta}.  Finally let $a\geq4$.  The map
\begin{align*}
 r\longmapsto r+2^{a-2}\bar r\pmod{2^a},
\end{align*}
where $r\bar r\equiv1\pmod{2^a}$, is a bijection of
$\mathbb A_{2^a}$ and changes $r^2$ by $2^{a-1}$ modulo $2^a$.
Indeed, on each class $r_0\pmod{2^{a-2}}$ it sends the four lifts
$r_0+2^{a-2}j$ to the same four lifts in a permuted order, since
$\bar r$ is odd.
By \eqref{e:ramprimepower},
\begin{align*}
 c_{2^a}(m)+c_{2^a}(m-2^{a-1})=0
\end{align*}
for every $m$.  Pairing the terms proves $C_{2^a}(n)=0$.  The lemma
now follows from \eqref{e:Cmultiplicative}.
\end{proof}

\subsection{Two exponential sums and their correlations}

We shall frequently use the following two exponential sums, whose basic properties we establish in this subsection.
\begin{definition}[Exponential sums along squares]
For a Dirichlet character $\chi$ modulo $q$, put
\begin{align}\label{e:mathfraksgeneral}
 \mathfrak s_\chi(a,q)
 =\frac1{\phi(q)}\sum_{r\in\mathbb A_q}
 \chi(r)e\left(-\frac{ar^2}{q}\right).
\end{align}
For the principal character we write
\begin{align}\label{e:mathfraks}
 \mathfrak s(a,q)
 =\frac1{\phi(q)}\sum_{r\in\mathbb A_q}
 e\left(-\frac{ar^2}{q}\right).
\end{align}
\end{definition}

\begin{definition}
Let $g=(q,y)$ and $(b,g)=1$.  Define
\begin{align}\label{e:upsilondef}
 \Upsilon(a,q,b,y)=\frac{\phi(y)}{\phi([y,q])}
 \sum_{\substack{r\in\mathbb A_q\\r\equiv b\pmod g}}
 e\left(-\frac{ar}{q}\right).
\end{align}
\end{definition}
We shall use the following progression form of Cohen's identity.
\begin{lemma}
\label{l:progressioncohen}
Let $g\mid q$ and $(b,g)=1$. Then
\begin{align}\label{e:progcohenid}
 \sum_{\substack{r\in\mathbb A_q\\r\equiv b\pmod g}}c_q(n-r)
 =
 \begin{cases}
  0,&(g,q/g)>1,\\
  \mu(q/g)c_{q/g}(n)c_g(n-b),&(g,q/g)=1.
 \end{cases}
\end{align}
\end{lemma}

\begin{proof}
This is \cite[Lemma~11]{laceyprogression}.
\end{proof}

Next, we evaluate the correlation between $\Upsilon$ and $\mathfrak s$ over the reduced residue classes $\mathbb A_q$.
\begin{lemma}\label{l:sumoverramanujan}
Let $g=(q,y)$ and $(b,g)=1$.  Then
\begin{align}\label{e:eqjustforLHS}
 &\frac{\phi(q)\phi([y,q])}{\phi(y)}
 \sum_{a\in\mathbb A_q}\Upsilon(a,q,b,y)
 \mathfrak s(a,q)e\left(\frac{ax}{q}\right)\notag\\
 &\qquad=\mu(q/g)\mathbf1_{(g,q/g)=1}
 \sum_{u\in\mathbb A_q}c_{q/g}(x-u^2)c_g(x-b-u^2).
\end{align}
In particular, the right-hand side vanishes unless $q/g$ is square-free;
if $(g,q/g)=1$, it factors over the prime powers dividing $q$.
\end{lemma}

\begin{proof}
Opening \eqref{e:upsilondef} and \eqref{e:mathfraks}, and then summing
over $a\in\mathbb A_q$, gives
\begin{align*}
 &\frac{\phi(q)\phi([y,q])}{\phi(y)}
 \sum_{a\in\mathbb A_q}\Upsilon(a,q,b,y)
 \mathfrak s(a,q)e(ax/q)\\
 &=\sum_{u\in\mathbb A_q}
   \sum_{\substack{r\in\mathbb A_q\\r\equiv b\pmod g}}
   c_q(x-u^2-r).
\end{align*}
Apply Lemma~\ref{l:progressioncohen} with $n=x-u^2$ to obtain
\eqref{e:eqjustforLHS}.  If $q=q_0q_1$, $(q_0,q_1)=1$, and
$g=g_0g_1$ with $g_i\mid q_i$, the Chinese remainder theorem and the
multiplicativity of Ramanujan sums split the sum over $u$ into the
product of the corresponding sums modulo $q_0$ and $q_1$.  This proves
the final assertion.
\end{proof}

The preceding identity gives the Euler product needed later. 

\begin{lemma}\label{l:motamakemajorprimesquareprog}
Let $Q>2$, let $k\geq1$ be an integer, and let
$x,b\in\mathbb Z$ satisfy $(b,k)=1$.
Then
\begin{align}\label{e:corupss}
 &\sum_{P^+(q)<Q}\sum_{a\in\mathbb A_q}
 \Upsilon(a,q,b,k)\mathfrak s(a,q)e(ax/q)\notag\\
 &=\mathcal T_2
 \prod_{\substack{2<p<Q\\p\mid k\\(\frac{x-b}{p})=1}}
 \left(1+\frac{p+1}{p-1}\right)
 \prod_{\substack{2<p<Q\\p\mid k\\(\frac{x-b}{p})\ne1}}0 
 \prod_{\substack{2<p<Q\\p\nmid k\\(\frac{x}{p})=1}}
 \left(1-\frac{p+1}{(p-1)^2}\right)
 \prod_{\substack{2<p<Q\\p\nmid k\\(\frac{x}{p})\ne1}}
 \left(1+\frac1{p-1}\right),
\end{align}
where, writing \(h=v_2(k)\),
\begin{align}\label{e:T2factor}
 \mathcal T_2=
 \begin{cases}
  1+(-1)^x,&h=0,\\
  2^{\min\{h,3\}}
  \mathbf1_{\{x-b\equiv1\;(\mathrm{mod}\;2^{\min\{h,3\}})\}},
  &h\geq1.
 \end{cases}
\end{align}
Thus the product is zero if some odd prime $p<Q$ dividing $k$
satisfies $\bigl(\frac{x-b}{p}\bigr)\ne1$.
\end{lemma}

\begin{proof}
Apply Lemma~\ref{l:sumoverramanujan} with $y=k$. At every odd prime, only the exponents \(0\) and \(1\) contribute.
At \(2\), the exponent \(1\) may always contribute, while the exponents
\(2\) and \(3\) may contribute only when \(4\mid k\) and \(8\mid k\),
respectively; all higher exponents give zero.

For fixed $k,x,b$, the summand is multiplicative as a function of $q$.
Indeed, if $q=q_0q_1$ with $(q_0,q_1)=1$ and
$g_i=(q_i,k)$, then $g=(q,k)=g_0g_1$.  By
Lemma~\ref{l:sumoverramanujan}, the Ramanujan-sum correlation factors
into its components modulo $q_0$ and $q_1$.  Moreover,
\[
 \frac{\phi(k)}{\phi(q)\phi([k,q])}
 =
 \prod_{p^\nu\parallel q}
 \frac{\phi(p^{v_p(k)})}
 {\phi(p^\nu)\phi(p^{\max\{\nu,v_p(k)\}})},
\]
where a factor with $v_p(k)=0$ is interpreted in the evident way.
Hence the complete summand factors over the prime powers dividing $q$,
and the sum over $P^+(q)<Q$ is an Euler product.

For an odd prime $p\nmid k$, the local denominator is
$\phi(p)^2$, while Lemma~\ref{l:cohenquadratic} and the factor
$\mu(p)=-1$ give the local coefficient
\begin{align*}
 \begin{cases}
  -(p+1)/(p-1)^2,&(\frac xp)=1,\\
  1/(p-1),&(\frac xp)\ne1.
 \end{cases}
\end{align*}
For $p\mid k$, the local denominator is $\phi(p)$ and the local
coefficient is
\begin{align*}
 \begin{cases}
  (p+1)/(p-1),&(\frac{x-b}{p})=1,\\
  -1,&(\frac{x-b}{p})\ne1.
 \end{cases}
\end{align*}
Adding $1$ to each coefficient gives the four odd-prime Euler products 
in \eqref{e:corupss}. At \(p=2\), all terms with \(\nu\geq4\) vanish.  If \(2\nmid k\),
only the term with \(\nu=1\) may remain, so the contribution at \(2\)
is \(1+(-1)^x\).  If \(2\mid k\), the contribution at \(2\) is
\[
 1+\sum_{\nu=1}^{\min\{v_2(k),3\}}c_{2^\nu}(x-b-1),
\]
which equals \(\mathcal T_2\) by \eqref{e:Ta}.  Multiplying the
contributions over all primes proves the lemma.
\end{proof}

\begin{corollary}\label{c:generalmajorprimesquareprog}
Let $Q>2$, let $k\geq1$ be an integer, and let $\ell,n,b\in\mathbb Z$ satisfy $(b,k)=1$. Then
\begin{align}\label{e:generalcorupss}
 &\sum_{P^+(q)<Q}\sum_{a\in\mathbb A_q}
 \Upsilon(a\ell,q,b,k)\mathfrak s(a,q)e(an/q)\notag\\
 &\qquad=
 \mathcal T_2(\ell,k,n;b)
 \prod_{2<p<Q}
 \frac{p\phi((p,k))}{(p-1)^2}
 \left|\left\{(u,v)\in\mathbb A_p^2:
 \begin{array}{l}
  u^2+\ell v\equiv n\pmod p,\\
  v\equiv b\pmod{(p,k)}
 \end{array}\right\}\right|,
\end{align}
where 
\begin{align}\label{e:generalT2}
 \mathcal T_2(\ell,k,n;b)
 :=
 1+\sum_{\nu=1}^{3}
 \frac{\phi(k)}
 {\phi(2^\nu)\phi([k,2^\nu])}
 \sum_{\substack{u,v\in\mathbb A_{2^\nu}\\
 v\equiv b\;(\mathrm{mod}\,(2^\nu,k))}}
 c_{2^\nu}(n-u^2-\ell v).
\end{align}
\end{corollary}


\begin{proof}
Opening the definitions of \(\Upsilon\) and \(\mathfrak s\), and
summing first over \(a\), gives
\begin{align}\label{e:fixedqcoefficient}
 &\sum_{a\in\mathbb A_q}
 \Upsilon(a\ell,q,b,k)\mathfrak s(a,q)e(an/q)=
 \frac{\phi(k)}{\phi(q)\phi([k,q])}
 \sum_{\substack{u,v\in\mathbb A_q\\
 v\equiv b\;(\mathrm{mod}\,(q,k))}}
 c_q(n-u^2-\ell v).
\end{align}
This expression is multiplicative in \(q\). Let \(p>2\), and denote by \(N_p\) the
cardinality appearing on the right-hand side of
\eqref{e:generalcorupss}.  The quadratic sum modulo \(p^\nu\) vanishes
for \(\nu\geq2\), so only the modulus \(p\) contributes.  Since
\[
 c_p(m)=p\mathbf1_{p\mid m}-1,
 \qquad
 \frac{\phi(k)}{\phi(p)\phi([k,p])}
 =\frac{\phi((p,k))}{(p-1)^2},
\]
and the number of pairs satisfying only
\(v\equiv b\pmod{(p,k)}\) is
\((p-1)^2/\phi((p,k))\), we obtain
\begin{align*}
 1+\frac{\phi(k)}{\phi(p)\phi([k,p])}
 \sum_{\substack{u,v\in\mathbb A_p\\
 v\equiv b\pmod{(p,k)}}}
 c_p(n-u^2-\ell v)
 &=
 1+\frac{\phi((p,k))}{(p-1)^2}
 \left(pN_p-\frac{(p-1)^2}{\phi((p,k))}\right)\\
 &=
 \frac{p\phi((p,k))}{(p-1)^2}N_p.
\end{align*}
At \(p=2\), the terms with modulus \(2^\nu\), \(\nu\geq4\), vanish,
while the sum of the remaining terms is
\(\mathcal T_2(\ell,k,n;b)\).  Multiplying the contributions over all
primes proves \eqref{e:generalcorupss}.
\end{proof}

\begin{lemma}\label{l:lowramanujan}
For $Q\geq3$ and every integer $n$,
\begin{align}
 \sum_{P^+(q)<Q}\frac1{q\phi(q)}
 \left|\sum_{r\in\mathbb A_q}c_q(n-r^2)\right|
 \ll\log Q.
\end{align}
\end{lemma}

\begin{proof}
Write $q=2^aq_1$ with $q_1$ odd.  By
Lemma~\ref{l:cohenquadratic}, the summand vanishes unless $a\leq3$ and
$q_1$ is square-free.  The contribution from $2^a$ is an absolute constant,
and hence
\begin{align*}
 \sum_{P^+(q)<Q}\frac1{q\phi(q)}
 \left|\sum_{r\in\mathbb A_q}c_q(n-r^2)\right|
 &\ll\prod_{2<p<Q}\left(1+\frac{p+1}{p(p-1)}\right)\\
 &\ll\exp\left(\sum_{p<Q}\frac1p+O(1)\right)
 \ll\log Q,
\end{align*}
by Mertens' theorem.
\end{proof}

\subsection{Pointwise Gauss-sum bounds}

We now record the pointwise bounds for the quadratic Gauss sums appearing in the major-arc multipliers.
\begin{lemma}\label{l:inftyboundsaq}
If $(a,q)=1$, then
\begin{align}\label{e:gaussboundknown}
 |\mathfrak s(a,q)|
 \leq q^{-1/2}
 \exp\left((2+o(1))\frac{(\log q)^{1/2}}{\log\log q}\right).
\end{align}
In particular, the right-hand side may be replaced by
$q^{-1/2}\exp(3(\log q)^{1/2}/\log\log q)$ for all sufficiently large
$q$.
\end{lemma}

\begin{proof}
The normalized sum is multiplicative in the modulus.  Let $p$ be odd.
Writing
\begin{align*}
 G(a;p^j)=\sum_{r\bmod p^j}e(-ar^2/p^j),
\end{align*}
the classical quadratic Gauss-sum evaluation gives
$|G(a;p^j)|=p^{j/2}$ and $G(a;p^j)=pG(a;p^{j-2})$ for $j\geq2$;
see, for example, \cite[Chapter~3]{iwaniec}.  Therefore
\begin{align}\label{e:primepowerestimate}
 \sum_{r\in\mathbb A_{p^j}}e(-ar^2/p^j)
 =G(a;p^j)-pG(a;p^{j-2})=0,\qquad j\geq2.
\end{align}
For $j=1$ its absolute value is at most $\sqrt p+1$.  The same
calculation at $2$ shows that the unit sum vanishes for $2^j$ with
$j\geq4$; the three remaining powers contribute only an absolute
constant.  Consequently $\mathfrak s(a,q)=0$ unless
$q=2^jq_0$, $0\leq j\leq3$, and $q_0$ is odd and square-free, and in
the latter case
\begin{align}\label{e:sproductbound}
 |\mathfrak s(a,q)|
 \ll q^{-1/2}\prod_{p\mid q_0}(1-p^{-1/2})^{-1}.
\end{align}
Taking logarithms, the additional factor is
\begin{align*}
 \exp\left(\sum_{p\mid q_0}p^{-1/2}+O\left(\sum_{p\mid q_0}p^{-1}\right)\right).
\end{align*}
The sum is maximal, subject to $\prod_{p\mid q_0}p\leq q$, by taking
the smallest primes.  The prime number theorem and partial summation
give
\begin{align*}
 \sum_{p\mid q_0}p^{-1/2}
 \leq(2+o(1))\frac{(\log q)^{1/2}}{\log\log q},\qquad
 \sum_{p\mid q_0}p^{-1}=O(\log\log\log q).
\end{align*}
Substitution in \eqref{e:sproductbound} proves the assertion.
\end{proof}

\begin{lemma}\label{l:Gcharbound}
Let \(\chi\) be a real Dirichlet character modulo \(q\), and let
\((a,q)=1\).  Then, for every \(\varepsilon>0\),
\begin{align}\label{e:Gcharbound}
 |\mathfrak s_\chi(a,q)|\ll_\varepsilon q^{-1/2+\varepsilon}.
\end{align}
Equivalently, \( |\mathfrak s_\chi(a,q)|\ll q^{-1/2+o(1)}\).
\end{lemma}

\begin{proof}
By the Chinese remainder theorem the numerator of \(\mathfrak s_\chi(a,q)\)
factors over the prime powers \(p^\nu\parallel q\).  Thus it is enough to
prove, uniformly for \((b,p)=1\),
\[
 \left|\sum_{r\in\mathbb A_{p^\nu}}\chi_p(r)e(-br^2/p^\nu)\right|
 \ll C p^{\nu/2},
\]
where \(\chi_p\) is the local component of \(\chi\) and \(C\) is absolute.

Let \(p\) be odd.  If \(\chi_p\) is principal, then
\[
 \sum_{r\in\mathbb A_{p^\nu}}e(-br^2/p^\nu)
 =
 \sum_{r\bmod p^\nu}e(-br^2/p^\nu)
 -
 \sum_{\substack{r\bmod p^\nu\\p\mid r}}e(-br^2/p^\nu).
\]
The first sum is a complete quadratic Gauss sum modulo \(p^\nu\), and the
second is \(p\) times a complete quadratic Gauss sum modulo \(p^{\nu-2}\)
when \(\nu\geq2\).  Hence both are \(O(p^{\nu/2})\).  If \(\chi_p\) is
non-principal, then \(\chi_p(r)=(r/p)\).  For \(\nu=1\), Weil's bound gives
\[
 \left|\sum_{r\bmod p}\left(\frac rp\right)e(-br^2/p)\right|
 \ll p^{1/2}.
\]
For \(\nu\geq2\), write \(r=u+p^{\nu-1}t\), with
\(u\bmod p^{\nu-1}\), \(p\nmid u\), and \(t\bmod p\).  Since \(\chi_p\) is
unchanged under this replacement,
\[
 (u+p^{\nu-1}t)^2\equiv u^2+2up^{\nu-1}t\pmod {p^\nu},
\]
and the inner sum over \(t\) is
\[
 \sum_{t\bmod p}e(-2but/p)=0.
\]
Thus the local sum vanishes for \(\nu\geq2\).

For \(p=2\), every real character modulo \(2^\nu\) is induced from a
character of conductor dividing \(8\).  For \(\nu\leq5\) the desired bound is
trivial.  For \(\nu\geq6\), write \(r=u+2^{\nu-3}t\), with
\(u\bmod 2^{\nu-3}\), \(u\) odd, and \(t\bmod 8\).  The character is
unchanged, and
\[
 (u+2^{\nu-3}t)^2\equiv u^2+u2^{\nu-2}t\pmod {2^\nu}.
\]
Therefore the inner sum over \(t\) is
\[
 \sum_{t\bmod 8}e(-but/4)=0,
\]
since \(bu\) is odd.  Hence the local sum is zero for \(\nu\geq6\), and the
bound follows also at \(p=2\).

Multiplying the local estimates gives
\[
 \left|\sum_{r\in\mathbb A_q}\chi(r)e(-ar^2/q)\right|
 \ll C^{\omega(q)}q^{1/2}.
\]
Dividing by \(\phi(q)\), we obtain
\begin{align}\label{e:actualboundfors}
     |\mathfrak s_\chi(a,q)|
 \ll C^{\omega(q)}q^{1/2}\phi(q)^{-1}
 \ll C^{\omega(q)}q^{-1/2}\frac q{\phi(q)}
 \ll_\varepsilon q^{-1/2+\varepsilon},
\end{align}
since \(C^{\omega(q)}q/\phi(q)\ll_\varepsilon q^\varepsilon\).  This proves
\eqref{e:Gcharbound}.
\end{proof}

We shall use the following refinement of Lemma~\ref{l:Gcharbound} in
Section~\ref{s:hilo}.

\begin{lemma}\label{l:exceptionalarithmetic}
Let \(x\in\mathbb Z\). Let $\chi_0$ be a primitive real character of conductor $q_0$, and for
each multiple $q$ of $q_0$ let $\chi_q$ denote the character modulo $q$
induced by $\chi_0$.  For \(0\leq j\leq5\), let
\(\chi_{2^jq_0}\) have the same meaning, and set
\[
 \sigma_0:=
 \max_{0\leq j\leq5}
 \sup_{a\in\mathbb A_{2^jq_0}}
 |\mathfrak s_{\chi_{2^jq_0}}(a,2^jq_0)|,
 \qquad
 \rho(d):=\prod_{p\mid d}(1-p^{-1/2})^{-1}.
\]
There exists $\mathcal B_0\geq1$, depending only on $q_0$, such that,
for every \(q=q_0d\),
\begin{align}\label{e:exceptionalarithmetic}
 \left|
 \sum_{a\in\mathbb A_q}
 \mathfrak s_{\chi_q}(a,q)e(ax/q)
 \right|
 &\ll
 \mathcal B_0\prod_{p\mid d}\frac{p+1}{p-1},\notag\\
 |\mathfrak s_{\chi_q}(a,q)|
 &\ll \sigma_0d^{-1/2}\rho(d),
 \qquad \textup{ for } a\in \mathbb{A}_q\\
  \sigma_0\mathcal B_0^{1/2}&\ll q_0^{-1/25}.
\end{align}
In particular,
\[
 (1+\log q_0)\sigma_0\mathcal B_0^{1/2}\ll1.
\]
Moreover, if
\(a\mapsto\mathfrak s_{\chi_q}(a,q)\) is not identically zero, then
\begin{align}\label{e:sqfreeconds}
    d=2^\nu d_1,\qquad
 0\leq\nu\leq5,\qquad
 \mu^2(d_1)=1,\qquad
 (d_1,2q_0)=1.
\end{align}
\end{lemma}
\begin{proof}
The first estimate in \eqref{e:exceptionalarithmetic} follows by factoring the complete sum over prime powers. The second estimate uses the vanishing of induced quadratic sums at higher prime powers. The final bound on \(\sigma_0\mathcal B_0^{1/2}\) comes from the restriction \(p\equiv1\pmod4\) on the odd primes in the conductor.

The Chinese remainder theorem factors both expressions in
\eqref{e:exceptionalarithmetic} over prime powers.  If $p\nmid q_0$,
the local factors in the first and second estimates are bounded,
respectively, by
\[
 \frac{p+1}{p-1}
 \qquad\text{and}\qquad
 p^{-1/2}(1-p^{-1/2})^{-1}.
\]
If \(p\mid q_0\) is odd and 
\(\chi_p(-1)=-1\), then the terms indexed by \(r\) and \(-r\) cancel.
Hence
\[
 \sum_{r\in\mathbb A_p}\chi_p(r)e(-br^2/p)=0.
\]
Thus a nonzero local quadratic sum requires \(\chi_p(-1)=1\), which for
the quadratic character modulo \(p\) is equivalent to
\(p\equiv1\pmod4\).  Weil's bound gives
\[
 \left|
 \sum_{r\in\mathbb A_p}\chi_p(r)e(-br^2/p)
 \right|
 \leq 2p^{1/2},
\]
so the normalized local factor in the second estimate of
\eqref{e:exceptionalarithmetic} is at most
\(2p^{1/2}/(p-1)\).
For the first estimate, Fourier inversion gives the local factor
\[
 \frac1{p-1}\sum_{r\in\mathbb A_p}
 \chi_p(r)c_p(x-r^2).
\]
Since \(\sum_{r\in\mathbb A_p}\chi_p(r)=0\) and
\[
 c_p(n)=
 \begin{cases}
  p-1,&p\mid n,\\
  -1,&p\nmid n,
 \end{cases}
\]
this equals
\[
 \frac{p}{p-1}
 \sum_{\substack{r\in\mathbb A_p\\r^2\equiv x\;(\mathrm{mod}\;p)}}
 \chi_p(r).
\]
The congruence \(r^2\equiv x\pmod p\) has at most two solutions, so this
factor has absolute value at most \(2p/(p-1)\).

At \(p=2\), every primitive real character has conductor dividing \(8\).
Moreover, the local quadratic sum vanishes once the exponent of \(2\)
in the induced modulus is sufficiently large, as shown in the proof of
Lemma~\ref{l:Gcharbound}.  Only finitely many \(2\)-adic factors remain,
and all of them are bounded by an absolute constant. We may therefore take
\[
 \mathcal B_0\asymp
 \prod_{\substack{p\mid q_0\\p>2}}\frac{2p}{p-1},
\]
which gives the first estimate in
\eqref{e:exceptionalarithmetic}.

Let \(p\) be odd.  If \(p\nmid q_0\), then
\(v_p(q)=v_p(d)\), and the local character modulo \(p^{v_p(q)}\) is
principal.  By \eqref{e:primepowerestimate}, the corresponding local
quadratic sum vanishes when \(v_p(d)\geq2\).  Hence
\[
 p\nmid q_0
 \quad\Longrightarrow\quad
 v_p(d)\in\{0,1\}.
\]

Now suppose that \(p\mid q_0\).  Since \(q_0\) is the conductor of a
primitive real character, \(v_p(q_0)=1\).  If \(p\mid d\), then
\(\nu:=v_p(q)\geq2\), and the local character modulo \(p^\nu\) is induced
from the quadratic character modulo \(p\).  Writing
\(r=u+p^{\nu-1}t\), we obtain
\begin{align*}
 &\sum_{\substack{r\;(\mathrm{mod}\;p^\nu)\\p\nmid r}}
 \chi_p(r)e(-br^2/p^\nu)=
 \sum_{\substack{u\;(\mathrm{mod}\;p^{\nu-1})\\p\nmid u}}
 \chi_p(u)e(-bu^2/p^\nu)
 \sum_{t\;(\mathrm{mod}\;p)}e(-2but/p)
 =0.
\end{align*}
Thus
\[
 p\mid q_0\quad\Longrightarrow\quad v_p(d)=0.
\]

At \(p=2\), the calculation in the proof of
Lemma~\ref{l:Gcharbound} shows that the local sum vanishes when
\(v_2(q)\geq6\).  Consequently, whenever the induced quadratic sum is
not identically zero, one may write
\[
 d=2^\nu d_1,\qquad
 0\leq\nu\leq5,\qquad
 \mu^2(d_1)=1,\qquad
 (d_1,2q_0)=1,
\]
which is \eqref{e:sqfreeconds}. Since \(q=(2^\nu q_0)d_1\) and
\((2^\nu q_0,d_1)=1\), the Chinese remainder theorem gives, for
suitable \(b_0\in\mathbb A_{2^\nu q_0}\) and
\(b_1\in\mathbb A_{d_1}\),
\[
 \mathfrak s_{\chi_q}(a,q)
 =
 \mathfrak s_{\chi_{2^\nu q_0}}(b_0,2^\nu q_0)
 \mathfrak s(b_1,d_1).
\]
For every \(p\mid d_1\), the local character is principal and
\(v_p(d_1)=1\).  Hence the corresponding normalized local sum
satisfies
\[
 \frac1{\phi(p)}
 \left|\sum_{r\in\mathbb A_p}e(-b_pr^2/p)\right|
 \leq\frac{\sqrt p+1}{p-1}
 =p^{-1/2}(1-p^{-1/2})^{-1}.
\]
Multiplying these estimates over \(p\mid d_1\) and using the definition
of \(\sigma_0\), we obtain
\[
 |\mathfrak s_{\chi_q}(a,q)|
 \leq
 \sigma_0d_1^{-1/2}
 \prod_{p\mid d_1}(1-p^{-1/2})^{-1}.
\]
Since \(d=2^\nu d_1\) with \(0\leq\nu\leq5\), the bounded \(2\)-adic
factor may be absorbed into the implied constant.  Therefore
\[
 |\mathfrak s_{\chi_q}(a,q)|
 \ll \sigma_0d^{-1/2}\rho(d),
\]
which proves the second estimate in
\eqref{e:exceptionalarithmetic}.

It remains to prove the last estimate in
\eqref{e:exceptionalarithmetic}.  The maximum over
\(0\leq j\leq5\) in the definition of \(\sigma_0\) affects only the
\(2\)-adic local factor and therefore changes the estimate by an
absolute constant.  If \(\sigma_0=0\), there is nothing to prove.
Otherwise, factorization over prime powers and pairing \(r\) with
\(-r\) show that, if \(\chi_p(-1)=-1\) for some odd prime
\(p\mid q_0\), then the corresponding local quadratic sum vanishes.
Thus a nonzero term requires $\chi_p(-1)=1$ for every odd prime $p\mid q_0$.
Since \(\chi_p\) is the quadratic character modulo \(p\), this is
equivalent to \(p\equiv1\pmod4\).
At such a prime, the local contribution to \(\sigma_0\) is at most
\(2p^{1/2}/(p-1)\), while the corresponding factor in
\(\mathcal B_0^{1/2}\) is \((2p/(p-1))^{1/2}\).  Multiplying these bounds
over \(p\mid q_0\) gives
\[
 \sigma_0\mathcal B_0^{1/2}
 \ll
 \prod_{\substack{p\mid q_0\\p>2}}
 \frac{2\sqrt2\,p}{(p-1)^{3/2}}.
\]
The factor corresponding to \(p=5\) is bounded absolutely, and for every
\(p\geq13\),
\[
 \frac{2\sqrt2\,p}{(p-1)^{3/2}}\leq p^{-1/25}.
\]
The \(2\)-adic factor is also bounded absolutely.  Since every remaining
odd prime divisor of \(q_0\) is at least \(13\), we conclude that
\[
 \sigma_0\mathcal B_0^{1/2}\ll q_0^{-1/25},
\]
as required.
\end{proof}

\section{Proof of the improving inequalities}
\label{s:improving}

We prove the fixed-scale improving estimates for the prime-square averages.
The argument is separated into an abstract high--low principle and the
number-theoretic decomposition proved later in Section~\ref{s:hilo}.
The low part is estimated from \(\ell^1\) to \(\ell^\infty\);
the high part is estimated on \(\ell^2\).  Balancing these two
bounds gives the endpoint restricted estimate, and interpolation with the
\(\ell^2\) bound gives the strong \(\ell^p\)-improving estimates for
\(p>3/2\).



We use
the logarithmic and Orlicz conventions in \eqref{e:Lognotation} and
\eqref{e:defPhi}; in particular, $\Phi_c=\Phi_{c,1/2}$.

\begin{theorem}\label{t:improvinggeneral}
Let $A_N$ and $B_N$ be positive convolution operators whose kernels are
supported in an interval of length $O(N)$.  Suppose that
\begin{align}
 |A_Nf(x)|&\leq(\log N)^aB_N|f|(x),\label{e:anbnsupbound}\\
 |B_Nf(x)|&\ll N^{-v}\|f\mathbf1_{[x-CN,x+CN]}\|_1,
 \qquad 0\leq v<1,\label{e:bnv}
\end{align}
and that $s\geq1-v$, $r>0$, and, whenever
$2\leq J\leq N^{(1-v)/s}/10$,
\begin{align}
 B_N=Lo_{J,N}+Hi_{J,N},\qquad
 \|Lo_{J,N}f\|_\infty&\ll N^{-1}J^s(\log J)^\gamma\|f\|_1,
 \label{e:lobn1}\\
 \|Hi_{J,N}f\|_2&\ll J^{-r}(\log J)^\beta\|f\|_2,
 \label{e:hibn1}
\end{align}
where $\gamma\geq\beta\geq0$.  Suppose also that, for
$2\leq Q$ and $\log Q\leq c_0(\log N)^c$,
\begin{align}
 A_N&=\operatorname{Low}_{Q,N}+\operatorname{High}_{Q,N}+E_{Q,N},
 \label{e:decompan}\\
 \|\operatorname{Low}_{Q,N}f\|_\infty
 &\ll N^{-1}Q^s(\log Q)^\alpha\|f\|_1,\label{e:lowgeneral}\\
 \|\operatorname{High}_{Q,N}f\|_2
 &\ll Q^{-r}\exp\left(c_1\frac{(\log Q)^\kappa}{\log\log Q}\right)\|f\|_2,
 \label{e:highgeneral}\\
 \|E_{Q,N}f\|_2&\ll_A Q^{-r}(\log Q)^{C_A}
 (1+\log N)^{-A}\|f\|_2
 \qquad(A>0),\label{e:errorgeneral}
\end{align}
where \(0<\kappa<1\) and \(a,\alpha,c_0,c_1,c>0\).  In
\eqref{e:errorgeneral}, the constant \(C_A\) may depend on \(A\).  Then, for
$F,G$ contained in an interval of length $N$,
\begin{align}
 \frac{|\langle B_N\mathbf1_F,\mathbf1_G\rangle|}{N}
 &\ll (|F||G|/N^2)^{1-\frac{s}{2(s+r)}}
 (\Log (|F||G|/N^2))^{\gamma+\frac{s(\beta-\gamma)}{s+r}},\label{e:endimprovbn}\\
 \frac{|\langle A_N\mathbf1_F,\mathbf1_G\rangle|}{N}
 &\ll (|F||G|/N^2)^{1-\frac{s}{2(s+r)}}\Phi_{c^*,\kappa}(|F||G|/N^2),
 \label{e:endimprovan}
\end{align}
where 
\begin{align}\label{e:cstar}
 c^*=\frac{c_1s}{2^\kappa(s+r)^{\kappa+1}}+o(1).
\end{align}
Here and below the $o(1)$ in a subscript means that the estimate holds
with $c^*$ replaced by any fixed larger constant.
\end{theorem}

\begin{proof}
Let $\delta \coloneqq |F||G|/N^2$. After enlarging the implied constant, we may assume that
$\delta\leq\delta_0$, where $\delta_0>0$ depends only on the fixed
parameters; otherwise we take $J=Q=2$.  From
\eqref{e:lobn1}, \eqref{e:hibn1}, the triangle inequality, and Cauchy--Schwarz,
\begin{align}\label{e:Bbalance}
 \frac{|\langle B_N\mathbf1_F,\mathbf1_G\rangle|}{N}
 &\ll N^{-1}\|Lo_{J,N}\mathbf{1}_F\|_\infty |G|+
 N^{-1}\|Hi_{J,N}\mathbf{1}_F\|_2 |G|^{1/2}\\
 &\ll\delta J^s(\log J)^\gamma
      +\delta^{1/2}J^{-r}(\log J)^\beta.
\end{align}
Choose a dyadic $J$ within a factor $2$ of
\begin{align}\label{e:Jchoice}
 \delta^{-1/(2(s+r))}(\Log\delta)^{(\beta-\gamma)/(s+r)}.
\end{align}
If this is at most $N^{(1-v)/s}/10$, substitution in
\eqref{e:Bbalance} gives \eqref{e:endimprovbn}.  Otherwise use
\eqref{e:bnv}.  Since $s\geq1-v$, the estimate 
$$\delta^{-1/(2(s+r))}(\Log\delta)^{(\beta-\gamma)/(s+r)} > N^{(1-v)/s}/10$$
implies that
\begin{align*}
 N^{-1-v}|F||G|
 &\ll \delta^{1-\frac{s}{2(s+r)}}
 (\Log\delta)^{\frac{s(\beta-\gamma)}{s+r}}.
\end{align*}
Because $\gamma\geq0$ and $\Log\delta\geq1$, this is bounded by the
right-hand side of \eqref{e:endimprovbn}.  

We turn to $A_N$.  If
$\Log\delta\geq c_3(\log N)^c$, where $c_3>0$ is sufficiently small, then
\eqref{e:anbnsupbound} and \eqref{e:endimprovbn} give
\begin{align*}
 N^{-1}|\langle A_N\mathbf1_F,\mathbf1_G\rangle|
 \ll (\log N)^a N^{-1}|\langle B_N\mathbf1_F,\mathbf1_G\rangle| &\ll \delta^{1-\frac{s}{2(s+r)}}(\Log\delta)^{\gamma+a/c+\frac{s(\beta-\gamma)}{s+r}}\\
 &\ll \delta^{1-\frac{s}{2(s+r)}}\Phi_{c^*,\kappa}(\delta).
\end{align*}
The last inequality follows because $\Phi_{c^*,\kappa}(t)$ dominates every
fixed power of $\Log t$ as $t\downarrow0$.

Suppose now that $\Log\delta<c_3(\log N)^c$.  Choose $Q$ dyadic.  From \eqref{e:decompan}--\eqref{e:errorgeneral},
\begin{align}\label{e:Abalance} \qquad
 \frac{|\langle A_N\mathbf1_F,\mathbf1_G\rangle|}{N} &\ll N^{-1}\|\operatorname{Low}_{Q,N}\mathbf1_F\|_\infty |G|   +N^{-1}\|\operatorname{High}_{Q,N}\mathbf1_F\|_2  |G|^{1/2} + N^{-1} \|E_{Q,N}\mathbf1_F\|_2 |G|^{1/2}\\
 &\ll \delta Q^s(\log Q)^\alpha
 +\delta^{1/2}Q^{-r}\exp\left(c_1\frac{(\log Q)^\kappa}{\log\log Q}\right)
 +\delta^{1/2}Q^{-r}(\log Q)^{C_A}
 (1+\log N)^{-A}.
\end{align}
Take
\begin{align}\label{e:Qchoice}
 \log Q=\frac{\Log \delta}{2(s+r)}+
 \frac{c_1}{s+r}\frac{\left((\Log \delta)/(2(s+r))\right)^\kappa}{\log \Log \delta}+O(\log \Log \delta).
\end{align}
The constant in the $O(\log \Log \delta)$ is chosen so that the two first terms
in \eqref{e:Abalance} are comparable.  Since $\Log \delta<c_3(\log N)^c$,
$c_3$ may be chosen so that $\log Q\leq c_0(\log N)^c$; thus the
decomposition is available.  Moreover,
\begin{align*}
 Q^s
 &=\delta^{-\frac{s}{2(s+r)}}
 \exp\left(\left(\frac{c_1s}{2^\kappa(s+r)^{\kappa+1}}+o(1)\right)
 \frac{(\Log \delta)^\kappa}{\log \Log \delta}\right).
\end{align*}
The powers of \(\log Q\) in the first two terms of
\eqref{e:Abalance} are absorbed by the \(o(1)\) in the exponential.  For
fixed \(A\), the factor \((\log Q)^{C_A}\) in the error term is dominated,
for sufficiently large \(Q\), by the subexponential factor in the high
term.  The bounded range of \(Q\) is absorbed into the implied constant.
Since \((1+\log N)^{-A}\leq1\), the error term is also absorbed by the
high term.  This proves \eqref{e:endimprovan} and \eqref{e:cstar}.
\end{proof}

The above proof does not use that $F$ and $G$ lie in a common interval,
except to ensure that $|F||G|\leq N^2$.  Hence
\eqref{e:endimprovbn} and \eqref{e:endimprovan} remain valid for
arbitrary finite sets $F,G\subset\mathbb Z$ satisfying
$|F||G|\leq N^2$.

\begin{corollary}\label{c:improvinggeneral}
Assume in addition that $A_N$ and $B_N$ are uniformly bounded on
$\ell^2$.  If
\begin{align*}
 \frac{2(r+s)}{2r+s}<p\leq2,
\end{align*}
then, for $f\in \ell^p(\mathbb{Z})$,
\begin{align}
 \|B_Nf\|_{p'}+\|A_Nf\|_{p'}
 \ll_p N^{1/p'-1/p}\|f\|_p.\label{e:lpimproving}
\end{align}
\end{corollary}

\begin{proof}
Set $p_0=\frac{2(r+s)}{2r+s}$ and assume that $f$ is supported on an interval of length $O(N)$ first. 
We  remove the logarithmic losses in
\eqref{e:endimprovbn} and \eqref{e:endimprovan}.  Indeed, for every
$\varepsilon>0$ and every $0<\delta\leq1$, $ \bigl(\log(2/\delta)\bigr)^{O(1)}
 \ll_\varepsilon \delta^{-\varepsilon}.$
 Since every fixed logarithmic or
subexponential factor occurring in \eqref{e:endimprovbn} and
\eqref{e:endimprovan} is bounded by $C_\varepsilon\delta^{-\varepsilon}$,
these estimates imply that, for every $p_0<\widetilde p<2$,
\begin{align}\label{e:restrictedimproving}
 |\langle T_N\mathbf1_F,\mathbf1_G\rangle|
 \ll_{\widetilde p}
 N^{1/\widetilde p'-1/\widetilde p}
 |F|^{1/\widetilde p}|G|^{1/\widetilde p},
 \qquad T_N\in\{A_N,B_N\}.
\end{align}
Thus $A_N$ and $B_N$ are uniformly of restricted weak type
$(\widetilde p,\widetilde p')$, with the required scaling in $N$.

Now fix $p_0<p<2$ and choose $\widetilde p$ with
$p_0<\widetilde p<p$. Interpolating \eqref{e:restrictedimproving} with the uniform
$\ell^2\to\ell^2$ estimate gives \eqref{e:lpimproving} when $f$ is
supported in an interval of length $O(N)$, see \cite[Chapter 1]{BB}.  The case $p=2$
is precisely the assumed uniform $\ell^2$ estimate.

Finally, partition $\mathbb Z$ into intervals $I$ of length comparable
to $N$.  Since the kernels are supported in $[-CN,CN]$,
$T_N(f\mathbf1_I)$ is supported in $I+[-CN,CN]$.  The intervals may
therefore be divided into finitely many families such that these
enlarged intervals are pairwise disjoint within each family.  For each
such family $\mathcal I$, the local estimate and $p'\geq p$ give
\begin{align*}
 \left\|\sum_{I\in\mathcal I}T_N(f\mathbf1_I)\right\|_{p'}
 &\ll_p N^{1/p'-1/p}
 \left(\sum_{I\in\mathcal I}\|f\mathbf1_I\|_p^{p'}\right)^{1/p'}\leq N^{1/p'-1/p}
 \left(\sum_{I\in\mathcal I}\|f\mathbf1_I\|_p^p\right)^{1/p}.
\end{align*}
Summing over the finitely many families proves
\eqref{e:lpimproving} for every $f\in\ell^p(\mathbb Z)$.
\end{proof}

We shall prove in Section~\ref{s:hilo} that the annular prime-square
averages admit the decomposition required below.  For dyadic $N$, put
\begin{align}\label{e:annularaverages}
 \widetilde A_{N,\Lambda}
 &=A_{N,\Lambda}-\frac12A_{N/2,\Lambda},&
 \widetilde A_{N,\mathbb N^2}
 &=A_{N,\mathbb N^2}-\frac12A_{N/2,\mathbb N^2}.
\end{align}
Thus
\begin{align}\label{e:annularreconstruction}
 A_{N,\Lambda}
 =\sum_{\substack{j\geq0\\2^{-j}N\geq2}}
 2^{-j}\widetilde A_{2^{-j}N,\Lambda},
\end{align}
For the square averages, telescoping gives
\[
 A_{N,\mathbb N^2}f
 =
 \sum_{\substack{j\geq0\\2^{-j}N\geq2}}
 2^{-j}\widetilde A_{2^{-j}N,\mathbb N^2}f
 +\frac1N f(\,\cdot-1).
\]

\begin{proposition}\label{p:decompositionofan}
Let \(N\) be dyadic, let \(Q\geq2\), and suppose that
\[
 \log Q\leq c_0(\log N)^{1/2},
\]
where \(c_0>0\) is a sufficiently small absolute constant.  Then
\begin{align}\label{e:primedecomposition}
 \widetilde A_{N,\Lambda}
 =
 \operatorname{Low}_{Q,N}
 +\operatorname{High}_{Q,N}+E_{Q,N},
\end{align}
where
\begin{align}
 \|\operatorname{Low}_{Q,N}f\|_\infty
 &\ll N^{-2}Q\log(2Q)\,\|f\|_1,
 \label{e:primelow}\\
 \|\operatorname{High}_{Q,N}f\|_2
 &\ll Q^{-1/2}
 \exp\left(
 (2+o(1))\frac{(\log Q)^{1/2}}{\log\log Q}
 \right)\|f\|_2,
 \label{e:primehigh}\\
 \|E_{Q,N}f\|_2
 &\ll_A Q^{-1/2}(\log(2Q))^{O_A(1)}
 (1+\log N)^{-A}\|f\|_2
 \label{e:primeerror}
\end{align}
for every \(A>0\).  More precisely, for some absolute \(C>0\),
the low term satisfies
\begin{align}\label{e:primelowlocal}
 |\operatorname{Low}_{Q,N}f(x)|
 \ll Q\log(2Q)\sup_{M\geq N}M^{-2}
 \sum_{|y-x|\leq CM^2}|f(y)|.
\end{align}
The maximal estimate needed in Section~\ref{s:sparse} is
\begin{align}\label{e:primehighmax}
 \left\|
 \sup_{\substack{M\geq N_0\\M\in2^{\mathbb N}}}
 |\operatorname{High}_{Q,M}f|
 \right\|_2
 &\ll Q^{-1/2}
 \exp\left(
 (2+o(1))\frac{(\log Q)^{1/2}}{\log\log Q}
 \right)\|f\|_2,
\end{align}
provided
\[
 \log Q\leq c_0(\log N_0)^{1/2}.
\]
 Here \(o(1)\) is taken as \(Q\to\infty\), uniformly in the displayed
ranges of \(N\) and \(N_0\); bounded values of \(Q\) are absorbed into
the implicit constants. 
The operators in
\eqref{e:primedecomposition} are defined in
Section~\ref{s:hilo}.
\end{proposition}

\begin{proof}[Proof of Theorem~\ref{t:improvingfixed}]
By positivity, it is enough to consider dyadic $N$.  Indeed, if
$N\leq N'<2N$ and $N'$ is dyadic, then
\[
 |A_{N,\Lambda}f|
 \leq 2A_{N',\Lambda}|f|,
\]
and the spatial scales $N^2$ and $(N')^2$ are comparable.
We first apply the high--low argument of
Theorem~\ref{t:improvinggeneral} to the annular operators
$\widetilde A_{N,\Lambda}$ and
$\widetilde A_{N,\mathbb N^2}$, with spatial scale $N^2$ and
\[
 r=\tfrac12,\qquad s=1,\qquad v=\tfrac12,\qquad
 \gamma=2,\qquad\beta=1,\qquad
 \alpha=1,\qquad\kappa=\tfrac12,\qquad c_1=2.
\]
The estimates required for \(\widetilde A_{N,\mathbb N^2}\) follow from
\cite[Lemma~2.2 and (2.5)]{LHF}, applied at the scales \(N\) and \(N/2\).
The estimates for \(\widetilde A_{N,\Lambda}\) are given by
Proposition~\ref{p:decompositionofan}, while
\[
 |\widetilde A_{N,\Lambda}f|
 \leq(\log N)\widetilde A_{N,\mathbb N^2}|f|.
\]
Both annular operators are uniformly bounded on \(\ell^2\).  For every
fixed \(\varepsilon>0\), the high estimate in
Proposition~\ref{p:decompositionofan} satisfies the hypothesis
\eqref{e:highgeneral} with \(c_1=2+\varepsilon\), for all sufficiently
large \(Q\).  Applying Theorem~\ref{t:improvinggeneral} and then absorbing
\(\varepsilon\) into the final \(o(1)\) gives the asserted endpoint estimate
for \(\widetilde A_{N,\Lambda}\).

It remains to pass from the annular averages to $A_{N,\Lambda}$.
Since $|I|\asymp N^2$, put $ \delta=\frac{|F||G|}{N^4}.$
The case $\delta\gtrsim1$ follows directly from the uniform $\ell^2$
bound, so we may assume that $\delta\leq1$.  Put
$M_j=2^{-j}N$, and let $j_0$ be the first index for which
$2^{4j_0}\delta>1$, with $j_0=\infty$ if there is no such index.

For \(j<j_0\), we have \(|F||G|\leq M_j^4\).  Applying the annular
endpoint estimate with \(c_1=2+\varepsilon\), and using
\[
 1-\frac{s}{2(s+r)}=\frac23,\qquad
 c_\varepsilon
 :=\frac{2+\varepsilon}{\sqrt2(3/2)^{3/2}}
 =\frac{4\sqrt3}{9}+O(\varepsilon),
\]
gives
\begin{align*}
 \frac{2^{-j}}{N^2}
 \bigl|\langle\widetilde A_{M_j,\Lambda}
 \mathbf1_F,\mathbf1_G\rangle\bigr|
 &\ll_\varepsilon
 2^{-j/3}\delta^{2/3}
 \Phi_{c_\varepsilon+o(1)}(2^{4j}\delta)\\
 &\ll_\varepsilon
 2^{-j/3}\delta^{2/3}
 \Phi_{c_\varepsilon+o(1)}(\delta).
\end{align*}
Here we used \(2^{4j}\delta\leq1\) for \(j<j_0\), so
\(\Phi_{c_\varepsilon+o(1)}(2^{4j}\delta)\leq
\Phi_{c_\varepsilon+o(1)}(\delta)\), after increasing the implicit
constant if necessary.

For \(j\geq j_0\), the uniform \(\ell^2\) estimate gives
\[
 \frac{2^{-j}}{N^2}
 \bigl|\langle\widetilde A_{M_j,\Lambda}
 \mathbf1_F,\mathbf1_G\rangle\bigr|
 \ll 2^{-j}\delta^{1/2}.
\]
Therefore
\begin{align*}
 \frac1{N^2}
 \bigl|\langle A_{N,\Lambda}\mathbf1_F,\mathbf1_G\rangle\bigr|
 &\ll_\varepsilon
 \delta^{2/3}\Phi_{c_\varepsilon+o(1)}(\delta)
 \sum_{j<j_0}2^{-j/3}
 +\delta^{1/2}\sum_{j\geq j_0}2^{-j}\\
 &\ll_\varepsilon
 \delta^{2/3}\Phi_{c_\varepsilon+o(1)}(\delta)
 +\delta^{3/4}.
\end{align*}
Since \(\delta^{3/4}\leq\delta^{2/3}\) and
\(\Phi_{c_\varepsilon+o(1)}(\delta)\geq1\), this proves the endpoint
estimate with Orlicz constant \(c_\varepsilon+o(1)\).  Letting
\(\varepsilon\downarrow0\) gives the stated constant
\(4\sqrt3/9+o(1)\) for \(A_{N,\Lambda}\).

This completes the proof of the endpoint case.  The estimates for
\(p>3/2\) follow in the same way from
Corollary~\ref{c:improvinggeneral}.
\end{proof}

\section{Proof of the sparse bounds}
\label{s:sparse}

We pass from the fixed-scale improving estimates to sparse bounds for the
prime-square maximal operator.  The new point is that, after linearizing
the maximal function, the scale depends on the spatial point.  Thus the
fixed-scale \(\ell^2\) estimate for the high part must be replaced by a
maximal \(\ell^2\) estimate.
We
then produce a sparse family through iteration.

 Notice that
there is no maximal-function loss for the low term, since
\begin{align}\label{e:supLinfty}
 \left\|\sup_{N>M}|Lo_{J,N}f|\right\|_\infty
 \leq \sup_{N>M}\|Lo_{J,N}f\|_\infty , \qquad M>0.
\end{align}
 We next estimate the low and high parts at admissible stopping times.

\begin{lemma}\label{l:localsparse}
Fix constants \(C\geq1\) and \(c>2C+1\).  Let \((B_M)\) be a family
of positive convolution operators satisfying
\begin{align}\label{e:bnv-local}
 |B_M f(x)|
 \ll M^{-v}\|f\mathbf1_{[x-CM,x+CM]}\|_1,
 \qquad 0\leq v<1.
\end{align}
Let \(s\geq1-v\), \(r>0\), and \(\gamma\geq\beta\geq0\), and suppose
that, whenever \(M\) is dyadic and \(2\leq J\leq M^{(1-v)/s}/10\),
\begin{align}\label{e:localdecomposition}
 B_M=Lo_{J,M}+Hi_{J,M}+Er_{J,M}.
\end{align}
Let \(I_0\) be an interval, and let
\(f=\mathbf1_F\) and \(g=\mathbf1_G\), where
\(F\subset cI_0\) and \(G\subset I_0\).  Suppose that
\[
 \tau:I_0\longrightarrow
 \{2^k:2^k\leq|I_0|\}
\]
satisfies
\begin{align}\label{e:admisslocal}
 \sup_{M\geq\tau(x)}M^{-1}
 \sum_{|y-x|\leq C_0M}f(y)
 \ll \langle f\rangle_{cI_0,1},
 \qquad x\in I_0.
\end{align}
Write \(B_\tau f(x)=B_{\tau(x)}f(x)\).

Assume that, for every \(J\geq2\) and every \(x\in I_0\) such that
\(\tau(x)\geq(10J)^{s/(1-v)}\),
\begin{align}
 |Lo_{J,\tau(x)}f(x)|
 &\ll J^s(\log J)^\gamma
 \sup_{M\geq\tau(x)}M^{-1}
 \sum_{|y-x|\leq CM}f(y),
 \label{e:localLow}
\end{align}
and that, for every \(h\in\ell^2(\mathbb Z)\),
\begin{align}
 \left\|
 \sup_{\substack{M\geq(10J)^{s/(1-v)}\\M\ {\rm dyadic}}}
 |Hi_{J,M}h|
 \right\|_2
 &\ll J^{-r}(\log J)^\beta\|h\|_2.
 \label{e:localHigh}
\end{align}
Finally, suppose that, for some fixed \(A>1/2\), uniformly in dyadic \(M\),
\begin{align}
 \|\widehat{Er}_{J,M}\|_{L^\infty(\mathbb T)}
 \ll_A
 J^{-r}(\log J)^\beta(1+\log M)^{-A}.
 \label{e:localError}
\end{align}
Then
\begin{align}\label{e:sparsesquares1}
 |I_0|^{-1}|\langle B_\tau f,g\rangle|
 \ll_c
 \left(\frac{|F||G|}{|I_0|^2}\right)^{
 1-\frac{s}{2(s+r)}}
 \left(\Log\frac{|F||G|}{|I_0|^2}\right)^{
 \gamma+\frac{s(\beta-\gamma)}{s+r}}.
\end{align}
\end{lemma}

\begin{proof}
Put $\delta_0=\frac{|F||G|}{|I_0|^2}$.  We may assume that
$\delta_0\leq1/4$, since otherwise the result follows by choosing the
implicit constant in \eqref{e:sparsesquares1} sufficiently large.
Choose a dyadic $J$ within a fixed factor of
\begin{align}\label{e:localJchoice}
 \delta_0^{-1/(2(s+r))}
 (\Log\delta_0)^{(\beta-\gamma)/(s+r)}
\end{align}
and put $T=(10J)^{s/(1-v)}$.

We first consider the set $\{\tau<T\}$. By \eqref{e:bnv-local} and admissibility,
\begin{align*}
 B_{\tau(x)}f(x)
 &\ll \tau(x)^{-v}
 \sum_{|y-x|\leq C\tau(x)}f(y)
 \ll \tau(x)^{1-v}\langle f\rangle_{cI_0,1}
 \ll T^{1-v}\langle f\rangle_{cI_0,1}.
\end{align*}
Since $T^{1-v}=(10J)^s$, it follows that
\begin{align}\label{e:smalltau}
 |I_0|^{-1}
 \big|\langle\mathbf1_{\{\tau<T\}}B_\tau f,g\rangle\big|
 \ll J^s\delta_0.
\end{align}

On $\{\tau\geq T\}$ we have
$J\leq\tau(x)^{(1-v)/s}/10$, so the fixed-scale decomposition
\eqref{e:localdecomposition} is available at the scale $M=\tau(x)$.
For each dyadic $M$ put
\[
 E_M:=\{x\in I_0:\tau(x)=M\}.
\]
Then, pointwise on $I_0$,
\[
 \mathbf1_{\{\tau\geq T\}}B_\tau f
 =
 \mathbf1_{\{\tau\geq T\}}Lo_{J,\tau}f
 +
 \mathbf1_{\{\tau\geq T\}}Hi_{J,\tau}f
 +
 \sum_{\substack{M\geq T\\M\text{ dyadic}}}
 \mathbf1_{E_M}Er_{J,M}f.
\]
The pointwise estimate \eqref{e:localLow} and admissibility give 
\begin{align*}
 &|I_0|^{-1}
 \big|\langle\mathbf1_{\{\tau\geq T\}}
 Lo_{J,\tau}f,g\rangle\big|\ll
 \frac{J^s(\log J)^\gamma}{|I_0|}
 \sum_{\substack{x\in G\\\tau(x)\geq T}}
 \sup_{M\geq\tau(x)}M^{-1}
 \sum_{|y-x|\leq CM}f(y)\ll
 \delta_0J^s(\log J)^\gamma.
\end{align*}
For the high-frequency part we  use the multifrequency estimate \eqref{e:localHigh}.
Thus Cauchy--Schwarz gives
\begin{align*}
 |I_0|^{-1}
 \big|\langle\mathbf1_{\{\tau\geq T\}}Hi_{J,\tau}f,g\rangle\big|
 &\leq |I_0|^{-1}\|g\|_2
 \left\|\sup_{M\geq T}|Hi_{J,M}f|\right\|_2  \\
 &\ll
 |I_0|^{-1}|G|^{1/2}J^{-r}(\log J)^\beta |F|^{1/2}  =\delta_0^{1/2}J^{-r}(\log J)^\beta .
\end{align*}
The error term is treated scale by scale.  By Cauchy--Schwarz and
Plancherel,
\begin{align*}
 |I_0|^{-1}
 \sum_{\substack{M\geq T\\ M\text{ dyadic}}}
 \big|\langle\mathbf1_{E_M}Er_{J,M}f,g\rangle\big| & \leq
 |I_0|^{-1}
 \sum_{\substack{M\geq T\\ M\text{ dyadic}}}
 |G\cap E_M|^{1/2}\|Er_{J,M}f\|_2 \\
 &\leq
 |I_0|^{-1}|F|^{1/2}
 \sum_{\substack{M\geq T\\ M\text{ dyadic}}}
 |G\cap E_M|^{1/2}\|\widehat{Er}_{J,M}\|_\infty .
\end{align*}
Using \eqref{e:localError} and Cauchy--Schwarz over the dyadic scales,
\begin{align*}
  &|I_0|^{-1}
 \sum_{\substack{M\geq T\\ M\text{ dyadic}}}
 \big|\langle\mathbf1_{E_M}Er_{J,M}f,g\rangle\big| \\
 &\ll_A
 J^{-r}(\log J)^{\beta}|I_0|^{-1}|F|^{1/2}
 \left(\sum_M |G\cap E_M|\right)^{1/2}
 \left(\sum_{M\text{ dyadic},M\geq T}(1+\log M)^{-2A}\right)^{1/2}\\
 &\ll_A J^{-r}(\log J)^{\beta}
 \delta_0^{1/2} (\log T)^{1/2-A},
\end{align*}
since the sets $G\cap E_M$ are disjoint subsets of $G$.

Combining the small-$\tau$, low, high, and error estimates gives
\begin{align*}
 |I_0|^{-1}|\langle B_\tau f,g\rangle|
 \ll
 \delta_0J^s(\log J)^\gamma
 +\delta_0^{1/2}J^{-r}(\log J)^\beta .
\end{align*}
Substituting \eqref{e:localJchoice} and using
$\log J\sim \Log\delta_0$ gives
\begin{align*}
 |I_0|^{-1}|\langle B_\tau f,g\rangle|
 \ll \delta_0^{1-\frac{s}{2(s+r)}}
 (\Log\delta_0)^{
 \gamma+\frac{s(\beta-\gamma)}{s+r}},
\end{align*}
which proves \eqref{e:sparsesquares1}.
\end{proof}

Choose a fixed odd integer \(c>2C+1\), where \(C\) is the constant in
\eqref{e:admisslocal}.  The proof of the following lemma is a standard argument in harmonic analysis. We include the details for completeness.

\begin{lemma}\label{l:sparserecursion}
Suppose that every admissible $\tau$ satisfies
\begin{align}\label{e:localform}
 |\langle T_\tau\mathbf1_F,\mathbf1_G\rangle|
 \ll |I_0|\,\Psi\left(\frac{|F|}{|I_0|},
                       \frac{|G|}{|I_0|}\right)
\end{align}
for $F\subset cI_0$ and $G\subset I_0$, where $(T_N)_{N\in2^{\mathbb N}}$ is a family of positive convolution operators whose kernel is supported in
$[-CN,CN]$ for some fixed $C>0$.  Then there is a sparse family $\mathcal S$ such that
\begin{align}\label{e:recursionconclusion}
 \langle\sup_{N\in2^{\mathbb N}}|T_N\mathbf1_F|,\mathbf1_G\rangle
 \ll\sum_{I\in\mathcal S}|I|\,
 \Psi\left(\frac{|F\cap cI|}{|I|},
            \frac{|G\cap I|}{|I|}\right).
\end{align}
\end{lemma}

Since the operators are positive and the input is nonnegative, the absolute
values in \(\sup_N|T_N\mathbf1_F|\) may be omitted.  We keep them only to
match the maximal-operator notation.

\begin{proof}
All scale suprema in this proof are over dyadic scales. Fix a dyadic interval \(I_0\), with
\(F\subset cI_0\) and \(G\subset I_0\), and let
\(\mathcal I(I_0)\) be the maximal dyadic subintervals \(I\subset I_0\)
for which
\begin{align}\label{e:badinterval}
 \langle\mathbf1_F\rangle_{cI,1}
 >K\langle\mathbf1_F\rangle_{cI_0,1}
 \quad\text{or}\quad
 \langle\mathbf1_G\rangle_{I,1}
 >K\langle\mathbf1_G\rangle_{I_0,1},
\end{align}
where \(K\) is sufficiently large.  The weak \((1,1)\) estimates for
the Hardy--Littlewood and dyadic maximal operators give
\[
 \sum_{I\in\mathcal I(I_0)}|I|
 \ll_c K^{-1}|I_0|.
\]
Thus \(K\) may be chosen so that
\begin{align}\label{e:packing}
 \sum_{I\in\mathcal I(I_0)}|I|\leq\tfrac14|I_0|.
\end{align}

For \(x\in I_0\), let \(I(x)\) be the interval in
\(\mathcal I(I_0)\) containing \(x\), if one exists, and put
\(I(x)=\varnothing\) otherwise.  Define
\[
 \mathcal S_{I_0}(x)=
 \begin{cases}
 \{N\in2^{\mathbb N}:|I(x)|<N\leq|I_0|\},
 &I(x)\neq\varnothing,\\
 \{N\in2^{\mathbb N}:N\leq|I_0|\},
 &I(x)=\varnothing.
 \end{cases}
\]
Choose \(\tau(x)\in\mathcal S_{I_0}(x)\) such that
\[
 T_{\tau(x)}\mathbf1_F(x)
 =
 \sup_{N\in\mathcal S_{I_0}(x)}T_N\mathbf1_F(x).
\]

We verify that \(\tau\) is admissible.  Let \(N\geq\tau(x)\) be
dyadic.  If \(N<|I_0|\), let \(I\subset I_0\) be the dyadic interval
containing \(x\) with \(N\leq|I|<2N\).  By maximality of the stopping
intervals, \eqref{e:badinterval} fails for \(I\), while
\(c>2C+1\) gives $[x-CN,x+CN]\subset cI.$
Consequently,
\[
 N^{-1}\sum_{|y-x|\leq CN}\mathbf1_F(y)
 \lesssim_c
 \langle\mathbf1_F\rangle_{cI,1}
 \lesssim
 \langle\mathbf1_F\rangle_{cI_0,1}.
\]
If \(N\geq|I_0|\), then
\[
 N^{-1}\sum_{|y-x|\leq CN}\mathbf1_F(y)
 \leq N^{-1}|F|
 \lesssim_c\langle\mathbf1_F\rangle_{cI_0,1}.
\]
Hence \(\tau\) satisfies \eqref{e:admisslocal}.

By \eqref{e:localform},
\begin{align}\label{e:parentterm}
 \big\langle
 \sup_{N\in\mathcal S_{I_0}(x)}|T_N\mathbf1_F|,
 \mathbf1_G
 \big\rangle
 \ll
 |I_0|\Psi\bigl(\langle\mathbf1_F\rangle_{cI_0,1},\langle\mathbf1_G\rangle_{I_0,1}\bigr).
\end{align}
For \(x\in I\in\mathcal I(I_0)\), the remaining scales are
\(N\leq |I|\).  If \(K_N\) denotes the kernel of \(T_N\), then
\(K_N(x,y)=0\) unless \(|x-y|\leq CN\).  Hence, for \(x\in I\) and
\(N\leq |I|\),
\[
 T_N\mathbf1_F(x)=T_N\mathbf1_{F\cap cI}(x).
\]
Therefore
\[
 \sup_{N\leq |I_0|}|T_N\mathbf1_F(x)|
 \leq
 \sup_{N\in\mathcal S_{I_0}(x)}|T_N\mathbf1_F(x)|
 +
 \sum_{I\in\mathcal I(I_0)}
 \mathbf1_I(x)
 \sup_{N\leq |I|}|T_N\mathbf1_{F\cap cI}(x)|.
\]
Pairing with \(\mathbf1_G\) and using \eqref{e:parentterm}, we obtain
\begin{align}\label{e:recursive}
 \big\langle\sup_{N\leq|I_0|}|T_N\mathbf1_F|,\mathbf1_G\big\rangle
 &\ll |I_0|\Psi\bigl(\langle\mathbf1_F\rangle_{cI_0,1},
                  \langle\mathbf1_G\rangle_{I_0,1}\bigr) +
 \sum_{I\in\mathcal I(I_0)}
 \big\langle\sup_{N\leq|I|}
 |T_N\mathbf1_{F\cap cI}|,\mathbf1_{G\cap I}\big\rangle .
\end{align}

Iterating \eqref{e:recursive} and using \eqref{e:packing} at every
generation gives a sparse collection \(\mathcal S\) of dyadic intervals
such that
\[
 \big\langle\sup_{N\leq|I_0|}|T_N\mathbf1_F|,\mathbf1_G\big\rangle
 \ll
 \sum_{I\in\mathcal S}
 |I|\Psi\bigl(\langle\mathbf1_F\rangle_{cI,1},\langle\mathbf1_G\rangle_{I,1}\bigr).
\]
The constants are uniform in the truncation \(N\leq |I_0|\).  Letting
the truncation increase and applying monotone convergence gives
\eqref{e:recursionconclusion}.
\end{proof}

We now apply the stopping-time argument of
Lemmas~\ref{l:localsparse} and \ref{l:sparserecursion} to the annular
prime-square averages.  At small stopping scales we use their domination
by the square averages, while at large scales we use the decomposition in
Proposition~\ref{p:decompositionofan}.

\begin{lemma}\label{l:primeendpointlocal}
Let $f=\mathbf1_F$ and $g=\mathbf1_G$, where
$F\subset cI_0$ and $G\subset I_0$, and put $\delta_0=\frac{|F||G|}{|I_0|^2}.$
 Let $\tau:I_0\to\{2^k:(2^k)^2\leq|I_0|\}$ be admissible at the
spatial scale $\tau(x)^2$, in the sense that
\begin{align}\label{e:admissprime}
 \sup_{M\geq\tau(x)}M^{-2}
 \sum_{|y-x|\leq CM^2}f(y)
 \ll\langle f\rangle_{cI_0,1}.
\end{align}
Under the hypotheses of Proposition~\ref{p:decompositionofan},
\begin{align}\label{e:primeendpointlocal}
 |I_0|^{-1}
 \big|\langle\widetilde A_{\tau,\Lambda}\mathbf1_F,
 \mathbf1_G\rangle\big|
 \ll\delta_0^{2/3}
 \Phi_{4\sqrt3/9+o(1)}(\delta_0).
\end{align}
\end{lemma}

\begin{proof}
Put \(L=\Log\delta_0^{-1}\) and
\[
 E_0=\{x\in I_0:\log\tau(x)\leq C_0L^2\},
 \qquad g_0=\mathbf1_{E_0}g.
\]
By positivity,
\[
 \mathbf1_{E_0}\widetilde A_{\tau,\Lambda}f
 \leq C_0L^2 A_{\tau,\mathbb N^2}f.
\]
The square averages satisfy the hypotheses of
Lemma~\ref{l:localsparse}, with
\[
 r=\tfrac12,\qquad s=1,\qquad v=\tfrac12,\qquad
 \gamma=2,\qquad\beta=1,
\]
by \cite[Lemmas~5.2--5.6]{LHF}.  Applying that lemma to \(f\), \(g_0\),
and the same admissible stopping time \(\tau\), and using
\(|\operatorname{supp}g_0|\leq|G|\), gives
\begin{align}\label{e:sparse1}
 |I_0|^{-1}
 \left|
 \left\langle
 \mathbf1_{E_0}\widetilde A_{\tau,\Lambda}f,g
 \right\rangle
 \right|
 &\ll
 L^2\delta_0^{2/3}L^{4/3}
 \ll \delta_0^{2/3}L^{10/3}.
\end{align}
For every fixed \(\varepsilon>0\), every fixed power of \(L\) is bounded by
\(\Phi_{\varepsilon}(\delta_0)\), after enlarging the implicit constant.
This gives the bound on  the right-hand side of
\eqref{e:primeendpointlocal}.

It remains to consider \(\{\log\tau>C_0L^2\}\).  Choose \(Q\) as in
\eqref{e:Qchoice}, with
\[
 s=1,\qquad r=\frac12,\qquad\kappa=\frac12,\qquad c_1=2+\varepsilon,
\]
and round it to a dyadic integer.  Then \(\log Q\ll L\).  By taking
\(C_0\) sufficiently large, we have $\log Q\leq c(\log\tau(x))^{1/2}$ whenever $\log\tau(x)>C_0L^2$.
Proposition~\ref{p:decompositionofan} therefore gives
\begin{align}\label{e:decomposproof}
     \widetilde A_{\tau(x),\Lambda}f(x)
 =
 \operatorname{Low}_{Q,\tau(x)}f(x)
 +\operatorname{High}_{Q,\tau(x)}f(x)
 +E_{Q,\tau(x)}f(x).
\end{align}
By \eqref{e:primelowlocal} and admissibility,
\[
 |\operatorname{Low}_{Q,\tau(x)}f(x)|
 \ll Q\log Q
 \sup_{M\geq\tau(x)}M^{-2}
 \sum_{|y-x|\leq CM^2}f(y)
 \ll Q\log Q\,\langle f\rangle_{cI_0,1}.
\]
Consequently,
\begin{align}\label{e:sparse21}
    \left|
 \left\langle
 \mathbf1_{\{\log\tau>C_0L^2\}}
 \operatorname{Low}_{Q,\tau}f,g
 \right\rangle
 \right|
 \ll |I_0|\delta_0Q\log Q.
\end{align}
For the high term, Cauchy--Schwarz and \eqref{e:primehighmax} give
\begin{align}\label{e:sparse22}
 \left|
 \left\langle
 \mathbf1_{\{\log\tau>C_0L^2\}}
 \operatorname{High}_{Q,\tau}f,g
 \right\rangle
 \right|&\leq
 \|g\|_2
 \left\|\sup_{M>\exp(C_0L^2)}
 |\operatorname{High}_{Q,M}f|\right\|_2\\
 &\ll
 |I_0|\delta_0^{1/2}Q^{-1/2}
 \exp\left((2+\varepsilon)\frac{(\log Q)^{1/2}}{\log\log Q}\right).
\end{align}
For the error term, \eqref{e:primeerror}  gives
\begin{align*}
 \left\|\sup_{\log M>C_0L^2}|E_{Q,M}f|\right\|_2^2
 &\leq\sum_{\log M>C_0L^2}\|E_{Q,M}f\|_2^2\\
 &\ll_A Q^{-1}(\log Q)^{2C_A}
 \sum_{\log M>C_0L^2}(1+\log M)^{-2A}\|f\|_2^2\\
 &\ll_A Q^{-1}(\log Q)^{2C_A}L^{2-4A}\|f\|_2^2,
\end{align*}
where the sums are over dyadic \(M\).  Hence
\begin{align}\label{e:sparse23}
    \left|
 \left\langle
 \mathbf1_{\{\log\tau>C_0L^2\}}E_{Q,\tau}f,g
 \right\rangle
 \right|
 \ll_A |I_0|\delta_0^{1/2}Q^{-1/2}
 (\log Q)^{C_A}L^{1-2A}.
\end{align}
For any fixed $A>1/2$, the remaining factor
$(\log Q)^{C_A}L^{1-2A}$ is a fixed power of $L$.  Since
$\log Q\asymp L$, it is absorbed by
\[
 \exp\left(
 \varepsilon\frac{(\log Q)^{1/2}}{\log\log Q}
 \right)
\]
for every fixed $\varepsilon>0$.  Thus the error is absorbed into the
high contribution after replacing $2+\varepsilon$ by $2+2\varepsilon$.

Combining the estimates \eqref{e:sparse1}, \eqref{e:decomposproof}, \eqref{e:sparse21}, \eqref{e:sparse22}, and \eqref{e:sparse23}, we obtain
\begin{align*}
 |I_0|^{-1}
 \big|\langle\widetilde A_{\tau,\Lambda}f,g\rangle\big|
 &\ll\delta_0^{2/3}\Phi_{\varepsilon}(\delta_0)
 +\delta_0Q\log Q+\delta_0^{1/2}Q^{-1/2}
 \exp\left((2+2\varepsilon)\frac{(\log Q)^{1/2}}{\log\log Q}\right).
\end{align*}
By the choice of \(Q\) in \eqref{e:Qchoice}, with \(c_1=2+2\varepsilon\),
the last two terms are bounded by $\delta_0^{2/3}
 \Phi_{4\sqrt3/9+O(\varepsilon)+o(1)}(\delta_0).$
The first term in \eqref{e:sparse1} is also absorbed into this bound after
fixing \(\varepsilon>0\).  Letting \(\varepsilon\downarrow0\) proves
\eqref{e:primeendpointlocal}.
\end{proof}

\begin{proof}[Proof of Theorem~\ref{t:sparsebound}]
For the square averages $A_{N,\mathbb N^2}$, the local estimates required
in Lemma~\ref{l:localsparse} are proved in
\cite[Lemmas~5.2--5.6]{LHF}: Lemma~5.2 controls the tails of the
smoothed kernel, Lemmas~5.3 and 5.4 give the two high contributions,
and Lemma~5.6 uses precisely the admissibility condition
\eqref{e:admisslocal}.  Their parameters are
\[
 r=\tfrac12,\qquad s=1,\qquad v=\tfrac12,\qquad
 \gamma=2,\qquad\beta=1.
\]
Lemmas~\ref{l:localsparse} and \ref{l:sparserecursion} therefore give
the asserted endpoint sparse form for $A_{N,\mathbb N^2}$.

For the annular prime-square averages,
Lemmas~\ref{l:primeendpointlocal} and \ref{l:sparserecursion},
applied at the spatial scale \(N^2\), give the corresponding
sparse bound over dyadic arithmetic cutoffs.
For $f\geq0$, positivity and annular reconstruction give
\[
 A_\Lambda^*f
 \leq4\sup_{\substack{N\in2^{\mathbb N}\\N\geq2}}
 |\widetilde A_{N,\Lambda}f|.
\]
We therefore obtain
\begin{align}\label{e:diagonalest}
 \left\langle\sup_N|A_{N,\Lambda}\mathbf1_F|,
 \mathbf1_G\right\rangle
 \ll\sum_{I\in\mathcal S}|I|
 \left(\frac{|F\cap cI||G\cap I|}{|I|^2}\right)^{2/3}
 \Phi_{4\sqrt3/9+o(1)}
 \left(\frac{|F\cap cI||G\cap I|}{|I|^2}\right).
\end{align}
It remains to obtain the open sparse region.  First,
\(\sup_N|A_{N,\Lambda}\mathbf1_F(x)|\leq C\) uniformly in \(x\), and
hence
\[
 \left|\left\langle
 \sup_N|A_{N,\Lambda}\mathbf1_F|,\mathbf1_G
 \right\rangle\right|\leq C|G|.
\]
This is the local restricted estimate corresponding to
\((1/p,1/q)=(0,1)\).  Next, by \cite[Theorem~1]{melhop},
\[
 \left\|\sup_N|A_{N,\Lambda}f|\right\|_u
 \ll_u\|f\|_u,\qquad u>1.
\]
Therefore H\"older's inequality gives
\[
 \left|\left\langle
 \sup_N|A_{N,\Lambda}\mathbf1_F|,\mathbf1_G
 \right\rangle\right|
 \ll_u |F|^{1/u}|G|^{1/u'}.
\]
After localization to an interval \(I\), this is the restricted local
estimate at
\[
 \left(\frac1p,\frac1q\right)
 =\left(\frac1u,\frac1{u'}\right),\qquad u>1.
\]
As \(u\) ranges over \((1,\infty)\), these points fill the open segment
between \((1,0)\) and \((0,1)\).  Interpolating them with
\eqref{e:diagonalest} gives every \((1/p,1/q)\) in the interior of the
triangle with vertices
\[
 (2/3,2/3),\qquad(1,0),\qquad(0,1).
\]
Applying Lemma~\ref{l:sparserecursion} to those interpolated local
estimates gives \eqref{e:introsparseopen}.
\end{proof}

\section{The high--low estimates}\label{s:hilo}

In this section we prove Proposition~\ref{p:decompositionofan}.  We work
with the annular averages
\begin{align}\label{e:tildedef2}
     \widetilde A_{N,\Lambda}f
 =A_{N,\Lambda}f-\tfrac12A_{N/2,\Lambda}f
 =\frac1N\sum_{N/2<n\leq N}\Lambda(n)f(\,\cdot-n^2),
\end{align}
where \(N\) is dyadic. For \(0<\beta\leq1\), put
\begin{align}\label{e:M}
 \widehat{\widetilde M}_{N,\beta}(\theta)
 =\frac1N\int_{N/2}^{N}t^{\beta-1}e(-t^2\theta)\,dt,
 \qquad \widetilde M_N=\widetilde M_{N,1}.
\end{align}
Splitting the integral at \(\|\theta\|^{-1/2}\) and integrating by parts
gives, uniformly in \(0<\beta\leq1\),
\begin{align}\label{e:MNless}
 |\widehat{\widetilde M}_{N,\beta}(\theta)|
 \ll N^{\beta-1}\min\{1,(N^2\|\theta\|)^{-1/2}\}+N^{-1}.
\end{align}
Fix an even smooth function \(0\leq\eta\leq1\), equal to one on
\([-1/2,1/2]\) and supported on \([-1,1]\).  For \(t\in\mathbb R\),
write $\eta_t(\theta)=\eta(16^t\theta).$
In denominator sums, \(q\sim2^s\) means
\(2^s\leq q<2^{s+1}\), with \(s\in\mathbb Z_{\geq0}\).
Set
\begin{align}\label{e:defR}
 R=\max\{\exp((\log N)^{1/10}),Q^8\}.
\end{align}
The assumption
\(\log Q\leq c_0(\log N)^{1/2}\), with \(c_0>0\) sufficiently small,
ensures that Page's theorem applies up to \(R^{11/10}\).  It also gives
\(QR^4=o(N^2)\), so the narrower cutoffs introduced below are contained
in the region where \(\eta_s=1\).

\subsection{The major-arc model and its error}

We use the following annular form of the standard prime-square
circle-method input.

\begin{lemma}\label{l:prime-square-input}
If \((a,q)=1\), \(q\leq R^{11/10}\), and
\(|\theta|\leq R/(qN^2)\), then
\begin{align}\label{e:usualnearzerosiegel}
 \widehat{\widetilde A}_{N,\Lambda}(a/q+\theta)
 &=\mathfrak s(a,q)\widehat{\widetilde M}_N(\theta)
 -\mathbf1_{q_0\mid q}\mathfrak s_{\chi_q}(a,q)
   \widehat{\widetilde M}_{N,\beta}(\theta)
 +O(R^{-5}).
\end{align}
Here \(q_0=q_0(N)\), \(\beta=\beta(N)\), and
\(\chi_0=\chi_0(N)\) are the exceptional conductor, zero, and primitive
character supplied by the Landau--Page theorem at scale \(N\), and
\(\chi_q\) is induced by \(\chi_0\).  These quantities may change when
\(N\) changes. The second term is omitted if no exceptional character
exists. If $\xi$ lies outside all major arcs
\[
 \left\|\xi-a/q\right\|\leq\frac{R}{qN^2},
 \qquad q\leq R,\quad(a,q)=1,
\]
where distance is taken modulo one, then
\begin{align}\label{e:prime-square-minor-input}
 |\widehat{\widetilde A}_{N,\Lambda}(\xi)|
 \ll R^{-1/4}(\log N)^{C_0}.
\end{align}
Here $R\geq1$, $N$ is sufficiently large, and $R^4\leq N/2$,
as is the case for the choice in \eqref{e:defR}.
\end{lemma}

\begin{proof}
Since $\Lambda$ is supported on prime powers, the condition $(n,q)>1$ in \eqref{e:tildedef2}
forces $n=p^\alpha$ with $p\mid q$.  The terms with $\alpha=1$
contribute \(O(qN^{-1}(\log N)^2)\), while all higher prime powers
contribute \(O(N^{-1/2}(\log N)^2)\).  Thus the total error from these
terms is
\[
 O\bigl(N^{-1/2}(\log N)^2+qN^{-1}(\log N)^2\bigr).
\]
Grouping the remaining terms into reduced residue classes modulo \(q\)
therefore gives
\begin{align*}
 \widehat{\widetilde A}_{N,\Lambda}(a/q+\theta)
 &=
 \frac1N\sum_{r\in\mathbb A_q}e(-ar^2/q)
 \sum_{\substack{N/2<n\leq N\\n\equiv r\;(\mathrm{mod}\;q)}}
 \Lambda(n)e(-n^2\theta)\\
 &\quad+
 O\bigl(N^{-1/2}(\log N)^2+qN^{-1}(\log N)^2\bigr).
\end{align*}
By \eqref{e:landaupage}, uniformly for \(N/2\leq x\leq N\),
\[
 \sum_{\substack{n\leq x\\n\equiv r\;(\mathrm{mod}\;q)}}\Lambda(n)
 =
 \frac{x}{\phi(q)}
 -\mathbf1_{q_0\mid q}
 \frac{\chi_q(r)x^\beta}{\beta\phi(q)}
 +O\bigl(xe^{-c\sqrt{\log N}}\bigr).
\]
Partial summation, using
\(
 |\frac{d}{dt}e(-t^2\theta)|\ll|\theta|t
\),
now yields
\begin{align*}
 \frac1N
 \sum_{\substack{N/2<n\leq N\\n\equiv r\;(\mathrm{mod}\;q)}}
 \Lambda(n)e(-n^2\theta)&=
 \frac1{N\phi(q)}\int_{N/2}^{N}e(-t^2\theta)\,dt
 -\mathbf1_{q_0\mid q}
 \frac{\chi_q(r)}{N\phi(q)}
 \int_{N/2}^{N}t^{\beta-1}e(-t^2\theta)\,dt\\
 &\qquad\quad+
 O\left(e^{-c\sqrt{\log N}}
 \bigl(1+|\theta|N^2\bigr)\right).
\end{align*}
Here \(t^{\beta-1}\) arises from differentiating \(t^\beta/\beta\).
Summing over \(r\in\mathbb A_q\) produces the two Gauss sums and gives
\begin{align*}
 \widehat{\widetilde A}_{N,\Lambda}(a/q+\theta)
 &=
 \mathfrak s(a,q)\widehat{\widetilde M}_N(\theta)
 -\mathbf1_{q_0\mid q}\mathfrak s_{\chi_q}(a,q)
  \widehat{\widetilde M}_{N,\beta}(\theta)\\
 &\quad+
 O\left(
 N^{-1/2}(\log N)^2+qN^{-1}(\log N)^2
 +qe^{-c\sqrt{\log N}}\bigl(1+|\theta|N^2\bigr)
 \right).
\end{align*}
Since \(q\leq R^{11/10}\) and
\(|\theta|N^2\leq R/q\), the error is
\[
 \ll N^{-1/2}(\log N)^2
 +R^{11/10}N^{-1}(\log N)^2
 +R^2e^{-c\sqrt{\log N}}
 \ll R^{-5},
\]
by \eqref{e:defR} and the sufficiently small choice of \(c_0\).
This proves \eqref{e:usualnearzerosiegel}.

For the minor arcs put $x=N/2$. Choose a reduced approximation $a/q$
with $q\leq x^{3/2}$, 
\begin{align}\label{e:balancedapp}
 |\xi-a/q|\leq(qx^{3/2})^{-1},
 \qquad |\xi-a/q|<q^{-2}.
\end{align}
If $q\leq x^{1/2}$, apply \cite[Theorem~1.1]{Xiumin2005}
directly to the annular sum. With $T=q(1+|\xi-a/q|x^2)$, being outside
the major arcs gives $T\geq R/4$, while \eqref{e:balancedapp} gives
$T\ll x^{1/2}$. Absorbing $\tau(q)^{3/2}\ll q^{1/12}$ in this range
yields the bound
\[
 |\widehat{\widetilde A}_{N,\Lambda}(\xi)|
 \ll (\log N)^{C_0}R^{-5/12}.
\]
For $x^{1/2}\leq q\leq x^{3/2}$, use Ghosh's quadratic prime
exponential-sum estimate \cite{Ghosh1981}, with its fixed exponent $1/16$.
Its normalized bound is $O(x^{-1/16})$, and $R^4\leq x$ implies
$x^{-1/16}\leq R^{-1/4}$. Combining the two cases proves
\eqref{e:prime-square-minor-input}.
\end{proof}

Define the major-arc model by
\begin{align}\label{e:VNdef}
 \widehat V_N(\xi)
 =\sum_{q\leq R}\sum_{a\in\mathbb A_q}
 \bigl(
 &\mathfrak s(a,q)\widehat{\widetilde M}_N(\xi-a/q) -\mathbf1_{q_0\mid q}\mathfrak s_{\chi_q}(a,q)
 \widehat{\widetilde M}_{N,\beta}(\xi-a/q)
 \bigr)\eta_s(\xi-a/q),
\end{align}
and put
\[
 E_{Q,N}=\widetilde A_{N,\Lambda}-V_N.
\]

\begin{lemma}[Approximation]\label{l:aaplemma}
For every \(A>0\),
\begin{align}\label{e:aaplemma}
 \|E_{Q,N}f\|_2
 \ll_A Q^{-1/2}(1+\log N)^{-A}\|f\|_2.
\end{align}
\end{lemma}

\begin{proof}
Fix \(\xi\in\mathbb T\). The localized products in
\eqref{e:VNdef} are understood to be periodized.
Distinct reduced fractions with \(q,q'\sim2^s\) are separated by
at least \((qq')^{-1}\geq2^{-2s-2}\), whereas
\(\eta_s(\xi-a/q)\) can be nonzero only when
\(\|\xi-a/q\|\leq16^{-s}\).
Thus only \(O(1)\) fractions contribute at each dyadic denominator
scale, each with \(O(1)\) active periodic lifts.

Consider the terms in \(\widehat V_N(\xi)\) for which
\[
 \|\xi-a/q\|>\frac{R}{qN^2}.
\]
Write \(\theta=\xi-a/q\). By \eqref{e:MNless}, uniformly in
\(0<\beta\leq1\) and \(j\in\mathbb Z\),
\[
 |\widehat{\widetilde M}_N(\theta+j)|
 +|\widehat{\widetilde M}_{N,\beta}(\theta+j)|
 \ll(q/R)^{1/2}+N^{-1}.
\]
Lemmas~\ref{l:inftyboundsaq} and \ref{l:Gcharbound}, with
\(\varepsilon=1/12\), bound both Gauss coefficients by
\(O(q^{-5/12})\). Each such fraction therefore contributes
\(O(R^{-1/2}q^{1/12}+N^{-1})\).
Summing over the \(O(\log(eR))\) denominator scales gives a total
contribution of size
\begin{align}\label{e:nondiagappaq}
    O\bigl((R^{-5/12}+N^{-1}) \log(eR)\bigr)=O(R^{-1/3}),
\end{align}
where we use \(R=N^{o(1)}\) in the permitted range of \(Q\).

If \(\xi\in\mathfrak M(R;N)\), there is a reduced fraction
\(a/q\), with \(q\leq R\), such that
\[
 \left\|\xi-\frac aq\right\|\leq\frac{R}{qN^2}.
\]
Such a fraction is unique by rational separation and
\(R^3=o(N^2)\).
Write \(\xi=a/q+\theta\pmod1\), with
\(|\theta|\leq R/(qN^2)\).
If \(q\sim2^s\), then
\[
 16^s|\theta|
 \leq q^4\frac{R}{qN^2}
 \leq\frac{R^4}{N^2}<\frac12,
\]
so \(\eta_s(\theta)=1\).
Lemma~\ref{l:prime-square-input} gives
\begin{align}\label{e:majortildeest}
 \widehat{\widetilde A}_{N,\Lambda}(\xi)
 &=
 \mathfrak s(a,q)\widehat{\widetilde M}_N(\theta)
 -\mathbf1_{q_0\mid q}\mathfrak s_{\chi_q}(a,q)
   \widehat{\widetilde M}_{N,\beta}(\theta)
 +O(R^{-5}).
\end{align}
The two main terms are the \(j=0\) terms in the contribution of
\(a/q\) to \(\widehat V_N(\xi)\).
The other active lifts satisfy \(|\theta+j|\geq1/2\) and contribute
\(O(N^{-1})\) in total by integration by parts in \eqref{e:M}.
Every other contributing fraction lies outside its corresponding
major arc. Thus \eqref{e:nondiagappaq} and
\eqref{e:majortildeest} give
\[
 |\widehat{\widetilde A}_{N,\Lambda}(\xi)-\widehat V_N(\xi)|
 \ll R^{-5}+N^{-1}+R^{-1/3}
 \ll R^{-1/3}.
\]

On the complementary minor arcs
\(\xi\in\mathfrak m(R;N)\),
\eqref{e:prime-square-minor-input} gives
\[
 |\widehat{\widetilde A}_{N,\Lambda}(\xi)|
 \ll R^{-1/4}(\log N)^{C_0}
 \ll R^{-1/6}.
\]
The last inequality follows from
\(R\geq\exp((\log N)^{1/10})\), for sufficiently large \(N\).
Every fraction contributing to \(\widehat V_N(\xi)\) now lies
outside its corresponding major arc, so \eqref{e:nondiagappaq} gives \(|\widehat V_N(\xi)|\ll R^{-1/3}\).
Consequently,
\[
 |\widehat{\widetilde A}_{N,\Lambda}(\xi)-\widehat V_N(\xi)| \leq |\widehat{\widetilde A}_{N,\Lambda}(\xi)|+|\widehat V_N(\xi)| 
 \ll R^{-1/6}.
\]

Combining the two cases, we obtain
\(\|\widehat E_{Q,N}\|_{L^\infty(\mathbb T)}\ll R^{-1/6}\).
Finally, \eqref{e:defR} gives
\[
 R^{-1/6}=R^{-1/16}R^{-5/48}
 \leq Q^{-1/2}
 \exp\left(-\frac5{48}(\log N)^{1/10}\right)
 \ll_A Q^{-1/2}(1+\log N)^{-A}.
\]
Plancherel proves \eqref{e:aaplemma}.
The remaining bounded scales are absorbed into the implied
constant, since \(Q\) is then bounded.
\end{proof}

\subsection{The decomposition}

For \(q\leq Q\), set
\begin{align}\label{e:mqdef}
 m_q=\frac{\log(N^2q/Q)}{\log16},
 \qquad
 16^{-m_q}=\frac{Q}{qN^2}.
\end{align}
If an exceptional primitive real character \(\chi_0=\chi_0(N)\) of
conductor \(q_0\) exists at the current scale \(N\), let
\(\mathcal B_0\) be the constant from
Lemma~\ref{l:exceptionalarithmetic}.  Put
\[
 D=\frac{Q}{\mathcal B_0},\qquad
 \ell_d=\frac{\log(\mathcal B_0dN^2/Q)}{\log16},\qquad
 16^{-\ell_d}=\frac{Q}{\mathcal B_0dN^2}
 \quad(1\le d\le D).
\]
In the exceptional sums below, \(q=q_0d\) and \(\chi_q\) denotes the
character modulo \(q\) induced by \(\chi_0\).  When \(D<1\), the
corresponding low sum is empty.  If no exceptional character exists,
all sums involving \(q_0\) are omitted.

Define
\begin{align}\label{e:lowpart1}
 \widehat{\operatorname{Low}}_{Q,N}(\xi)
 &=\sum_{q\leq Q}\sum_{a\in\mathbb A_q}
 \mathfrak s(a,q)\widehat{\widetilde M}_N(\xi-a/q)
 \eta_{m_q}(\xi-a/q)\notag\\
 &\quad-\sum_{\substack{q=q_0d\leq R\\d\leq D}}
 \sum_{a\in\mathbb A_q}
 \mathfrak s_{\chi_q}(a,q)
 \widehat{\widetilde M}_{N,\beta}(\xi-a/q)
 \eta_{\ell_d}(\xi-a/q).
\end{align}
The high multiplier is the complementary part of \(V_N\) in \eqref{e:VNdef}:
\begin{align}\label{e:highpart}
 \widehat{\operatorname{High}}_{Q,N}(\xi)
 &=\sum_{Q<q\leq R}\sum_{a\in\mathbb A_q}
 \mathfrak s(a,q)\widehat{\widetilde M}_N(\xi-a/q)
 \eta_s(\xi-a/q)\notag\\
 &\quad+\sum_{q\leq Q}\sum_{a\in\mathbb A_q}
 \mathfrak s(a,q)\widehat{\widetilde M}_N(\xi-a/q)
 (\eta_s-\eta_{m_q})(\xi-a/q)\notag\\
 &\quad-\sum_{\substack{q=q_0d\leq R\\d>D}}
 \sum_{a\in\mathbb A_q}
 \mathfrak s_{\chi_q}(a,q)
 \widehat{\widetilde M}_{N,\beta}(\xi-a/q)
 \eta_s(\xi-a/q)\notag\\
 &\quad-\sum_{\substack{q=q_0d\leq R\\d\leq D}}
 \sum_{a\in\mathbb A_q}
 \mathfrak s_{\chi_q}(a,q)
 \widehat{\widetilde M}_{N,\beta}(\xi-a/q)
 (\eta_s-\eta_{\ell_d})(\xi-a/q).
\end{align}
For each ordinary term with \(q\leq Q\), the low part contains the
cutoff \(\eta_{m_q}\), while the high part contains the complementary
cutoff \(\eta_s-\eta_{m_q}\).  Similarly, for each exceptional term
\(q=q_0d\) with \(d\leq D\), the two cutoffs are
\(\eta_{\ell_d}\) and \(\eta_s-\eta_{\ell_d}\).  Thus
\[
 \eta_{m_q}+(\eta_s-\eta_{m_q})=\eta_s,
 \qquad
 \eta_{\ell_d}+(\eta_s-\eta_{\ell_d})=\eta_s.
\]
The ordinary terms with \(q>Q\) and the exceptional terms with \(d>D\)
occur entirely in the high part, with cutoff \(\eta_s\).  Consequently,
\(\operatorname{Low}_{Q,N}+\operatorname{High}_{Q,N}=V_N\).
Since \(E_{Q,N}=\widetilde A_{N,\Lambda}-V_N\), we obtain
\begin{align}\label{e:decompdecomp}
 \widetilde A_{N,\Lambda}
 =\operatorname{Low}_{Q,N}
 +\operatorname{High}_{Q,N}+E_{Q,N}.
\end{align}

\subsection{The low estimate}
We prove the \(\ell^1\to\ell^\infty\) bound for the low part.
\begin{proposition}\label{p:lowcontribution}
For nonnegative \(f\),
\begin{align}\label{e:lowlocal}
 |\operatorname{Low}_{Q,N}f(x)|
 \ll Q\log(2Q)
 \sup_{M\geq N}M^{-2}
 \sum_{|y-x|\leq CM^2}f(y).
\end{align}
Consequently,
\[
 \|\operatorname{Low}_{Q,N}f\|_\infty
 \ll N^{-2}Q\log(2Q)\|f\|_1.
\]
\end{proposition}

\begin{proof}
Fourier inversion and Lemma~\ref{l:lowramanujan} give, uniformly in \(x\),
\[
 \sum_{q\leq Q}\frac1q
 \left|
 \sum_{a\in\mathbb A_q}\mathfrak s(a,q)e(ax/q)
 \right|
 \ll\log(2Q).
\]
The cutoff \(\eta_{m_q}\) has spatial scale
\(qN^2/Q\leq N^2\). By scaling and the rapid decay of \(\eta^\vee\), for every \(K\geq1\),
\[
 |\eta_{m_q}^\vee(x)|
 \ll_K 16^{-m_q}(1+16^{-m_q}|x|)^{-K}
 =
 \frac{Q}{qN^2}
 \left(1+\frac{Q|x|}{qN^2}\right)^{-K}.
\]
The measure corresponding to \(\widehat{\widetilde M}_N\) has total
mass \(O(1)\) and is supported in an interval \(J_N\) of length
\(O(N^2)\).  Since \(qN^2/Q\leq N^2\), convolution with this measure
gives
\[
 |(\widetilde M_N*\eta_{m_q}^\vee)(x)|
 \ll
 \frac{Q}{qN^2}
 \left(1+\frac{\operatorname{dist}(x,J_N)}{N^2}\right)^{-10}
 \ll
 \frac{Q}{qN^2}
 \left(1+\frac{|x|}{N^2}\right)^{-10},
\]
where the last inequality uses \(J_N\subset[-CN^2,CN^2]\).
Consequently, Fourier inversion and Lemma~\ref{l:lowramanujan} give
\begin{align*}
 |\sum_{q\leq Q}\sum_{a\in\mathbb A_q}
 \mathfrak s(a,q)e(ax/q)(\widetilde M_N\ast 
 \eta_{m_q}^{\vee})(x)|
 &\ll
 \frac{Q}{N^2}\left(1+\frac{|x|}{N^2}\right)^{-10}
 \sum_{q\leq Q}\frac1q
 \left|\sum_{a\in\mathbb A_q}
 \mathfrak s(a,q)e(ax/q)\right|\\
 &\ll
 \frac{Q\log(2Q)}{N^2}
 \left(1+\frac{|x|}{N^2}\right)^{-10}.
\end{align*}
We now estimate the exceptional part of the low kernel.  By Fourier
inversion, this kernel is
\[
 -\sum_{\substack{q=q_0d\leq R\\d\leq D}}
 \left(
 \sum_{a\in\mathbb A_q}
 \mathfrak s_{\chi_q}(a,q)e(ax/q)
 \right)
 (\widetilde M_{N,\beta}*\eta_{\ell_d}^{\vee})(x).
\]
The measure corresponding to
\(\widehat{\widetilde M}_{N,\beta}\) is supported in
\([N^2/4,N^2]\), and its total mass satisfies
\[
 \frac1N\int_{N/2}^{N}t^{\beta-1}\,dt
 \ll N^{\beta-1}\leq1.
\]
Moreover, since
\(16^{-\ell_d}=Q/(\mathcal B_0dN^2)\), rapid decay of
\(\eta^\vee\) gives
\[
 |\eta_{\ell_d}^{\vee}(x)|
 \ll
 \frac{Q}{\mathcal B_0dN^2}
 \left(
 1+\frac{Q|x|}{\mathcal B_0dN^2}
 \right)^{-10}.
\]
For \(d\leq D=Q/\mathcal B_0\), the corresponding spatial scale
\(\mathcal B_0dN^2/Q\) is at most \(N^2\).  Convolving with
\(\widetilde M_{N,\beta}\) therefore yields
\[
 |(\widetilde M_{N,\beta}*\eta_{\ell_d}^{\vee})(x)|
 \ll
 \frac{Q}{\mathcal B_0dN^2}
 \left(1+\frac{|x|}{N^2}\right)^{-10},
\]
after enlarging the absolute constant to account for the support
\([N^2/4,N^2]\).

If the arithmetic sum corresponding to \(q=q_0d\) is nonzero,
Lemma~\ref{l:exceptionalarithmetic} gives
\[
 \left|
 \sum_{a\in\mathbb A_{q_0d}}
 \mathfrak s_{\chi_{q_0d}}(a,q_0d)e(ax/(q_0d))
 \right|
 \ll
 \mathcal B_0\prod_{p\mid d}\frac{p+1}{p-1}.
\]
Consequently,
\[
 |\sum_{\substack{q=q_0d\leq R\\d\leq D}}
 \left(
 \sum_{a\in\mathbb A_q}
 \mathfrak s_{\chi_q}(a,q)e(ax/q)
 \right)
 (\widetilde M_{N,\beta}*\eta_{\ell_d}^{\vee})(x)|
 \ll
 \frac{Q}{N^2}
 \left(1+\frac{|x|}{N^2}\right)^{-10}
 \sum_{\substack{d\leq D\\
       \mathfrak s_{\chi_{q_0d}}\not\equiv0}}
 \frac1d\prod_{p\mid d}\frac{p+1}{p-1}.
\]
Here we have discarded the additional restriction \(q_0d\leq R\).

It remains to estimate the sum over \(d\).  By
Lemma~\ref{l:exceptionalarithmetic}, every contributing \(d\) can be
written as
\[
 d=2^\nu d_1,\qquad
 0\leq\nu\leq5,\qquad
 \mu^2(d_1)=1,\qquad
 (d_1,2q_0)=1.
\]
The finitely many possible powers of \(2\) contribute only an absolute
constant.  Hence
\begin{align*}
 \sum_{\substack{d\leq D\\
       \mathfrak s_{\chi_{q_0d}}\not\equiv0}}
 \frac1d\prod_{p\mid d}\frac{p+1}{p-1}
 &\ll
 \sum_{\substack{d_1\leq D\\
       \mu^2(d_1)=1\\(d_1,2q_0)=1}}
 \frac1{d_1}\prod_{p\mid d_1}\frac{p+1}{p-1}\leq
 \prod_{\substack{p\leq D\\p\nmid2q_0}}
 \left(1+\frac{p+1}{p(p-1)}\right).
\end{align*}
Since
\[
 \frac{p+1}{p(p-1)}=\frac1p+O(p^{-2}),
\]
Mertens' theorem gives
\[
 \prod_{\substack{p\leq D\\p\nmid2q_0}}
 \left(1+\frac{p+1}{p(p-1)}\right)
 \ll
 \exp\left(\sum_{p\leq D}\frac1p+O(1)\right)
 \ll\log(2D)\leq\log(2Q).
\]
We have therefore proved
\[
 |\sum_{\substack{q=q_0d\leq R\\d\leq D}}
 \left(
 \sum_{a\in\mathbb A_q}
 \mathfrak s_{\chi_q}(a,q)e(ax/q)
 \right)
 (\widetilde M_{N,\beta}*\eta_{\ell_d}^{\vee})(x)|
 \ll
 \frac{Q\log(2Q)}{N^2}
 \left(1+\frac{|x|}{N^2}\right)^{-10}.
\]
Combining this with the ordinary low-kernel estimate gives the same
bound for the full low kernel.

Finally, for \(f\geq0\), split the convolution into the regions
\[
 |x-y|\leq CN^2
 \quad\text{and}\quad
 2^jN^2<|x-y|\leq2^{j+1}N^2,\qquad 2^j\geq C.
\]
On the \(j\)-th region the kernel is
\(O(Q\log(2Q)N^{-2}2^{-10j})\).  Taking
\(M=2^{j/2}N\), up to a fixed factor, gives
\[
 N^{-2}2^{-10j}
 \sum_{|x-y|\leq C2^jN^2}f(y)
 \ll
 2^{-9j}
 \sup_{M\geq N}M^{-2}
 \sum_{|x-y|\leq CM^2}f(y).
\]
Summing the geometric series proves \eqref{e:lowlocal}.  Since
\[
 \sup_{M\geq N}M^{-2}
 \sum_{|x-y|\leq CM^2}f(y)
 \leq N^{-2}\|f\|_1,
\]
the global \(\ell^1\to\ell^\infty\) estimate follows.
\end{proof}

\subsection{The high estimate}

We use the following standard multifrequency estimate for the terms in
\eqref{e:highpart} with cutoff $\eta_s$.  We include its proof for
completeness, including the uniformity in $0<\beta\leq1$ and in the
coefficients.  The two terms containing
\(\eta_s-\eta_{m_q}\) and \(\eta_s-\eta_{\ell_d}\) are treated
separately by square-function estimates.

\begin{lemma}[Multifrequency estimate]\label{l:notBourgain}
Let \(0<\beta\leq1\), and let \(\mathcal R_s\) be any subset of the
reduced fractions with denominators \(q\sim2^s\).  Then
\begin{align}\label{e:notBourgain}
\left\|
 \sup_{\substack{N\in2^{\mathbb N}\\N\gtrsim2^s}}
 \left|
 \left(
 \sum_{a/q\in\mathcal R_s}
 c_{a,q}\widehat{\widetilde M}_{N,\beta}(\xi-a/q)
 \eta_s(\xi-a/q)\widehat f(\xi)
 \right)^\vee
 \right|
 \right\|_2
 \ll
 (1+s)\sup_{a/q\in\mathcal R_s}|c_{a,q}|\,\|f\|_2.
\end{align}
\end{lemma}

\begin{proof}
Write
\(\mathcal R_s=\{\theta_1,\ldots,\theta_J\}\), where
\(\theta_j=a_j/q_j\).  Distinct reduced fractions with denominators
\(q_j\sim2^s\) satisfy
\[
 \|\theta_j-\theta_k\|_{\mathbb T}\gg2^{-2s},
 \qquad
 J\ll2^{2s},
 \qquad
 \log(2+J)\ll1+s.
\]
Put $ \widehat f_j(\theta)
 =\eta_s(\theta)\widehat f(\theta+\theta_j).$ 
The expression inside the supremum in
\eqref{e:notBourgain} equals
\[
 \sum_{j=1}^J c_{a_j,q_j}e(\theta_jx)
 \bigl(\widetilde M_{N,\beta}*f_j\bigr)(x).
\]
For \(\beta=1\), apply the separated-frequency maximal inequality
\cite[Theorem~1.1]{lacey93illinois}, followed by transference to
\(\ell^2(\mathbb Z)\), to the family $ g_j=c_{a_j,q_j}f_j.$
Since the coefficients are independent of \(x\) and \(N\), linearity gives
\begin{align}\label{e:estmulti}
 \left\|\sup_{N\gtrsim2^s}
 \left|\sum_{j=1}^Jc_{a_j,q_j}e(\theta_jx)
 \bigl(\widetilde M_N*f_j\bigr)(x)\right|\right\|_2&\ll
 \log(2+J)
 \left(\sum_{j=1}^J|c_{a_j,q_j}|^2\|f_j\|_2^2\right)^{1/2}\\
 &\leq
 \log(2+J)\sup_j|c_{a_j,q_j}|
 \left(\sum_{j=1}^J\|f_j\|_2^2\right)^{1/2}.
\end{align}
This is bounded by the right-hand side of
\eqref{e:notBourgain}, using $\sum_{j=1}^J\|f_j\|_2^2\ll\|f\|_2^2.$

The estimate \eqref{e:estmulti} also holds uniformly for
\(0<\beta\leq1\).  Indeed, on \(N/2\leq t\leq N\), we can write $t^{\beta-1}=N^{\beta-1}(t/N)^{\beta-1},$
and the functions \(u\mapsto u^{\beta-1}\) have uniformly bounded size
and variation on \([1/2,1]\). Abel summation in the defining integral for
\(\widehat{\widetilde M}_{N,\beta}\) therefore reduces the required
multifrequency estimate to the case
\(\widetilde M_N=\widetilde M_{N,1}\), with an additional factor
\(N^{\beta-1}\), which is bounded by $1$.  
\end{proof}

\begin{corollary}
\label{c:orlicz-prefix}
Let \(0<\beta\leq1\), \(s\geq0\), and
\(|c_{a,q}|\leq b\).  Then
\begin{align}\label{e:orlicz-prefix}
 &\left\|
 \sup_{\substack{N\in2^{\mathbb N}\\N\geq2^s}}
 \sup_K
 \left|
 \left(
 \sum_{\substack{q\sim2^s\\q\leq K}}
 \sum_{a\in\mathbb A_q}
 c_{a,q}
 \widehat{\widetilde M}_{N,\beta}(\xi-a/q)
 \eta_s(\xi-a/q)\widehat f(\xi)
 \right)^\vee
 \right|
 \right\|_2
 \ll (1+s)^2b\|f\|_2.
\end{align}
The same estimate holds if each summand is multiplied by a fixed
multiplier of modulus at most one supported where
\(\eta_s(\xi-a/q)\neq0\).
\end{corollary}

\begin{proof}
Order the denominators \(q\sim2^s\) increasingly.  Every initial
segment is a disjoint union of at most one dyadic block of each
possible length.  If \(\mathcal D_j\) denotes the family of blocks of
one fixed length, the Rademacher--Menshov decomposition gives
\[
 \sup_K\sup_{N\geq2^s}
 \left|\sum_{\substack{q\sim2^s\\q\leq K}}T_{N,q}f\right|
 \leq
 \sum_j
 \left(
 \sum_{I\in\mathcal D_j}
 \sup_{N\geq2^s}
 \left|\sum_{q\in I}T_{N,q}f\right|^2
 \right)^{1/2},
\]
where \(T_{N,q}\) denotes the sum over \(a\in\mathbb A_q\) in
\eqref{e:orlicz-prefix}.

For each block \(I\), Lemma~\ref{l:notBourgain} gives
\[
 \left\|
 \sup_{N\geq2^s}
 \left|\sum_{q\in I}T_{N,q}f\right|
 \right\|_2
 \ll (1+s)b\|f_I\|_2,
\]
where \(\widehat f_I\) is the restriction of \(\widehat f\) to the
frequency supports corresponding to the denominators in \(I\).
These supports are disjoint, and hence, for each \(j\),
\[
 \sum_{I\in\mathcal D_j}\|f_I\|_2^2\leq\|f\|_2^2.
\]
There are \(O(1+s)\) possible block lengths, proving
\eqref{e:orlicz-prefix}.  The same argument applies after inserting
the additional multipliers.
\end{proof}
\begin{proposition}\label{p:highupperbound}
Let $N_0\geq2$ and suppose that $\log Q\leq c_0(\log N_0)^{1/2}.$
Then, for every dyadic $N\geq N_0$,
\begin{align}
 \|\operatorname{High}_{Q,N}f\|_2
 &\ll Q^{-1/2}
 \exp\left(
 (2+o(1))\frac{\sqrt{\log Q}}{\log\log Q}
 \right)\|f\|_2,
 \label{e:highcontribution}\\
 \left\|
 \sup_{\substack{M\geq N_0\\M\in2^{\mathbb N}}}
 |\operatorname{High}_{Q,M}f|
 \right\|_2
 &\ll Q^{-1/2}
 \exp\left(
 (2+o(1))\frac{\sqrt{\log Q}}{\log\log Q}
 \right)\|f\|_2.
 \label{e:highcontributionmaximal}
\end{align}
\end{proposition}

Here \(o(1)\) is taken as \(Q\to\infty\), uniformly in the allowed
range of \(N_0\).  The bounded values of \(Q\) are absorbed into the
implicit constant.

\begin{proof}
It is enough to prove \eqref{e:highcontributionmaximal}, since
\eqref{e:highcontribution} follows by taking \(N_0=N\).  At scale \(M\),
write
\[
 R(M)=\max\{\exp((\log M)^{1/10}),Q^8\}.
\]
The allowed range of \(Q\) gives \(R(M)=M^{o(1)}\), uniformly for
\(M\geq N_0\).  Hence every  denominator \(q\leq R(M)\) satisfies
\(M\gtrsim q\).  In particular, \(M\gtrsim2^s\) whenever \(q\sim2^s\),
as required in Lemma~\ref{l:notBourgain}.

Consider first the terms in \eqref{e:highpart} of the form
\[
 \mathfrak s(a,q)\widehat{\widetilde M}_M(\xi-a/q)
 \eta_s(\xi-a/q),
 \qquad q\sim2^s,\quad q>Q.
\]
At a fixed denominator scale, the condition \(q\leq R(M)\)
introduces an upper cutoff in the sum over \(q\), depending on \(M\).
Apply Corollary~\ref{c:orlicz-prefix} with
\(K=\lfloor R(M)\rfloor\), setting the coefficients with \(q\leq Q\)
equal to zero. Together with Lemma~\ref{l:inftyboundsaq}, this bounds
the maximal norm of these terms by
\[
 \sum_{2^s\gtrsim Q}(1+s)^2(2^s)^{-1/2}
 \exp\left(
 (2+o(1))\frac{\sqrt{\log(2^s)}}{\log\log(2^s)}
 \right)\|f\|_2.
\]
The additional fixed logarithmic factor is absorbed into the
\(o(1)\) in the following estimate.
The summands decrease geometrically, up to the displayed
subexponential factor, so this is
\[
 \ll
 Q^{-1/2}
 \exp\left(
 (2+o(1))\frac{\sqrt{\log Q}}{\log\log Q}
 \right)\|f\|_2.
\]
We next consider the terms containing
\[
 \mathfrak s(a,q)\widehat{\widetilde M}_M(\xi-a/q)
 \bigl(\eta_s-\eta_{m_q(M)}\bigr)(\xi-a/q),
 \qquad q\sim2^s,\quad q\leq Q,
\]
where $16^{-m_q(M)}=\frac{Q}{qM^2}.$
Fix \(q,a\), put \(\theta=\xi-a/q\), and let \(M_1\geq N_0\) be the
first dyadic scale for which
\((\eta_s-\eta_{m_q(M_1)})(\theta)\neq0\).  If no such scale exists,
there is nothing to estimate.  Otherwise,
\[
 |\theta|\gg\frac{Q}{qM_1^2}.
\]
For every dyadic \(M\geq M_1\), \eqref{e:MNless} gives
\[
 |\widehat{\widetilde M}_M(\theta)|
 \ll
 (M^2|\theta|)^{-1/2}+M^{-1}
 \ll
 \left(\frac qQ\right)^{1/2}\frac{M_1}{M}.
\]
Here the \(M^{-1}\) term is absorbed because
\(|\theta|\leq1/2\).  Consequently,
\[
 \sum_{\substack{M\geq M_1\\M\in2^{\mathbb N}}}
 |\widehat{\widetilde M}_M(\theta)|^2
 \ll\frac qQ.
\]
At a fixed denominator scale, the rational frequency supports are
disjoint.  Bounding the supremum by the square function and applying
Plancherel therefore gives
\begin{align*}
 &\left\|
 \sup_{\substack{M\geq N_0\\M\in2^{\mathbb N}}}
 \left|
 \left(
 \sum_{q\sim2^s}\sum_{a\in\mathbb A_q}
 \mathfrak s(a,q)\widehat{\widetilde M}_M(\xi-a/q)
 (\eta_s-\eta_{m_q(M)})(\xi-a/q)\widehat f(\xi)
 \right)^\vee
 \right|
 \right\|_2\\
 &\qquad\ll
 \sup_{\substack{q\sim2^s\\a\in\mathbb A_q}}
 |\mathfrak s(a,q)|
 \left(\frac qQ\right)^{1/2}\|f\|_2\ll
 Q^{-1/2}
 \exp\left(
 (2+o(1))\frac{\sqrt{\log(2^s)}}{\log\log(2^s)}
 \right)\|f\|_2.
\end{align*}
Summing over \(2^s\lesssim Q\) gives
\[
 \ll
 Q^{-1/2}(1+\log Q)
 \exp\left(
 (2+o(1))\frac{\sqrt{\log Q}}{\log\log Q}
 \right)\|f\|_2,
\]
which has the required form after absorbing the logarithmic factor into
the \(o(1)\).

It remains to estimate the exceptional terms.  The exceptional
character may depend on \(M\).  For every primitive real character
\(\chi\) that occurs at some dyadic scale \(M\geq N_0\), let
\[
 \mathcal M(\chi)
 =
 \{M\geq N_0:\chi\text{ is the exceptional character at scale }M\}.
\]
Write \(q_\chi,\beta_\chi,\sigma_\chi\), and \(\mathcal B_\chi\) for
the corresponding conductor, exceptional zero, and constants from
Lemma~\ref{l:exceptionalarithmetic}, and put
\[
 D_\chi=\frac{Q}{\mathcal B_\chi}.
\]
These quantities are fixed when \(M\) is restricted to
\(\mathcal M(\chi)\). Let \(\operatorname{High}^{\mathrm{exc},\chi}_{Q,M}\) denote the sum of
the two exceptional terms in \eqref{e:highpart} when
\(M\in\mathcal M(\chi)\), and set it equal to zero otherwise.
  We first consider the terms
\[
 -\mathfrak s_{\chi_{q_\chi d}}(a,q_\chi d)
 \widehat{\widetilde M}_{M,\beta_\chi}
       (\xi-a/(q_\chi d))
 \eta_s(\xi-a/(q_\chi d)),
 \qquad d>D_\chi.
\]
For the fixed character \(\chi\), the condition
\(q_\chi d\leq R(M)\) restricts the sum to
\(d\leq R(M)/q_\chi\), with an upper limit depending on \(M\).
Apply Corollary~\ref{c:orlicz-prefix} with
\(K=\lfloor R(M)\rfloor\), retaining only denominators
\(q=q_\chi d\) with \(d>D_\chi\), and restrict the time supremum
to \(M\in\mathcal M(\chi)\). At the denominator scale
\(q_\chi d\sim2^s\), the resulting factor is
\[
 (1+s)^2\ll(1+\log(q_\chi d))^2.
\]
By Lemma~\ref{l:exceptionalarithmetic},
\begin{align*}
 (1+\log(q_\chi d))^2
 |\mathfrak s_{\chi_{q_\chi d}}(a,q_\chi d)|&\ll
 (1+\log(q_\chi d))^2
 \sigma_\chi d^{-1/2}\rho(d)\\
 &=
 Q^{-1/2}\sigma_\chi\mathcal B_\chi^{1/2}
 (1+\log(q_\chi d))^2\rho(d)
 \left(\frac{D_\chi}{d}\right)^{1/2}.
\end{align*}
Grouping \(d\) into dyadic ranges and using the maximal-order
estimate for \(\rho\), the maximal norm of these terms is
\begin{align*}
 &\ll
 Q^{-1/2}\sigma_\chi\mathcal B_\chi^{1/2}
 \sum_{2^k\gtrsim D_\chi}
 (1+\log q_\chi+k)^2
 \sup_{d\sim2^k}\rho(d)
 \left(\frac{D_\chi}{2^k}\right)^{1/2}\|f\|_2\\
 &\ll
 Q^{-1/2}\sigma_\chi\mathcal B_\chi^{1/2}
 (1+\log q_\chi+\log Q)^{O(1)}
 \exp\left(
 (2+o(1))\frac{\sqrt{\log Q}}{\log\log Q}
 \right)\|f\|_2.
\end{align*}
Indeed, since \(D_\chi\leq Q\), the ranges \(2^k\leq Q\)
contribute only fixed logarithmic powers. For \(2^k>Q\),
the factor \((D_\chi/2^k)^{1/2}\) gives geometric decay
after including the displayed logarithmic and subexponential
factors.

We finally consider the terms
\[
 -\mathfrak s_{\chi_{q_\chi d}}(a,q_\chi d)
 \widehat{\widetilde M}_{M,\beta_\chi}
       (\xi-a/(q_\chi d))
 \bigl(\eta_s-\eta_{\ell_{d,\chi}(M)}\bigr)
       (\xi-a/(q_\chi d)),
 \qquad d\leq D_\chi,
\]
where
\[
 16^{-\ell_{d,\chi}(M)}
 =\frac{Q}{\mathcal B_\chi dM^2}.
\]
Fix \(d,a\), put \(\theta=\xi-a/(q_\chi d)\), and let
\(M_1\in\mathcal M(\chi)\) be the first scale for which the cutoff
difference is nonzero.  If there is no such scale, the corresponding
term vanishes.  Otherwise,
\[
 |\theta|\gg\frac{Q}{\mathcal B_\chi dM_1^2}.
\]
Using \eqref{e:MNless}, \(M^{\beta_\chi-1}\leq1\), and
\(M_1\geq N_0\), we obtain
\begin{align*}
 \sum_{\substack{M\in\mathcal M(\chi)\\M\geq M_1}}
 |\widehat{\widetilde M}_{M,\beta_\chi}(\theta)|^2
 &\leq
 \sum_{\substack{M\geq M_1\\M\in2^{\mathbb N}}}
 |\widehat{\widetilde M}_{M,\beta_\chi}(\theta)|^2\ll
 \frac{\mathcal B_\chi d}{Q}+M_1^{-2}
 \ll\frac{\mathcal B_\chi d}{Q}.
\end{align*}
In the last inequality we used
\(\mathcal B_\chi d\geq1\) and \(Q\ll N_0^2\).

The square-function and Plancherel argument therefore bounds the
maximal norm at each dyadic \(d\)-scale by
\begin{align*}
 &\sup_{\substack{d\sim2^k\\a\in\mathbb A_{q_\chi d}}}
 |\mathfrak s_{\chi_{q_\chi d}}(a,q_\chi d)|
 \left(\frac{\mathcal B_\chi d}{Q}\right)^{1/2}\|f\|_2\ll
 Q^{-1/2}\sigma_\chi\mathcal B_\chi^{1/2}
 \sup_{d\sim2^k}\rho(d)\,\|f\|_2.
\end{align*}
Since \(d\leq D_\chi\leq Q\), summing the dyadic \(d\)-scales gives
\[
 \ll
 Q^{-1/2}\sigma_\chi\mathcal B_\chi^{1/2}
 (1+\log Q)
 \exp\left(
 (2+o(1))\frac{\sqrt{\log Q}}{\log\log Q}
 \right)\|f\|_2.
\]
Combining the two exceptional contributions, we have proved
\begin{align}\label{e:fixedcharacterhigh}
 \left\|
 \sup_{M\in\mathcal M(\chi)}
 |\operatorname{High}^{\mathrm{exc},\chi}_{Q,M}f|
 \right\|_2
 &\ll
 Q^{-1/2}\sigma_\chi\mathcal B_\chi^{1/2}
 (1+\log q_\chi+\log Q)^{O(1)}
 \exp\left(
 (2+o(1))\frac{\sqrt{\log Q}}{\log\log Q}
 \right)\|f\|_2.
\end{align}
At each scale \(M\), at most one exceptional character occurs.
Consequently, pointwise,
\begin{align}\label{e:pontchar}
 \sup_{M\geq N_0}
 |\operatorname{High}^{\mathrm{exc}}_{Q,M}f|^2
 \leq
 \sum_\chi
 \sup_{M\in\mathcal M(\chi)}
 |\operatorname{High}^{\mathrm{exc},\chi}_{Q,M}f|^2.
\end{align}
If \(q_\chi\sim2^j\), then, after decreasing the absolute constant in
the definition of the Page region, the zero of \(\chi\) lies in the
exceptional region for the modulus bound \(2^{j+1}\).  The
Landau--Page theorem therefore permits at most one such character in
each dyadic conductor range.  By the stronger estimate in
Lemma~\ref{l:exceptionalarithmetic},
\begin{align}\qquad\label{e:sumoverchi}
 &\sum_\chi
 \left(
 \sigma_\chi\mathcal B_\chi^{1/2}
 (1+\log q_\chi+\log Q)^{O(1)}
 \right)^2\ll
 \sum_{j\geq0}
 2^{-2j/25}(1+j+\log Q)^{O(1)}
 \ll(1+\log Q)^{O(1)}.
\end{align}
Combining \eqref{e:pontchar} with
\eqref{e:fixedcharacterhigh}, and then using \eqref{e:sumoverchi} gives the required exceptional contribution.  The
remaining logarithmic factor is absorbed into the \(o(1)\), proving
\eqref{e:highcontributionmaximal}.
\end{proof}

Proposition~\ref{p:lowcontribution},
Proposition~\ref{p:highupperbound}, Lemma~\ref{l:aaplemma}, and the
exact identity \eqref{e:decompdecomp} prove
Proposition~\ref{p:decompositionofan}.  This completes the proofs of
Theorems~\ref{t:improvingfixed} and \ref{t:sparsebound}.

\section{Quantitative maximal estimates near \(L^1\)}
\label{s:orlicz}

This section proves the distributional estimate underlying
Theorem~\ref{t:pet}\textup{(ii)}.
For \(C>0\) and \(t\geq0\), define
\begin{align}\label{e:orlicz-gauge}
 \Psi_C(t)
 =
 \int_0^t
 \exp\left(
 C\sqrt{\log(e+u)\log(e+\log(e+u))}
 \right)\,du.
\end{align}
Then \(\Psi_C\) is an increasing convex function and
\begin{align}\label{e:orlicz-gauge-comparison}
 \Psi_C(t)
 \asymp_C
 t\exp\left(
 C\sqrt{\log(e+t)\log(e+\log(e+t))}
 \right).
\end{align}
In particular, \(\Psi_C(2t)\ll_C\Psi_C(t)\) and
\(\Psi_C(t)\geq t\).
For \(N\geq2\), define the unweighted prime-square averages by
\[
 A_{N,\mathbb P^2}f(x)
 =
 \frac1{\pi(N)}
 \sum_{p\leq N}f(x-p^2).
\]

\begin{theorem}[Orlicz maximal estimate]
\label{t:orlicz-maximal}
 There are absolute constants \(C,C_0\geq1\) such that the
following holds. Choose \(A_N\) from \(A_{N,\Lambda}\),
\(A_{N,\mathbb N^2}\), and \(A_{N,\mathbb P^2}\), and put
\(A^*f=\sup_{N\geq2}|A_Nf|\).
Then, for every finitely supported \(f:\mathbb Z\to\mathbb C\) and
every \(\lambda>0\),
\begin{align}\label{e:orlicz-truncated}
 \bigl|\{x\in\mathbb Z:A^*f(x)>\lambda\}\bigr|
 \ll
 \sum_{\substack{x\in\mathbb Z\\
 |f(x)|>\lambda/C_0}}
 \Psi_C\left(\frac{|f(x)|}{\lambda}\right).
\end{align}
Consequently,
\begin{align}\label{e:orlicz-modular}
 \bigl|\{x\in\mathbb Z:A^*f(x)>\lambda\}\bigr|
 \ll
 \sum_{x\in\mathbb Z}
 \Psi_C\left(\frac{|f(x)|}{\lambda}\right).
\end{align}
In particular, for every finite \(F\subset\mathbb Z\) and \(t\geq1\),
\begin{align}\label{e:orlicz-restricted}
 \bigl|\{x\in\mathbb Z:A^*\mathbf1_F(x)>t^{-1}\}\bigr|
 \ll |F|\Psi_C(t).
\end{align}
\end{theorem}

The truncated estimate gives the following quantitative
\(\ell^p(\mathbb Z)\) bound.

\begin{corollary}\label{c:orlicz-p}
For \(1<p\leq2\) and \(f\in\ell^p(\mathbb Z)\),
\begin{align}\label{e:orlicz-p}
 \|A^*f\|_{\ell^p}
 \leq
 \exp\left(
 \frac{C}{p-1}\log\frac{e}{p-1}
 \right)\|f\|_{\ell^p},
\end{align}
after enlarging the absolute constant \(C\). 
\end{corollary}

\begin{proof}
By \eqref{e:orlicz-truncated} and the layer-cake formula,
\begin{align*}
 \|A^*f\|_{\ell^p}^p
 &=p\int_0^\infty
 \lambda^{p-1}
 \bigl|\{A^*f>\lambda\}\bigr|\,d\lambda\\
 &\ll
 p\sum_x
 \int_0^{C_0|f(x)|}
 \lambda^{p-1}
 \Psi_C\left(\frac{|f(x)|}{\lambda}\right)\,d\lambda\\
 &=p\|f\|_{\ell^p}^p
 \int_{1/C_0}^\infty
 \Psi_C(u)u^{-p-1}\,du.
\end{align*}
Put \(\varepsilon=p-1\).  By
\eqref{e:orlicz-gauge-comparison} and the substitution \(u=e^v\),
\begin{align*}
 \int_{1/C_0}^\infty
 \Psi_C(u)u^{-p-1}\,du
 &\ll
 1+\int_0^\infty
 \exp\left(
 -\varepsilon v+
 C\sqrt{(1+v)\log(e+v)}
 \right)\,dv.
\end{align*}
For \(v\geq0\),
\[
 C\sqrt{(1+v)\log(e+v)}
 \leq
 \frac{\varepsilon v}{2}
 +\frac{C_1}{\varepsilon}\log\frac{e}{\varepsilon},
\]
with an absolute \(C_1>0\).  Hence
\[
 \int_{1/C_0}^\infty
 \Psi_C(u)u^{-p-1}\,du
 \ll
 \varepsilon^{-1}
 \exp\left(
 \frac{C_1}{\varepsilon}\log\frac{e}{\varepsilon}
 \right).
\]
Taking \(p\)-th roots and enlarging \(C\) proves
\eqref{e:orlicz-p} for finitely supported \(f\).
Truncation and Fatou's lemma extend the estimate to every
\(f\in\ell^p(\mathbb Z)\).
\end{proof}
Next, we bound the maximal operator associated with a
Ionescu--Wainger model, then compare it with the annular
arithmetic averages in \(\ell^2(\mathbb Z)\).
Balancing these estimates gives the restricted bound.

\begin{proposition}
\label{p:IWinput7}
Let \(r,k\in\mathbb N\), with \(r\geq2\) and \(k\geq4\),
and let \(Q\geq2^k\). Put
\begin{equation}\label{e:IWparameters7}
 \begin{gathered}
 p_r=\frac{2r}{2r-1},
 \qquad B=C_*^{\,k^2Q^{1/(k-1/2)}},\\
 \varepsilon=(2^{20}rB^{2r+8})^{-1},
 \qquad C_{r,k}=(C_*r\log(2+k))^{C_*k},
 \end{gathered}
\end{equation}
where \(C_*>1\) is a sufficiently large absolute constant.
Then there exist an integer \(P\leq3^Q\) and a
Ionescu--Wainger set
\(\Sigma=\Sigma_{\leq k}\subset P^{-1}\mathbb Z/\mathbb Z\)
with the following properties.

\begin{enumerate}[label=\textup{(\roman*)}]
\item
Every rational point of denominator at most \(Q\) belongs to
\(\Sigma\):
\begin{align}\label{e:IWcontainment7}
 \bigcup_{1\leq q\leq Q}
 \left\{\frac aq\pmod1:a\in\mathbb Z\right\}
 \subseteq\Sigma.
\end{align}
Every element of \(\Sigma\), when written in reduced form, has
denominator at most \(B\). With the choice of \(\varepsilon\) in \eqref{e:IWparameters7}, the
corresponding major-arc parameter set is
\((r,\frac12)\)-good.

\item
Let \(\nu\) be a complex measure on \(\mathbb Z_P\) with
\(\|\nu\|_{\mathrm{TV}}\leq1\), and let \(m\) be a multiplier on
\(\mathbb R\) supported in \([-\varepsilon,\varepsilon]\).  Then, for
\(p_r\leq p\leq2r\),
\begin{align}\label{e:IWmultiplier7}
 \left\|
 T_{\displaystyle
   \sum_{\alpha\in\Sigma}
   \widehat\nu(\alpha)m(\,\cdot-\alpha)}
 \right\|_{\ell^p(\mathbb Z)\to\ell^p(\mathbb Z)}
 \leq
  C_{r,k}
 \|T_m\|_{L^p(\mathbb R)\to L^p(\mathbb R)},
\end{align}
where the translates of \(m\) are periodized on \(\mathbb T\).

\item
Put $\Omega=[-\varepsilon,\varepsilon]\times\Sigma
 \subset\mathbb R\times\mathbb Q/\mathbb Z.$
If \(f\in\ell^p(\mathbb Z)\) has Fourier support in
\(\pi(\Omega)\), then
\begin{align}\label{e:IWsampling7}
  C_{r,k}^{-1}\|f\|_{\ell^p(\mathbb Z)}
 \leq
 \|\mathcal S_\Omega^{-1}f\|_{L^p(\mathbb A_{\mathbb Z})}
 \leq
  C_{r,k}\|f\|_{\ell^p(\mathbb Z)}
\end{align}
for \(p_r\leq p\leq2r\).

\item
Let \(V\) be a finite-dimensional normed space and
\(0<\delta\leq(100P)^{-1}\).  If
\(G:\mathbb R\times\mathbb Z_P\to V\) has real Fourier support in
\([-\delta,\delta]\), then, for \(1\leq p\leq\infty\),
\begin{align}\label{e:IWdiagonal7}
 \left\|
 \bigl(G(n,n\bmod P)\bigr)_{n\in\mathbb Z}
 \right\|_{\ell^p(\mathbb Z;V)}
 \asymp
 \|G\|_{L^p(\mathbb R\times\mathbb Z_P;V)},
\end{align}
with absolute implied constants.
\end{enumerate}
\end{proposition}

\begin{proof}
Apply \cite[Lemma~5.1]{Tao2021} with $\rho=\frac{2}{k-\frac12}.$
Then
\[
 \left\lfloor\frac2\rho\right\rfloor+1=k,
 \qquad
 Q^{\rho/2}=Q^{1/(k-1/2)}.
\]
The lemma produces a set of pairwise coprime integers defining
\(\Sigma_{\leq k}\), together with an integer \(P\leq3^Q\), such that
\eqref{e:IWcontainment7} holds.  Its constituent coprime integers are
at most \(B^{1/k}\), after enlarging \(C_*\), and therefore every
denominator occurring in \(\Sigma_{\leq k}\) is at most \(B\).
The choice of \(\varepsilon\) in \eqref{e:IWparameters7} satisfies
the goodness condition in \cite[(5.1)]{Tao2021}.  This proves \textup{(i)}.

For \textup{(ii)}, define the adelic multiplier
\[
 M(\theta,\alpha)=m(\theta)\widehat\nu(\alpha),
 \qquad
 (\theta,\alpha)\in\mathbb R\times P^{-1}\mathbb Z/\mathbb Z.
\]
Extend \(M\) by zero when \(\alpha\notin P^{-1}\mathbb Z/\mathbb Z\). The associated operator factors into its real and finite-group parts.
Indeed, if \(F\in L^p(\mathbb R\times\mathbb Z_P)\), then
\[
 T_MF(x,z)
 =
 \int_{\mathbb Z_P}
 T_m\bigl(F(\,\cdot\,,z-y)\bigr)(x)\,d\nu(y).
\]
Minkowski's inequality and translation invariance on \(\mathbb Z_P\)
therefore give
\begin{align}\label{e:rzptorzp}
 \|T_MF\|_{L^p(\mathbb R\times\mathbb Z_P)}
 &\leq
 \int_{\mathbb Z_P}
 \left\|
 T_m\bigl(F(\,\cdot\,,\,\cdot-y)\bigr)
 \right\|_{L^p(\mathbb R\times\mathbb Z_P)}
 \,d|\nu|(y)\leq
 \|T_m\|_{L^p(\mathbb R)\to L^p(\mathbb R)}
 \|\nu\|_{\mathrm{TV}}
 \|F\|_{L^p(\mathbb R\times\mathbb Z_P)}.
\end{align}
For a general \(F\in L^p(\mathbb A_{\mathbb Z})\), define
\[
 F_P(x,z)
 =
 \int_{\widehat{\mathbb Z}}F(x,z+Pu)\,du,
 \qquad z\in\mathbb Z_P,
\]
where \(du\) is normalized Haar measure and any integer
representative of \(z\) may be used.
Jensen's inequality gives
\[
 \|F_P\|_{L^p(\mathbb R\times\mathbb Z_P)}
 \leq
 \|F\|_{L^p(\mathbb A_{\mathbb Z})}.
\]
This averaging preserves the Fourier coefficients at
\(\alpha\in\mathbb T[P]\) and removes the remaining arithmetic
frequencies. Since \(M\) vanishes outside these frequencies,
\(T_MF\) is the lift of \(T_MF_P\), and hence
\[
 \|T_MF\|_{L^p(\mathbb A_{\mathbb Z})}
 =
 \|T_MF_P\|_{L^p(\mathbb R\times\mathbb Z_P)}.
\]
Using \eqref{e:rzptorzp} and \(\|\nu\|_{\mathrm{TV}}\leq1\),
we obtain
\begin{align}\label{e:adelictoeuclidean}
 \|T_M\|_{L^p(\mathbb A_{\mathbb Z})\to
 L^p(\mathbb A_{\mathbb Z})}
 \leq
 \|T_m\|_{L^p(\mathbb R)\to L^p(\mathbb R)}.
\end{align}
The multiplier on \(\mathbb T\) obtained from \(M\) and the rational
frequency set \(\Sigma\) is
\[
 M_\Sigma^{\mathbb Z}(\xi)
 =
 \sum_{\alpha\in\Sigma}
 M(\xi-\alpha,\alpha)
 =
 \sum_{\alpha\in\Sigma}
 \widehat\nu(\alpha)m(\xi-\alpha).
\]
Hence the adelic Ionescu--Wainger multiplier theorem
\cite[Theorem~1.7]{Tao2021} and \eqref{e:adelictoeuclidean} give
\begin{align*}
 \|T_{M_\Sigma^{\mathbb Z}}\|_{\ell^p(\mathbb Z)\to\ell^p(\mathbb Z)}
 &\leq
 (r\log(2+k))^{O(k)}
 \|T_M\|_{L^p(\mathbb A_{\mathbb Z})\to
 L^p(\mathbb A_{\mathbb Z})}\leq
  C_{r,k}
 \|T_m\|_{L^p(\mathbb R)\to L^p(\mathbb R)}.
\end{align*}
This is \eqref{e:IWmultiplier7}.

Estimate \eqref{e:IWsampling7} is the major-arc sampling theorem
\cite[Theorem~1.6]{Tao2021}, again with its constant bounded by
\( C_{r,k}\).  This proves \textup{(iii)}.

For \textup{(iv)}, expand \(G\) in the finite variable:
\[
 G(x,z)
 =
 \sum_{a\in\mathbb Z_P}G_a(x)e(az/P),
 \qquad
 \operatorname{supp}\widehat{G_a}
 \subseteq[-\delta,\delta].
\]
Its restriction to the diagonal \(n\mapsto(n,n\bmod P)\) is therefore
\[
 G(n,n\bmod P)
 =
 \sum_{a\in\mathbb Z_P}G_a(n)e(an/P).
\]
Thus the Fourier transform of the \(a\)-th summand is supported in $\frac aP+[-\delta,\delta]\subset\mathbb T.$
If \(a\neq a'\), then
\[
 \left\|\frac aP-\frac{a'}P\right\|_{\mathbb T}
 \geq\frac1P>2\delta.
\]
Hence these intervals are pairwise disjoint.  Equivalently, the map
\[
 \pi:
 [-\delta,\delta]\times\mathbb T[P]\longrightarrow\mathbb T,
 \qquad
 \pi(\theta,\alpha)=\theta+\alpha\pmod1,
\]
is injective.
Put
\[
 \Omega_P=
 \left[-\frac1{100P},\frac1{100P}\right]
 \times\mathbb T[P].
\]
Since \(\delta\leq(100P)^{-1}\), the Fourier support of \(G\) is
contained in \(\Omega_P\), and its sampled function is precisely $\mathcal S_{\Omega_P}G(n)=G(n,n\bmod P).$
The quantitative Shannon sampling theorem
\cite[Proposition~1.5]{Tao2021}, applied with \(Q=P\),
\(c=1/100\), and target space \(V\), now gives
\[
 \left\|
 \bigl(G(n,n\bmod P)\bigr)_{n\in\mathbb Z}
 \right\|_{\ell^p(\mathbb Z;V)}
 \asymp
 \|G\|_{L^p(\mathbb R\times\mathbb Z_P;V)},
\]
with constants independent of \(P\), \(\delta\), and \(G\).  This is
\eqref{e:IWdiagonal7}.
\end{proof}

The following standard continuous estimates are included for
completeness and to record their uniformity in \(\beta\).

\begin{lemma}\label{l:continuousestimates7}
Uniformly in \(0<\beta\leq1\), \(J\geq1\), and
\(|\sigma_j|\leq1\), one has, for \(1<p\leq2\),
\begin{align}
 \left\|
 T_{\sum_{j=1}^J\sigma_j
   (\widehat{\widetilde M}_{2^j,\beta}
    -\widehat{\widetilde M}_{2^{j-1},\beta})}f
 \right\|_{L^p(\mathbb R)}
 &\ll
 \frac1{p-1}\|f\|_{L^p(\mathbb R)},
 \label{e:continuousincrements7}\\
 \left\|
 \sup_{N\geq2}
 |\widetilde M_{N,\beta}*f|
 \right\|_{L^p(\mathbb R)}
 &\ll
 \frac1{p-1}\|f\|_{L^p(\mathbb R)}.
 \label{e:continuousmaximal7}
\end{align}
Here \(\widetilde M_{N,\beta}\) is the continuous annular measure
defined in \eqref{e:M}.  Since \eqref{e:continuousmaximal7} holds with
the supremum over all \(N\geq2\), it also holds when the supremum is
restricted to dyadic \(N\).
\end{lemma}

\begin{proof}
The change of variables \(x=t^2\) in \eqref{e:M} gives
\begin{align}\label{e:continuouskernel7}
 \widetilde M_{N,\beta}(x)
 =
 \frac{x^{\beta/2-1}}{2N}
 \mathbf1_{[N^2/4,N^2]}(x).
\end{align}
On the Fourier side, setting \(t=2^jv\) in \eqref{e:M} gives
\[
 \widehat{\widetilde M}_{2^j,\beta}(\theta)
 =
 2^{j(\beta-1)}F_\beta(4^j\theta),
 \qquad
 F_\beta(u)=
 \int_{1/2}^1v^{\beta-1}e(-uv^2)\,dv.
\]
Integration by parts, together with subtraction at \(u=0\), gives
uniformly in \(\beta\)
\begin{align}\label{e:Fbetabounds7}
 |F_\beta(u)|
 &\ll\min\{1,|u|^{-1}\},\qquad 
 |F_\beta(u)-F_\beta(u/4)|
 \ll\min\{|u|,|u|^{-1}\}.
\end{align}
For each \(j\), adding and subtracting
\(2^{(j-1)(\beta-1)}F_\beta(4^j\theta)\) gives
\begin{align*}
 &\left|
 2^{j(\beta-1)}F_\beta(4^j\theta)
 -
 2^{(j-1)(\beta-1)}F_\beta(4^{j-1}\theta)
 \right|\\
 &\leq
 \left|2^{j(\beta-1)}-2^{(j-1)(\beta-1)}\right|
 |F_\beta(4^j\theta)|+
 2^{(j-1)(\beta-1)}
 \left|F_\beta(4^j\theta)-F_\beta(4^{j-1}\theta)\right|.
\end{align*}
Since \(2^{(j-1)(\beta-1)}\leq1\), it follows that
\begin{align}\label{e:l2prep}
 &\sum_{j=1}^J
 \left|
 \widehat{\widetilde M}_{2^j,\beta}(\theta)
 -
 \widehat{\widetilde M}_{2^{j-1},\beta}(\theta)
 \right|\ll
 \sum_{j=1}^J
 \left|
 2^{j(\beta-1)}-2^{(j-1)(\beta-1)}
 \right|
 +
 \sum_{j=1}^J
 \min\{|4^j\theta|,|4^j\theta|^{-1}\} \ll 1. \qquad 
\end{align}
Let \(K\) be the convolution kernel of the operator in
\eqref{e:continuousincrements7}.  Estimate \eqref{e:l2prep}
and Plancherel's theorem give
\begin{align}\label{e:continuousincrementsL27}
 \|K*f\|_{L^2(\mathbb R)}
 \ll \|f\|_{L^2(\mathbb R)}.
\end{align}
Collecting equal kernels, we have
\[
 K
 =
 -\sigma_1\widetilde M_{1,\beta}
 +\sum_{j=1}^{J-1}
   (\sigma_j-\sigma_{j+1})\widetilde M_{2^j,\beta}
 +\sigma_J\widetilde M_{2^J,\beta}.
\]
Thus every \(\widetilde M_{2^j,\beta}\) occurs with a coefficient of
modulus at most \(2\).  If \(4^{j-1}<x<4^j\), then
\eqref{e:continuouskernel7} and \(0<\beta\leq1\) give
\begin{align*}
 |K(x)|
 &\ll
 \frac{x^{\beta/2-1}}{2^j}
 =
 \frac1x\frac{x^{\beta/2}}{2^j}
 \ll \frac1x, \qquad 
 |K'(x)|
 \ll
 \frac{x^{\beta/2-2}}{2^j}
 =
 \frac1{x^2}\frac{x^{\beta/2}}{2^j}
 \ll \frac1{x^2},
\end{align*}
At an endpoint \(x=4^j\), the jump of \(K\) is bounded by
\[
 \left|K(4^j+)-K(4^j-)\right|
 \ll
 \frac{(4^j)^{\beta/2-1}}{2^j}
 \ll 4^{-j}.
\]
Consequently,
\begin{align}\label{e:Hormander7}
 \int_{|x|>2|h|}
 |K(x-h)-K(x)|\,dx
 &\ll
 |h|\int_{2|h|}^{\infty}\frac{dx}{x^2}
 +
 |h|\sum_{4^j\gtrsim|h|}4^{-j}
 \ll1.
\end{align}
Thus \(K\) is a Calderón--Zygmund kernel, with its \(L^2\) and kernel
constants bounded uniformly in \(J,\beta\), and \((\sigma_j)\).
The Calderón--Zygmund weak endpoint estimate gives
\[
 \bigl|\{x\in\mathbb R:|K*f(x)|>\lambda\}\bigr|
 \ll \lambda^{-1}\|f\|_{L^1(\mathbb R)}.
\]
For \(1<p\leq3/2\), Marcinkiewicz interpolation with
\eqref{e:continuousincrementsL27} gives
\[
 \|K*f\|_{L^p(\mathbb R)}
 \ll (p-1)^{-1}\|f\|_{L^p(\mathbb R)}.
\]
For \(3/2\leq p\leq2\), interpolation between the \(L^{3/2}\) and
\(L^2\) estimates gives
\[
 \|K*f\|_{L^p(\mathbb R)}
 \ll\|f\|_{L^p(\mathbb R)}
 \ll(p-1)^{-1}\|f\|_{L^p(\mathbb R)}.
\]
This proves \eqref{e:continuousincrements7}.

By \eqref{e:continuouskernel7}, we have $0\leq\widetilde M_{N,\beta}(x)
 \ll N^{-2}\mathbf1_{[0,N^2]}(x),$
and hence
\[
 \sup_{N\geq1}
 |\widetilde M_{N,\beta}*f(x)|
 \ll M_{\mathrm{HL}}^+f(x).
\]
The Hardy--Littlewood estimate recalled in Section~\ref{ss:notaion}
proves \eqref{e:continuousmaximal7}.
\end{proof}

\begin{proposition}
\label{p:orlicz-low}
Let \(r,k,Q,P,\Sigma_{\leq k},B,\varepsilon,p_r\), and
\(C_{r,k}\) be as in Proposition~\ref{p:IWinput7}.
Let \(\nu\) be a complex measure on \(\mathbb Z_P\) with
\(\|\nu\|_{\mathrm{TV}}\leq1\), and let \(0<\beta\leq1\).  Define
\(L_N^{\nu,\beta}\) by
\begin{align}\label{e:deforliczlow7}
 \widehat{L_N^{\nu,\beta}}(\xi)
 =
 \sum_{\alpha\in\Sigma_{\leq k}}
 \widehat\nu(\alpha)
 \widehat{\widetilde M}_{N,\beta}(\xi-\alpha)
 \eta\left(\frac{\xi-\alpha}{\varepsilon}\right).
\end{align}
Then
\begin{align}\label{e:orlicz-low}
 \left\|
 \sup_{\substack{N\geq1\\N\in2^{\mathbb N}}}
 |L_N^{\nu,\beta}f|
 \right\|_{\ell^{p_r}(\mathbb Z)}
 \ll
 rC_{r,k}^2
 \log\bigl(e+Q+r\log B\bigr)
 \|f\|_{\ell^{p_r}(\mathbb Z)}.
\end{align}
The implied constant is independent of \(\nu\) and \(\beta\).
\end{proposition}

\begin{proof}
With \(\varepsilon\) as in \eqref{e:IWparameters7}, set
\[
 \delta=\min\left\{\frac{\varepsilon}{16},\frac1{100P}\right\},
 \qquad
 J_*=\left\lceil\frac12\log_2\delta^{-1}\right\rceil,
 \qquad
 N_*=2^{J_*}.
\]
Thus \(1\leq N_*\sqrt\delta<2\).
For the scales \(N\leq N_*\), write
\(L_j=L_{2^j}^{\nu,\beta}\), \(0\leq j\leq J_*\).
For arbitrary \(|\sigma_j|\leq1\),
Proposition~\ref{p:IWinput7}, applied to the continuous multiplier
\[
 \sum_{j=1}^{J_*}\sigma_j
 \left(
 \widehat{\widetilde M}_{2^j,\beta}
 -
 \widehat{\widetilde M}_{2^{j-1},\beta}
 \right)\eta(\,\cdot\,/\varepsilon),
\]
and Lemma~\ref{l:continuousestimates7} give
\begin{align}\label{e:randomizedincrements7}
 \left\|
 \sum_{j=1}^{J_*}\sigma_j(L_j-L_{j-1})f
 \right\|_{\ell^{p_r}(\mathbb Z)}
 &\ll
 C_{r,k}(p_r-1)^{-1}\|f\|_{\ell^{p_r}(\mathbb Z)}
 \ll
 rC_{r,k}\|f\|_{\ell^{p_r}(\mathbb Z)}.
\end{align}
Here we also used that the inverse Fourier transform of
\(\eta(\,\cdot\,/\varepsilon)\) has uniformly bounded
\(L^1(\mathbb R)\)-norm. The same multiplier estimate, without taking
differences, gives
\[
 \|L_0f\|_{\ell^{p_r}(\mathbb Z)}
 \ll C_{r,k}\|f\|_{\ell^{p_r}(\mathbb Z)}.
\]
The Rademacher--Menshov inequality
\cite[(2.7)]{MirekSteinZorinKranich} gives
\begin{align*}
 &\left\|
 \max_{1\leq m\leq J_*}
 \left|\sum_{j=1}^m(L_j-L_{j-1})f\right|
 \right\|_{\ell^{p_r}(\mathbb Z)}\ll
 \sum_{\ell\leq\lceil\log_2J_*\rceil}
 \left\|
 \left(
 \sum_{I\in\mathcal D_\ell}
 \left|\sum_{j\in I}(L_j-L_{j-1})f\right|^2
 \right)^{1/2}
 \right\|_{\ell^{p_r}(\mathbb Z)}.
\end{align*}
Here \(\mathcal D_\ell\) consists of the nonempty intersections
of \(\{1,\ldots,J_*\}\) with the blocks
\((m2^\ell,(m+1)2^\ell]\), \(m\geq0\).
For fixed \(\ell\), let
\((\epsilon_I)_{I\in\mathcal D_\ell}\) be independent Rademacher
variables and set \(\sigma_j=\epsilon_I\) for the unique
\(I\in\mathcal D_\ell\) containing \(j\).  Khintchine's inequality,
Fubini's theorem, and \eqref{e:randomizedincrements7} give
\begin{align*}
 \left\|
 \left(
 \sum_{I\in\mathcal D_\ell}
 \left|\sum_{j\in I}(L_j-L_{j-1})f\right|^2
 \right)^{1/2}
 \right\|_{\ell^{p_r}(\mathbb Z)}
 &\ll
 \left(
 \mathbb E_\epsilon
 \left\|
 \sum_{j=1}^{J_*}\sigma_j(L_j-L_{j-1})f
 \right\|_{\ell^{p_r}(\mathbb Z)}^{p_r}
 \right)^{1/p_r}\\
 &\ll rC_{r,k}\|f\|_{\ell^{p_r}(\mathbb Z)}.
\end{align*}
 Summing over \(\ell\) and including
the estimate for \(L_0\) yields
\begin{align}\label{e:boundedscales7}
 \left\|
 \max_{0\leq j\leq J_*}|L_jf|
 \right\|_{\ell^{p_r}(\mathbb Z)}
 \ll
 rC_{r,k}\log(2+J_*)
 \|f\|_{\ell^{p_r}(\mathbb Z)}.
\end{align}

\medskip

For \(N\geq N_*\), introduce the narrower cutoff
\(\kappa(\theta)=\eta(\theta/\delta)\), and define
\[
 F(x,y)=
 \sum_{\alpha\in\Sigma_{\leq k}}
 \int_{\mathbb R}
 \eta\left(\frac{\theta}{\varepsilon}\right)
 \widehat f(\alpha+\theta)e(x\theta+y\alpha)\,d\theta.
\]
 Estimate  \eqref{e:IWsampling7}, followed by
\eqref{e:IWmultiplier7} with $m(\theta)=\eta(\theta/\varepsilon)$ and $\nu=\delta_0,$
gives
\begin{align}\label{e:Fadelicbound7}
 \|F\|_{L^{p_r}(\mathbb R\times\mathbb Z_P)}
 &\ll
 C_{r,k}
 \left\|
 \left(
 \sum_{\alpha\in\Sigma_{\leq k}}
 \eta\left(\frac{\xi-\alpha}{\varepsilon}\right)
 \widehat f(\xi)
 \right)^\vee
 \right\|_{\ell^{p_r}(\mathbb Z)}\ll
 C_{r,k}^2\|f\|_{\ell^{p_r}(\mathbb Z)}.
\end{align}
For a finite set \(\mathcal N\) of dyadic \(N\geq N_*\), define
\[
 G_N
 =
 \bigl(
   (\widetilde M_{N,\beta}*_{\mathbb R}\kappa^\vee)
   \otimes\nu
 \bigr)*_{\mathbb R\times\mathbb Z_P}F,
 \qquad N\in\mathcal N.
\]
Equivalently, $ \widehat G_N(\theta,\alpha)
 =
 \widehat{\widetilde M}_{N,\beta}(\theta)
 \kappa(\theta)\widehat\nu(\alpha)
 \widehat F(\theta,\alpha).$
Since convolution by \(\kappa^\vee\) and by \(\nu\) has uniformly
bounded \(L^{p_r}\)-operator norm, Lemma~\ref{l:continuousestimates7},
applied in the real variable, gives
\[
 \left\|
 \sup_{N\in\mathcal N}|G_N|
 \right\|_{L^{p_r}(\mathbb R\times\mathbb Z_P)}
 \ll
 r\|F\|_{L^{p_r}(\mathbb R\times\mathbb Z_P)}.
\]
The real Fourier support of every \(G_N\) is contained in
\([-\delta,\delta]\).
Apply \eqref{e:IWdiagonal7} to the
\(\ell^\infty(\mathcal N)\)-valued function
\((G_N)_{N\in\mathcal N}\).
Together with \eqref{e:Fadelicbound7}, this gives
\begin{align}\label{e:largescalesmain7}
 \left\|
 \sup_{N\in\mathcal N}
 |G_N(n,n\bmod P)|
 \right\|_{\ell^{p_r}(\mathbb Z)}
 &\ll
 r\|F\|_{L^{p_r}(\mathbb R\times\mathbb Z_P)}\ll
 rC_{r,k}^2\|f\|_{\ell^{p_r}(\mathbb Z)}.
\end{align}
Because \(\delta\leq\varepsilon/16\),
\(\eta(\theta/\varepsilon)=1\) on \(\operatorname{supp}\kappa\).
Consequently,
\[
 G_N(n,n\bmod P)
 =
 \left(
 \sum_{\alpha\in\Sigma_{\leq k}}
 \widehat\nu(\alpha)
 \widehat{\widetilde M}_{N,\beta}(\xi-\alpha)
 \kappa(\xi-\alpha)\widehat f(\xi)
 \right)^\vee(n).
\]
Thus \eqref{e:largescalesmain7} controls the part of
\(L_N^{\nu,\beta}f\) carrying the common cutoff \(\kappa\).  Letting
\(\mathcal N\) increase to all dyadic \(N\geq N_*\) gives the same
bound for the full supremum of these terms.

It remains to estimate the difference between the original cutoff and
\(\kappa\).  Since
\(\eta(\theta/\varepsilon)\kappa(\theta)=\kappa(\theta)\), this
difference has continuous multiplier
\[
 R_N(\theta)
 =
 \widehat{\widetilde M}_{N,\beta}(\theta)
 \eta\left(\frac{\theta}{\varepsilon}\right)
 \bigl(1-\kappa(\theta)\bigr).
\]
The measures \(\widetilde M_{N,\beta}\) have uniformly bounded mass,
and the inverse Fourier transforms of the cutoff functions have
uniformly bounded \(L^1(\mathbb R)\)-norm.  Hence
\[
 \|T_{R_N}\|_{L^1(\mathbb R)\to L^1(\mathbb R)}
 \ll1.
\]
Moreover, \(1-\kappa(\theta)\neq0\) implies
\(|\theta|\geq\delta/2\).  The Fourier estimate
\eqref{e:MNless} and Plancherel's theorem therefore give
\[
 \|T_{R_N}\|_{L^2(\mathbb R)\to L^2(\mathbb R)}
 \ll (N\sqrt{\delta})^{-1}.
\]
Since $ \frac1{p_r}
 =
 \left(1-\frac1r\right)
 +\frac1{2r},$
interpolation yields
\[
 \|T_{R_N}\|_{L^{p_r}(\mathbb R)\to L^{p_r}(\mathbb R)}
 \ll (N\sqrt{\delta})^{-1/r}.
\]
Applying the Ionescu--Wainger multiplier estimate in
Proposition~\ref{p:IWinput7}, we obtain
\[
 \left\|
 \left(
 \sum_{\alpha\in\Sigma_{\leq k}}
 \widehat\nu(\alpha)R_N(\xi-\alpha)\widehat f(\xi)
 \right)^\vee
 \right\|_{\ell^{p_r}(\mathbb Z)}
 \ll
 C_{r,k}(N\sqrt{\delta})^{-1/r}
 \|f\|_{\ell^{p_r}(\mathbb Z)}.
\]
Taking the supremum over \(N\) by summing these estimates over the
dyadic scales and using \(N_*\sqrt{\delta}\asymp1\), we get
\begin{align}\label{e:cutoffremainder7}
 &\left\|
 \sup_{\substack{N\geq N_*\\N\in2^{\mathbb N}}}
 \left|
 \left(
 \sum_{\alpha\in\Sigma_{\leq k}}
 \widehat\nu(\alpha)R_N(\xi-\alpha)\widehat f(\xi)
 \right)^\vee
 \right|
 \right\|_{\ell^{p_r}(\mathbb Z)}
\ll
 C_{r,k}
 \sum_{\ell\geq0}
 \bigl(2^\ell N_*\sqrt{\delta}\bigr)^{-1/r}
 \|f\|_{\ell^{p_r}(\mathbb Z)}
 \ll
 rC_{r,k}\|f\|_{\ell^{p_r}(\mathbb Z)}.
\end{align}
Finally, \(P\leq3^Q\) and \eqref{e:IWparameters7} imply
\[
 \log\delta^{-1}\ll Q+r\log B,
 \qquad
 \log(2+J_*)
 \ll\log\bigl(e+Q+r\log B\bigr).
\]
Combining \eqref{e:boundedscales7},
\eqref{e:largescalesmain7}, and
\eqref{e:cutoffremainder7} proves
\[
 \left\|
 \sup_{\substack{N\geq1\\N\in2^{\mathbb N}}}
 |L_N^{\nu,\beta}f|
 \right\|_{\ell^{p_r}(\mathbb Z)}
 \ll
 rC_{r,k}^2
 \log\bigl(e+Q+r\log B\bigr)
 \|f\|_{\ell^{p_r}(\mathbb Z)},
\]
as required.
\end{proof}

We now compare the rational-frequency operator in
Proposition~\ref{p:orlicz-low} with the circle-method model
\(V_N\) defined in \eqref{e:VNdef}.  Let \(N_0\) be the least dyadic
integer such that
\begin{equation}\label{e:orlicz-N0}
 N_0\geq
 \exp\left(C\bigl((\log Q)^2+r\log B\bigr)\right),
\end{equation}
where \(C\) is a sufficiently large absolute constant.  For every
dyadic \(N\geq N_0\), the definition of \(R=R(N)\) gives
\begin{align}\label{e:orlicz-parameter-range}
 R\leq N^{1/8},\qquad
 N_0\geq2B,\qquad
 \frac{\log(eQ)}{N_0\sqrt\varepsilon}\leq Q^{-1}.
\end{align}
In particular, the hypotheses of Lemma~\ref{l:aaplemma} hold throughout
this range.

\begin{proposition}\label{p:orlicz-hilo}
Using the assumptions and notation of
Propositions~\ref{p:IWinput7} and \ref{p:orlicz-low}, put
\begin{equation}\label{e:orlicz-D}
 D(r,k,Q)=rC_{r,k}^2\log(e+Q+r\log B).
\end{equation}
For dyadic \(N\geq N_0\), there are decompositions
\[
 \widetilde A_{N,\Lambda}=L_N+H_N,
 \qquad
 \widetilde A_{N,\mathbb N^2}=L_N^\square+H_N^\square,
\]
such that
\begin{align}
 \left\|\sup_{N\geq N_0}|L_Nf|\right\|_{\ell^{p_r}(\mathbb Z)}
 +
 \left\|\sup_{N\geq N_0}|L_N^\square f|\right\|_{\ell^{p_r}(\mathbb Z)}
 &\ll D(r,k,Q)\|f\|_{\ell^{p_r}(\mathbb Z)},
 \label{e:orlicz-low-comparison}\\
 \left\|\sup_{N\geq N_0}|H_Nf|\right\|_{\ell^2(\mathbb Z)}
 +
 \left\|\sup_{N\geq N_0}|H_N^\square f|\right\|_{\ell^2(\mathbb Z)}
 &\ll Q^{-1/4}\|f\|_{\ell^2(\mathbb Z)}.
 \label{e:orlicz-high-comparison}
\end{align}
\end{proposition}

\begin{proof}
All suprema and sums in \(N\) below are over dyadic scales.
All rational fractions are written in reduced form, and
\(s\) is chosen so that \(2^s\leq q<2^{s+1}\).
We write \(R=R(N)\).
By Lemma~\ref{l:aaplemma}, with \(A=2\),
\[
 \|\widetilde A_{N,\Lambda}-V_N\|
 _{\ell^2(\mathbb Z)\to\ell^2(\mathbb Z)}
 \ll Q^{-1/2}(1+\log N)^{-2}.
\]
Summing over dyadic \(N\geq N_0\), we obtain
\begin{align}\label{e:orlicz-circle-error}
 \left\|
 \sup_{N\geq N_0}
 |(\widetilde A_{N,\Lambda}-V_N)f|
 \right\|_{\ell^2(\mathbb Z)}
 &\leq
 \sum_{N\geq N_0}
 \|(\widetilde A_{N,\Lambda}-V_N)f\|_{\ell^2(\mathbb Z)}
\ll Q^{-1/2}\|f\|_{\ell^2(\mathbb Z)}.
\end{align}

Recall from Section~\ref{s:hilo} that \(\mathcal M(\chi)\) is the
set of scales at which \(\chi\) is exceptional, while \(q_\chi\)
and \(\beta_\chi\) denote its conductor and exceptional zero.
Since \(R(N)\geq Q^8\), all exceptional zeros with \(q_\chi\leq Q\)
arising for \(N\geq N_0\) lie in a common Page region for the fixed
conductor bound \(Q\). Thus \eqref{e:landaupagezero} permits at most
one such character. Denote it, when it exists, by \(\chi_*\).
By \eqref{e:IWcontainment7} and
\(\Sigma_{\leq k}\subset P^{-1}\mathbb Z/\mathbb Z\), its conductor
divides \(P\).

For a primitive real character \(\chi\) with \(q_\chi\mid P\), put
\[
 \nu_P^\chi
 =
 \frac1{\phi(P)}
 \sum_{u\in\mathbb A_P}\chi(u)\delta_{u^2}.
\]
This measure has total variation at most one.
For \(q\mid P\), reduction from \(\mathbb A_P\) to
\(\mathbb A_q\) has equal-size fibers. If \(q_\chi\mid q\), the
character descends to the character \(\chi_q\) modulo \(q\);
otherwise, it has mean zero on each fiber. Consequently,
\begin{align}\label{e:orlicz-character-coefficient}
 \widehat{\nu_P^\chi}(a/q)
 =
 \begin{cases}
  \mathfrak s_{\chi_q}(a,q),&q_\chi\mid q,\\
  0,&q_\chi\nmid q.
 \end{cases}
\end{align}
In particular,
\(\widehat{\nu_P^1}(a/q)=\mathfrak s(a,q)\).
Define
\begin{align}\label{e:orlicz-LNdef}
 L_N
 &=
 L_N^{\nu_P^1,1}
 -
 \mathbf1_{\mathcal M(\chi_*)}(N)
 L_N^{\nu_P^{\chi_*},\beta_{\chi_*}},
 \qquad
 H_N=\widetilde A_{N,\Lambda}-L_N,
\end{align}
with the second term omitted when \(\chi_*\) does not exist.
Proposition~\ref{p:orlicz-low}, applied to the two measures above,
gives
\begin{align}\label{e:orlicz-low-prime}
 \left\|\sup_{N\geq N_0}|L_Nf|\right\|_{\ell^{p_r}(\mathbb Z)}
 \ll
 D(r,k,Q)\|f\|_{\ell^{p_r}(\mathbb Z)}.
\end{align}

We now estimate \(H_N\) in \(\ell^2(\mathbb Z)\).
The exceptional character may vary with \(N\), but every character
of conductor at most \(Q\) has already been included in
\eqref{e:orlicz-LNdef}. Let \(E_N^{\mathrm{exc}}\) denote the
contribution to \(V_N\) from exceptional characters of conductor
greater than \(Q\):
\[
 \widehat{E_N^{\mathrm{exc}}f}(\xi)
 =
 -
 \sum_{\substack{\chi\\q_\chi>Q}}
 \mathbf1_{\mathcal M(\chi)}(N)
 \sum_{\substack{q\leq R(N)\\q_\chi\mid q}}
 \sum_{a\in\mathbb A_q}
 \mathfrak s_{\chi_q}(a,q)
 \widehat{\widetilde M}_{N,\beta_\chi}(\xi-a/q)
 \eta_s(\xi-a/q)\widehat f(\xi).
\]
For each fixed \(N\), at most one term in the outer sum is nonzero.
Moreover, every denominator in this expression exceeds \(Q\),
since \(q_\chi\mid q\).

To describe the remaining terms, put
\[
 \mathcal C_N(a,q;\theta)
 =
 \mathfrak s(a,q)\widehat{\widetilde M}_N(\theta)
 -
 \mathbf1_{\mathcal M(\chi_*)}(N)
 \mathbf1_{\{q_{\chi_*}\mid q\}}
 \mathfrak s_{\chi_{*,q}}(a,q)
 \widehat{\widetilde M}_{N,\beta_{\chi_*}}(\theta).
\]
Here \(\chi_{*,q}\) is induced by \(\chi_*\) when
\(q_{\chi_*}\mid q\); the entire second term is taken to be zero
if \(\chi_*\) does not exist or \(q_{\chi_*}\nmid q\).
Define
\begin{align*}
 \widehat{V_N^{>Q}f}(\xi)
 &=
 \sum_{Q<q\leq R(N)}\sum_{a\in\mathbb A_q}
 \mathcal C_N(a,q;\xi-a/q)
 \eta_s(\xi-a/q)\widehat f(\xi),\\
 \widehat{L_N^{>Q}f}(\xi)
 &=
 \sum_{\substack{a/q\in\Sigma_{\leq k}\\q>Q}}
 \mathcal C_N(a,q;\xi-a/q)
 \eta\left(\frac{\xi-a/q}{\varepsilon}\right)\widehat f(\xi),\\
 \widehat{D_Nf}(\xi)
 &=
 \sum_{q\leq Q}\sum_{a\in\mathbb A_q}
 \mathcal C_N(a,q;\xi-a/q)
 \left(
 \eta_s(\xi-a/q)
 -
 \eta\left(\frac{\xi-a/q}{\varepsilon}\right)
 \right)\widehat f(\xi).
\end{align*}
Thus \(V_N^{>Q}\) is the part of \(V_N-E_N^{\mathrm{exc}}\)
with \(q>Q\), and \(L_N^{>Q}\) is the corresponding denominator
restriction of \(L_N\). For \(q\leq Q\), every fraction belongs to
\(\Sigma_{\leq k}\) by \eqref{e:IWcontainment7}, and the arithmetic
coefficients agree by \eqref{e:orlicz-character-coefficient}.
Their difference is therefore precisely \(D_N\).
It follows from \eqref{e:VNdef} and \eqref{e:orlicz-LNdef} that
\begin{align}\label{e:orlicz-Hdecomposition}
 H_N
 =
 (\widetilde A_{N,\Lambda}-V_N)
 +E_N^{\mathrm{exc}}
 +V_N^{>Q}-L_N^{>Q}+D_N.
\end{align}

By Lemmas~\ref{l:inftyboundsaq} and \ref{l:Gcharbound},
\begin{equation}\label{e:orlicz-coefficient-bound}
 |\mathfrak s(a,q)|+|\mathfrak s_{\chi_q}(a,q)|
 \ll q^{-3/8}.
\end{equation}
Fix an exceptional character \(\chi\) with \(q_\chi>Q\).
For \(N\in\mathcal M(\chi)\), the operator
\(E_N^{\mathrm{exc}}\) consists of its contribution, whose
coefficients vanish unless \(q_\chi\mid q\).
Corollary~\ref{c:orlicz-prefix}, together with
\eqref{e:orlicz-coefficient-bound}, gives
\begin{align*}
 \left\|
 \sup_{\substack{N\geq N_0\\N\in\mathcal M(\chi)}}
 |E_N^{\mathrm{exc}}f|
 \right\|_{\ell^2(\mathbb Z)}
 &\ll
 \sum_{2^s\geq q_\chi/2}
 (1+s)^2\,2^{-3s/8}\|f\|_{\ell^2(\mathbb Z)}\ll
 q_\chi^{-3/8} \log(eq_\chi)^2\|f\|_{\ell^2(\mathbb Z)}.
\end{align*}
At each \(N\), at most one exceptional character occurs.
The Landau--Page argument preceding \eqref{e:sumoverchi} also
permits at most one in each dyadic conductor range. Hence
\begin{align}\label{e:orlicz-large-conductors}
 \left\|
 \sup_{N\geq N_0}|E_N^{\mathrm{exc}}f|
 \right\|_{\ell^2(\mathbb Z)}^2
 &\leq
 \sum_{\substack{\chi\\q_\chi>Q}}
 \left\|
 \sup_{\substack{N\geq N_0\\N\in\mathcal M(\chi)}}
 |E_N^{\mathrm{exc}}f|
 \right\|_{\ell^2(\mathbb Z)}^2
 \notag\\
 &\ll
 \sum_{2^{j+1}>Q}
 2^{-3j/4}(1+j)^4\|f\|_{\ell^2(\mathbb Z)}^2
 \notag\\
 &\ll
 Q^{-3/4}\log(eQ)^4\|f\|_{\ell^2(\mathbb Z)}^2
 \ll Q^{-1/2}\|f\|_{\ell^2(\mathbb Z)}^2.
\end{align}

Both \(V_N^{>Q}\) and \(L_N^{>Q}\) contain ordinary terms and,
possibly, \(\chi_*\)-terms with \(q>Q\).
For \(V_N^{>Q}\), Corollary~\ref{c:orlicz-prefix} allows the
\(N\)-dependent upper bound \(q\leq R(N)\).
For \(L_N^{>Q}\), every denominator is at most \(B\), and the
choice of \(\varepsilon\) in \eqref{e:IWparameters7} gives
\[
 \operatorname{supp}\eta(\,\cdot\,/\varepsilon)
 \subseteq\{\theta:\eta_s(\theta)=1\},
 \qquad q\sim2^s,\quad q\leq B.
\]
Since \(N_0\geq2B\) by \eqref{e:orlicz-parameter-range}, the
filtered form of the same corollary applies.
For the ordinary terms use \(\beta=1\); for the exceptional terms
use \(\beta=\beta_{\chi_*}\) and restrict \(N\) to
\(\mathcal M(\chi_*)\).
The coefficient bound \eqref{e:orlicz-coefficient-bound} and
\[
 \sum_{2^s>Q/2}(1+s)^2\,2^{-3s/8}
 \ll Q^{-3/8}\log(eQ)^2
 \ll Q^{-1/4}
\]
therefore give
\begin{align}
 \left\|
 \sup_{N\geq N_0}|V_N^{>Q}f|
 \right\|_{\ell^2(\mathbb Z)}
 &\ll Q^{-1/4}\|f\|_{\ell^2(\mathbb Z)},
 \label{e:orlicz-large-denominators}\\
 \left\|
 \sup_{N\geq N_0}|L_N^{>Q}f|
 \right\|_{\ell^2(\mathbb Z)}
 &\ll Q^{-1/4}\|f\|_{\ell^2(\mathbb Z)}.
 \label{e:orlicz-extra-frequencies}
\end{align}

It remains to estimate the cutoff difference \(D_N\).
By \eqref{e:MNless} and
\eqref{e:orlicz-coefficient-bound}, uniformly for \(q\leq Q\)
and \(0<|\theta|\leq1/2\),
\[
 |\mathcal C_N(a,q;\theta)|
 \ll \frac1{N\sqrt{|\theta|}}.
\]
Consequently, summing the dyadic geometric series gives
\[
 \sum_{N\geq N_0}|\mathcal C_N(a,q;\theta)|^2
 \ll
 \frac1{|\theta|}\sum_{N\geq N_0}\frac1{N^2}
 \ll\frac1{N_0^2|\theta|}.
\]
The cutoff difference vanishes for
\(|\theta|<\varepsilon/2\), and its support is contained in that
of \(\eta_s\).
Let \(D_{N,s}\) denote the part of \(D_N\) with \(q\sim2^s\).
The rational-frequency supports are disjoint at this scale.
The pointwise bound of the supremum by the square function and
Plancherel's theorem therefore give
\begin{align*}
 &\left\|
 \sup_{N\geq N_0}|D_{N,s}f|
 \right\|_{\ell^2(\mathbb Z)}^2
 \leq
 \sum_{N\geq N_0}\|D_{N,s}f\|_{\ell^2(\mathbb Z)}^2\\
 &=
 \int_{\mathbb T}
 \sum_{\substack{q\sim2^s\\q\leq Q}}\sum_{a\in\mathbb A_q}
 \left|
 \eta_s(\xi-a/q)
 -
 \eta\left(\frac{\xi-a/q}{\varepsilon}\right)
 \right|^2
 \sum_{N\geq N_0}
 |\mathcal C_N(a,q;\xi-a/q)|^2
 |\widehat f(\xi)|^2\,d\xi\\
 &\ll
 \frac1{N_0^2\varepsilon}\|f\|_{\ell^2(\mathbb Z)}^2.
\end{align*}
Taking square roots and summing over the
\(O(\log(eQ))\) denominator scales, we obtain
\begin{equation}\label{e:orlicz-cutoff-difference}
 \left\|
 \sup_{N\geq N_0}|D_Nf|
 \right\|_{\ell^2(\mathbb Z)}
 \ll
 \frac{\log(eQ)}{N_0\sqrt\varepsilon}
 \|f\|_{\ell^2(\mathbb Z)}
 \ll Q^{-1}\|f\|_{\ell^2(\mathbb Z)},
\end{equation}
where the last inequality is
\eqref{e:orlicz-parameter-range}.
This estimate includes the possible \(\chi_*\)-contribution
through \(\mathcal C_N\).
Combining \eqref{e:orlicz-Hdecomposition} with
\eqref{e:orlicz-circle-error},
\eqref{e:orlicz-large-conductors},
\eqref{e:orlicz-large-denominators},
\eqref{e:orlicz-extra-frequencies}, and
\eqref{e:orlicz-cutoff-difference}, we conclude that
\[
 \left\|
 \sup_{N\geq N_0}|H_Nf|
 \right\|_{\ell^2(\mathbb Z)}
 \ll Q^{-1/4}\|f\|_{\ell^2(\mathbb Z)}.
\]

\medskip

For the square averages, define
\[
 \nu_P^\square
 =
 \frac1P\sum_{u\in\mathbb Z_P}\delta_{u^2},
 \qquad
 L_N^\square=L_N^{\nu_P^\square,1},
 \qquad
 H_N^\square=\widetilde A_{N,\mathbb N^2}-L_N^\square.
\]
The measure \(\nu_P^\square\) has mass one, so
Proposition~\ref{p:orlicz-low} gives
\[
 \left\|
 \sup_{N\geq N_0}|L_N^\square f|
 \right\|_{\ell^{p_r}(\mathbb Z)}
 \ll D(r,k,Q)\|f\|_{\ell^{p_r}(\mathbb Z)}.
\]
For a reduced fraction \(a/q\), with \(q\mid P\), reduction
modulo \(q\) gives
\[
 \widehat{\nu_P^\square}(a/q)
 =
 S(a,q)
 :=
 \frac1q\sum_{u\in\mathbb Z_q}e(-au^2/q).
\]
By the classical quadratic Gauss-sum bound
\cite[Section~3]{iwaniec}, we have
\(|S(a,q)|\leq\sqrt{2/q}\).
Define the square circle-method model by
\[
 \widehat V_N^\square(\xi)
 =
 \sum_{q\leq R}\sum_{a\in\mathbb A_q}
 S(a,q)\widehat{\widetilde M}_N(\xi-a/q)
 \eta_s(\xi-a/q).
\]
If \(\theta=\xi-a/q\), with \(q\leq R\) and
\(|\theta|\leq2R/(qN^2)\), summation over residue classes modulo
\(q\), followed by partial summation, gives
\begin{align*}
 \widehat{\widetilde A}_{N,\mathbb N^2}(\xi)
 &=
 S(a,q)\widehat{\widetilde M}_N(\theta)
 +O\left(\frac qN(1+N^2|\theta|)\right)=
 S(a,q)\widehat{\widetilde M}_N(\theta)+O(R/N).
\end{align*}

Outside these arcs, apply \cite[Theorem~3.1]{LHF} to the initial
square averages at \(M=N\) and \(M=N/2\).
Dirichlet's theorem supplies a reduced fraction \(a/q\) with
\(q\leq4M\) and \(|2\xi-a/q|\leq(4Mq)^{-1}\).
If \(d\) is the reduced denominator of \(a/(2q)\), then
\(q\leq d\leq2q\).
The Gauss factor and the continuous integral bound the main
term in that theorem by
\[
 \ll
 d^{-1/2}
 \min\left\{1,
 \bigl(M^2|\xi-a/(2q)|\bigr)^{-1/2}\right\}.
\]
If \(d>R\), this is \(O(R^{-1/2})\).
If \(d\leq R\), being outside the arcs implies
\(|\xi-a/(2q)|>2R/(dN^2)\), which gives the same bound since
\(M\geq N/2\).
The approximation error is
\(O(M^{-1}\sqrt q)=O(N^{-1/2})\).
Using
\(\widetilde A_{N,\mathbb N^2}
 =A_{N,\mathbb N^2}-\frac12A_{N/2,\mathbb N^2}\),
we therefore obtain
\[
 |\widehat{\widetilde A}_{N,\mathbb N^2}(\xi)|
 \ll R^{-1/2}+N^{-1/2}
\]
outside the arcs.

By \eqref{e:orlicz-parameter-range}, \(R\leq N^{1/8}\).
Hence the arcs above are disjoint and their corresponding
cutoffs \(\eta_s\) equal one.
Every term of \(V_N^\square\) outside its corresponding arc
satisfies, by \eqref{e:MNless},
\[
 |S(a,q)\widehat{\widetilde M}_N(\theta)\eta_s(\theta)|
 \ll
 q^{-1/2}\left(\frac1{N\sqrt{|\theta|}}+\frac1N\right)
 \ll R^{-1/2}.
\]
At each \(\xi\), at most one rational frequency contributes at
each dyadic denominator scale, so there are
\(O(\log(eR))\) such terms.
Combining the major-arc approximation with the bounds outside
the arcs and applying Plancherel's theorem gives
\begin{align}\label{e:orlicz-square-error}
 \|\widetilde A_{N,\mathbb N^2}-V_N^\square\|
 _{\ell^2(\mathbb Z)\to\ell^2(\mathbb Z)}
 &\ll
 R/N+R^{-1/2}\log(eR)+N^{-1/2}
 \notag\\
 &\ll R^{-1/3}
 \ll Q^{-1/2}(1+\log N)^{-2}.
\end{align}
The last two inequalities follow from \eqref{e:defR} and
\eqref{e:orlicz-parameter-range}.
Summing over dyadic \(N\geq N_0\), as in
\eqref{e:orlicz-circle-error}, yields
\[
 \left\|
 \sup_{N\geq N_0}
 |(\widetilde A_{N,\mathbb N^2}-V_N^\square)f|
 \right\|_{\ell^2(\mathbb Z)}
 \ll Q^{-1/2}\|f\|_{\ell^2(\mathbb Z)}.
\]

Finally, compare \(V_N^\square\) with \(L_N^\square\).
For the terms with \(q>Q\), Corollary~\ref{c:orlicz-prefix}
and \(|S(a,q)|\leq\sqrt{2/q}\) give, for each of the two sums,
maximal \(\ell^2(\mathbb Z)\)-operator norm at most
\[
 \ll
 \sum_{2^s>Q/2}(1+s)^2\,2^{-s/2}
 \ll Q^{-1/2}\log(eQ)^2
 \ll Q^{-1/4}.
\]
For the circle-method sum, allow the upper limit
\(q\leq R(N)\) to depend on \(N\).
For the Ionescu--Wainger sum, use the filtered form: its
denominators are at most \(B\), its cutoff is supported where
\(\eta_s=1\), and \(N_0\geq2B\).

For \(q\leq Q\), the arithmetic coefficients agree, and their
difference has multiplier
\[
 \sum_{q\leq Q}\sum_{a\in\mathbb A_q}
 S(a,q)\widehat{\widetilde M}_N(\xi-a/q)
 \left(
 \eta_s(\xi-a/q)
 -
 \eta\left(\frac{\xi-a/q}{\varepsilon}\right)
 \right).
\]
By \eqref{e:MNless}, the coefficients
\(S(a,q)\widehat{\widetilde M}_N(\theta)\) satisfy the same
dyadic square-sum bound used above for \(\mathcal C_N\).
The cutoff vanishes for \(|\theta|<\varepsilon/2\), and the
supports are disjoint at each denominator scale.
The preceding Plancherel calculation therefore bounds the
maximal \(\ell^2(\mathbb Z)\)-operator norm of this difference by
\[
 \ll
 \frac{\log(eQ)}{N_0\sqrt\varepsilon}
 \ll Q^{-1},
\]
using \eqref{e:orlicz-parameter-range}.

Combining these three bounds for
\[
 H_N^\square
 =
 (\widetilde A_{N,\mathbb N^2}-V_N^\square)
 +(V_N^\square-L_N^\square)
\]
gives
\[
 \left\|
 \sup_{N\geq N_0}|H_N^\square f|
 \right\|_{\ell^2(\mathbb Z)}
 \ll Q^{-1/4}\|f\|_{\ell^2(\mathbb Z)}.
\]
Together with \eqref{e:orlicz-low-prime} and the low estimate
for \(L_N^\square\), this proves
\eqref{e:orlicz-low-comparison} and
\eqref{e:orlicz-high-comparison}.
\end{proof}

We now optimize the parameters in
Proposition~\ref{p:orlicz-hilo} and deduce the endpoint estimates.

\begin{proof}[Proof of Theorem~\ref{t:orlicz-maximal}]
Let \(\widetilde A^*\) denote either of the dyadic annular maximal
operators
\[
 f\longmapsto
 \sup_{N\in2^{\mathbb N}}|\widetilde A_{N,\Lambda}f|,
 \qquad
 f\longmapsto
 \sup_{N\in2^{\mathbb N}}|\widetilde A_{N,\mathbb N^2}f|.
\]
For either choice of \(\widetilde A^*\), write
\(\widetilde A_N=L_N+H_N\) for the corresponding decomposition in
Proposition~\ref{p:orlicz-hilo}.  We first work with dyadic scales.
Chebyshev's estimate \(\sum_{n\leq N}\Lambda(n)\ll N\) gives
\(\|\widetilde A_N\mathbf1_F\|_{\ell^1(\mathbb Z)}\ll|F|\)
for prime squares; the same bound is immediate for squares.
There are \(O(1+\log N_0)\) scales below \(N_0\).

For \(N\geq N_0\), if
\(|\widetilde A_N\mathbf1_F|>t^{-1}\), then either
\(|L_N\mathbf1_F|>(2t)^{-1}\) or
\(|H_N\mathbf1_F|>(2t)^{-1}\).
The union bound and Chebyshev's inequality, with exponents
\(1\), \(p_r\), and \(2\), followed by
\eqref{e:orlicz-low-comparison} and
\eqref{e:orlicz-high-comparison}, therefore give, for \(t\geq1\),
\begin{align}
 |\{\widetilde A^*\mathbf1_F>t^{-1}\}|
 &\leq
 t\sum_{N<N_0}
 \|\widetilde A_N\mathbf1_F\|_{\ell^1(\mathbb Z)}
 +
 (2t)^{p_r}
 \left\|\sup_{N\geq N_0}|L_N\mathbf1_F|
 \right\|_{\ell^{p_r}(\mathbb Z)}^{p_r}
+
 4t^2
 \left\|\sup_{N\geq N_0}|H_N\mathbf1_F|
 \right\|_{\ell^2(\mathbb Z)}^2
 \notag\\
 &\ll
 |F|\left(
 t(1+\log N_0)+t^{p_r}D(r,k,Q)^{p_r}+t^2Q^{-1/2}
 \right)
 \notag\\
 &=t|F|\left(
 1+\log N_0+
 t^{1/(2r-1)}D(r,k,Q)^{p_r}
 +tQ^{-1/2}
 \right).
 \label{e:orlicz-optimize}
\end{align}
Here we used
\(\|\mathbf1_F\|_{\ell^p(\mathbb Z)}^p=|F|\)
and \(p_r-1=(2r-1)^{-1}\).
For sufficiently large \(t\), put
\[
 v=\log(e+t),\qquad
 Q=2^{\lceil2\log_2(e+t)\rceil},\qquad
 r=\lceil v\rceil,\qquad
 k=\left\lfloor
 \sqrt{\frac{v}{\log(e+v)}}
 \right\rfloor.
\]
Then \(k\geq4\), \(Q\geq2^k\), $t^{1/(2r-1)}\ll1,$ and $tQ^{-1/2}\ll1.$
Moreover, the definitions of \(B\) and \(C_{r,k}\) in
\eqref{e:IWparameters7} give
\[
 \log B\ll k^2\exp\left(\frac{Cv}{k}\right),
 \qquad
 \log C_{r,k}\ll k\log(e+v).
\]
Together with the definitions of \(N_0\) and \(D(r,k,Q)\) in
\eqref{e:orlicz-N0} and \eqref{e:orlicz-D}, these bounds give
\[
 \log(1+\log N_0)+\log D(r,k,Q)
 \ll
 k\log(e+v)+\frac{v}{k}
 \ll\sqrt{v\log(e+v)}.
\]
Thus the expression in parentheses in
\eqref{e:orlicz-optimize} is at most $\exp\left(C\sqrt{v\log(e+v)}\right).$
For bounded \(t\), the same conclusion follows by taking fixed
admissible values of \(r,k,Q\).  By \eqref{e:orlicz-gauge-comparison},
\begin{align}\label{e:orlicz-restricted-annular}
 |\{\widetilde A^*\mathbf1_F>t^{-1}\}|
 \ll |F|\Psi_C(t).
\end{align}
By \eqref{e:annularreconstruction} and its square analogue,
\[
 \sup_{\substack{N\in2^{\mathbb N}\\N\geq2}}|A_Nf|
 \leq
 2\widetilde A^*f+\frac12|f(\,\cdot-1)|.
\]
The last term is needed only for squares and satisfies the
required restricted estimate, since \(\Psi_C(t)\geq t\geq1\).
If \(N\leq N'<2N\) and \(N'\) is dyadic, positivity gives
\(|A_Nf|\leq2A_{N'}|f|\).
Together with \eqref{e:orlicz-restricted-annular}, these bounds
prove \eqref{e:orlicz-restricted} for the square and
von Mangoldt averages.

It remains to pass to the unweighted prime-square averages.  For
\(f\geq0\) and \(2^m\leq N<2^{m+1}\), the prime number theorem gives
\begin{align*}
 A_{N,\mathbb P^2}f
 &\ll
 \frac{m+1}{2^m}
 \sum_{j=0}^m\frac1{j+1}
 \sum_{\substack{2^j\leq p<2^{j+1}\\p\in\mathbb P}}
 (\log p)f(\,\cdot-p^2)\ll
 A_\Lambda^*f\,
 (m+1)\sum_{j=0}^m\frac{2^{j-m}}{j+1}
 \ll A_\Lambda^*f.
\end{align*}
Thus \eqref{e:orlicz-restricted} holds for all three maximal
operators.

We pass to general functions by the standard dyadic layer-cake
decomposition; see \cite[Chapter~1, \S2]{BB}.
Let \(C_\infty\geq1\) bound the masses of the three families of
kernels.  For finitely supported \(f\), the part with
\(|f|\leq\lambda/(2C_\infty)\) contributes at most \(\lambda/2\).
Put
\[
 E_j=
 \left\{
 \frac{2^j\lambda}{2C_\infty}<|f|
 \leq\frac{2^{j+1}\lambda}{2C_\infty}
 \right\},
 \qquad j\geq0.
\]
Since \(\sum_{j\geq0}(j+1)^{-2}<2\), positivity and
\eqref{e:orlicz-restricted} give
\begin{align*}
 |\{A^*f>\lambda\}|
 &\leq
 \sum_{j\geq0}
 \left|\left\{
 A^*\mathbf1_{E_j}>
 \frac{2^{-j}}{4(j+1)^2}
 \right\}\right|\ll
 \sum_{j\geq0}|E_j|
 \Psi_{C_1}\bigl(4\cdot2^j(j+1)^2\bigr)\\
 &\ll
 \sum_{j\geq0}|E_j|
 \Psi_{C_2}\left(\frac{2^j}{2C_\infty}\right)\leq
 \sum_{\substack{x\in\mathbb Z\\
 |f(x)|>\lambda/(2C_\infty)}}
 \Psi_{C_2}\left(\frac{|f(x)|}{\lambda}\right).
\end{align*}
Here \(C_1\) is the constant in \eqref{e:orlicz-restricted}.
The third inequality follows from
\eqref{e:orlicz-gauge-comparison}, with \(C_2>C_1\)
sufficiently large, since
\(\log(j+1)=o(\sqrt{j\log(e+j)})\).
The last inequality uses monotonicity of \(\Psi_{C_2}\).
Renaming \(C_2\) as \(C\) proves
\eqref{e:orlicz-truncated} with \(C_0=2C_\infty\);
enlarging the sum gives \eqref{e:orlicz-modular}.
\end{proof}

Let \((X,\mu,T)\) be an invertible measure-preserving system.
For any of the three families in Theorem~\ref{t:orlicz-maximal},
use the ergodic notation of Section~\ref{ss:transference}.

\begin{corollary}\label{c:orlicz-ergodic}
For \(1<p\leq2\),
\begin{align}\label{e:orlicz-ergodic-lp}
 \|A_T^*f\|_{L^p(X)}
 \leq
 \exp\left(
 \frac{C}{p-1}\log\frac{e}{p-1}
 \right)\|f\|_{L^p(X)}.
\end{align}
Moreover,
\begin{align}\label{e:orlicz-ergodic}
 \mu\{A_T^*f>\lambda\}
 \ll
 \int_X\Psi_C\left(\frac{|f|}{\lambda}\right)\,d\mu,
 \qquad \lambda>0.
\end{align}
If $ \int_X\Psi_C(|f|/a)\,d\mu<\infty$
for some \(a>0\), then the corresponding averages converge almost
everywhere and in \(L^1(X)\).
\end{corollary}

\begin{proof}
Calder\'on's transference principle
\cite[Theorem~1]{Calderon1968} gives
\eqref{e:orlicz-ergodic-lp} from \eqref{e:orlicz-p}. 
The modular transference argument in
Section~\ref{ss:transference}, applied to
\eqref{e:orlicz-modular}, gives \eqref{e:orlicz-ergodic}.

For bounded \(g\), \cite[Corollary~2(ii)]{MehlhopSlomian},
applied to the polynomial \(n^2\), gives almost everywhere
convergence of the square averages and the logarithmically
weighted prime-square averages. The latter are normalized by
\(\sum_{p\leq N}\log p\); the prime number theorem permits
normalization by \(N\), and Abel summation gives convergence
of the unweighted prime averages. The additional prime powers
in \(A_{N,\Lambda,T}g\) contribute at most
\[
 \frac{\|g\|_{L^\infty(X)}}{N}
 \sum_{\substack{p^j\leq N\\j\geq2}}\log p
 \ll
 \|g\|_{L^\infty(X)}N^{-1/2}\log^2(2N)
 =o(1).
\]
Thus all three families converge almost everywhere for bounded
functions.

Suppose that \(\int_X\Psi_C(|f|/a)\,d\mu<\infty\) for some
\(a>0\), and set \(g_m=f\mathbf1_{\{|f|\leq m\}}\).
The doubling property from \eqref{e:orlicz-gauge-comparison}
gives
\(\Psi_C(|f|/\lambda)\ll_{a,\lambda,C}\Psi_C(|f|/a)\)
for every \(\lambda>0\).
Since the averages of each \(g_m\) converge almost everywhere,
\eqref{e:orlicz-ergodic} yields
\begin{align*}
 &\mu\left\{
 \limsup_{N,M\to\infty}
 |A_{N,T}f-A_{M,T}f|>2\lambda
 \right\}\ll
 \mu\{A_T^*(f-g_m)>\lambda\}
 \ll
 \int_X\Psi_C\left(\frac{|f-g_m|}{\lambda}\right)\,d\mu.
\end{align*}
Using dominated convergence, pointwise convergence  $f-g_m\to0 $, and  \(|f-g_m|\leq|f|\), we conclude that 
$$\mu\left\{
 \limsup_{N,M\to\infty}
 |A_{N,T}f-A_{M,T}f|>2\lambda
 \right\} =0.$$
Taking \(\lambda\) through the positive rationals proves
almost everywhere convergence.

Convergence in \(L^1(X)\) follows first for bounded functions
by dominated convergence, and then for all \(L^1(X)\) functions
by density and the uniform \(L^1(X)\)-operator bound.
\end{proof}

\begin{proof}[Proof of Theorem~\ref{t:pet}\textup{(ii)}]
For the unweighted prime-square averages,
\[A_T^*f=\sup_{m\geq1}
\left|m^{-1}\sum_{k=1}^m f\circ T^{p_k^2}\right|,\]
so \eqref{e:orlicz-ergodic-lp} gives the claim.
\end{proof}

\section{Persistent universal \(L^1\)-badness of prime squares}
\label{s:l1}

We now study the endpoint \(L^1\) behaviour of averages along the prime
squares.  We use LaVictoire's criterion \cite{L09}, which reduces the
ergodic conclusion to distributional properties of subsets of finite
cyclic groups \(\mathbb Z/q\mathbb Z\).  For each modulus \(q\), we take
\(\Lambda_q\) to be the set of squares of units modulo \(q\).  The Chinese remainder theorem gives its multiplicative structure.
We quote Granville--Kurlberg \cite{GK08} for the Poisson-spacing
conclusion, but verify the two spacing conditions needed for our purpose
directly. Dirichlet's theorem shows that every
element of \(\Lambda_q\) is represented by a prime square, and the prime
number theorem in arithmetic progressions gives the required limiting
estimate.

Throughout this section, $p_1<p_2<\cdots$ denotes the sequence of primes and $n_k=p_k^2.$
We prove that $(n_k)$ is persistently universally $L^1$-bad, in the standard
sense that for every nonatomic ergodic measure-preserving system and every
infinite $\mathcal S\subset\mathbb N$, there is an $f\in L^1$ for which the
averages 
$$N^{-1}\sum_{k\leq N} T^{n_k}f, \qquad N\in\mathcal S, $$ 
diverge on a
set of positive measure.  We use the following criterion of LaVictoire.
For $\Lambda\subset\mathbb Z/q\mathbb Z$, write
\[
 \mathbb P_q(\Lambda)=\frac{|\Lambda|}{q},\qquad
 s_q(\Lambda)=\mathbb P_q(\Lambda)^{-1},\qquad
 \Lambda^\gamma
 =\Lambda+(0,\gamma s_q(\Lambda))\subset\mathbb Z/q\mathbb Z,
\]
where intervals and distances are understood cyclically.  When
\(\Lambda=\Lambda_q\), we write simply \(s_q=s_q(\Lambda_q)\).

The following form of LaVictoire's criterion combines
\cite[Theorem~3.1]{L09} with the Conze transference principle and
Sawyer's Banach principle, as explained in \cite[Section~1]{L09}.

\begin{theorem}\label{t:l1badblackbox}
Let $(n_k)$ be an increasing sequence of natural numbers.  Suppose that, for
every $q>1$, there is a set $\Lambda_q\subset\{n_k\bmod q:k\geq1\}$
which is multiplicative under the Chinese remainder theorem: if $(q,r)=1$,
then $\Lambda_{qr}$ corresponds to $\Lambda_q\times\Lambda_r$. Suppose also that there are two sequences \((u_j)\) and \((q_j)\) of
odd integers whose combined collection is pairwise coprime and such that
\begin{align}
 &\mathbb P_{q_j}(\Lambda_{q_j})\longrightarrow0,
 \qquad
 \inf_j\mathbb P_{u_j}(\Lambda_{u_j})>0,
 \label{e:l1-density}\\
 &\gamma-\liminf_{j\to\infty}
       \mathbb P_{q_j}(\Lambda_{q_j}^{\gamma})=o(\gamma)
       \quad(\gamma\downarrow0),\label{e:l1-thickening}
\end{align}
and that there is an $\alpha>0$ such that, for every sufficiently small
$\gamma>0$,
\begin{align}
 \liminf_{j\to\infty}\frac1{|\Lambda_{q_j}|^2}
 \big|\{(u,v)\in\Lambda_{q_j}^2:
   \operatorname{dist}_{q_j}(u-v,\Lambda_{q_j})>
       \gamma s_{q_j}\}\big|>5\alpha.             \label{e:l1-difference}
\end{align}
Finally, suppose that some $\beta>0$ has the following property. If \(Q\) is any product of distinct terms from the two sequences
\((u_j)\) and \((q_j)\), then for every \(a\in\Lambda_Q\),
\begin{align}
 \liminf_{N\to\infty}\frac1N
 \big|\{1\leq k\leq N:n_k\equiv a\pmod Q\}\big|
       >\frac{\beta}{|\Lambda_Q|}.                 \label{e:l1-equidistribution}
\end{align}
Then $(n_k)$ is persistently universally $L^1$-bad.
\end{theorem}


For comparison, we record the following special case of the theorem of
Granville and Kurlberg \cite[Theorem 1]{GK08}.

\begin{theorem}[Granville--Kurlberg \cite{GK08}]\label{t:grankul}
For each odd prime $p$, let $\Omega_p\subset\mathbb Z/p\mathbb Z$, put
$r_p=|\Omega_p|/p$ and $s_p=r_p^{-1}$, and define $\Omega_q$ from the
$\Omega_p$, for squarefree $q$, by the Chinese remainder theorem.  Suppose
that, for every fixed $k\geq1$ and all distinct
$0,h_1,\ldots,h_{k-1}\pmod p$,
\begin{align}
 N_k(\mathbf h,\Omega_p)
 &:=\big|\{t\in \mathbb{Z}/p\mathbb{Z}:t,t+h_1,\ldots,t+h_{k-1}\in\Omega_p\}\big|=r_p^kp\big(1+O_k((1-r_p)p^{-c})\big)             \label{e:GK-hyp}
\end{align}
for some $c>0$, and suppose that $s_p=p^{o(1)}$.  Then the normalized spacings
of $\Omega_q$ are Poisson distributed as $s_q=q/|\Omega_q|\to\infty$.
In particular, if the elements of $\Omega_q$ are listed cyclically and
$F_q(\theta)$ is the proportion of normalized gaps exceeding $\theta$, then
$F_q(\theta)\to e^{-\theta}$ for every $\theta>0$.
\end{theorem}

We verify the correlation hypothesis for the nonzero quadratic
residues.

\begin{lemma}\label{l:quadratic-GK}
Let $p$ be an odd prime and let
\[
 \Lambda_p=\{x^2\pmod p:x\not\equiv0\pmod p\}.
\]
For every fixed $k\geq1$ and distinct
$0,h_1,\ldots,h_{k-1}\pmod p$,
\[
 N_k(\mathbf h,\Lambda_p)=\frac{p}{2^k}+O_k(p^{1/2})
 =\left(\frac{p-1}{2p}\right)^kp\big(1+O_k(p^{-1/2})\big).
\]
Consequently \eqref{e:GK-hyp} holds, for example with $c=1/3$.
\end{lemma}

\begin{proof}
Let $\chi_p$ be the quadratic character modulo $p$, extended by
$\chi_p(0)=0$, and put $h_0=0$.  For every $x\pmod p$,
\[
 \mathbf1_{\Lambda_p}(x)
   =\frac12\big(1+\chi_p(x)-\mathbf1_{\{0\}}(x)\big).
\]
Therefore
\begin{align*}
 N_k(\mathbf h,\Lambda_p)
 &=2^{-k}\sum_{t\bmod p}\prod_{i=0}^{k-1}
     \big(1+\chi_p(t+h_i)-\mathbf1_{\{0\}}(t+h_i)\big).
\end{align*}
The term obtained by choosing $1$ from every factor is $p/2^k$.
Every term containing at least one factor
$\mathbf1_{\{0\}}(t+h_i)$ is $O_k(1)$: one such factor fixes
$t=-h_i$, and there are at most $3^k$ terms in the expansion.
The remaining nonconstant terms have the form
\[
 \sum_{t\bmod p}\chi_p\!\left(\prod_{i\in I}(t+h_i)\right),
 \qquad \varnothing\ne I\subset\{0,\ldots,k-1\}.
\]
Because the $h_i$ are distinct, the polynomial inside the character is
squarefree.  For $|I|=1$ the sum is $0$, and for $|I|\geq2$ Weil's bound for
multiplicative character sums gives absolute value at most
$(|I|-1)p^{1/2}$, see \cite[Section 11.1]{iwaniec}.  Summing over the $2^k-1$ nonempty sets $I$ proves
\[
 N_k(\mathbf h,\Lambda_p)=2^{-k}p+O_k(p^{1/2}).
\]
Since
\[
 \left(\frac{p-1}{2p}\right)^kp
   =2^{-k}p+O_k(1),
\]
the second formula follows.  Finally $1-r_p=(p+1)/(2p)\sim 1$, so the
relative error $O_k(p^{-1/2})$ is
$O_k((1-r_p)p^{-1/3})$.
\end{proof}

\begin{proof}[Proof of Theorem \ref{t:pet}\textup{(i)}]
For every \(q\geq2\), define
\[
 \Lambda_q=\{b^2\pmod q:(b,q)=1\}.
\]
If \((q,r)=1\), the Chinese remainder theorem identifies
\(\Lambda_{qr}\) with \(\Lambda_q\times\Lambda_r\).  Moreover,
\[
 \Lambda_q\subseteq\{p_k^2\pmod q:k\geq1\}.
\]
Indeed, if \(a\in\Lambda_q\), choose a unit \(b\pmod q\) such that
\(b^2\equiv a\pmod q\).  Dirichlet's theorem gives a prime
\(p\equiv b\pmod q\), and hence \(p^2\equiv a\pmod q\).

Choose distinct odd primes \(u_j\) and pairwise disjoint sets
\(B_j\) of odd primes, with \(|B_j|=j\), such that
\(u_i\notin B_j\) for all \(i,j\geq1\) and
\[
 p>3^j\quad(p\in B_j),\qquad
 \sum_{p\in B_j}\frac{64}{\sqrt p}\leq1.
\]
This can be done successively, since arbitrarily large unused
primes are available at each stage. Put $q_j=\prod_{p\in B_j}p$. The combined collection
$(u_j),(q_j)$ is odd and pairwise coprime.

For every odd squarefree $q$, the Chinese remainder theorem and the
local square-root count give
\begin{align}
 |\Lambda_q|=\prod_{p\mid q}\frac{p-1}{2}
   =\frac{\phi(q)}{2^{\omega(q)}},
 \qquad
 s_q=\frac{q}{|\Lambda_q|}
   =\prod_{p\mid q}\frac{2p}{p-1}.                \label{e:l1-size}
\end{align}
Consequently
\[
 \mathbb P_{u_j}(\Lambda_{u_j})\geq\frac13,
 \qquad
 3^{-j}\leq\mathbb P_{q_j}(\Lambda_{q_j})\leq2^{-j},
 \qquad
 2^j\leq s_{q_j}\leq3^j.
\]
This proves \eqref{e:l1-density} and $s_{q_j}\to\infty$.

We verify the two remaining spacing conditions directly by finite
moment estimates. This avoids using the strict finite-fraction
inequality in \cite[Lemma~3.2]{L09}. For $\Lambda=\Lambda_q$ and
$m=|\Lambda|$, define
\[
 R_q(t)=|\{a\in\Lambda:a+t\in\Lambda\}|,
 \qquad
 T_q(t)=|\{(a,b)\in\Lambda^2:a-b-t\in\Lambda\}|.
\]
For an odd prime $p$, write $m_p=(p-1)/2$. The quadratic-character
estimate used in Lemma~\ref{l:quadratic-GK}, with two distinct shifts,
gives the following convenient uniform bound:
\begin{align}\label{e:rpup}
    R_p(t)\leq\left(1+\frac{54}{\sqrt p}\right)\frac{m_p^2}{p}
 \qquad(t\not\equiv0\pmod p).
\end{align}
Indeed the absolute error is at most $6\sqrt p$, while $m_p\geq p/3$.
Also $R_p(0)=m_p$ and
\[
 T_p(t)=\sum_{b\in\Lambda_p}R_p(-b-t).
\]
At most one term in this sum has zero shift. Applying \eqref{e:rpup} to the other terms and using
$m_p\leq(9/\sqrt p)m_p^3/p$ gives
\[
 T_p(t)\leq\left(1+\frac{63}{\sqrt p}\right)\frac{m_p^3}{p}
 \qquad(t\in\mathbb Z/p\mathbb Z).
\]
Both counts multiply under the Chinese remainder theorem. For
\(q=q_j\), our choice of \(B_j\) gives
\[\prod_{p\mid q}(1+64/\sqrt p)
\leq\exp(\sum_{p\mid q}64/\sqrt p)\leq e<3.\]
Estimate \eqref{e:rpup} therefore yields
\begin{align}\label{e:tqrqbound}
 R_q(t)\leq3\frac{m^2}{q}
 \quad\hbox{if }t\not\equiv0\pmod p\text{ for every }p\mid q,
 \qquad
 T_q(t)\leq3\frac{m^3}{q}\quad\hbox{for every }t.
\end{align}
Fix \(0<\gamma\leq1\). There are between \(\gamma q/m-1\) and
\(\gamma q/m\) integers \(d\) with \(0<d<\gamma q/m\).
For such distinct pairs \(d,e\), we have
\(0<|d-e|<\gamma q/m\leq3^j<p\) for every \(p\in B_j\).
The bound for \(R_q\) therefore gives
\(|(\Lambda+d)\cap(\Lambda+e)|\leq3m^2/q\).
The union bound and the first two terms of inclusion--exclusion give
\[
 \left(\frac{\gamma q}{m}-1\right)m
 -\frac32\left(\frac{\gamma q}{m}\right)^2\frac{m^2}{q}
 \leq|\Lambda^\gamma|
 \leq\frac{\gamma q}{m}\,m.
\]
Dividing by \(q\), we obtain
\[
 \gamma-\frac mq-\frac32\gamma^2
 \leq\mathbb P_{q_j}(\Lambda_{q_j}^{\gamma})\leq\gamma.
\]
Since \(m/q\leq2^{-j}\to0\), it follows that
\[
 0\leq\gamma-\liminf_{j\to\infty}
 \mathbb P_{q_j}(\Lambda_{q_j}^{\gamma})
 \leq\frac32\gamma^2.
\]
This proves \eqref{e:l1-thickening}.

If \(\operatorname{dist}_q(a-b,\Lambda)\leq\gamma q/m\), then
\(a-b-w\equiv t\pmod q\) for some \(w\in\Lambda\) and some
integer \(t\) with \(|t|\leq\gamma q/m\).
Summing \eqref{e:tqrqbound} over these at most
\(2\gamma q/m+1\) integers gives
\[
 \frac1{m^2}
 |\{(a,b)\in\Lambda^2:
       \operatorname{dist}_q(a-b,\Lambda)\leq\gamma q/m\}|
 \leq
 \frac1{m^2}\left(\frac{2\gamma q}{m}+1\right)\frac{3m^3}{q}
 =6\gamma+\frac{3m}{q}.
\]
Since \(m/q\leq2^{-j}\to0\), taking complements and then the
lower limit yields
\[
 \liminf_{j\to\infty}\frac1{|\Lambda_{q_j}|^2}
 |\{(a,b)\in\Lambda_{q_j}^2:
      \operatorname{dist}_{q_j}(a-b,\Lambda_{q_j})>
                      \gamma s_{q_j}\}|
 \geq1-6\gamma.
\]
For \(0<\gamma<1/16\), we have
\(1-6\gamma>5/8>5(1/160)\).
Thus \eqref{e:l1-difference} holds with
\(\gamma_0=1/16\) and \(\alpha=1/160\).

Finally, let \(Q\) be a product of distinct terms from \((u_j)\) and
\((q_j)\), and let \(a\in\Lambda_Q\).  The integer \(Q\) is odd and
square-free.  By the Chinese
remainder theorem, the congruence $b^2\equiv a\pmod Q$ has exactly
$2^{\omega(Q)}$ unit solutions $b\pmod Q$.  For each such $b$, the prime
number theorem in the fixed arithmetic progression $b\pmod Q$ gives
\[
 |\{k\leq N:p_k\equiv b\pmod Q\}|
    =\frac{N}{\phi(Q)}+o_Q(N).
\]
  Distinct roots give disjoint
progressions, apart from the finitely many primes dividing $Q$, and hence
\begin{align*}
 \lim_{N\to\infty}\frac1N
 |\{k\leq N:p_k^2\equiv a\pmod Q\}|
 &=\frac{2^{\omega(Q)}}{\phi(Q)}
  =\frac1{|\Lambda_Q|},
\end{align*}
where the last equality is \eqref{e:l1-size}.  Therefore
\eqref{e:l1-equidistribution} holds with any fixed $0<\beta<1$, for example
$\beta=1/2$.  All hypotheses of Theorem \ref{t:l1badblackbox} are now
verified, and $(p_k^2)$ is persistently universally $L^1$-bad.
\end{proof}

\section{Quantitative Goldbach-type estimate: reductions}
\label{s:reductionsGB}


We now prove that almost every integer $n\equiv12\pmod{24}$ is the sum
of a prime square and a prime of the form $x^2+y^2+1$, with $(x,y)=1$.
The proof combines an analytic estimate for prime squares with a
lower-bound sieve for the primitive sum-of-two-squares condition and an
upper-bound sieve for the remaining large prime factors congruent to
$3\pmod4$.

With \(\mathcal P\) as defined in \eqref{e:defcalP}, put
\begin{align*}
 \mathcal N
 =\{n<N^2:n\equiv12\pmod{24},\
 n\in\mathcal P+\mathbb P^2\}.
\end{align*}
To prove Theorem~\ref{t:repmotthm}, it is enough to show that, for every
\(A>0\),
\begin{align}\label{e:exceptionalgoal}
 |\{n<N^2:n\equiv12\pmod{24}\}\setminus\mathcal N|
 \ll_A N^2(\log N)^{-A}.
\end{align}
We may restrict throughout to
\(n\gg_A N^2(\log N)^{-A}\), since the number of smaller integers is
already of the required order.
For such \(n\), let
\begin{align}
 \mathbb P_{3,n}
 &=\{p\in\mathbb P:p\equiv3\pmod4,\ p\nmid n-1\},
 \label{e:seta}\\
 \mathcal A_n
 &=\{r-1<N^2:r\equiv11\pmod{24},\ r\in\mathbb P,\
 n-r\in\mathbb P^2\}.
 \label{e:seta2}
\end{align}
We use the sieve notation \(S(\mathcal F,\mathfrak P,z)\) introduced in
\eqref{e:sievenotation}.

\begin{lemma}\label{l:fromsiftedtop}
For every $\varepsilon>0$,
\begin{align}\label{e:fromsiftedtop}
 |\{r\in\mathcal P:n-r\in\mathbb P^2\}|
 \geq S(\mathcal A_n,\mathbb P_{3,n},N^2)-O_\varepsilon(N^\varepsilon).
\end{align}
\end{lemma}

\begin{proof}
Let $r-1$ be counted by $S(\mathcal A_n,\mathbb P_{3,n},N^2)$ and write
$n-r=q^2$ with $q$ prime. If no prime $p\equiv3\pmod4$ divides $r-1$,
then $r-1\equiv2\pmod4$ and the primitive sum-of-two-squares criterion
shows that $r-1=x^2+y^2$ with $(x,y)=1$. Hence $r\in\mathcal P$.

Otherwise, let $p\equiv3\pmod4$ divide $r-1$. Since $p\notin\mathbb
P_{3,n}$, we have $p\mid n-1$. Therefore
\begin{align*}
 p\mid(n-1)-(r-1)=q^2,
\end{align*}
so $p=q$. Thus every exceptional $r$ is of the form $r=n-q^2$, where
$q$ is a prime divisor of $n-1$. There are at most
$\omega(n-1)\ll_\varepsilon N^\varepsilon$ such primes.
\end{proof}
We now truncate the lower sieve. Fix
\begin{align}\label{e:alpharange}
 1<\alpha<2,\qquad z=N^{1/\alpha}.
\end{align}

\begin{lemma}\label{l:truncationT}
For every \(\varepsilon>0\),
\begin{align}\label{e:truncationT}
 S(\mathcal A_n,\mathbb P_{3,n},z)
 -
 |\{r\in\mathcal P:n-r\in\mathbb P^2\}|
 =
 |T_n|+O_\varepsilon(N^\varepsilon),
\end{align}
where
\begin{align}\label{e:Tfactorization}
    T_n=
 \left\{
 r-1=2up_1p_2:
 \begin{array}{l}
 r\in\mathbb P,\quad n-r\in\mathbb P^2,\\
 z\leq p_2\leq p_1,\quad p_i\in\mathbb P_{3,n},\\
 p\mid u\Longrightarrow p\equiv1\pmod4
 \end{array}
 \right\}.
\end{align}
\end{lemma}

\begin{proof}
Let \(r-1\) be counted by the difference on the left of
\eqref{e:truncationT}.  Then \(r-1\) has no prime divisor
\(p<z\) belonging to \(\mathbb P_{3,n}\), but \(r\notin\mathcal P\).
By the primitive sum-of-two-squares criterion, \(r-1\) has a prime
divisor congruent to \(3\pmod4\).

Write \(n-r=q^2\), with \(q\in\mathbb P\).  If a prime
\(p\equiv3\pmod4\) divides both \(r-1\) and \(n-1\), then
\[
 p\mid(n-1)-(r-1)=q^2,
\]
and hence \(p=q\).  Thus
\(q\mid n-1\), and \(q\) uniquely determines \(r=n-q^2\).  The number
of these representations is therefore at most
\(\omega(n-1)\ll_\varepsilon N^\varepsilon\). 
After removing them, every prime divisor of \(r-1\) congruent to
\(3\pmod4\) belongs to \(\mathbb P_{3,n}\), and hence is at least \(z\).

Since \(r\equiv11\pmod{24}\), one has \(r-1\equiv2\pmod8\).  The
number of prime factors of \(r-1\) congruent to \(3\pmod4\), counted
with multiplicity, is even.  There are therefore at least two such
factors.  There cannot be four, since
\[
 z^4=N^{4/\alpha}>N^2>r-1.
\]
Consequently,
\[
 r-1=2up_1p_2,
\]
where \(z\leq p_2\leq p_1\), \(p_i\in\mathbb P_{3,n}\), and every
prime divisor of \(u\) is congruent to \(1\pmod4\).  This gives \(T_n\).
The converse implication is immediate from the primitive
sum-of-two-squares criterion, completing the proof.
\end{proof}

For the factorization in \eqref{e:Tfactorization}, put $\ell=2up_2.$
Since $p_1p_2u<N^2/2$ and $p_2\leq p_1$,
\begin{align*}
 u<\frac{N^2}{2z^2},\qquad
 z\leq p_2<\left(\frac{N^2}{2u}\right)^{1/2}.
\end{align*}
Accordingly, define
\begin{align}\label{e:defL}
 L=\left\{
 2up_2:
 \begin{array}{l}
  u\in\mathbb N,\quad
  u<\dfrac{N^2}{2z^2},\\[2mm]
  p_2\in\mathbb P,\quad p_2\equiv3\pmod4,\quad
  z\leq p_2<\left(\dfrac{N^2}{2u}\right)^{1/2},\\[2mm]
  p\mid u\Longrightarrow p\equiv1\pmod4
  \quad\text{for every prime }p
 \end{array}
 \right\},
\end{align}
and put
\begin{align}\label{e:defLn}
 L_n=\{\ell\in L:(\ell,n-1)=1\}.
\end{align}
Thus, if \(\ell=2up_2\in L_n\), then
\(p_2\nmid n-1\), and hence \(p_2\in\mathbb P_{3,n}\).
The ordering \(p_1\geq p_2\) makes the choice of \(\ell=2up_2\)
canonical.
For every $\ell\in L$, define
\begin{align}\label{e:defMn}
 M_n(\ell)
 =\{\ell p+1\colon \ell p+1<N^2:p\equiv3\pmod4,\ p\in\mathbb P,\
 n-1-\ell p=m^2\text{ for some }m\in\mathbb P,\ m>3\}.
\end{align}
and let
\begin{align}\label{e:defWln}
 W_{\ell,n}=\{p\in\mathbb P:p\nmid12\ell n\},\qquad
 z_\ell=(N^2/\ell)^{1/4}.
\end{align}

\begin{lemma}\label{l:upperreduction}
For every $\varepsilon>0$,
\begin{align}\label{e:upperreduction}
 |T_n|
 \leq\sum_{\ell\in L_n}
 S(M_n(\ell),W_{\ell,n},z_\ell)+O_\varepsilon(N^\varepsilon).
\end{align}
\end{lemma}

\begin{proof}
Let \(r-1\in T_n\).  By \eqref{e:Tfactorization}, it has the canonical
factorization $r-1=2up_1p_2=\ell p_1$ and $\ell=2up_2\in L$,
where \(z\leq p_2\leq p_1\).  Write \(n-r=q^2\), with
\(q\in\mathbb P\).  Since \(n\equiv12\pmod{24}\) and
\(r\equiv11\pmod{24}\), we have $q^2=n-r\equiv1\pmod{24},$
so \(q>3\).  Moreover,
\[
 n-1-\ell p_1=n-r=q^2.
\]

We first remove the cases in which \(\ell\notin L_n\).  If
\((\ell,n-1)>1\), choose a prime \(s\mid(\ell,n-1)\). So $s\mid q^2$, 
and hence \(s=q\).  In particular, \(q\mid n-1\), and \(q\) determines
\(r=n-q^2\).  The number of such \(r\) is therefore at most $\omega(n-1)\ll_\varepsilon N^\varepsilon.$

For every remaining \(r-1\in T_n\), we have \(\ell\in L_n\).  The
relations
\[
 r=\ell p_1+1<N^2,
 \qquad
 n-1-\ell p_1=q^2,
 \qquad q>3,
\]
show that \(r\in M_n(\ell)\).  Furthermore,
\[
 z_\ell=\left(\frac{N^2}{\ell}\right)^{1/4}
 \leq N^{1/2}<z\leq p_1<r.
\]
Since \(r\) is prime, no prime below \(z_\ell\) divides \(r\).
Consequently, \(r\) is counted by
\(S(M_n(\ell),W_{\ell,n},z_\ell)\).

Thus, apart from \(O_\varepsilon(N^\varepsilon)\) elements, every
\(r-1\in T_n\) gives an \(\ell\in L_n\) for which \(r\) is counted by
the corresponding sifted set.  Summing over \(\ell\in L_n\) proves
\eqref{e:upperreduction}.
\end{proof}

Combining Lemmas~\ref{l:fromsiftedtop}, \ref{l:truncationT}, and
\ref{l:upperreduction}, we obtain the reduction used below.

\begin{corollary}\label{c:mastersievereduction}
For every $\varepsilon>0$,
\begin{align}\label{e:mastersievereduction}
 |\{r\in\mathcal P:n-r\in\mathbb P^2\}|
 &\geq S(\mathcal A_n,\mathbb P_{3,n},z)-\sum_{\ell\in L_n}
 S(M_n(\ell),W_{\ell,n},z_\ell)
 -O_\varepsilon(N^\varepsilon).
\end{align}
\end{corollary}
It remains to obtain a lower bound for
\(S(\mathcal A_n,\mathbb P_{3,n},z)\) and an upper bound for 
$$ \sum_{\ell\in L_n}S(M_n(\ell),W_{\ell,n},z_\ell).$$
The common analytic input for these two sieve estimates is the following
Bombieri--Vinogradov approximation, proved in
Section~\ref{s:app2trial}.
Fix a nonzero real function
\(\psi\in C_c^\infty((-1/8,1/8))\), and put
\[
 \eta=\frac{\psi*\widetilde\psi}{(\psi*\widetilde\psi)(0)},
 \qquad
 \widetilde\psi(t)=\psi(-t).
\]
Thus \(\eta(0)=1\), \(\operatorname{supp}\eta\subset[-1/4,1/4]\), and
\(\eta^\vee\geq0\).

\begin{lemma}[Bombieri--Vinogradov approximation]
\label{l:circmethod}
Let \(A,\delta>0\) and \(L_0<N^{2-\delta}\).  There are sufficiently
large constants \(B=B(A,\delta)\) and \(C=C(A,\delta)\) such that the
following holds.  Put
\[
 R=(\log N)^C,\qquad
 \eta_N(\theta)=\eta\left(\frac{N^2\theta}{R}\right),
\]
and suppose that $V_\ell\leq (N^2/\ell)^{1/2}(\log N)^{-B}$ for $\ell<L_0.$
Let \(|d_{k,\ell}|\leq1\) and
\((a_{k,\ell},k)=1\).
Define 
\begin{align}\quad \quad\label{e:defGR9}
 G_R(\ell,k,n)
 &:=
 \mathcal T_2(\ell,k,n;a_{k,\ell})
 \prod_{2<p\leq R}
 \frac{p\phi((p,k))}{(p-1)^2}
 \left|\left\{(u,v)\in\mathbb A_p^2:
 \begin{array}{l}
  u^2+\ell v\equiv n\pmod p,\\
  v\equiv a_{k,\ell}\pmod{(p,k)}
 \end{array}\right\}\right|,
\end{align}
where $ \mathcal T_2(\ell,k,n;a_{k,\ell}) $ is defined in \eqref{e:generalT2} and the second congruence is understood to be trivial when \(p\nmid k\).
Also define
\begin{align}\label{e:defHR9}
 H_{\ell,R}(n)
 :=
 \ell\int_{\mathbb R}\eta_N(\theta)^2
 \left(\int_2^N\frac{e(\theta u^2)}{\log u}\,du\right)
 \left(\int_{\substack{2\leq v\leq N^2\\\ell v<N^2}}
 \frac{e(\ell\theta v)}{\log v}\,dv\right)
 e(-n\theta)\,d\theta.
\end{align}
Then
\begin{align}\label{e:circmethod}
 \sum_{n<N^2}\left|
 \sum_{\ell<L_0}\sum_{k<V_\ell}
 d_{k,\ell}r_{k,\ell}(n)\right|^2
 \ll_{A,\delta}N^4(\log N)^{-A},
\end{align}
where
\begin{align}\label{e:unweightedremainder}
 r_{k,\ell}(n)
 &=
 \sum_{\substack{m^2+\ell t=n\\
 t\equiv a_{k,\ell}\pmod k}}
 \one_{\mathbb P}(m)\one_{\mathbb P}(t)
 -
 \frac{G_R(\ell,k,n)}{\ell\phi(k)}H_{\ell,R}(n).
\end{align}
\end{lemma}
We now specify the residue classes used in
Lemma~\ref{l:circmethod}.  In the lower sieve, for square-free \(d\)
supported on primes greater than \(3\), let
\(a_{24d,1}\pmod{24d}\) be determined by
\begin{align}\label{e:lowerresidue9}
 a_{24d,1}&\equiv11\pmod{24},
 &a_{24d,1}&\equiv1\pmod d,
\end{align}
Then \((a_{24d,1},24d)=1\).

For the upper sieve, write \(\ell=2v\), where \((v,6)=1\).  If \(d\)
is square-free and \((d,12\ell)=1\), the Chinese remainder theorem
determines \(a_{12d,\ell}\pmod{12d}\) by
\begin{align}\label{e:upperresidue9}
 va_{12d,\ell}&\equiv5\pmod{12},
 &\ell a_{12d,\ell}&\equiv-1\pmod d.
\end{align}
Since \(v\) is invertible modulo \(12\) and \(\ell\) is invertible
modulo \(d\), one has \((a_{12d,\ell},12d)=1\).  Moreover,
\[
 \ell a_{12d,\ell}\equiv10\pmod{12},
 \qquad
 \ell a_{12d,\ell}\equiv-1\pmod d.
\]
As \((12,d)=1\), these two congruences determine a unique class
\(b_{12d}\pmod{12d}\), depending only on \(d\), such that
\[
 \ell a_{12d,\ell}\equiv b_{12d}\pmod{12d}.
\]
In particular, \((b_{12d},12d)=2\).

\medskip

Fix \(A_0,K>0\).  Choose \(C=C(A_0,K)\) sufficiently large for
Lemma~\ref{l:circmethod} with logarithmic saving \(K\), and put
\begin{align}\label{e:alphachoice8}
 \alpha=\frac98,\qquad
 z=N^{1/\alpha}=N^{8/9},\qquad
 R=(\log N)^C.
\end{align}
To simplify notation, we suppress the dependence on \(R\) and set
\begin{align}\label{e:lowerXomega}
 X_n&=\frac{G_R(1,24,n)}{\phi(24)}H_{1,R}(n),\notag\\
 \frac{\omega_n(d)}d
 &=\frac{\phi(24)G_R(1,24d,n)}
 {\phi(24d)G_R(1,24,n)},\qquad
 \Theta_n(z)=
 \prod_{\substack{3<p<z\\p\equiv3\;(\mathrm{mod}\;4)}}
 \left(1-\frac{\omega_n(p)}p\right).
\end{align}
The following propositions are proved in
Section~\ref{s:sieveGB} by combining Lemma~\ref{l:circmethod} with the
lower and upper Rosser--Iwaniec sieves.  We use the sieve functions
\(f_\kappa\) and \(F_\kappa\) with the normalization fixed in
Section~\ref{ss:notaion}.  Their estimates are uniform for
\[
 n\equiv12\pmod{24},
 \qquad
 N^2(\log N)^{-A_0}\leq n<N^2.
\]

\begin{proposition}[Lower sieve]\label{p:sieveneedlower8}
For $N^2(\log N)^{-A_0}\leq n<N^2$, $n\equiv12\pmod{24}$,
\begin{align}\label{e:sieveneedlower8}
 S(\mathcal A_n,\mathbb P_{3,n},z)
 \geq X_n\Theta_n(z)\bigl(f_{1/2}(9/8)+o(1)\bigr)
 +\mathcal E_1(n),
\end{align}
where
\begin{align}\label{e:E1section8}
\sum_{n<N^2}|\mathcal E_1(n)|^2
\ll_{A_0,K} N^4(\log N)^{-K}.
\end{align}
\end{proposition}
For \(\alpha>1\), define
\begin{equation}\label{e:Wdef8}
 \mathcal W(\alpha)
 :=
 \frac{\alpha}{4\sqrt2}
 \int_2^{2\alpha}
 \frac{t-2+(t-1)\log(t-1)}
      {t^2(t-1)\sqrt{1-t/(2\alpha)}}\,dt .
\end{equation}
The integral is finite, since the singularity at \(t=2\alpha\) is integrable.
\begin{proposition}[Upper sieve]\label{p:sieveneedupper8}
For $N^2(\log N)^{-A_0}\leq n<N^2$, $n\equiv12\pmod{24}$,
\begin{align}\label{e:sieveneedupper8}
 \sum_{\ell\in L_n}S(M_n(\ell),W_{\ell,n},z_\ell)
 &\leq
 \bigl(c_{9/8}F_1(2)\mathcal W(9/8)+o(1)\bigr)
 X_n\Theta_n(z)+\mathcal E_2(n),
\end{align}
where
\begin{align}\label{e:calphadef8}
 c_\alpha=\frac{16e^{-\gamma/2}}{\sqrt{\pi\alpha}}
\end{align}
and
\begin{align}\label{e:E2section8}
\sum_{n<N^2}|\mathcal E_2(n)|^2
\ll_{A_0,K} N^4(\log N)^{-K}.
\end{align}
\end{proposition}

Corollary~\ref{c:mastersievereduction}, Proposition~\ref{p:sieveneedlower8}, and Proposition~\ref{p:sieveneedupper8} reduce the problem to comparing the two
sieve constants. Put
\begin{align}\label{e:Deltadef}
 \Delta
 =f_{1/2}(9/8)-c_{9/8}F_1(2)\mathcal W(9/8).
\end{align}

We shall also need the following two lemmas.  Their proofs are given in
Section~\ref{s:sieveGB}.
\begin{lemma}\label{l:maintermsize8}
Let \(A_0>0\).  There is \(C_{A_0}>0\) such that, for every fixed
sufficiently large \(C\), with \(R=(\log N)^C\), one has
\begin{align}\label{e:maintermlower8}
 X_n\Theta_n(z)\gg_{A_0,C}
 N(\log N)^{-C_{A_0}}
\end{align}
whenever $n\equiv12\pmod{24}$ and $N^2(\log N)^{-A_0}\leq n<N^2$.
\end{lemma}

\begin{lemma}\label{l:numericalcomparison}
With the normalization of \cite[(1.8)--(1.10)]{I80},
\begin{align}\label{e:Deltanumerical}
 \Delta>0.
\end{align}
\end{lemma}

We now complete the proof of Theorem~\ref{t:repmotthm}.

\begin{proof}[Proof of Theorem~\ref{t:repmotthm}]
Fix \(A_0>0\), and let \(C_{A_0}\) be given by
Lemma~\ref{l:maintermsize8}.  Set $ K=A_0+2C_{A_0}+10.$
Choose the constants in Lemma~\ref{l:circmethod} and in the two sieve
propositions for this value of \(K\), enlarging \(C\), if necessary,
so that Lemma~\ref{l:maintermsize8} also applies.

By Lemma~\ref{l:numericalcomparison}, $\Delta>0.$
Hence, for all sufficiently large \(N\), Corollary
\ref{c:mastersievereduction} and
Propositions~\ref{p:sieveneedlower8} and
\ref{p:sieveneedupper8} give
\begin{align}\label{e:finalcomparison8}
 |\{r\in\mathcal P:n-r\in\mathbb P^2\}|
 &\geq
 \frac{\Delta}2X_n\Theta_n(z)
 +\mathcal E_1(n)-\mathcal E_2(n)
 -O_\varepsilon(N^\varepsilon).
\end{align}
By \eqref{e:maintermlower8}, there is \(c_{A_0}>0\) such that
\[
 \frac{\Delta}2X_n\Theta_n(z)
 \geq
 4c_{A_0}N(\log N)^{-C_{A_0}}
\]
throughout the stated range of \(n\).  Taking \(\varepsilon=1/2\), we
also have, for sufficiently large \(N\),
\[
 O_\varepsilon(N^\varepsilon)
 \leq c_{A_0}N(\log N)^{-C_{A_0}}.
\]
Discard the integers \(n<N^2\) for which
\[
 |\mathcal E_1(n)|+|\mathcal E_2(n)|
 >
 c_{A_0}N(\log N)^{-C_{A_0}}.
\]
By Cauchy--Schwarz, Chebyshev's inequality, \eqref{e:E1section8}, and \eqref{e:E2section8}, their number is at most
\begin{align*}
 \frac{
 N^4(\log N)^{-K}
 }{
 N^2(\log N)^{-2C_{A_0}}
 }
 \ll_{A_0}
 N^2(\log N)^{-A_0}.
\end{align*}
For every remaining \(n\), the right-hand side of
\eqref{e:finalcomparison8} is at least
\[
 2c_{A_0}N(\log N)^{-C_{A_0}}>0.
\]
Thus \(n\in\mathcal P+\mathbb P^2\).

Finally, the integers below \(N^2(\log N)^{-A_0}\), together with the
discarded exceptional set, contribute
\(O_{A_0}(N^2(\log N)^{-A_0})\).  This proves
\eqref{e:exceptionalgoal}.
\end{proof}

\section{Quantitative Goldbach-type estimate: the circle method}
\label{s:app2trial}

This section proves Lemma~\ref{l:circmethod}.  Fourier inversion reduces
the problem to estimates for a prime-square exponential sum and a prime
sum in an arithmetic progression. 
Let
\begin{align*}
 P_U(\theta)
 &=
 \sum_{m\leq U}\one_{\mathbb P}(m)e(m^2\theta),\\
 F_{k,\ell,W}(\theta)
 &=
 \sum_{\substack{t\leq W\\
 \ell t<N^2\\
 t\equiv a_{k,\ell}\pmod k}}
 \one_{\mathbb P}(t)e(\ell t\theta),
\end{align*}
where \(2\leq U\leq N\) and \(2\leq W\leq N^2\).  The restriction
\(\ell t<N^2\) does not change any Fourier coefficient indexed by
\(n<N^2\).

Choose \(C=C(A,\delta)\) sufficiently large, put
\(R=(\log N)^C\), and define
\begin{align*}
 \mathfrak M
 &=
 \bigcup_{q\leq R}\ \bigcup_{a\in\mathbb A_q}
 \left\{\theta\in\mathbb T:
 \left\|\theta-\frac aq\right\|_{\mathbb T}
 \leq\frac{R}{N^2}\right\},
 \qquad
 \mathfrak m=\mathbb T\setminus\mathfrak M.
\end{align*}
Distinct reduced fractions with denominators at most \(R\) are
separated by at least \(R^{-2}\).  Hence the arcs in \(\mathfrak M\)
are disjoint for sufficiently large \(N\).  Moreover, they contain
the usual major arcs
\[
 \left\{\theta:
 \left\|\theta-\frac aq\right\|_{\mathbb T}
 \leq\frac{R}{qN^2}\right\},
 \qquad q\leq R.
\]
Thus \(\mathfrak m\) is contained in the minor-arc set occurring in
\cite[Theorem~1]{Xiumin2005}.  Partial summation therefore gives, for
every \(D>0\),
\begin{align}\label{e:prime-square-minor}
 \sup_{\theta\in\mathfrak m}|P_U(\theta)|
 \ll_D \frac{U}{\log U}(\log N)^{-D}
\end{align}
whenever \(U\geq N(\log N)^{-2D}\).  Below this range we use
\(|P_U(\theta)|\leq U\).

\begin{lemma}
\label{l:hybridminorarc}
Under the hypotheses of Lemma~\ref{l:circmethod}, uniformly for
\(2\leq U\leq N\) and \(2\leq W\leq N^2\),
\begin{align}\label{e:BVinput9}
 \sum_{n<N^2}\left|
 \int_{\mathfrak m}P_U(\theta)
 \sum_{\ell<L_0}\sum_{k<V_\ell}
 d_{k,\ell}F_{k,\ell,W}(\theta)e(-n\theta)\,d\theta
 \right|^2
 \ll_{A,\delta}N^4(\log N)^{-A}.
\end{align}
\end{lemma}

\begin{proof}
Bessel's inequality gives
\begin{align}\label{e:minorBessel}
 &\sum_{n<N^2}\left|
 \int_{\mathfrak m}P_U(\theta)
 \sum_{\ell<L_0}\sum_{k<V_\ell}
 d_{k,\ell}F_{k,\ell,W}(\theta)e(-n\theta)\,d\theta
 \right|^2\leq
 \sup_{\theta\in\mathfrak m}|P_U(\theta)|^2
 \int_{\mathbb T}\left|
 \sum_{\ell<L_0}\sum_{k<V_\ell}
 d_{k,\ell}F_{k,\ell,W}(\theta)\right|^2d\theta.
\end{align}
Opening the final integral and using orthogonality, we obtain
\begin{align*}
 &\int_{\mathbb T}\left|
 \sum_{\ell<L_0}\sum_{k<V_\ell}
 d_{k,\ell}F_{k,\ell,W}(\theta)\right|^2d\theta=
 \sum_{\ell_1,\ell_2<L_0}
 \sum_{\substack{k_1<V_{\ell_1}\\k_2<V_{\ell_2}}}
 d_{k_1,\ell_1}\overline{d_{k_2,\ell_2}}
 \sum_{\substack{t_i\in\mathbb P,\ t_i\leq W\\
 \ell_it_i<N^2\\
 t_i\equiv a_{k_i,\ell_i}\;(\mathrm{mod}\;k_i)\\
 \ell_1t_1=\ell_2t_2}}1.
\end{align*}
We estimate separately the terms with \(t_1=t_2\) and
\(t_1\neq t_2\).

If \(t_1=t_2\), then \(\ell_1=\ell_2=\ell\).  When the two
congruences are compatible, they restrict \(t_1\) to one residue class
modulo \([k_1,k_2]\).  Since \(|d_{k,\ell}|\leq1\), the contribution
of these terms is at most
\begin{align}\label{e:t1t2diag}
    \sum_{\ell<L_0}\sum_{k_1,k_2<V_\ell}
 \left(\frac{N^2}{\ell[k_1,k_2]}+1\right)
 \ll N^2(\log N)^4,
\end{align}
where we used
\[
 \sum_{k_1,k_2\leq V}\frac1{[k_1,k_2]}
 \ll(\log(2V))^3,
 \qquad
 V_\ell^2\leq\frac{N^2}{\ell}.
\]

Suppose now that \(t_1\neq t_2\).  Since \(t_1\) and \(t_2\) are
distinct primes, the relation $ \ell_1t_1=\ell_2t_2$
implies
\[
 \ell_1=st_2,\qquad \ell_2=st_1
\]
for some \(s\geq1\).  Put
\[
 w_\ell(t)=
 \sum_{k<V_\ell}
 \one_{t\equiv a_{k,\ell}\,(\mathrm{mod}\,k)}.
\]
Using again \(|d_{k,\ell}|\leq1\), the off-diagonal contribution is
bounded by
\[
 \sum_{\substack{s\geq1,\ t_1,t_2\in\mathbb P\\
 t_1\neq t_2,\ t_i\leq W\\
 st_1,st_2<L_0,\ st_1t_2<N^2}}
 w_{st_2}(t_1)w_{st_1}(t_2).
\]
For every dyadic interval \(t\sim T\), expanding the square gives
\begin{align}\label{e:wellsquare9}
 \sum_{t\sim T}w_\ell(t)^2
 &\leq
 \sum_{k_1,k_2<V_\ell}
 \left(\frac{T}{[k_1,k_2]}+1\right)
 \ll T(\log N)^3+V_\ell^2.
\end{align}
Indeed, for each \(k_1,k_2\), the two congruences are either
incompatible or determine one residue class modulo \([k_1,k_2]\).

Restrict \(s,t_1,t_2\) to dyadic intervals
\(s\sim S\) and \(t_i\sim T_i\).  For every nonempty such range,
\(ST_1T_2\ll N^2\).  After dropping the primality conditions,
Cauchy--Schwarz gives
\begin{align}\label{e:CSwst}
 &\sum_{\substack{s\sim S,\ t_i\sim T_i\\
 st_1,st_2<L_0,\ st_1t_2<N^2}}
 w_{st_2}(t_1)w_{st_1}(t_2)\leq
 \left(\sum_{\substack{s\sim S,\ t_2\sim T_2\\st_2<L_0}}
       \sum_{t_1\sim T_1}w_{st_2}(t_1)^2\right)^{1/2}
 \left(\sum_{\substack{s\sim S,\ t_1\sim T_1\\st_1<L_0}}
       \sum_{t_2\sim T_2}w_{st_1}(t_2)^2\right)^{1/2}.
\end{align}
By \eqref{e:wellsquare9} and \(V_\ell^2\leq N^2/\ell\),
\begin{align*}
 \sum_{\substack{s\sim S,\ t_2\sim T_2\\st_2<L_0}}
 \sum_{t_1\sim T_1}w_{st_2}(t_1)^2
 &\ll ST_1T_2(\log N)^3
   +N^2\sum_{s\sim S}\frac1s
        \sum_{t_2\sim T_2}\frac1{t_2}\ll N^2(\log N)^3.
\end{align*}
The same estimate holds with \(t_1\) and \(t_2\) interchanged.
Therefore \eqref{e:CSwst} is \(O(N^2(\log N)^3)\).

There are \(O((\log N)^3)\) choices for \(S,T_1,T_2\).
Combining \eqref{e:t1t2diag} and \eqref{e:CSwst}, we obtain
\begin{align}\label{e:remainingL2}
 \int_{\mathbb T}\left|
 \sum_{\ell<L_0}\sum_{k<V_\ell}
 d_{k,\ell}F_{k,\ell,W}(\theta)\right|^2d\theta
 &\ll N^2(\log N)^4
 +\sum_{S,T_1,T_2}N^2(\log N)^3\ll N^2(\log N)^6.
\end{align}

If \(U\geq N(\log N)^{-2D}\), then
\eqref{e:prime-square-minor}, \eqref{e:minorBessel}, and
\eqref{e:remainingL2} give \(O(N^4(\log N)^{4-2D})\).
For smaller \(U\), the bound \(|P_U|\leq U\) gives
\(O(N^4(\log N)^{6-4D})\).  Taking \(D\) sufficiently large proves
\eqref{e:BVinput9}.
\end{proof}
We now rescale the function \(\eta\) fixed before
Lemma~\ref{l:circmethod}.  On \(\mathbb R\), put
\[
 \eta_N(\theta)
 =\eta\left(\frac{N^2\theta}{R}\right),
\]
and use the same notation for its periodization on \(\mathbb T\).
Moreover, \(|\eta_N|\leq1\) and the Euclidean inverse
Fourier transform of \(\eta_N\) is nonnegative.

Since $\eta$ is smooth and even, $\eta'(0)=0$, and hence
$\eta(t)=1+O(t^2)$ for $|t|\leq1$.  Together with $|\eta|\leq1$, this
gives
\begin{align}\label{e:weyl1bound}
     |1-\eta_N(\theta)^2|
 \ll
 \min\left\{1,\left(\frac{N^2|\theta|}{R}\right)^2\right\}.
\end{align}
For \(2\leq U\leq N\) and \(2\leq W\leq N^2\), put
\begin{align}\label{e:euclideankernel}
 I_U(\gamma)
 &=
 \int_2^U\frac{e(\gamma u^2)}{\log u}\,du, \qquad 
 J_{\ell,W}(\gamma)
 =
 \int_{\substack{2\leq v\leq W\\\ell v<N^2}}
 \frac{e(\ell\gamma v)}{\log v}\,dv,
\end{align}
and define
\begin{align}\label{e:defHR92}
 H_{\ell,R}(n;U,W)
 :=
 \ell\int_{\mathbb R}
 \eta_N(\gamma)^2I_U(\gamma)J_{\ell,W}(\gamma)
 e(-n\gamma)\,d\gamma.
\end{align}
Thus
\[
 H_{\ell,R}(n)=H_{\ell,R}(n;N,N^2).
\]
\begin{lemma}\label{l:positiveH9}
Uniformly in \(\ell,U,W\) and \(n<N^2\),
\[
 H_{\ell,R}(n;U,W)\geq0,
 \qquad
 H_{\ell,R}(n;U,W)
 =
 H_\ell(n;U,W)+O(NR^{-1/2}),
\]
where
\[
 H_\ell(n;U,W)
 =
 \int_{\substack{2\leq u\leq U\\
 2\leq(n-u^2)/\ell\leq W}}
 \frac{du}{\log u\log((n-u^2)/\ell)}.
\]
\end{lemma}

\begin{proof}
Fourier inversion gives
\[
 \int_{\mathbb R}I_U(\gamma)J_{\ell,W}(\gamma)e(-n\gamma)\,d\gamma
 =\frac1\ell H_\ell(n;U,W).
\]
The inverse Fourier transforms of \(I_U\) and \(J_{\ell,W}\) are
nonnegative measures.  Moreover,
\[
 (\eta_N^2)^\vee
 =\eta_N^\vee*\eta_N^\vee\geq0.
\]
Hence the inverse Fourier transform of
\(\eta_N^2I_UJ_{\ell,W}\) is nonnegative, proving
\(H_{\ell,R}(n;U,W)\geq0\).

Van der Corput's estimate and integration by parts give
\begin{align}\label{e:euclideanbounds}
    |I_U(\gamma)|\ll\min\{N,|\gamma|^{-1/2}\},
 \qquad
 |J_{\ell,W}(\gamma)|
 \ll\ell^{-1}\min\{N^2,|\gamma|^{-1}\}.
\end{align}
By the definitions of \(H_{\ell,R}\) and \(H_\ell\),
\begin{align*}
 H_{\ell,R}(n;U,W)-H_\ell(n;U,W)
 &=
 \ell\int_{\mathbb R}
 \bigl(\eta_N(\gamma)^2-1\bigr)
 I_U(\gamma)J_{\ell,W}(\gamma)e(-n\gamma)\,d\gamma.
\end{align*}
Using \eqref{e:weyl1bound}, \eqref{e:euclideanbounds}, 
and splitting the integral at \(N^{-2}\) and \(R/N^2\), we obtain
\begin{align}\label{e:lrtol}
 |H_{\ell,R}(n;U,W)-H_\ell(n;U,W)|
 &\ll
 N^3\int_0^{N^{-2}}
 \left(\frac{N^2\gamma}{R}\right)^2d\gamma\\
 &\quad+
 \int_{N^{-2}}^{R/N^2}
 \left(\frac{N^2\gamma}{R}\right)^2\gamma^{-3/2}\,d\gamma
 +\int_{R/N^2}^{\infty}\gamma^{-3/2}\,d\gamma\\
 &\ll \frac{N}{R^2}+\frac{N}{R^{1/2}}
 \ll NR^{-1/2}.
\end{align}
\end{proof}

We next approximate the two exponential sums on \(\mathfrak M\).  If
\(\theta=a/q+\gamma\in\mathfrak M\), define
\begin{align*}
 P_U^{(0)}(\theta)
 &=
 \mathfrak s(-a,q)I_U(\gamma),\\
 F_{k,\ell,W}^{(0)}(\theta)
 &=
 \frac{\Upsilon(-a\ell,q,a_{k,\ell},k)}{\phi(k)}
 J_{\ell,W}(\gamma).
\end{align*}

\begin{lemma}
\label{l:majorapproximation9}
Uniformly for \(2\leq U\leq N\) and \(2\leq W\leq N^2\),
\begin{align}\label{e:majorapproximation9}
 \sum_{n<N^2}\left|
 \sum_{\ell<L_0}\sum_{k<V_\ell}d_{k,\ell}
 \int_{\mathfrak M}
 \bigl(P_UF_{k,\ell,W}
 -P_U^{(0)}F_{k,\ell,W}^{(0)}\bigr)(\theta)
 e(-n\theta)\,d\theta
 \right|^2
 \ll_{A,\delta}N^4(\log N)^{-A}.
\end{align}
\end{lemma}

\begin{proof}
Fix \(\theta=a/q+\gamma\in\mathfrak M\), put \(X_\ell=N^2/\ell\),
and write \(h=[k,q]\). Define
\[
 \Delta(X;h)=
 \sup_{\substack{(b,h)=1\\2\leq y\leq X}}
 \left|\pi(y;h,b)-\frac{\operatorname{Li}(y)}{\phi(h)}\right|,
 \qquad
 E(y;h,b)=\pi(y;h,b)-\frac{\operatorname{Li}(y)}{\phi(h)}.
\]
Both \(F_{k,\ell,W}\) and \(J_{\ell,W}\) are unchanged when \(W\)
is replaced by \(\min\{W,X_\ell\}\), so assume \(W\leq X_\ell\).
Every prime in the progression is coprime to \(k\). Those not
coprime to \(h\) must therefore divide \(q\), and their total
contribution has absolute value at most \(q\). For the remaining
primes, split into the reduced classes
\[
 C_{k,q}
 =\{b\in\mathbb A_h:b\equiv a_{k,\ell}\pmod k\}.
\]
Their contribution is
\[
 \sum_{b\in C_{k,q}}e(a\ell b/q)
 \sum_{\substack{t\leq W,\ \ell t<N^2\\t\equiv b\pmod h}}
 \one_{\mathbb P}(t)e(\ell\gamma t).
\]
Since \(W\leq X_\ell\), the strict cutoff can remove only the term
\(t=W=X_\ell\). Thus
\begin{align*}
 F_{k,\ell,W}(\theta)
 &=
 \sum_{\substack{t\leq W\\
 t\equiv a_{k,\ell}\;(\mathrm{mod}\;k)}}
 \one_{\mathbb P}(t)e(\ell\theta t)+O(1),\qquad 
 J_{\ell,W}(\gamma)
=
 \int_2^W\frac{e(\ell\gamma v)}{\log v}\,dv.
\end{align*}
The \(O(1)\) correction is zero unless \(W=X_\ell\) is a prime
in the progression, in which case it equals
\(-e(\ell\theta X_\ell)\).
For each \(b\in C_{k,q}\), Stieltjes integration gives
\begin{align*}
 \sum_{\substack{t\leq W\\t\equiv b\;(\mathrm{mod}\;h)}}
 \one_{\mathbb P}(t)e(\ell\gamma t)
 &=
 \frac1{\phi(h)}\int_2^W\frac{e(\ell\gamma v)}{\log v}\,dv
 +\int_{2^-}^W e(\ell\gamma v)\,dE(v;h,b).
\end{align*}
Partial summation gives
\begin{align*}
 \left|
 \int_{2^-}^W e(\ell\gamma v)\,dE(v;h,b)
 \right|
 &\leq
 |E(W;h,b)|
 +2\pi|\gamma|\ell\int_2^W|E(v;h,b)|\,dv\\
 &\leq
 \bigl(1+2\pi|\gamma|\ell(W-2)\bigr)\Delta(X_\ell;h)\ll
 \bigl(1+|\gamma|\ell X_\ell\bigr)\Delta(X_\ell;h).
\end{align*}
Since \(|\gamma|\leq R/N^2\), \(\ell X_\ell=N^2\), and
\(|C_{k,q}|\leq q\leq R\), summing these errors gives
\(O(R^2\Delta(X_\ell;[k,q]))\). Also,
\[
 \frac1{\phi([k,q])}
 \sum_{b\in C_{k,q}}e(a\ell b/q)
 =
 \frac{\Upsilon(-a\ell,q,a_{k,\ell},k)}{\phi(k)}.
\]
Restoring both the strict endpoint and the primes dividing \(q\),
we obtain an estimate for the full original sum:
\begin{align}\label{e:fklwapp}
 F_{k,\ell,W}(\theta)
 &=F_{k,\ell,W}^{(0)}(\theta)+E_{k,\ell,W}(\theta),\notag\\
 \sup_{2\leq W'\leq X_\ell}|E_{k,\ell,W'}(\theta)|
 &\ll R^2\Delta(X_\ell;[k,q])+q+1.
\end{align}

We use the  Bombieri--Vinogradov theorem: for every
\(T>0\), there is \(L=L(T)\) such that
\begin{align}\label{e:weightedBV9}
 \sum_{h\leq X^{1/2}(\log X)^{-L}}\Delta(X;h)
 \ll_T X(\log X)^{-T}.
\end{align}
For fixed \(q\) and \(h\), every positive integer \(k\) satisfying
\([k,q]=h\) has \(h/k\in\{1,\ldots,q\}\).
The map \(k\mapsto h/k\) is injective on this set, so each \(h\)
has at most \(q\) preimages. We obtain
\begin{align*}
 \sum_{1\leq k<V_{\ell}}
 \sup_{2\leq W'\leq X_\ell}
 |E_{k,\ell,W'}(a/q+\gamma)|
 &\ll
 R^3\sum_{h\leq RV_\ell}\Delta(X_\ell;h)
 +(q+1)V_{\ell}.
\end{align*}
Here \([k,q]\leq kq\leq RV_\ell\). Since
\(N^\delta<X_\ell\leq N^2\), we have $ \delta\log N<\log X_\ell\leq2\log N.$
For sufficiently large \(B=B(T,C,\delta)\) and \(N\),
\[
 RV_\ell
 \leq X_\ell^{1/2}(\log N)^{-B+C}
 \leq X_\ell^{1/2}(\log X_\ell)^{-L}.
\]
Thus \eqref{e:weightedBV9}, with a sufficiently strong logarithmic
saving, controls the first term uniformly in the denominator and
all cutoffs. Taking \(T=K+3C+1\) in \eqref{e:weightedBV9} and choosing
\(B=B(K,C,\delta)\) sufficiently large, we obtain
\begin{align*}
 R^3\sum_{1\leq\ell<L_0}
 \sum_{h\leq RV_\ell}\Delta(X_\ell;h)
 &\ll_{K,C,\delta}
 (\log N)^{-K-1}\sum_{1\leq\ell<L_0}X_\ell\\
 &=N^2(\log N)^{-K-1}
 \sum_{1\leq\ell<L_0}\frac1\ell \ll N^2(\log N)^{-K}.
\end{align*}
The finite error in \eqref{e:fklwapp} satisfies, uniformly for
\(q\leq R\),
\begin{align*}
 (q+1)\sum_{1\leq\ell<L_0}\sum_{1\leq k<V_\ell}1
 &\leq2RN\sum_{1\leq\ell<L_0}\ell^{-1/2}\\
 &\ll RN\sqrt{L_0}
 \ll N^{2-\delta/2}(\log N)^C \ll_{K,C,\delta}N^2(\log N)^{-K},
\end{align*}
where we used \(V_\ell\leq N/\sqrt\ell\),
\(L_0<N^{2-\delta}\), and \(R=(\log N)^C\).
Together with \eqref{e:fklwapp} and the estimate obtained from
\eqref{e:weightedBV9}, this proves
\begin{align}\label{e:majorlinearerror9}
 \sup_{\theta\in\mathfrak M}
 \sum_{1\leq\ell<L_0}\sum_{1\leq k<V_\ell}
 \sup_{2\leq W'\leq X_\ell}
 |E_{k,\ell,W'}(\theta)|
 \ll_{K,C,\delta}N^2(\log N)^{-K}.
\end{align}

We next approximate \(P_U\).  If
\(U<N(\log N)^{-2D}\), then
\[
 |P_U(\theta)|+|P_U^{(0)}(\theta)|
 \ll U\ll N(\log N)^{-D},
\]
so \eqref{e:majorsquare9} follows trivially.  We may therefore assume
that \(U\geq N(\log N)^{-2D}\), and hence
\(\log U\sim\log N\).
  Separating the primes
dividing \(q\) and then splitting the remaining primes into reduced
residue classes modulo \(q\), we have
\[
 P_U(\theta)
 =
 \sum_{r\in\mathbb A_q}e(ar^2/q)
 \sum_{\substack{m\leq U\\m\equiv r\;(\mathrm{mod}\;q)}}
 \one_{\mathbb P}(m)e(\gamma m^2)
 +O(\log q).
\]
Since \(q\leq R=(\log N)^C\), the Siegel--Walfisz theorem gives,
with an arbitrarily large logarithmic saving,
\[
 \max_{\substack{r\in\mathbb A_q\\2\leq y\leq U}}
 \left|
 \pi(y;q,r)-\frac{\operatorname{Li}(y)}{\phi(q)}
 \right|
 \ll_D N(\log N)^{-D-2C-3}.
\]
Partial summation and
\(|\gamma|U^2\leq R\) therefore give
\[
 \sum_{\substack{m\leq U\\m\equiv r\;(\mathrm{mod}\;q)}}
 \one_{\mathbb P}(m)e(\gamma m^2)
 =
 \frac1{\phi(q)}
 \int_2^U\frac{e(\gamma u^2)}{\log u}\,du
 +O_D\bigl(NR(\log N)^{-D-2C-3}\bigr).
\]
Summing this formula over \(r\in\mathbb A_q\), using \(q\leq R\), and
absorbing the primes dividing \(q\), we obtain
\begin{align}\label{e:majorsquare9}
 P_U(a/q+\gamma)
 &=
 \frac1{\phi(q)}
 \sum_{r\in\mathbb A_q}e(ar^2/q)
 \int_2^U\frac{e(\gamma u^2)}{\log u}\,du
 +O_D\bigl(N(\log N)^{-D}\bigr)\notag\\
 &=
 P_U^{(0)}(a/q+\gamma)
 +O_D\bigl(N(\log N)^{-D}\bigr).
\end{align}

Since the arcs are disjoint, and since there are
\(\ll R^2\) reduced fractions with denominator at most \(R\), we conclude that  $ \operatorname{meas}(\mathfrak M)\ll\frac{R^3}{N^2}.$
We have
\[
 P_UF_{k,\ell,W}-P_U^{(0)}F_{k,\ell,W}^{(0)}
 =
 P_U^{(0)}E_{k,\ell,W}
 +(P_U-P_U^{(0)})F_{k,\ell,W}.
\]
For the first term, Bessel's inequality gives
\begin{align*}
 \sum_{n<N^2}\left|
 \int_{\mathfrak M}P_U^{(0)}(\theta)
 \sum_{\ell<L_0}\sum_{k<V_\ell}
 d_{k,\ell}E_{k,\ell,W}(\theta)e(-n\theta)\,d\theta
 \right|^2&\leq
 \int_{\mathfrak M}|P_U^{(0)}(\theta)|^2
 \left|\sum_{\ell<L_0}\sum_{k<V_\ell}
 d_{k,\ell}E_{k,\ell,W}(\theta)\right|^2d\theta\\
 &\ll
 N^2\frac{R^3}{N^2}
 \bigl(N^2(\log N)^{-K}\bigr)^2
 =
 N^4R^3(\log N)^{-2K},
\end{align*}
where we used \eqref{e:majorlinearerror9} and
\(|P_U^{(0)}|\ll N\).
For the second term, Bessel's inequality, \eqref{e:majorsquare9}, and
\eqref{e:remainingL2} give
\begin{align*}
 &\sum_{n<N^2}\left|
 \int_{\mathfrak M}(P_U-P_U^{(0)})(\theta)
 \sum_{\ell<L_0}\sum_{k<V_\ell}
 d_{k,\ell}F_{k,\ell,W}(\theta)e(-n\theta)\,d\theta
 \right|^2\\
 &\qquad\leq
 \|P_U-P_U^{(0)}\|_{L^\infty(\mathfrak M)}^2
 \int_{\mathbb T}\left|
 \sum_{\ell<L_0}\sum_{k<V_\ell}
 d_{k,\ell}F_{k,\ell,W}(\theta)\right|^2d\theta\ll
 N^2(\log N)^{-2D}\,
 N^2(\log N)^6.
\end{align*}
Choosing \(K\) and \(D\) so that
\(2K>A+3C\) and \(2D>A+6\) completes the proof.
\end{proof}

We  identify the arithmetic factor.  For \(q\geq1\), put
\begin{align}\label{e:defGq9}
 \mathfrak G_q(\ell,k,n)
 =
 \sum_{a\in\mathbb A_q}
 \Upsilon(a\ell,q,a_{k,\ell},k)
 \mathfrak s(a,q)e(an/q).
\end{align}

\begin{lemma}
\label{l:arithmeticcoefficients9}
The series
\[
 \sum_{\substack{q\geq1\\P^+(q)\leq R}}
 \mathfrak G_q(\ell,k,n)
\]
converges absolutely and satisfies
\begin{align}\label{e:localGR92}
 \sum_{\substack{q\geq1\\P^+(q)\leq R}}
 \mathfrak G_q(\ell,k,n)
 &=
 \prod_{p\leq R}
 \left(1+\sum_{\nu\geq1}
 \mathfrak G_{p^\nu}(\ell,k,n)\right)
 =G_R(\ell,k,n),
\end{align}
where \(G_R\) is defined in \eqref{e:defGR9}.
Moreover, for some absolute \(C_1>0\),
\begin{align}\label{e:arithmeticcoefficients9}
 |\mathfrak s(a,q)|
 &\ll q^{-1/2}\tau(q)^{C_1},\notag\\
 |\Upsilon(a\ell,q,a_{k,\ell},k)|
 &\ll q^{-1}(q,k\ell)\tau(q)^{C_1},\\
 |\mathfrak G_q(\ell,k,n)|
 &\ll q^{-1}(q,k\ell)\tau(q)^{C_1},\notag
\end{align}
and
\begin{align}\label{e:harmonicsum9}
 \sum_{\ell<L_0}\sum_{k<V_\ell}
 \frac{(q,k\ell)}{\ell\phi(k)}
 \ll \tau(q)^{C_1}(\log N)^2.
\end{align}
\end{lemma}

\begin{proof}
The first two estimates in \eqref{e:arithmeticcoefficients9} follow
from \eqref{e:actualboundfors} and
\cite[Lemma~3.8]{laceyprogression}.  Equation
\eqref{e:fixedqcoefficient}, with \(b=a_{k,\ell}\), shows that
\(\mathfrak G_q(\ell,k,n)\) is multiplicative in \(q\).
Comparing the Euler products in
\eqref{e:generalcorupss}, we see that, for every odd prime \(p\),
\[
 \mathfrak G_{p^\nu}(\ell,k,n)=0\quad(\nu\geq2),
 \qquad
 |\mathfrak G_p(\ell,k,n)|
 \ll p^{-1}(p,k\ell).
\]
At \(p=2\), only finitely many exponents contribute and their total is
bounded absolutely.  It follows by multiplicativity that
\[
 |\mathfrak G_q(\ell,k,n)|
 \ll q^{-1}(q,k\ell)\tau(q)^{C_1}.
\]
This proves the third estimate in
\eqref{e:arithmeticcoefficients9}.  Moreover,
\(\mathfrak G_{p^\nu}(\ell,k,n)=0\) for every odd \(p\) and
\(\nu\geq2\), while at \(p=2\) only finitely many exponents contribute.
Thus, among the integers \(q\) satisfying \(P^+(q)\leq R\), only
finitely many have \(\mathfrak G_q(\ell,k,n)\neq0\).  Consequently,
\[
 \sum_{\substack{q\geq1\\P^+(q)\leq R}}
 |\mathfrak G_q(\ell,k,n)|<\infty.
\]
 Choose \(Q>R\) so that the interval \((R,Q)\) contains no prime.
Then \(p<Q\) is equivalent to \(p\leq R\).  Taking
\(b=a_{k,\ell}\) in
Corollary~\ref{c:generalmajorprimesquareprog} therefore gives
\eqref{e:localGR92}.

For \eqref{e:harmonicsum9}, we use
\[
 (q,k\ell)
 \leq
 \sum_{\substack{d_1d_2\mid q\\d_1\mid k,\ d_2\mid\ell}}
 \phi(d_1d_2).
\]
Consequently,
\begin{align*}
 \sum_{\ell<L_0}\sum_{k<V_\ell}
 \frac{(q,k\ell)}{\ell\phi(k)}
 &\leq
 \sum_{d_1d_2\mid q}\phi(d_1d_2)
 \sum_{\substack{\ell<L_0\\d_2\mid\ell}}\frac1\ell
 \sum_{\substack{k<V_\ell\\d_1\mid k}}\frac1{\phi(k)}\\
 &\ll
 (\log N)^2
 \sum_{d_1d_2\mid q}
 \frac{\phi(d_1d_2)}{d_2\phi(d_1)}\\
 &\leq
 (\log N)^2\sum_{d_1d_2\mid q}1
 \leq \tau(q)^2(\log N)^2,
\end{align*}
where we used
\(\phi(d_1d_2)\leq d_2\phi(d_1)\). After enlarging \(C_1\), if necessary, this proves
\eqref{e:harmonicsum9}.
\end{proof}

We complete the major arcs.

\begin{lemma}\label{l:majorcompletion9}
Under the hypotheses of Lemma~\ref{l:circmethod}, uniformly for
\(2\leq U\leq N\) and \(2\leq W\leq N^2\),
\begin{align}\label{e:majorcompletion9}
& \sum_{n<N^2}\left|
 \sum_{\ell<L_0}\sum_{k<V_\ell}d_{k,\ell}
 \left(
 \int_{\mathfrak M}
 P_U^{(0)}(\theta)F_{k,\ell,W}^{(0)}(\theta)e(-n\theta)\,d\theta
 -
 \frac{G_R(\ell,k,n)}{\ell\phi(k)}
 H_{\ell,R}(n;U,W)
 \right)\right|^2\\
&\hspace{7cm} \ll_{A,\delta}N^4(\log N)^{-A}.
\end{align}
\end{lemma}

\begin{proof}
After replacing \(a\) by \(-a\), the integral in
\eqref{e:majorcompletion9} equals
\begin{align}\label{e:hardmajor9}
 \frac1{\phi(k)}
 \sum_{q\leq R}\sum_{a\in\mathbb A_q}
 &\Upsilon(a\ell,q,a_{k,\ell},k)
 \mathfrak s(a,q)e(an/q)
 \int_{|\gamma|\leq R/N^2}
 I_U(\gamma)J_{\ell,W}(\gamma)e(-n\gamma)\,d\gamma.
\end{align}

We first replace the sharp cutoff in \eqref{e:hardmajor9} by
\(\eta_N(\gamma)^2\).  At the denominator scale
\(q\sim2^s\), Plancherel and disjointness of the corresponding
frequency intervals give
\begin{align*}
 &\left\|
 \left(
 \sum_{q\sim2^s}\sum_{a\in\mathbb A_q}
 \bigl(\one_{|\gamma|\leq R/N^2}-\eta_N(\gamma)^2\bigr)
 \mathfrak s(a,q)I_U(\gamma)
 \sum_{\ell<L_0}\sum_{k<V_\ell}
 \frac{d_{k,\ell}}{\phi(k)}
 \Upsilon(a\ell,q,a_{k,\ell},k)J_{\ell,W}(\gamma)
 \right)^\vee
 \right\|_2^2\\
 &\quad\ll
 \sum_{q\sim2^s}\sum_{a\in\mathbb A_q}
 q^{-3}\tau(q)^{O(1)}(\log N)^{O(1)}
 \int_{\mathbb R}
 \left|\one_{|\gamma|\leq R/N^2}-\eta_N(\gamma)^2\right|^2
 \min\{N^2,|\gamma|^{-1}\}
 \min\{N^4,|\gamma|^{-2}\}\,d\gamma\\
 &\quad\ll
 \sum_{q\sim2^s}\sum_{a\in\mathbb A_q}
 q^{-3}\tau(q)^{O(1)}(\log N)^{O(1)}
 \left(
 \frac{N^{14}}{R^4}\int_0^{N^{-2}}\gamma^4\,d\gamma
 +\frac{N^8}{R^4}\int_{N^{-2}}^{R/N^2}\gamma\,d\gamma
 \right)\\
 &\quad\ll
 N^4R^{-2}(\log N)^{O(1)}
 \sum_{q\sim2^s}\phi(q)q^{-3}\tau(q)^{O(1)}\\
 &\quad\ll
 N^4R^{-2}2^{-s+o(s)}(\log N)^{O(1)}.
\end{align*}
Here \(\gamma=\theta-a/q\); in the first inequality we used
\eqref{e:arithmeticcoefficients9}, \eqref{e:harmonicsum9}, and the
bounds for \(I_U\) and \(J_{\ell,W}\).  Since the frequency supports corresponding to distinct reduced
fractions \(a/q\), \(q\leq R\), are disjoint, we may sum the preceding
squared estimates over \(2^s\leq R\).  Their total is
\[
 \ll
 N^4R^{-2}(\log N)^{O(1)}
 \sum_{2^s\leq R}2^{-s+o(s)}
 \ll
 N^4R^{-2}(\log N)^{O(1)}
 \ll_A N^4(\log N)^{-A},
\]
provided \(C\) is sufficiently large.

Put $ Y=\frac{N}{100R^{1/2}}.$
Distinct reduced fractions with denominators at most \(Y\) are
separated by at least \(Y^{-2}\), whereas \(\eta_N\) is supported in
an interval of width \(O(R/N^2)\).  Thus all the translates
\(\eta_N(\,\cdot-a/q\,)\), \(q\leq Y\), have disjoint supports.
We may therefore extend the major-arc sum to the additional \(R\)-smooth
moduli
\[
 R<q\leq Y,\qquad P^+(q)\leq R,
\]
provided that we estimate the \(L^2\)-contribution of the added terms.
This contribution is negligible by the following argument.

For \(q\sim2^s\),  the Fourier
transform in \(n\) of the added terms is
\[
 \sum_{\substack{q\sim2^s\\P^+(q)\leq R}}\sum_{a\in\mathbb A_q}
 \eta_N(\gamma)^2\mathfrak s(a,q)I_U(\gamma)
 \sum_{\ell<L_0}\sum_{k<V_\ell}
 \frac{d_{k,\ell}}{\phi(k)}
 \Upsilon(a\ell,q,a_{k,\ell},k)J_{\ell,W}(\gamma),
 \qquad \gamma=\theta-a/q.
\]
The translated supports are disjoint at this denominator scale.
Therefore, by Plancherel, \eqref{e:arithmeticcoefficients9}, and
\eqref{e:harmonicsum9}, its squared \(\ell^2\)-norm is
\begin{align*}
 &\ll
 \sum_{\substack{q\sim2^s\\P^+(q)\leq R}}\sum_{a\in\mathbb A_q}
 q^{-3}\tau(q)^{O(1)}(\log N)^{O(1)}
 \int_{\mathbb R}|\eta_N(\gamma)|^4
 \min\{N^2,|\gamma|^{-1}\}
 \min\{N^4,|\gamma|^{-2}\}\,d\gamma\\
 &\ll
 N^4(\log N)^{O(1)}
 \sum_{q\sim2^s}\phi(q)q^{-3}\tau(q)^{O(1)}
 \ll N^4\,2^{-s+o(s)}(\log N)^{O(1)}.
\end{align*}
Taking square roots, summing over the dyadic scales \(2^s>R\), and
then squaring gives
\begin{align*}
 N^4(\log N)^{O(1)}
 \left(\sum_{2^s>R}2^{-s/2+o(s)}\right)^2
 \ll N^4R^{-1+o(1)}(\log N)^{O(1)},
\end{align*}
which shows the contribution of these terms are negligible, as claimed above. 
After replacing
\(\mathbf1_{\{|\gamma|\leq R/N^2\}}\) by \(\eta_N^2\) and extending
the sum from \(q\leq R\) to $q\leq Y$ and $P^+(q)\leq R$.
Fourier inversion and \eqref{e:defHR92} give
\begin{align}\label{e:completedmajor9}
 &\sum_{\ell<L_0}\sum_{k<V_\ell}d_{k,\ell}
 \int_{\mathfrak M}
 P_U^{(0)}(\theta)F_{k,\ell,W}^{(0)}(\theta)
 e(-n\theta)\,d\theta\notag\\
 &\qquad=
 \sum_{\ell<L_0}\sum_{k<V_\ell}
 d_{k,\ell}\frac{H_{\ell,R}(n;U,W)}{\ell\phi(k)}
 \sum_{\substack{q\leq Y\\P^+(q)\leq R}}
 \mathfrak G_q(\ell,k,n)
 +\mathcal E(n),
\end{align}
where
\[
 \sum_{n<N^2}|\mathcal E(n)|^2
 \ll_{A,\delta}N^4(\log N)^{-A}.
\]

For the remaining terms with \(q>Y\) and \(P^+(q)\leq R\),
\eqref{e:arithmeticcoefficients9}, \eqref{e:harmonicsum9}, and
\(|H_{\ell,R}(n;U,W)|\ll N\) give
\begin{align*}
 \left|
 \sum_{\ell<L_0}\sum_{k<V_\ell}
 d_{k,\ell}\frac{H_{\ell,R}(n;U,W)}{\ell\phi(k)}
 \sum_{\substack{q>Y\\P^+(q)\leq R}}
 \mathfrak G_q(\ell,k,n)
 \right|&\leq
 \sum_{\substack{q>Y\\P^+(q)\leq R}}
 \sum_{\ell<L_0}\sum_{k<V_\ell}
 \frac{|H_{\ell,R}(n;U,W)|}{\ell\phi(k)}
 |\mathfrak G_q(\ell,k,n)|\\
 &\ll
 N(\log N)^{O(1)}
 \sum_{\substack{q>Y\\P^+(q)\leq R}}
 \frac{\tau(q)^{O(1)}}q.
\end{align*}
Rankin's argument gives
\begin{align}\label{e:arbitlogbound}
 \sum_{\substack{q>Y\\P^+(q)\leq R}}
 \frac{\tau(q)^{O(1)}}q
 &\leq
 Y^{-1/\log(2R)}
 \prod_{p\leq R}
 \left(
 \sum_{\nu\geq0}
 \frac{(\nu+1)^{O(1)}}
 {p^{\nu(1-1/\log(2R))}}
 \right)\\
 &\ll
 Y^{-1/\log(2R)}
 \exp\left(
 O\left(\sum_{p\leq R}
 p^{-1+1/\log(2R)}\right)\right)\\
 &\ll
 Y^{-1/\log(2R)}
 \exp\left(O\left(\sum_{p\leq R}\frac1p\right)\right)\\
 &\ll
 Y^{-1/\log(2R)}(\log R)^{O(1)}
 \ll_A(\log N)^{-A}.
\end{align}
After increasing \(C\), its square summed over \(n<N^2\) is
\(O_A(N^4(\log N)^{-A})\).  The completed arithmetic sum is now
exactly \(G_R(\ell,k,n)\), by \eqref{e:localGR92}, \eqref{e:completedmajor9}, and \eqref{e:arbitlogbound}.  This proves
\eqref{e:majorcompletion9}.
\end{proof}

We finish the proof of the stated approximation.

\begin{proof}[Proof of Lemma~\ref{l:circmethod}]
For \(n<N^2\), Fourier inversion gives
\begin{align}\label{e:fourierinversion9}
 \sum_{\substack{m\leq U,\ t\leq W\\
 m^2+\ell t=n\\
 t\equiv a_{k,\ell}\pmod k}}
 \one_{\mathbb P}(m)\one_{\mathbb P}(t)
 =
 \int_{\mathbb T}
 P_U(\theta)F_{k,\ell,W}(\theta)e(-n\theta)\,d\theta.
\end{align}
After multiplying by \(d_{k,\ell}\) and summing over \(\ell\) and
\(k\), introduce \(E_1(n),E_2(n),E_3(n)\) through the following
successive equalities:
\begin{align}\label{e:e123def}
 &\sum_{\ell<L_0}\sum_{k<V_\ell}d_{k,\ell}
 \sum_{\substack{m\leq U,\ t\leq W\\
 m^2+\ell t=n\\
 t\equiv a_{k,\ell}\pmod k}}
 \one_{\mathbb P}(m)\one_{\mathbb P}(t)\\
 &\quad=
 \sum_{\ell<L_0}\sum_{k<V_\ell}d_{k,\ell}
 \int_{\mathfrak M}
 P_U(\theta)F_{k,\ell,W}(\theta)e(-n\theta)\,d\theta
 +E_1(n),
 \qquad\text{by Lemma~\ref{l:hybridminorarc}},\\
 &\quad=
 \sum_{\ell<L_0}\sum_{k<V_\ell}d_{k,\ell}
 \int_{\mathfrak M}
 P_U^{(0)}(\theta)F_{k,\ell,W}^{(0)}(\theta)
 e(-n\theta)\,d\theta
 +E_1(n)+E_2(n),
 \qquad\text{by Lemma~\ref{l:majorapproximation9}},\\
 &\quad=
 \sum_{\ell<L_0}\sum_{k<V_\ell}d_{k,\ell}
 \frac{G_R(\ell,k,n)}{\ell\phi(k)}
 H_{\ell,R}(n;U,W)
 +E_1(n)+E_2(n)+E_3(n),
 \qquad\text{by Lemma~\ref{l:majorcompletion9}}.
\end{align}
The three cited lemmas give, uniformly in \(U\) and \(W\),
\begin{align}\label{e:threeerrors9}
 \sum_{n<N^2}
 \bigl(
 |E_1(n)|^2+|E_2(n)|^2+|E_3(n)|^2
 \bigr)
 \ll_{A,\delta}N^4(\log N)^{-A}.
\end{align}
Estimate \eqref{e:e123def} therefore implies
\begin{align}\label{e:cutoff9}
 &\sum_{n<N^2}\Bigg|
 \sum_{\ell<L_0}\sum_{k<V_\ell}d_{k,\ell}
 \left(
 \sum_{\substack{m\leq U,\ t\leq W\\
 m^2+\ell t=n\\
 t\equiv a_{k,\ell}\pmod k}}
 \one_{\mathbb P}(m)\one_{\mathbb P}(t)
 -
 \frac{G_R(\ell,k,n)}{\ell\phi(k)}
 H_{\ell,R}(n;U,W)
 \right)\Bigg|^2
 \ll_{A,\delta}N^4(\log N)^{-A}.
\end{align}
Taking \(U=N\) and \(W=N^2\) proves
\eqref{e:circmethod}.  Indeed, when \(n<N^2\), the equation
\(m^2+\ell t=n\) forces \(m<N\) and \(\ell t<N^2\), while $ H_{\ell,R}(n;N,N^2)=H_{\ell,R}(n).$
\end{proof}

\section{Quantitative Goldbach-type estimate: sieve arguments}
\label{s:sieveGB}

We prove Propositions~\ref{p:sieveneedlower8} and
\ref{p:sieveneedupper8}.  Recall that \(A_0,K>0\) and the parameters
\(\alpha,z,R\) were fixed in \eqref{e:alphachoice8}, with \(C\) chosen
so that Lemma~\ref{l:circmethod} gives the logarithmic saving \(K\).
Throughout this section,
\begin{align}\label{e:nrange10}
     n\equiv12\pmod{24},
 \qquad
 N^2(\log N)^{-A_0}\leq n<N^2,
\end{align}
and all constants implicit in \(o(1)\) are uniform in this range. Put
\[
 \mathbb P_3
 =\{p\in\mathbb P:p>3,\ p\equiv3\pmod4\},
 \qquad
 W_\ell=\{p\in\mathbb P:p\nmid12\ell\}.
\]
The following observation replaces the \(n\)-dependent sets of sieve
primes by fixed sets, allowing Lemma~\ref{l:circmethod} to be applied
with coefficients independent of \(n\).

\begin{lemma}\label{l:fixedsievesets10}
For every \(\varepsilon>0\),
\begin{align}
 S(\mathcal A_n,\mathbb P_{3,n},z)
 &\geq S(\mathcal A_n,\mathbb P_3,z),
 \label{e:fixedlower10}\\
 \sum_{\ell\in L_n}S(M_n(\ell),W_{\ell,n},z_\ell)
 &\leq
 \sum_{\ell\in L}S(M_n(\ell),W_\ell,z_\ell)
 +O_\varepsilon(N^\varepsilon).
 \label{e:fixedupper10}
\end{align}
\end{lemma}

\begin{proof}
Since \(r-1\equiv10\pmod{24}\), the prime \(3\) never divides an
element of \(\mathcal A_n\).  Thus the first inequality follows from
\(\mathbb P_{3,n}\setminus\{3\}\subseteq\mathbb P_3\).  For the
second, suppose that an integer \(r=\ell t+1\in M_n(\ell)\) is counted
after sifting by \(W_{\ell,n}\), but is removed after sifting by
\(W_\ell\).  Then some prime \(p<z_\ell\) divides both \(n\) and
\(r\).  Since \(n-r=m^2\) for a prime \(m\), it follows that \(m=p\),
and hence
\[
 r=n-p^2.
\]
For each \(p\mid n\), the number of representations
\(r-1=\ell t\), with \(\ell\in L\) and \(t\) prime, is at most
\(\tau(r-1)\).  Consequently, the total number of such pairs
\((\ell,r)\) is
\[
 \ll\sum_{\substack{p\mid n\\n-p^2>1}}
 \tau(n-p^2-1)\ll_\varepsilon N^\varepsilon.
\]
Enlarging the remaining sum from \(L_n\) to \(L\) can only increase it,
which proves \eqref{e:fixedupper10}.
\end{proof}

For an odd prime \(p\), set
\[
 \rho_p(a)=|\{x\in\mathbb F_p^\times:x^2=a\}|.
\]
Thus \(\rho_p(a)=1+(\frac ap)\) when \(p\nmid a\), and
\(\rho_p(a)=0\) when \(p\mid a\).

We now apply the residue classes fixed in
\eqref{e:lowerresidue9} and \eqref{e:upperresidue9} to the two sieve
sequences.  For the lower sieve, let \(d\) be square-free and supported
on \(\mathbb P_3\).  Since \((d,24)=1\), the congruences in
\eqref{e:lowerresidue9} give
\[
 d\mid r-1,\quad r\equiv11\pmod{24}
 \quad\Longleftrightarrow\quad
 r\equiv a_{24d,1}\pmod{24d}.
\]
Consequently,
\begin{align*}
 |\mathcal A_{n,d}|
 &=
 \sum_{\substack{m^2+t=n\\
 t\equiv a_{24d,1}\;(\mathrm{mod}\;24d)}}
 \one_{\mathbb P}(m)\one_{\mathbb P}(t)=
 \frac{G_R(1,24d,n)}{\phi(24d)}H_{1,R}(n)
 +r_{24d,1}(n),
\end{align*}
where the second equality follows from the definition of the remainder
in Lemma~\ref{l:circmethod}. Hence, recalling the definitions in \eqref{e:lowerXomega},
\begin{align}\label{e:lowersequence10}
 |\mathcal A_{n,d}|
 =
 X_n\frac{\omega_n(d)}d+r_{24d,1}(n).
\end{align}
Indeed,
\[
 X_n\frac{\omega_n(d)}d
 =
 \frac{G_R(1,24,n)}{\phi(24)}H_{1,R}(n)
 \frac{\phi(24)G_R(1,24d,n)}
      {\phi(24d)G_R(1,24,n)}
 =
 \frac{G_R(1,24d,n)}{\phi(24d)}H_{1,R}(n).
\]

\begin{lemma}\label{l:lowerlocaldensity}
For \(p>3\), \(p\equiv3\pmod4\),
\begin{align}\label{e:lowerlocaldensity10}
 \frac{\omega_n(p)}p
 =
 \begin{cases}
 \displaystyle
 \frac{\rho_p(n-1)}{p-1-\rho_p(n)},&p\leq R,\\[7pt]
 \displaystyle\frac1{p-1},&p>R.
 \end{cases}
\end{align}
In particular, for \(p\geq7\),
\(0\leq\omega_n(p)/p\leq2/(p-3)<1\).
\end{lemma}

\begin{proof}
For \(p\leq R\), the unrestricted local equation is
\(x^2+t=n\), with \(x,t\in\mathbb F_p^\times\).  It has
\(p-1-\rho_p(n)\) solutions.  The additional condition
\(p\mid t-1\) fixes \(t=1\) and leaves \(\rho_p(n-1)\) choices for
\(x\).  Taking the ratio of these two local factors proves the first
case.  If \(p>R\), the prime \(p\) does not occur in \(G_R\), and hence
\[
 \frac{\phi(24)G_R(1,24p,n)}
 {\phi(24p)G_R(1,24,n)}=\frac1{p-1}.
\]
The final claim follows from \(\rho_p(a)\leq2\).
\end{proof}

For the upper sieve, write \(\ell=2v\) and recall the class
\(a_{12d,\ell}\) from \eqref{e:upperresidue9}.  Since the prime \(m\)
is greater than \(3\), one has \(m^2\equiv1\pmod{24}\).  Thus
\[
 m^2+\ell t=n-1,\qquad n\equiv12\pmod{24}
\]
implies
\[
 \ell t\equiv10\pmod{24},
 \qquad\text{or equivalently}\qquad
 vt\equiv5\pmod{12}.
\]
Moreover,
\[
 d\mid\ell t+1
 \quad\Longleftrightarrow\quad
 \ell t\equiv-1\pmod d.
\]
By \eqref{e:upperresidue9}, these two conditions are equivalent to
\[
 t\equiv a_{12d,\ell}\pmod{12d}.
\]
It follows that
\begin{align}\label{e:upperexact10}
 |M_n(\ell)_d|
 =
 \sum_{\substack{m^2+\ell t=n-1\\
 t\equiv a_{12d,\ell}\;(\mathrm{mod}\;12d)}}
 \one_{\mathbb P}(m)\one_{\mathbb P}(t).
\end{align}

\begin{lemma}\label{l:upperlocaldensity}
There is a multiplicative function
\(\omega_n(\ell,\cdot)\), defined on square-free integers supported on
\(W_\ell\), such that
\begin{align}\label{e:upperfactorization10}
 \frac{G_R(\ell,12d,n-1)}{\phi(12d)}
 =
 \frac{G_R(\ell,12,n-1)}{\phi(12)}
 \frac{\omega_n(\ell,d)}d .
\end{align}
For \(p\nmid12\ell\), its  factors are
\begin{align}\label{e:upperlocaldensity10}
 \frac{\omega_n(\ell,p)}p
 =
 \begin{cases}
 \displaystyle
 \frac{\rho_p(n)}{p-1-\rho_p(n-1)},&p\leq R,\\[7pt]
 \displaystyle\frac1{p-1},&p>R.
 \end{cases}
\end{align}
In particular, \(0\leq\omega_n(\ell,p)/p<1\).
\end{lemma}

\begin{proof}
Suppose first that \(p\leq R\).  The unrestricted local equation
\[
 x^2+\ell t\equiv n-1\pmod p,
 \qquad x,t\in\mathbb F_p^\times,
\]
has \(p-1-\rho_p(n-1)\) solutions.  The condition
\(p\mid\ell t+1\) fixes \(t=-\ell^{-1}\), and the equation then becomes
\(x^2\equiv n\pmod p\), which has \(\rho_p(n)\) solutions.  The ratio
of these local factors is therefore
\[
 \frac{\rho_p(n)}{p-1-\rho_p(n-1)}.
\]
For \(p>R\), the prime \(p\) does not occur in the product defining
\(G_R\), so the factor
\(p\) to the modulus contributes only
\[
 \frac{\phi(12)}{\phi(12p)}=\frac1{p-1}.
\]
Multiplicativity gives \eqref{e:upperfactorization10}. If the local factor of \(G_R(\ell,12,n-1)\) at some prime \(p\le R\)
vanishes, then the corresponding local factor of
\(G_R(\ell,12d,n-1)\) also vanishes for every square-free \(d\) supported
on \(W_\ell\).  In this case we define \(\omega_n(\ell,p)\) by
\eqref{e:upperlocaldensity10}; the identity
\eqref{e:upperfactorization10} is then examined locally, and both sides
vanish whenever the common obstructing factor appears.

For \(p\geq7\), \eqref{e:upperlocaldensity10} is at most
\(2/(p-3)<1\).  The case \(p=5\) follows directly from the possible
values of \(\rho_5(n)\) and \(\rho_5(n-1)\), and gives a bound of
\(1/2\).
\end{proof}
Put
\begin{align}\label{e:upperXomega10}
 X_n(\ell)
 &=
 \frac{G_R(\ell,12,n-1)}{\ell\phi(12)}
 H_{\ell,R}(n-1),\qquad
 \Theta_{n,\ell}(y)
 =
 \prod_{\substack{p<y\\p\in W_\ell}}
 \left(1-\frac{\omega_n(\ell,p)}p\right).
\end{align}
Applying Lemma~\ref{l:circmethod} to \eqref{e:upperexact10} and using
\eqref{e:upperfactorization10}, we obtain
\begin{align}\label{e:uppersequence10}
 |M_n(\ell)_d|
 =
 X_n(\ell)\frac{\omega_n(\ell,d)}d
 +r_{12d,\ell}(n-1).
\end{align}

For \(p>R\), both local densities equal \(1/(p-1)\),
independently of \(n\). The uniform interval forms of Mertens'
theorem give, for some fixed \(K_0\geq2\) and all \(2\leq w<y\),
\begin{align}\label{e:outerdimensions10}
 \prod_{\substack{w\leq p<y\\p>3,\ p\equiv3\;(\mathrm{mod}\;4)}}
 \left(1-\frac1{p-1}\right)^{-1}
 &\leq
 \left(1+\frac{K_0}{\log w}\right)
 \left(\frac{\log y}{\log w}\right)^{1/2},\notag\\
 \prod_{\substack{w\leq p<y\\p\in W_\ell}}
 \left(1-\frac1{p-1}\right)^{-1}
 &\leq
 \left(1+\frac{K_0}{\log w}\right)
 \frac{\log y}{\log w}.
\end{align}
The second bound is uniform in \(\ell\), since restricting the
product over odd primes only decreases it. Restricting both
products to \(p>R\) gives outer dimensions \(1/2\) and \(1\),
respectively. We retain the factors at \(p\leq R\) by pre-sieving.


\begin{lemma}\label{l:presieving}
Suppose that
\begin{align}\label{e:circestsieve}
 |\mathcal F_d|=Xg(d)+r_d,
\end{align}
where \(X\geq0\), and \(g\) is multiplicative on the square-free
integers supported on a fixed set of primes \(\mathfrak P\),
with \(g(1)=1\) and \(0\leq g(p)<1\).
Fix \(C,A>0\), set \(R=(\log N)^C\),
\(D_0=\exp((\log\log N)^3)\), and \(D_1=D/D_0\), and assume
\(R+1\leq z\leq D_1\leq N^A\) for sufficiently large \(N\).
Suppose that fixed constants \(\kappa_0,\kappa>0\) and
\(K_0\geq2\) satisfy, for every \(2\leq w<y\leq z\),
\begin{align}\label{e:tworangesievedimension}
 \prod_{\substack{w\leq p<y\\p\in\mathfrak P,\ p\leq R}}
 (1-g(p))^{-1}
 &\leq
 \left(1+\frac{K_0}{\log w}\right)
 \left(\frac{\log y}{\log w}\right)^{\kappa_0},\notag\\
 \prod_{\substack{w\leq p<y\\p\in\mathfrak P,\ p>R}}
 (1-g(p))^{-1}
 &\leq
 \left(1+\frac{K_0}{\log w}\right)
 \left(\frac{\log y}{\log w}\right)^\kappa.
\end{align}
All \(o(1)\) terms below are uniform for fixed
\(C,A,\kappa_0,\kappa,K_0\).
Then
\begin{align}\label{e:presievedbounds10}
 S(\mathcal F,\mathfrak P,z)
 &\geq XV(z)\bigl(f_\kappa(s)+o(1)\bigr)
 +\sum_{m<D}c_m^-r_m,\notag\\
 S(\mathcal F,\mathfrak P,z)
 &\leq XV(z)\bigl(F_\kappa(s)+o(1)\bigr)
 +\sum_{m<D}c_m^+r_m,
\end{align}
where
\[
 V(z)=\prod_{\substack{p<z\\p\in\mathfrak P}}(1-g(p)),
 \qquad
 s=\frac{\log D_1}{\log z},
 \qquad
 |c_m^\pm|\leq1.
\]
The coefficients \(c_m^\pm\) depend only on the fixed prime sets and
the sieve levels, and not on the values of \(g\).  
The dimension
\(\kappa_0\) is used only in the preliminary sieve over \(p\leq R\);
the sieve functions in \eqref{e:presievedbounds10} are determined by
the outer dimension \(\kappa\).
\end{lemma}

\begin{proof}
Put
\[
 V_0(R)=\prod_{\substack{p\leq R\\ p\in\mathfrak P}}(1-g(p)),
 \quad
s_0=\frac{\log D_0}{\log(R+1)},\quad
 \mathfrak P_0=\{p\in\mathfrak P:p\leq R\},
 \quad
 \mathfrak P_1=\{p\in\mathfrak P:R<p<z\}.
\] 
Choose Rosser--Iwaniec weights satisfying
\begin{align}
 \operatorname{supp}\lambda_0^\pm
 \subseteq
 \{e<D_0:\mu^2(e)=1,\ p\mid e\Rightarrow p\in\mathfrak P_0\},\\
 \operatorname{supp}\lambda_1^\pm
 \subseteq
 \{d<D_1:\mu^2(d)=1,\ p\mid d\Rightarrow p\in\mathfrak P_1\}.
\end{align}
The preliminary Rosser--Iwaniec weights are chosen independently
of \(g\), with \(|\lambda_{0,e}^{\pm}|\leq1\) and the usual
pointwise sieve inequalities.
Apply the coefficient form of the fundamental lemma
\cite[(3.4)--(3.5)]{I80} to \(\mathfrak P_0\), level \(D_0\),
and sifting limit \(R+1\). This includes every prime in
\(\mathfrak P_0\), and the corresponding sieve product is
\(V_0(R)\). We obtain
\begin{align}\label{e:preliminarycoefficients10}
 &V_0(R)
 \bigl(f_{\kappa_0}(s_0)
       -O(Q_{\kappa_0}(s_0)(\log D_0)^{O(1)})\bigr)
 \leq
 \sum_e\lambda_{0,e}^-g(e)
 \leq V_0(R),\notag\\
 &V_0(R)
 \leq
 \sum_e\lambda_{0,e}^+g(e)
 \leq
 V_0(R)
 \bigl(F_{\kappa_0}(s_0)
       +O(Q_{\kappa_0}(s_0)(\log D_0)^{O(1)})\bigr).
\end{align}
The sieve functions satisfy
\[
 |F_{\kappa_0}(s)-1|+|f_{\kappa_0}(s)-1|
 \ll e^{-s},
\]
while \cite[(1.6)]{I80} gives
\[
 Q_{\kappa_0}(s)
 \leq
 \exp\{-s\log s+s\log\log(3s)+O_{\kappa_0}(s)\}.
\]
Since \(D_0=\exp((\log\log N)^3)\) and \(R=(\log N)^C\),
we have
\(s_0=\log D_0/\log(R+1)\gg_C(\log\log N)^2\).
Consequently, for every fixed \(K>0\),
\[
 |F_{\kappa_0}(s_0)-1|+|f_{\kappa_0}(s_0)-1|
 +Q_{\kappa_0}(s_0)(\log D_0)^{O(1)}
 \ll e^{-s_0}\ll_K(\log N)^{-K}.
\]
Thus, for either sign,
\begin{equation}\label{e:preliminarymain10}
 \sum_e\lambda_{0,e}^{\pm}g(e)
 =
 V_0(R)\bigl(1+O_K((\log N)^{-K})\bigr).
\end{equation}

For the lower bound, write
\[
 P_0=\prod_{p\in\mathfrak P_0}p,
 \qquad
 P_1=\prod_{p\in\mathfrak P_1}p.
\]
The outer lower weights satisfy the pointwise Rosser inequality
\[
 \sum_{d\mid a}\lambda_{1,d}^-
 \leq 1_{(a,P_1)=1}.
\]
Hence
\[
\begin{aligned}
 S(\mathcal F,\mathfrak P,z)
 =\sum_{a\in\mathcal F}
   1_{(a,P_0)=1}1_{(a,P_1)=1} &\geq
   \sum_{a\in\mathcal F}1_{(a,P_0)=1}
   \sum_{d\mid a}\lambda_{1,d}^-  \\
 &=\sum_d\lambda_{1,d}^-
   \sum_{\substack{a\in\mathcal F\\ d\mid a}}
   1_{(a,P_0)=1}
 =\sum_d\lambda_{1,d}^-S(\mathcal F_d,\mathfrak P_0,2R).
\end{aligned}
\]
Here the parameter \(2R\) is harmless because \(\mathfrak P_0\) already contains only primes
\(p\leq R\).

For the lower bound, if \(\lambda_{1,d}^-\geq0\) we insert the preliminary
lower weights into \(S(\mathcal F_d,\mathfrak P_0,2R)\), while if
\(\lambda_{1,d}^-<0\) we insert the preliminary upper weights.  In either
case, the remainders that occur are those in the approximation
\[
 |(\mathcal F_d)_e|
 =
 |\mathcal F_{ed}|
 =
 Xg(ed)+r_{ed},
\]
where $ \mathcal F_{ed}:=\{a\in\mathcal F:ed\mid a\}.$
Thus the coefficient of \(r_{ed}\) is
\[
 c_{ed}^-=
 \begin{cases}
  \lambda_{1,d}^-\lambda_{0,e}^-,&\lambda_{1,d}^-\geq0,\\
  \lambda_{1,d}^-\lambda_{0,e}^+,&\lambda_{1,d}^-<0.
 \end{cases}
\]
For the upper bound the choices are reversed, and
\[
 c_{ed}^+=
 \begin{cases}
  \lambda_{1,d}^+\lambda_{0,e}^+,&\lambda_{1,d}^+\geq0,\\
  \lambda_{1,d}^+\lambda_{0,e}^-,&\lambda_{1,d}^+<0.
 \end{cases}
\]
The two prime ranges are disjoint, so every \(m=ed\) has a unique such
factorization with \(e\) supported on \(\mathfrak P_0\) and \(d\) supported
on \(\mathfrak P_1\).  Defining \(c_m^\pm:=c_{ed}^\pm\), we therefore have
\[
 |c_m^\pm|\leq1,\qquad m=ed<D_0D_1=D.
\]
Put \(V_1(z)=\prod_{p\in\mathfrak P_1}(1-g(p))\), so that
\(V(z)=V_0(R)V_1(z)\).
By \eqref{e:preliminarymain10}, the preliminary error is
\[
 \ll_K
 X V_0(R)(\log N)^{-K}
 \sum_d|\lambda_{1,d}^{\pm}|g(d).
\]
The weights satisfy \(|\lambda_{1,d}^{\pm}|\leq1\) and are
supported on square-free integers whose prime factors lie in
\(\mathfrak P_1\). Hence
\[
 \sum_d|\lambda_{1,d}^{\pm}|g(d)
 \leq
 \sum_{\substack{d\ {\rm square\text{-}free}\\
                 p\mid d\Rightarrow p\in\mathfrak P_1}}g(d)
 =
 \prod_{p\in\mathfrak P_1}(1+g(p)).
\]
The equality follows by expanding the product: choosing \(g(p)\)
from the factors corresponding to the primes dividing \(d\)
produces \(g(d)\), by multiplicativity, and each such square-free
integer occurs exactly once.
Since \(0\leq g(p)<1\), we also have
\(1+g(p)\leq(1-g(p))^{-1}\).
Thus the preliminary error is
\[
 \ll_K
 X V_0(R)(\log N)^{-K}V_1(z)^{-1}
 =
 X V(z)(\log N)^{-K}V_1(z)^{-2}.
\]
Taking \(y=z\) and letting \(w\downarrow R\) in
\eqref{e:tworangesievedimension} gives
\[
 V_1(z)^{-1}
 \ll \left(\frac{\log z}{\log R}\right)^\kappa
 \ll(\log N)^\kappa,
\]
where we used \(z\leq N^A\) and \(R\geq2\).
The preliminary error is therefore
\(O_K(XV(z)(\log N)^{-K+2\kappa})\), which is
\(o(XV(z))\) upon choosing \(K>2\kappa\).
By the dimension-\(\kappa\) hypothesis
\eqref{e:tworangesievedimension}, the Rosser--Iwaniec coefficient
bounds for the outer weights give
\begin{align}\label{e:outercoefficients10}
 \sum_d\lambda_{1,d}^-g(d)
 &\geq V_1(z)\bigl(f_\kappa(s)+o(1)\bigr),\notag\\
 \sum_d\lambda_{1,d}^+g(d)
 &\leq V_1(z)\bigl(F_\kappa(s)+o(1)\bigr),
 \qquad
 s=\frac{\log D_1}{\log z}.
\end{align}

For the lower bound, the contribution of the main terms is
\begin{align*}
 &X\sum_{\lambda_{1,d}^-\geq0}
 \lambda_{1,d}^-g(d)
 \sum_e\lambda_{0,e}^-g(e)
 +X\sum_{\lambda_{1,d}^-<0}
 \lambda_{1,d}^-g(d)
 \sum_e\lambda_{0,e}^+g(e)\\
 &\qquad=
 XV_0(R)\sum_d\lambda_{1,d}^-g(d)+o(XV(z))\\
 &\qquad\geq
 XV_0(R)V_1(z)\bigl(f_\kappa(s)+o(1)\bigr)
 =
 XV(z)\bigl(f_\kappa(s)+o(1)\bigr),
\end{align*}
where the first equality follows from
\eqref{e:preliminarymain10} and the preceding estimate for the
preliminary error.

Similarly, the contribution of the main terms to the upper bound is
\begin{align*}
 &X\sum_{\lambda_{1,d}^+\geq0}
 \lambda_{1,d}^+g(d)
 \sum_e\lambda_{0,e}^+g(e)
 +X\sum_{\lambda_{1,d}^+<0}
 \lambda_{1,d}^+g(d)
 \sum_e\lambda_{0,e}^-g(e)\\
 &\qquad=
 XV_0(R)\sum_d\lambda_{1,d}^+g(d)+o(XV(z))\\
 &\qquad\leq
 XV_0(R)V_1(z)\bigl(F_\kappa(s)+o(1)\bigr)
 =
 XV(z)\bigl(F_\kappa(s)+o(1)\bigr).
\end{align*}
Combining these estimates with the remainder coefficients
\(c_m^\pm\) defined above proves \eqref{e:presievedbounds10}.
\end{proof}

The preliminary dimension bound
\eqref{e:tworangesievedimension} holds uniformly in both
applications with \(\kappa_0=2\). Indeed,
Lemmas~\ref{l:lowerlocaldensity} and
\ref{l:upperlocaldensity} give
\(g(p)\leq2/(p-3)=2/p+O(p^{-2})\) for \(p\geq7\), while
\(g(5)\leq1/2\) when that prime occurs. Mertens' uniform interval
estimate supplies the required bound. This preliminary
dimension affects only the error; the final sieve functions
retain the outer dimensions \(1/2\) and \(1\).

The separation of the small primes from the remaining sieve primes is
similar to the two-stage sieve in
\cite[Section~2, especially (14)--(16)]{MatomakiAlmostPrimes},
where preliminary beta-sieve weights are combined with linear-sieve
weights on the complementary range of primes.

\begin{proof}[Proof of Proposition~\ref{p:sieveneedlower8}]
By \eqref{e:fixedlower10}, it is enough to bound
\(S(\mathcal A_n,\mathbb P_3,z)\) from below.  Apply
Lemma~\ref{l:presieving} to \eqref{e:lowersequence10}, with
\[
 D=\frac{N}{24(\log N)^B}.
\]
The moduli \(24d\), \(d<D\), lie in the range of
Lemma~\ref{l:circmethod}. Since \(z=N^{8/9}\) and \(D/D_0=N^{1-o(1)}\), the level
conditions in Lemma~\ref{l:presieving} hold with \(A=2\).
  By \eqref{e:outerdimensions10}, the outer
dimension is \(1/2\), and
\[
 \frac{\log(D/D_0)}{\log z}=\frac98+o(1).
\]
Consequently,
\begin{align}\label{e:lowerbound10}
 S(\mathcal A_n,\mathbb P_{3,n},z)
 \geq X_n\Theta_n(z)
 \bigl(f_{1/2}(9/8)+o(1)\bigr)+\mathcal E_1(n).
\end{align}
Here \(\Theta_n(z)\) is defined in \eqref{e:lowerXomega}, and
\[
 \mathcal E_1(n)=\sum_{d<D}c_d^-r_{24d,1}(n),
\]
where Lemma~\ref{l:presieving} supplies coefficients \(c_d^-\)
independent of \(n\), with \(|c_d^-|\leq1\).
Lemma~\ref{l:circmethod}, with \(\ell=1\), gives
\[
 \sum_{n<N^2}|\mathcal E_1(n)|^2
 \ll_{A_0,K}N^4(\log N)^{-K}.
\]
This proves the proposition.
\end{proof}
We next estimate the sum of the upper-sieve main terms.  The calculation
at primes \(p\leq R\) is exact, while all factors at \(p>R\) are
independent of \(n\).  Recall that \(H_1\) denotes the \(\ell=1\) case of
the archimedean main term from Lemma~\ref{l:positiveH9}, namely
\[
 H_1(n)= H_1(n;N,N^2) =
 \int_{\substack{2\leq y\leq N\\ 2\leq n-y^2\leq N^2}}
 \frac{dy}{\log y\,\log(n-y^2)} .
\]
For \(n\) satisfying \eqref{e:nrange10}, the interval $\sqrt n/4\leq y\leq \sqrt n/2$
is contained in the domain of integration.  On this interval we have $y\sim \sqrt n$ and $n-y^2\sim n,$
and hence
\[
 \log y\ll \log N,\qquad \log(n-y^2)\ll \log N.
\]
Therefore
\begin{align}\label{e:H1elementary10}
 H_1(n)
 &\geq
 \int_{\sqrt n/4}^{\sqrt n/2}
 \frac{dy}{\log y\,\log(n-y^2)}
 \gg \frac{\sqrt n}{(\log N)^2}
 \gg_{A_0} N(\log N)^{-A_0/2-2},
\end{align}
uniformly for \(n\) in the range \eqref{e:nrange10}.

\begin{lemma}
\label{l:wumatomakisum}
Let \(\alpha=9/8\).  Uniformly in the range \eqref{e:nrange10},
\begin{align}\label{e:wumatomakisum10}
 \sum_{\ell\in L}X_n(\ell)\Theta_{n,\ell}(z_\ell)
 \leq
 \bigl(c_\alpha\mathcal W(\alpha)+o(1)\bigr)
 X_n\Theta_n(z),
 \qquad
 c_\alpha=\frac{16e^{-\gamma/2}}{\sqrt{\pi\alpha}}.
\end{align}
Here \(\mathcal W(\alpha)\) is the integral defined in
\eqref{e:Wdef8}.
\end{lemma}

\begin{proof}
The proof contains a few steps. 

\emph{Step 1: The parameters attached to \(\ell\).}
Write \(\ell=2up_2\) as in \eqref{e:defL}.  Thus every prime divisor
of \(u\) is \(1\pmod4\), while \(p_2\equiv3\pmod4\), and
\[
 u<\frac{N^2}{2z^2},\qquad
 z\leq p_2<\left(\frac{N^2}{2u}\right)^{1/2}.
\]
Recall also that the sifting level attached to \(M_n(\ell)\) is
\[
 z_\ell=\left(\frac{N^2}{\ell}\right)^{1/4}
       =\left(\frac{N^2}{2up_2}\right)^{1/4}.
\]
Since \(z=N^{1/\alpha}\) and \(\alpha=9/8\), we have \(z=N^{8/9}\), and
hence \(u<N^2/(2z^2)=N^{2/9}/2\).  Using the upper bound for \(p_2\), we get
\[
 2u^5p_2
 <
 2u^5\left(\frac{N^2}{2u}\right)^{1/2}
 =
 \sqrt2\,Nu^{9/2}
 <
 N^2.
\]
This is equivalent to \(u^4<N^2/(2up_2)=z_\ell^4\), and therefore
\(u<z_\ell\).  In particular, every prime divisor of \(u\) lies below
\(z_\ell\). Finally, since \(p_2\geq z\) and \(\ell=2up_2\geq p_2\), we have
\[
 z_\ell
 =\left(\frac{N^2}{\ell}\right)^{1/4}
 \leq\left(\frac{N^2}{p_2}\right)^{1/4}
 \leq\left(\frac{N^2}{z}\right)^{1/4}
 =N^{5/18}<N^{8/9}=z.
\]
Thus \(p_2\geq z>z_\ell\), and, since \(R=(\log N)^C\), also
\(p_2>R\) for all sufficiently large \(N\).

\emph{Step 2: Factors at primes \(p\leq R\).}
We compare the factors at primes \(p\leq R\) in
\[
 G_R(1,24,n)\Theta_n(z)
 \quad\text{and}\quad
 G_R(\ell,12,n-1)\Theta_{n,\ell}(z_\ell),
 \qquad \ell=2up_2.
\]
For \(3<p\leq R\), put
\[
 U_p=
 \frac{p\bigl(p-1-\rho_p(n)-\rho_p(n-1)\bigr)}
 {(p-1)^2}.
\]
We compute the factors in the two sieve products using
Lemmas~\ref{l:lowerlocaldensity} and
\ref{l:upperlocaldensity}.

Suppose first that \(p\nmid\ell\).  Since \(y\) is uniquely determined
by \(x\) in
\[
 x^2+\ell y\equiv n-1\pmod p,
\]
and \(y=0\) precisely when \(x^2\equiv n-1\pmod p\), we have
\[
 \left|\{(x,y)\in\mathbb A_p^2:
 x^2+\ell y\equiv n-1\pmod p\}\right|
 =p-1-\rho_p(n-1).
\]
Thus the factor of \(G_R(\ell,12,n-1)\) at \(p\) is
\[
 \frac{p\bigl(p-1-\rho_p(n-1)\bigr)}{(p-1)^2}.
\]
Lemma~\ref{l:upperlocaldensity} gives
\[
 \frac{\omega_n(\ell,p)}p
 =
 \frac{\rho_p(n)}
 {p-1-\rho_p(n-1)}.
\]
Consequently, the factor of
\(G_R(\ell,12,n-1)\Theta_{n,\ell}(z_\ell)\) at \(p\) is
\begin{align*}
 &\frac{p\bigl(p-1-\rho_p(n-1)\bigr)}{(p-1)^2}
 \left(
 1-\frac{\rho_p(n)}
 {p-1-\rho_p(n-1)}
 \right)=
 \frac{p\bigl(p-1-\rho_p(n)-\rho_p(n-1)\bigr)}
 {(p-1)^2}
 =U_p.
\end{align*}

If \(p\mid\ell\), then \(p\leq R<z\leq p_2\), so \(p\mid u\) and hence
\(p\equiv1\pmod4\).  In this case
\[
 x^2+\ell y\equiv n-1\pmod p
 \quad\Longleftrightarrow\quad
 x^2\equiv n-1\pmod p,
\]
and therefore
\[
 \left|\{(x,y)\in\mathbb A_p^2:
 x^2+\ell y\equiv n-1\pmod p\}\right|
 =(p-1)\rho_p(n-1).
\]
Since \(p\notin W_\ell\), there is no factor from
\(\Theta_{n,\ell}(z_\ell)\).  The factor contributed by \(p\) is
therefore
\[
 \frac{p}{(p-1)^2}(p-1)\rho_p(n-1)
 =
 \frac{p\rho_p(n-1)}{p-1}.
\]

For the lower-sieve product, the congruence
\[
 x^2+y\equiv n\pmod p
\]
has \(p-1-\rho_p(n)\) solutions in \(\mathbb A_p^2\).  Hence the factor
of \(G_R(1,24,n)\) at \(p\) is
\[
 \frac{p\bigl(p-1-\rho_p(n)\bigr)}{(p-1)^2}.
\]
If \(p\equiv1\pmod4\), then \(p\) is not a lower-sieve prime, so this
is also its factor in \(G_R(1,24,n)\Theta_n(z)\).  If
\(p\equiv3\pmod4\), Lemma~\ref{l:lowerlocaldensity} gives
\[
 \frac{\omega_n(p)}p
 =
 \frac{\rho_p(n-1)}{p-1-\rho_p(n)}.
\]
The factor in the full lower-sieve product is therefore
\begin{align*}
 &\frac{p\bigl(p-1-\rho_p(n)\bigr)}{(p-1)^2}
 \left(
 1-\frac{\rho_p(n-1)}{p-1-\rho_p(n)}
 \right)=
 \frac{p\bigl(p-1-\rho_p(n)-\rho_p(n-1)\bigr)}
 {(p-1)^2}
 =U_p.
\end{align*}
Thus the factor contributed by \(p\) to
\(G_R(1,24,n)\Theta_n(z)\) is
\[
 \begin{cases}
 \displaystyle
 \frac{p\bigl(p-1-\rho_p(n)\bigr)}{(p-1)^2},
     &p\equiv1\pmod4,\\[7pt]
 U_p,&p\equiv3\pmod4.
 \end{cases}
\]

To compare these products, define the non-negative multiplicative
function \(h_R\), supported on integers composed of primes
\(1\pmod4\), by
\begin{align}\label{e:defhR10}
 h_R(p^j)=
 \begin{cases}
 \displaystyle
 \frac{\rho_p(n-1)(p-1)}
 {p-1-\rho_p(n)-\rho_p(n-1)},&p\leq R,\\[8pt]
 \displaystyle\frac{p-1}{p-2},&p>R,
 \end{cases}
 \qquad j\geq1.
\end{align}
The second case will account for the factors with \(p>R\) when the
sum over \(u\) is estimated below.  The denominator in the first case
is positive.  Indeed, \(\rho_p(n),\rho_p(n-1)\in\{0,2\}\), so it is at
least \(p-5>0\) for \(p\geq7\).  For \(p=5\), two consecutive residues
cannot both be non-zero quadratic residues, and hence
\(\rho_5(n)+\rho_5(n-1)\leq2\).

For \(p\leq R\), \(p\equiv1\pmod4\), we have
\begin{align}\label{e:smallprimecancel10}
 U_p\left(1+\frac{h_R(p)}{p-1}\right)
 &=
 \frac{p\bigl(p-1-\rho_p(n)\bigr)}{(p-1)^2},
 \qquad
 U_ph_R(p)=\frac{p\rho_p(n-1)}{p-1}.
\end{align}
Thus the first identity gives the factor in
\(G_R(1,24,n)\Theta_n(z)\), while the second gives the factor in
\(G_R(\ell,12,n-1)\Theta_{n,\ell}(z_\ell)\) when \(p\mid u\).
Set
\begin{align}\label{e:BRCR10}
 B_R(n)&=\prod_{3<p\leq R}U_p,\qquad
 C_R(n)=
 B_R(n)
 \prod_{\substack{3<p\leq R\\p\equiv1\;(\mathrm{mod}\;4)}}
 \left(1+\frac{h_R(p)}{p-1}\right).
\end{align}
A direct calculation at \(2\) and \(3\) gives the same factor
\(8\cdot3=24\) in \(G_R(1,24,n)\) and
\(G_R(\ell,12,n-1)\). This and Lemma~\ref{l:lowerlocaldensity} imply that the product of the factors at \(p\leq R\) in
\(G_R(1,24,n)\Theta_n(z)\) is
\begin{align}\label{e:smallprimelower10}
    24
 \prod_{\substack{3<p\leq R\\p\equiv3\;(\mathrm{mod}\;4)}}U_p
 \prod_{\substack{3<p\leq R\\p\equiv1\;(\mathrm{mod}\;4)}}
 \frac{p\bigl(p-1-\rho_p(n)\bigr)}{(p-1)^2}
 =24C_R(n),
\end{align}
Similarly, Lemma~\ref{l:upperlocaldensity} gives the product of the factors at \(p\leq R\) in
\(G_R(\ell,12,n-1)\Theta_{n,\ell}(z_\ell)\) to be
\begin{align}\label{e:smallprimeupper10}
     24
 \prod_{\substack{3<p\leq R\\p\nmid u}}U_p
 \prod_{\substack{3<p\leq R\\p\mid u}}
 \frac{p\rho_p(n-1)}{p-1}
 =
 24B_R(n)
 \prod_{\substack{p\mid u\\p\leq R}}h_R(p).
\end{align}
Thus \(C_R(n)\) is precisely the small-prime part of the lower-sieve
main term, while \(B_R(n)\prod_{p\mid u,\,p\leq R}h_R(p)\) is the
small-prime part of the upper-sieve term indexed by
\(\ell=2up_2\).

For primes beyond the pre-sieving range, let \(q>R\) be an outer sieve
prime.  Since \(q\nmid\ell\), the condition that \(q\) divide the sifted
integer \(\ell y+1\) is equivalent to \(y\equiv-\ell^{-1}\pmod q\).  Thus
it selects one reduced residue class for the prime variable \(y\) modulo
\(q\).  Hence the corresponding local density is $ \frac{1}{\varphi(q)}=\frac{1}{q-1}.
$
This factor is independent of \(n\).  Mertens' theorem gives, uniformly
for \(y\geq N^\varepsilon\),
\begin{align}\label{e:mertensratios10}
 \prod_{R<p<y}\left(1-\frac1{p-1}\right)
 &=(1+o(1))\frac{\log R}{\log y},\notag\\
 \prod_{\substack{R<p<y\\p\equiv3\; (\mathrm{mod}\;4)}}
 \left(1-\frac1{p-1}\right)
 &=(1+o(1))\left(\frac{\log R}{\log y}\right)^{1/2}.
\end{align}
We now compute the lower normalization.  By definition,
\[
 X_n\Theta_n(z)
 =
 \frac{G_R(1,24,n)}{\varphi(24)}H_{1,R}(n)
 \prod_{\substack{3<p<z\\ p\equiv3\pmod4}}
 \left(1-\frac{\omega_n(p)}p\right).
\]
We split the product at \(R\).  The factors with \(p\leq R\) are part of
the pre-sieved local computation. By \eqref{e:smallprimelower10}, the factors at \(p\leq R\) give
\[
 \frac{G_R(1,24,n)}{\phi(24)}
 \prod_{\substack{3<p\leq R\\p\equiv3\;(\mathrm{mod}\;4)}}
 \left(1-\frac{\omega_n(p)}p\right)
 =
 \frac{24}{\phi(24)}C_R(n)
 =
 3C_R(n).
\]  
Combining this with \eqref{e:mertensratios10} for \(p>R\), we obtain
\begin{align}\label{e:lowernormalization10}
 X_n\Theta_n(z)
 =
 \bigl(3+o(1)\bigr)
 C_R(n)H_{1,R}(n)
 \left(\frac{\log R}{\log z}\right)^{1/2}.
\end{align}
For \(q>R\), Lemma~\ref{l:upperlocaldensity} gives $ \frac{\omega_n(\ell,q)}q=\frac1{q-1}.$
Since \(p_2\geq z>z_\ell\), the primes \(q<z_\ell\) that divide
\(\ell=2up_2\) are precisely the prime divisors of \(u\).  Hence
\[
 \prod_{\substack{R<q<z_\ell\\q\in W_\ell}}
 \left(1-\frac{\omega_n(\ell,q)}q\right)
 =
 \prod_{R<q<z_\ell}
 \left(1-\frac1{q-1}\right)
 \prod_{\substack{q\mid u\\q>R}}
 \left(1-\frac1{q-1}\right)^{-1},
\]
where we also used \(u<z_\ell\), proved above.  For \(q>R\),
\[
 \left(1-\frac1{q-1}\right)^{-1}
 =\frac{q-1}{q-2}
 =h_R(q).
\]
Moreover, \(h_R(q^j)=h_R(q)\) for every \(j\geq1\), and therefore
\[
 h_R(u)=\prod_{q\mid u}h_R(q).
\]
Consequently,
\[
 \prod_{\substack{q\mid u\\q\leq R}}h_R(q)
 \prod_{\substack{q\mid u\\q>R}}
 \left(1-\frac1{q-1}\right)^{-1}
 =h_R(u).
\]
Combining this identity with \eqref{e:mertensratios10} gives, uniformly
for \(\ell=2up_2\in L\),
\begin{align}\label{e:uppernormalization10}
 X_n(\ell)\Theta_{n,\ell}(z_\ell)
 =
 \bigl(6+o(1)\bigr)
 B_R(n)h_R(u)\frac{H_{\ell,R}(n-1)}{\ell}
 \frac{\log R}{\log z_\ell}.
\end{align}

\emph{Step 3: The mean value of \(h_R\).}
We prove that, uniformly in \(n\),
\begin{align}\label{e:meanvaluehR10}
 \sum_{u\leq y}h_R(u)
 &=
 \left(\frac{e^{-\gamma/2}}{\sqrt\pi}+o(1)\right)
 \frac{C_R(n)}{B_R(n)}
 \frac{y}{\sqrt{\log R\log y}}
\end{align}
whenever $\log y\geq\frac{\log N}{\log\log N}.$
We shall also prove the uniform upper bound
\begin{align}\label{e:meanvaluehRupper10}
 \sum_{u\leq y}h_R(u)
 &\ll
 \frac{C_R(n)}{B_R(n)}
 \frac{y}{\sqrt{\log R\log y}}
 \qquad (y\geq R^2).
\end{align}

We first establish the corresponding asymptotic for a fixed
multiplicative function.  Let \(h_\infty\) be supported on integers
composed of primes congruent to \(1\pmod4\), and define
\[
 h_\infty(p^j)=\frac{p-1}{p-2}
 \qquad
 (p\equiv1\pmod4,\ j\geq1).
\]
Then
\begin{align*}
 \sum_{p\leq x}h_\infty(p)\log p
 &=
 \sum_{\substack{p\leq x\\p\equiv1\;(\mathrm{mod}\;4)}}
 \log p
 +O\left(\sum_{p\leq x}\frac{\log p}{p}\right)=
 \frac{x}{2}+O_K\bigl(x(\log x)^{-K}\bigr)
\end{align*}
for every fixed \(K>0\).  Moreover,
\[
 h_\infty(p^j)=\frac{p-1}{p-2}\leq j+1=\tau_2(p^j),
\]
and hence \(h_\infty(n)\leq\tau_2(n)\).  Since the prime-sum
estimate holds with an arbitrarily large logarithmic saving,
\cite[Theorem~1]{GranvilleKoukoulopoulos}, applied with
\(\alpha=1/2\) and \(k=2\), gives
\begin{align}\label{e:meanvaluehinfty10}
 \sum_{u\leq y}h_\infty(u)
 &=
 \left(\frac{c_\infty}{\sqrt\pi}
 +O\left(\frac1{\log y}\right)\right)
 \frac{y}{\sqrt{\log y}},
\end{align}
where
\begin{align}\label{e:defcinfty10}
 c_\infty
 =
 \prod_p
 \left(1+\sum_{j\geq1}\frac{h_\infty(p^j)}{p^j}\right)
 \left(1-\frac1p\right)^{1/2}.
\end{align}

We next compare \(h_R\) with \(h_\infty\), so we can replace $h_{\infty}$ in \eqref{e:meanvaluehinfty10} with $h_R$.  For this purpose, we
claim that there is a multiplicative function \(a_{R,n}\), supported
on integers composed of primes \(p\leq R\), \(p\equiv1\pmod4\), such
that
\begin{align}\label{e:hRconvolution10}
 h_R=h_\infty*a_{R,n},
\end{align}
and, for every fixed \(J\geq0\),
\begin{align}\label{e:aRnmoments10}
 \sum_{d\geq1}
 \frac{|a_{R,n}(d)|}{d}(1+\log d)^J
 \ll_J(\log R)^{O_J(1)}.
\end{align}
Here $*$ is the Dirichlet's convolution. We shall also need the estimate
\begin{align}\label{e:aRnshifted10}
 \sum_{d\geq1}
 \frac{|a_{R,n}(d)|}
 {d^{\,1-1/\log(2R)}}
 \ll(\log R)^{O(1)}.
\end{align}
To prove these claims, define \(a_{R,n}\) by
\begin{align}\label{e:defaRn10}
 \sum_{d\geq1}\frac{a_{R,n}(d)}{d^s}
 =
 \prod_{\substack{p\leq R\\p\equiv1\;(\mathrm{mod}\;4)}}
 \frac{
  1+\sum_{j\geq1}h_R(p^j)p^{-js}
 }{
  1+\sum_{j\geq1}h_\infty(p^j)p^{-js}
 },
 \qquad \Re s>1.
\end{align}
This immediately gives \eqref{e:hRconvolution10}.  Since both
\(h_R\) and \(h_\infty\) are constant on the positive powers of a
prime, the local quotient in \eqref{e:defaRn10} is $ \frac{1+(h_R(p)-1)p^{-s}}
      {1+(h_\infty(p)-1)p^{-s}}.$
It follows that
\begin{align}\label{e:aRnlocal10}
 a_{R,n}(p^j)
 &=
 \bigl(h_R(p)-h_\infty(p)\bigr)
 \bigl(1-h_\infty(p)\bigr)^{j-1}=
 \bigl(h_R(p)-h_\infty(p)\bigr)
 \left(-\frac1{p-2}\right)^{j-1}
\end{align}
for \(p\leq R\), \(p\equiv1\pmod4\), and \(j\geq1\).
By \eqref{e:defhR10}, one has \(0\leq h_R(p)\leq6\):
if \(\rho_p(n-1)=0\), then \(h_R(p)=0\), while otherwise this follows
from the lower bounds for the denominator established immediately
after \eqref{e:defhR10}.  Since
\(h_\infty(p)=(p-1)/(p-2)\leq4/3\), \eqref{e:aRnlocal10} gives
\begin{align}\label{e:aRnlocalbound10}
 |a_{R,n}(p^j)|\ll(p-2)^{1-j}.
\end{align}
We now prove \eqref{e:aRnshifted10}.  Uniformly for $ 1-\frac1{\log(2R)}\leq\sigma\leq1,$
equation \eqref{e:aRnlocalbound10} gives
\[
 \sum_{j\geq1}
 \frac{|a_{R,n}(p^j)|}{p^{j\sigma}}
 \ll p^{-\sigma}.
\]
Therefore
\begin{align*}
 \sum_{d\geq1}\frac{|a_{R,n}(d)|}{d^\sigma}
 &\leq
 \prod_{\substack{p\leq R\\p\equiv1\;(\mathrm{mod}\;4)}}
 \left(1+O(p^{-\sigma})\right)\ll
 \exp\left(O\left(\sum_{p\leq R}p^{-\sigma}\right)\right)
 \ll(\log R)^{O(1)},
\end{align*}
because $p^{1/\log(2R)}\ll1$ for $p\leq R$.
Taking \(\sigma=1-1/\log(2R)\) proves
\eqref{e:aRnshifted10}.

We next prove the logarithmic-moment estimate
\eqref{e:aRnmoments10}.  Put
\[
 A_{R,n}^{*}(s)
 :=
 \sum_{d\geq1}\frac{|a_{R,n}(d)|}{d^s}
 =
 \prod_{\substack{p\leq R\\p\equiv1\;(\mathrm{mod}\;4)}}
 \left(
 1+\sum_{\nu\geq1}
 \frac{|a_{R,n}(p^\nu)|}{p^{\nu s}}
 \right).
\]
By \eqref{e:aRnlocalbound10}, $ A_{R,n}^{*}(1)\ll(\log R)^{O(1)}.$
Moreover, for every fixed \(j\geq1\), logarithmic
differentiation of the Euler product gives
\begin{align*}
 \left|
 \left(\frac{d}{ds}\right)^j
 \log A_{R,n}^{*}(s)\bigg|_{s=1}
 \right|
 &\ll_j
 \sum_{p\leq R}\sum_{\nu\geq1}
 \frac{|a_{R,n}(p^\nu)|
       (\nu\log p)^j}{p^\nu}\ll_j
 \sum_{p\leq R}\frac{(\log p)^j}{p}
 \ll_j(\log R)^j,
\end{align*}
where the sum over \(\nu\) is convergent uniformly by
\eqref{e:aRnlocalbound10}.  Repeated differentiation of
\(A_{R,n}^{*}=\exp(\log A_{R,n}^{*})\) therefore gives
\[
 \left|
 \left(\frac{d}{ds}\right)^j
 A_{R,n}^{*}(s)\bigg|_{s=1}
 \right|
 \ll_j(\log R)^{O_j(1)}.
\]
On the other hand,
\[
 (-1)^j
 \left(\frac{d}{ds}\right)^j
 A_{R,n}^{*}(s)\bigg|_{s=1}
 =
 \sum_{d\geq1}
 \frac{|a_{R,n}(d)|}{d}(\log d)^j.
\]
It follows that
\[
 \sum_{d\geq1}
 \frac{|a_{R,n}(d)|}{d}(\log d)^j
 \ll_j(\log R)^{O_j(1)}
 \qquad (j\geq0).
\]
Expanding \((1+\log d)^J\) and summing these estimates for
\(0\leq j\leq J\) proves \eqref{e:aRnmoments10}.

We next claim that the contribution of large \(d\) in
\eqref{e:hRconvolution10} is negligible.  Put $D_y=\exp(\sqrt{\log y}).$
Rankin's inequality and \eqref{e:aRnshifted10} give
\begin{align}\label{e:aRntail10}
 \sum_{d>D_y}\frac{|a_{R,n}(d)|}{d}
 &\leq
 D_y^{-1/\log(2R)}
 \sum_{d\geq1}
 \frac{|a_{R,n}(d)|}
 {d^{\,1-1/\log(2R)}}\notag\\
 &\ll
 \exp\left(
 -\frac{\sqrt{\log y}}{\log(2R)}
 \right)(\log R)^{O(1)}
 \ll_K(\log N)^{-K}
\end{align}
for every fixed \(K>0\).  The last estimate is uniform when
\(\log y\geq\log N/\log\log N\).

We can now transfer \eqref{e:meanvaluehinfty10} to \(h_R\).  By
\eqref{e:hRconvolution10},
\begin{align}\label{e:hRconvolutionsum10}
 \sum_{u\leq y}h_R(u)
 =
 \sum_{d\leq y}a_{R,n}(d)
 \sum_{v\leq y/d}h_\infty(v).
\end{align}
For \(d\leq D_y\), we have $\log(y/d)\asymp\log y$
and
\begin{align}\label{e:logyd10}
 \frac1{\sqrt{\log(y/d)}}
 =
 \frac1{\sqrt{\log y}}
 +O\left(\frac{\log d}{(\log y)^{3/2}}\right).
\end{align}
The factor \(\log d\) in \eqref{e:aRnmoments10} is needed precisely
to sum the error in \eqref{e:logyd10}.  Indeed,
\begin{align*}
 &\sum_{d\leq D_y}\frac{|a_{R,n}(d)|}{d}
 \left|
 \frac1{\sqrt{\log(y/d)}}-
 \frac1{\sqrt{\log y}}
 \right|\ll
 \frac1{(\log y)^{3/2}}
 \sum_{d\geq1}\frac{|a_{R,n}(d)|}{d}\log d
 \ll
 \frac{(\log R)^{O(1)}}{(\log y)^{3/2}}.
\end{align*}
The error term in \eqref{e:meanvaluehinfty10} gives the same bound,
by \eqref{e:aRnmoments10} with \(J=0\).

For \(d>D_y\), we use \(h_\infty\leq\tau_2\), and hence
\[
 \sum_{v\leq x}h_\infty(v)\ll x\log(2x).
\]
Together with \eqref{e:aRntail10}, this shows that the contribution of
\(d>D_y\) to \eqref{e:hRconvolutionsum10} is $O_K(y\log y\,(\log N)^{-K})$.
Choosing \(K\) sufficiently large, and using
\[
 R=(\log N)^C,\qquad
 \log y\geq\frac{\log N}{\log\log N},
\]
we conclude that
\begin{align}\label{e:hRmeanintermediate10}
 \sum_{u\leq y}h_R(u)
 =
 \left(
 \frac{c_\infty}{\sqrt\pi}
 \sum_{d\geq1}\frac{a_{R,n}(d)}d
 +o\left(
 \frac{C_R(n)}{B_R(n)\sqrt{\log R}}
 \right)
 \right)
 \frac{y}{\sqrt{\log y}}.
\end{align}
Here we also used \(C_R(n)/B_R(n)\geq1\), which follows from
\eqref{e:BRCR10} and the non-negativity of \(h_R\).

It remains to evaluate the coefficient in
\eqref{e:hRmeanintermediate10}.  By
\eqref{e:aRnmoments10} with \(J=0\), the series $ \sum_d\frac{a_{R,n}(d)}d$
converges absolutely.  We may therefore set \(s=1\) in
\eqref{e:defaRn10}, obtaining
\begin{align*}
 \sum_{d\geq1}\frac{a_{R,n}(d)}d
 &=
 \prod_{\substack{p\leq R\\p\equiv1\;(\mathrm{mod}\;4)}}
 \frac{
  1+\displaystyle\sum_{j\geq1}h_R(p^j)p^{-j}
 }{
  1+\displaystyle\sum_{j\geq1}h_\infty(p^j)p^{-j}
 }.
\end{align*}
Since \(h_R(p^j)=h_R(p)\) and
\(h_\infty(p^j)=h_\infty(p)\) for every \(j\geq1\),
\[
 \sum_{j\geq1}\frac{h_R(p^j)}{p^j}
 =
 h_R(p)\sum_{j\geq1}\frac1{p^j}
 =
 \frac{h_R(p)}{p-1},
\]
and similarly $ \sum_{j\geq1}\frac{h_\infty(p^j)}{p^j}
 =
 \frac{h_\infty(p)}{p-1}.$
Consequently,
\begin{align*}
 \sum_{d\geq1}\frac{a_{R,n}(d)}d
 &=
 \prod_{\substack{p\leq R\\p\equiv1\;(\mathrm{mod}\;4)}}
 \frac{1+h_R(p)/(p-1)}
      {1+h_\infty(p)/(p-1)}=
 \frac{
 \displaystyle
 \prod_{\substack{p\leq R\\p\equiv1\;(\mathrm{mod}\;4)}}
 \left(1+\frac{h_R(p)}{p-1}\right)}
 {\displaystyle
 \prod_{\substack{p\leq R\\p\equiv1\;(\mathrm{mod}\;4)}}
 \left(1+\frac{h_\infty(p)}{p-1}\right)}.
\end{align*}
By \eqref{e:BRCR10}, the product in the numerator is
\[
 \prod_{\substack{p\leq R\\p\equiv1\;(\mathrm{mod}\;4)}}
 \left(1+\frac{h_R(p)}{p-1}\right)
 =
 \frac{C_R(n)}{B_R(n)}.
\]
Multiplying by \(c_\infty\), we conclude that
\begin{align}\label{e:aRncoefficient10}
 c_\infty\sum_{d\geq1}\frac{a_{R,n}(d)}d
 &=
 \frac{C_R(n)}{B_R(n)}
 \frac{c_\infty}
 {\displaystyle
  \prod_{\substack{p\leq R\\p\equiv1\;(\mathrm{mod}\;4)}}
  \left(1+\frac{h_\infty(p)}{p-1}\right)}.
\end{align}
By \eqref{e:defcinfty10},
\[
 \prod_{\substack{p\leq R\\p\equiv1\;(\mathrm{mod}\;4)}}
 \left(1+\frac{h_\infty(p)}{p-1}\right)
 =
 \bigl(c_\infty+o(1)\bigr)
 \prod_{p\leq R}\left(1-\frac1p\right)^{-1/2}.
\]
Mertens' theorem therefore gives
\begin{align}\label{e:hinftyproduct10}
 \prod_{\substack{p\leq R\\p\equiv1\;(\mathrm{mod}\;4)}}
 \left(1+\frac{h_\infty(p)}{p-1}\right)
 =
 \bigl(e^{\gamma/2}c_\infty+o(1)\bigr)\sqrt{\log R}.
\end{align}
Combining \eqref{e:aRncoefficient10} and
\eqref{e:hinftyproduct10}, we obtain
\[
 c_\infty\sum_{d\geq1}\frac{a_{R,n}(d)}d
 =
 \left(e^{-\gamma/2}+o(1)\right)
 \frac{C_R(n)}{B_R(n)\sqrt{\log R}}.
\]
Substitution into \eqref{e:hRmeanintermediate10} proves
\eqref{e:meanvaluehR10}.

It remains to prove the upper bound
\eqref{e:meanvaluehRupper10}.  Shiu's theorem for non-negative
multiplicative functions \cite[Theorem~1]{Shiu}, applied with modulus
\(1\), gives
\[
 \sum_{u\leq y}h_R(u)
 \ll
 \frac{y}{\log y}
 \prod_{p\leq y}
 \left(1+\sum_{j\geq1}\frac{h_R(p^j)}{p^j}\right).
\]
It therefore suffices to show that
\begin{align}\label{e:hReulerupper10}
 \prod_{p\leq y}
 \left(1+\sum_{j\geq1}\frac{h_R(p^j)}{p^j}\right)
 \ll
 \frac{C_R(n)}{B_R(n)}
 \left(\frac{\log y}{\log R}\right)^{1/2}.
\end{align}
For \(y\geq R^2\), \eqref{e:defhR10} and \eqref{e:BRCR10} give
\begin{align*}
 \prod_{p\leq y}
 \left(1+\sum_{j\geq1}\frac{h_R(p^j)}{p^j}\right)
 &=
 \frac{C_R(n)}{B_R(n)}
 \prod_{\substack{R<p\leq y\\p\equiv1\;(\mathrm{mod}\;4)}}
 \left(1+\frac1{p-2}\right)\ll
 \frac{C_R(n)}{B_R(n)}
 \left(\frac{\log y}{\log R}\right)^{1/2},
\end{align*}
by Mertens' theorem in the progression \(1\pmod4\).  This proves
\eqref{e:hReulerupper10}, and hence
\eqref{e:meanvaluehRupper10}.

\emph{Step 4: The weighted sum over \(\ell\).}
Next, we prove an estimate for a weighted sum involving $h_R$. 
For every fixed \(A_1>0\), uniformly for
\(N^2(\log N)^{-A_1}\leq m\leq N^2\), one has
\begin{align}\label{e:Wsum10}
 &\sum_{\substack{\ell=2up_2\in L\\ \ell<m/2}}
 \frac{h_R(u)}{\ell(\log(m/\ell))^2}\leq
 \left(\frac{2\sqrt2e^{-\gamma/2}}{\sqrt\pi}
 \mathcal W(\alpha)+o(1)\right)
 \frac{C_R(n)}{B_R(n)}
 \frac1{\sqrt{\log R}(\log m)^{3/2}}.
\end{align}
Fix $N^2(\log N)^{-A_1}\leq m\leq N^2.$
We use partial summation twice.  Put
\begin{align}\label{e:rvdef}
     r=\frac{\log u}{\log m},
 \qquad
 v=\frac{\log p_2}{\log m}.
\end{align}
Since $\log m=2\log N+O_{A_1}(\log\log N)$ and $\alpha=\frac98,$
the conditions defining \(L\),
\[
 u<\frac{N^2}{2z^2},\qquad
 z\leq p_2<\left(\frac{N^2}{2u}\right)^{1/2},
 \qquad z=N^{1/\alpha},
\]
become, up to \(o(1)\)-changes in the endpoints,
\[
 0\leq r\leq1-\frac1\alpha=\frac19,
 \qquad
 \frac1{2\alpha}=\frac49
 \leq v\leq\frac{1-r}{2}.
\]
Moreover, \(\ell=2up_2\ll N^{10/9}\), whereas
\(m\geq N^2(\log N)^{-A_1}\).  Hence \(\ell<m/2\) is automatic for all
sufficiently large \(N\).
Writing out the condition \(\ell=2up_2\in L\), the left-hand side of
\eqref{e:Wsum10} is
\begin{align}\label{e:sumweighthr}
     \sum_{\substack{\ell=2up_2\in L}}
 \frac{h_R(u)}{\ell(\log(m/\ell))^2}
 =
 \frac12
 \sum_{\substack{u<N^2/(2z^2)\\ p\mid u\Rightarrow p\equiv1\pmod4}}
 \frac{h_R(u)}u
 \sum_{\substack{z\leq p_2<(N^2/(2u))^{1/2}\\ p_2\equiv3\pmod4}}
 \frac1{p_2(\log(m/(2up_2)))^2}. \qquad 
\end{align}
Put
\[
 \mathcal H_R(x)=
 \sum_{\substack{u\leq x\\p\mid u\Rightarrow p\equiv1\;(\mathrm{mod}\;4)}}
 h_R(u),
 \qquad
 \pi_3(y)=|\{p\leq y:p\equiv3\pmod4\}|.
\]
By partial summation, \eqref{e:sumweighthr} equals
\begin{align*}
 \frac12\int_{1^-}^{N^2/(2z^2)}\frac1x
 \left(
 \int_{z^-}^{(N^2/(2x))^{1/2}}
 \frac{d\pi_3(y)}
 {y(\log(m/(2xy)))^2}
 \right)d\mathcal H_R(x).
\end{align*}
Applying \eqref{e:meanvaluehR10} to the outer integral and the prime
number theorem in the progression \(3\pmod4\) to the inner integral
gives
\begin{align*}
 &\sum_{\substack{\ell=2up_2\in L\\ \ell<m/2}}
 \frac{h_R(u)}{\ell(\log(m/\ell))^2}\\
 &\quad\leq
 \left(\frac{e^{-\gamma/2}}{4\sqrt\pi}+o(1)\right)
 \frac{C_R(n)}{B_R(n)\sqrt{\log R}}
 \int_1^{N^2/(2z^2)}
 \frac{dx}{x\sqrt{\log x}}
 \int_z^{(N^2/(2x))^{1/2}}
 \frac{dy}
 {y\log y\,(\log(m/(2xy)))^2}.
\end{align*}
With \eqref{e:rvdef} in mind, 
we get
\[
 \sum_{\substack{\ell=2up_2\in L\\ \ell<m/2}}
 \frac{h_R(u)}{\ell(\log(m/\ell))^2} \leq \left(\frac{e^{-\gamma/2}}{4\sqrt\pi}+o(1)\right)
 \frac{C_R(n)}{B_R(n)}
 \frac{I_\alpha}
 {\sqrt{\log R}(\log m)^{3/2}},
\]
where
\[
 I_\alpha=
 \int_0^{1-1/\alpha}\frac{dr}{\sqrt r}
 \int_{1/(2\alpha)}^{(1-r)/2}
 \frac{dv}{v(1-r-v)^2}.
\]
The factor \(1/4\) comes from the prefactor \(1/2\) in
\(\ell=2up_2\), together with the density \(1/2\) of primes
\(p_2\equiv3\pmod4\). For fixed \(r\), put
\(t=2\alpha(1-r)\) and \(v=(1-r)y\).  Since
\[
 \int_{1/t}^{1/2}\frac{dy}{y(1-y)^2}
 =
 \frac{t-2+(t-1)\log(t-1)}{t-1},
\]
we get
\[
 I_\alpha
 =
 2\alpha\int_2^{2\alpha}
 \frac{t-2+(t-1)\log(t-1)}
 {t^2(t-1)(1-t/(2\alpha))^{1/2}}\,dt
 =
 8\sqrt2\,\mathcal W(\alpha).
\]
This gives \eqref{e:Wsum10} on the range
\(\log u\geq\log N/\log\log N\).

Next consider \(R^2\leq u<\exp(\log N/\log\log N)\).  Repeating the
preceding partial-summation argument, but using the upper bound
\eqref{e:meanvaluehRupper10} in place of the asymptotic
\eqref{e:meanvaluehR10}, shows that this range contributes
\[
 \ll
 \frac{C_R(n)}{B_R(n)}
 \frac1{\sqrt{\log R}(\log m)^{3/2}}
 \int_0^{1/(2\log\log N)+o(1)}
 \frac{dr}{\sqrt r}
 \int_{1/(2\alpha)}^{(1-r)/2}
 \frac{dv}{v(1-r-v)^2}.
\]
The inner integral is \(O_\alpha(1)\), since \(v\geq1/(2\alpha)\) and
\(1-r-v\geq(1-r)/2\gg_\alpha1\) on this range.  Hence the last display is
\[
 \ll
 \frac{C_R(n)}{B_R(n)}
 \frac1{\sqrt{\log R}(\log m)^{3/2}}
 \int_0^{1/(2\log\log N)+o(1)}
 \frac{dr}{\sqrt r}
 \ll
 (\log\log N)^{-1/2}
 \frac{C_R(n)}{B_R(n)}
 \frac1{\sqrt{\log R}(\log m)^{3/2}}.
\]
This is \(O((\log\log N)^{-1/2})\) times the right-hand side of
\eqref{e:Wsum10}.  The remaining range \(u<R^2\) is even smaller; indeed
\[
 \sum_{u<R^2}\frac{h_R(u)}u\ll(\log R)^{O(1)},
\]
and the corresponding \(p_2\)-sum is \(\ll(\log N)^{-2}\).  Since
\(R=(\log N)^C\), this contribution is also
\[
 o\left(\frac{C_R(n)/B_R(n)}
 {\sqrt{\log R}(\log m)^{3/2}}\right).
\]
This completes the proof of the claim.

\emph{Step 5: Replacing \(H_{\ell,R}\) by \(H_\ell\).}
We now return to \eqref{e:uppernormalization10}.  First note that
\begin{align}\label{e:hroverl}
 \sum_{\ell=2up_2\in L}\frac{h_R(u)}{\ell}
 \leq
 \frac12
 \sum_{u<N^{2/9}}\frac{h_R(u)}u
 \sum_{z\leq p_2<N}\frac1{p_2}
 \ll(\log N)^{O(1)}.
\end{align}
Indeed \(h_R\) is divisor-bounded uniformly in \(n\), so
\(\sum_{u<N^{2/9}}h_R(u)/u\ll(\log N)^{O(1)}\), while
\(\sum_{z\leq p_2<N}p_2^{-1}\ll\log\log N\).  On the other hand,
\eqref{e:lrtol} gives, uniformly in \(\ell\),
\[
 |H_{\ell,R}(n-1)-H_\ell(n-1)|\ll NR^{-1/2}.
\]
Therefore
\begin{align}\label{e:smoothingerror10}
 \sum_{\ell\in L}
 \frac{h_R(u)}\ell
 |H_{\ell,R}(n-1)-H_\ell(n-1)|
 \ll NR^{-1/2}(\log N)^{O(1)}.
\end{align}
We now insert this error into \eqref{e:uppernormalization10}.  As discussed
above, \(z_\ell=(N^2/\ell)^{1/4}\gg N^{2/9}\).  Thus
\(\log z_\ell\gg\log N\) uniformly in \(\ell\in L\).  Therefore the total
contribution of this error term is
\[
 \ll
 B_R(n)\log R
 \sum_{\ell\in L}
 \frac{h_R(u)}{\ell\log z_\ell}
 |H_{\ell,R}(n-1)-H_\ell(n-1)|
 \ll
 B_R(n)NR^{-1/2}(\log N)^{O(1)}
 \frac{\log R}{\log N}.
\]
By \eqref{e:lowernormalization10} and \eqref{e:H1elementary10},
\[
 X_n\Theta_n(z)
 \gg
 C_R(n)N(\log N)^{-O_A(1)}
 \left(\frac{\log R}{\log N}\right)^{1/2}.
\]
Since \(B_R(n)\leq C_R(n)\) and \(R=(\log N)^C\), the contribution of the
error from replacing \(H_{\ell,R}\) by \(H_\ell\) is \(o(X_n\Theta_n(z))\)
after choosing \(C\) sufficiently large.

\emph{Step 6: Integrating the weighted \(\ell\)-sum.}
We now expand the integral defining \(H_\ell(n-1)\).  Put
\(m=m(y):=n-1-y^2\).  Then
\begin{align*}
 \sum_{\ell\in L}
 \frac{h_R(u)H_\ell(n-1)}
 {\ell\log(N^2/\ell)}
 &=
 \int_{\substack{2\leq y\leq N\\ m(y)>0}}
 \frac{dy}{\log y}
 \sum_{\substack{\ell=2up_2\in L\\ \ell<m(y)/2}}
 \frac{h_R(u)}
 {\ell\log(N^2/\ell)\log(m(y)/\ell)}  \\
 &\leq
 \int_{\substack{2\leq y\leq N\\ m(y)>0}}
 \frac{dy}{\log y}
 \sum_{\substack{\ell=2up_2\in L\\ \ell<m(y)/2}}
 \frac{h_R(u)}
 {\ell(\log(m(y)/\ell))^2},
\end{align*}
since \(m(y)\leq N^2\).  Using \eqref{e:hroverl} and
\[
 \bigl|\{y:0<n-1-y^2<N^2(\log N)^{-A_1}\}\bigr|
 \ll
 \frac{N^2(\log N)^{-A_1}}{\sqrt n}
 \ll_{A_0}N(\log N)^{-A_1+A_0/2},
\]
the contribution of \(m(y)<N^2(\log N)^{-A_1}\) is
\[
 O_{A_0}\bigl(N(\log N)^{-A_1+A_0/2+O(1)}\bigr),
\]
which is negligible after choosing \(A_1\) sufficiently large. On the retained range \(m(y)\geq N^2(\log N)^{-A_1}\), we have $ \log m(y)=2\log N+O_{A_1}(\log\log N).$
Thus \eqref{e:Wsum10} gives, uniformly on this range,
\[
 \sum_{\substack{\ell=2up_2\in L\\2\ell\leq m(y)}}
 \frac{h_R(u)}{\ell(\log(m(y)/\ell))^2}
 \leq
 \left(\frac{e^{-\gamma/2}}{\sqrt\pi}
 \mathcal W(\alpha)+o(1)\right)
 \frac{C_R(n)}{B_R(n)}
 \frac{1}{\sqrt{\log R}(\log N)^{3/2}}.
\]
Here the factor \(2\sqrt2\) in \eqref{e:Wsum10} cancels with
the factor \(2^{3/2}\) arising from \((\log m(y))^{3/2}\).

For these \(y\), we have
\(2\leq n-y^2=m(y)+1\leq N^2\), so the retained range is
contained in the integration domain defining \(H_1(n)\).
Since \(\log(n-y^2)\leq2\log N\), it follows that
\[
 \int_{\substack{2\leq y\leq N\\
                 m(y)\geq N^2(\log N)^{-A_1}}}
 \frac{dy}{\log y}
 \leq
 2\log N
 \int_{\substack{2\leq y\leq N\\2\leq n-y^2\leq N^2}}
 \frac{dy}{\log y\,\log(n-y^2)}
 =
 2\log N\,H_1(n).
\]
Integrating the preceding bound therefore gives the required
estimate on the retained range.

To include the discarded range, recall that
\(C_R(n)/B_R(n)\geq1\) by \eqref{e:BRCR10}, while
\eqref{e:H1elementary10} gives
\(H_1(n)\gg_{A_0}N(\log N)^{-A_0/2-2}\).
Since \(R=(\log N)^C\), choosing \(A_1\) sufficiently large
in terms of \(A_0\) makes the previously estimated error satisfy
\[
 N(\log N)^{-A_1+A_0/2+O(1)}
 =
 o\left(
 \frac{C_R(n)}{B_R(n)}
 \frac{H_1(n)}{\sqrt{\log R\log N}}
 \right).
\]
Combining the two ranges, we obtain
\begin{align}\label{e:integratedW10}
 \sum_{\ell\in L}
 \frac{h_R(u)H_\ell(n-1)}
 {\ell\log(N^2/\ell)}
 &\leq
 \left(\frac{2e^{-\gamma/2}}{\sqrt\pi}
 \mathcal W(\alpha)+o(1)\right)
 \frac{C_R(n)}{B_R(n)}
 \frac{H_1(n)}{\sqrt{\log R\log N}}.
\end{align}

\emph{Step 7: Concluding the proof.}
We now sum \eqref{e:uppernormalization10}.  The smoothing error
\eqref{e:smoothingerror10} allows us to replace \(H_{\ell,R}\) by
\(H_\ell\) at a cost \(o(X_n\Theta_n(z))\).  Since
\(\log z_\ell=\frac14\log(N^2/\ell)\), we get
\begin{align*}
 \sum_{\ell\in L}X_n(\ell)\Theta_{n,\ell}(z_\ell)
 &\leq
 \bigl(6+o(1)\bigr)
 B_R(n)\log R
 \sum_{\ell\in L}
 \frac{h_R(u)H_\ell(n-1)}{\ell\log z_\ell}
 +o(X_n\Theta_n(z))\\
 &=
 \bigl(24+o(1)\bigr)
 B_R(n)\log R
 \sum_{\ell\in L}
 \frac{h_R(u)H_\ell(n-1)}
 {\ell\log(N^2/\ell)}
 +o(X_n\Theta_n(z)).
\end{align*}
Using \eqref{e:integratedW10}, this gives
\begin{align}\label{e:upperWfinal10}
 \sum_{\ell\in L}X_n(\ell)\Theta_{n,\ell}(z_\ell)
 \leq
 \left(\frac{48e^{-\gamma/2}}{\sqrt\pi}
 \mathcal W(\alpha)+o(1)\right)
 C_R(n)H_1(n)
 \left(\frac{\log R}{\log N}\right)^{1/2}.
\end{align}
Lemma~\ref{l:positiveH9} and \eqref{e:H1elementary10} give
\[
 H_{1,R}(n)=(1+o(1))H_1(n).
\]
Since \(\log z=\alpha^{-1}\log N\),
\eqref{e:lowernormalization10} becomes
\[
 X_n\Theta_n(z)
 =
 \bigl(3\sqrt\alpha+o(1)\bigr)
 C_R(n)H_1(n)
 \left(\frac{\log R}{\log N}\right)^{1/2}.
\]
Dividing \eqref{e:upperWfinal10} by this expression proves
\eqref{e:wumatomakisum10}, since
\[
 \frac{48e^{-\gamma/2}/\sqrt\pi}{3\sqrt\alpha}
 =
 \frac{16e^{-\gamma/2}}{\sqrt{\pi\alpha}}.
\]
\end{proof}

\begin{proof}[Proof of Proposition~\ref{p:sieveneedupper8}]
By Lemma~\ref{l:fixedsievesets10}, it is enough to estimate the fixed
sum over \(\ell\in L\) and \(W_\ell\).  For each such \(\ell\), apply
Lemma~\ref{l:presieving} to \eqref{e:uppersequence10}, with
\[
 D_\ell=\frac{(N^2/\ell)^{1/2}}
 {12(\log N)^B},
 \qquad z_\ell=(N^2/\ell)^{1/4}.
\]
Since \(\ell\ll N^{10/9}\), we have
\(z_\ell\gg N^{2/9}\). Together with \(D_0=N^{o(1)}\) and
\(D_\ell=z_\ell^2/(12(\log N)^B)\), this verifies the level
conditions in Lemma~\ref{l:presieving} uniformly in \(\ell\),
with \(A=2\). The outer dimension is \(1\), by \eqref{e:outerdimensions10}, and
\[
 \frac{\log(D_\ell/D_0)}{\log z_\ell}=2+o(1)
\]
uniformly in \(\ell\in L\).   It follows that
\begin{align*}
 S(M_n(\ell),W_\ell,z_\ell)
 &\leq X_n(\ell)\Theta_{n,\ell}(z_\ell)
 \bigl(F_1(2)+o(1)\bigr)+\sum_{d<D_\ell}c_{d,\ell}^+
 r_{12d,\ell}(n-1),
\end{align*}
where \(|c_{d,\ell}^+|\leq1\) and the coefficients are independent of
\(n\).  Summing over \(\ell\in L\) and applying
Lemma~\ref{l:wumatomakisum} to the main terms gives
\begin{align}\label{e:upperbound10}
 \sum_{\ell\in L_n}S(M_n(\ell),W_{\ell,n},z_\ell)
 &\leq
 \bigl(c_\alpha F_1(2)\mathcal W(\alpha)+o(1)\bigr)
 X_n\Theta_n(z)+\mathcal E_2(n),
\end{align}
where
\[
 \mathcal E_2(n)=
 \sum_{\ell\in L}\sum_{d<D_\ell}
 c_{d,\ell}^+r_{12d,\ell}(n-1)+O(N^{1/4}).
\]
The bounds defining \(L\) and \eqref{e:alphachoice8} give
\(\ell\ll N^{10/9}\).  Thus Lemma~\ref{l:circmethod} applies with
\(L_0=N^{6/5}\) and, for example, \(\delta=3/4\), and yields
\[
 \sum_{n<N^2}|\mathcal E_2(n)|^2
 \ll_{A_0,K}N^4(\log N)^{-K}.
\]
Here we used \eqref{e:fixedupper10} with \(\varepsilon=1/4\).  The
added term contributes \(O(N^{5/2})\) after squaring and summing over
\(n<N^2\), and is therefore absorbed by the preceding mean-square
bound.  This proves the proposition.
\end{proof}

We finish by proving the size of the common main term and the numerical
inequality used in Section~\ref{s:reductionsGB}.

\begin{proof}[Proof of Lemma~\ref{l:maintermsize8}]
For \(p>3\), every factor occurring in \(C_R(n)\) is positive and is
at least \(1-c/p\) for an absolute constant \(c\).  Mertens' theorem
therefore gives
\begin{align}\label{e:CRlower10}
 C_R(n)\gg(\log R)^{-c}.
\end{align}
Lemma~\ref{l:positiveH9} gives
\(H_{1,R}(n)=H_1(n)+O(NR^{-1/2})\). Taking \(C\) sufficiently large in terms of \(A_0\), \eqref{e:H1elementary10} therefore implies
\(H_{1,R}(n)\geq H_1(n)/2\).

Combining \eqref{e:lowernormalization10}, \eqref{e:CRlower10}, and
\eqref{e:H1elementary10}, and recalling that
\(\log R=C\log\log N\) and \(\log z=(8/9)\log N\), we obtain
\[
 X_n\Theta_n(z)
 \gg_{A_0,K} N(\log N)^{-C_{A_0}}
\]
for some \(C_{A_0}>0\), uniformly in \eqref{e:nrange10}.  This proves
Lemma~\ref{l:maintermsize8}.
\end{proof}

\begin{proof}[Proof of Lemma~\ref{l:numericalcomparison}]
With the normalization of \cite[(1.8)--(1.10)]{I80},
\[
 f_{1/2}(s)=
 \left(\frac{e^\gamma}{\pi s}\right)^{1/2}
 \int_1^s\frac{dt}{\sqrt{t(t-1)}},
 \qquad 1<s\leq3,
 \qquad F_1(2)=e^\gamma.
\]
At \(s=9/8\), the substitution \(t=\cosh^2 u\) gives
\[
 f_{1/2}(9/8)
 =\left(\frac{8e^\gamma}{9\pi}\right)^{1/2}\log2.
\]
It remains to bound \(\mathcal W(9/8)\). The change of variables
\(r=1-4t/9\) in \eqref{e:Wdef8} gives
\[
 8\sqrt2\,\mathcal W(9/8)
 =\int_0^{1/9}\frac1{\sqrt r}
   \int_{4/9}^{(1-r)/2}\frac{dv}{v(1-r-v)^2}\,dr.
\]
On the inner interval, both \(v\) and \(1-r-v\) are at least \(4/9\).
The inner integrand is therefore at most \(729/64\), and its interval
has length \((1/9-r)/2\). Consequently,
\[
 8\sqrt2\,\mathcal W(9/8)
 \leq\frac{729}{128}\int_0^{1/9}\frac{1/9-r}{\sqrt r}\,dr
 =\frac9{32}.
\]
Using the stated value of \(c_{9/8}\) and \(F_1(2)=e^\gamma\), we obtain
\[
 \Delta\geq\frac{e^{\gamma/2}}{\sqrt\pi}
   \left(\frac{2\sqrt2}{3}\log2-\frac38\right)>0.
\]
The last inequality follows from \(\log2\geq1/2\) and \(\sqrt2>9/8\).
This proves the numerical comparison.
\end{proof}

\bibliography{BibPrimeSquares}

\end{document}